\documentclass[reqno,10pt]{amsart}
\usepackage{fullpage}
\usepackage{rotating}

\usepackage{amssymb}

\usepackage{amsfonts}
\usepackage{amscd}

\usepackage{enumerate}
\usepackage{graphicx}
\usepackage{subfigure}
\usepackage[english]{babel}
\usepackage{amsthm}
\usepackage[T1]{fontenc}
\usepackage[latin1]{inputenc} 
\usepackage{lmodern} 
\usepackage{url}

\usepackage{hyperref}
\hypersetup{ 
  pdfauthor   = {Reynald Lercier, Christophe Ritzenthaler}
  pdftitle    = {Hyperelliptic curves and their invariants: geometric, arithmetic and algorithmic aspects},
  pdfsubject  = {},
  pdfkeywords = {},
  backref=true, pagebackref=true, hyperindex=true, colorlinks=true,
  breaklinks=true, urlcolor=blue, linkcolor=blue, citecolor=blue, bookmarks=true,
  bookmarksopen=true}
\usepackage{multirow}

\usepackage[ruled,vlined,boxed,linesnumbered]{algorithm2e}

\usepackage{breqn}
\newcounter{myequation}[equation]

\makeatletter
\gdef\th@break{\normalfont\slshape
  \def\@begintheorem##1##2{\item[%
       \rlap{\vbox{\hbox{\hskip \labelsep\theorem@headerfont ##1\ ##2}%
                   \hbox{\strut}}}]}%
\def\@opargbegintheorem##1##2##3{%
  \item[\rlap{\vbox{\hbox{\hskip \labelsep \theorem@headerfont
                     ##1\ ##2\ ##3}%
                    \hbox{\strut}}}]}}
\makeatother

\usepackage[all]{xy}

\theoremstyle{plain}
\newtheorem{theorem}{Theorem}[section]

\newtheorem{proposition}[theorem]{Proposition}
\newtheorem{lemma}[theorem]{Lemma}
\newtheorem{corollary}[theorem]{Corollary}

\theoremstyle{definition}
\newtheorem{definition}[theorem]{Definition}

\theoremstyle{remark}
\newtheorem{remark}[theorem]{Remark}
\newtheorem{example}[theorem]{Example}

\makeatletter

\newcommand{\stratum}[1]{%
  \subsubsection*{{\underline{\textsc{stratum #1}}}}
}

\DeclareMathOperator{\Aut}{Aut}

\DeclareMathOperator{\Div}{div}

\DeclareMathOperator{\Gal}{Gal}
\DeclareMathOperator{\GL}{GL}

\DeclareMathOperator{\Jac}{Jac}

\DeclareMathOperator{\PGL}{PGL}

\DeclareMathOperator{\SL}{SL}

\DeclareMathOperator{\Spec}{Spec}

\newcommand{\Id}{\textrm{id}}

\newcommand{\Hm}{\mathsf{H}}
\newcommand{\Mm}{\mathsf{M}}

\newcommand{\Xm}{\mathsf{X}}

\newcommand{\bb}{\mathcal{B}}
\newcommand{\cc}{\mathcal{C}}

\newcommand{\ff}{\mathcal{F}}
\newcommand{\ggm}{\mathcal{G}}
\newcommand{\hh}{\mathcal{H}}
\newcommand{\ii}{\mathcal{I}}
\newcommand{\jj}{\mathcal{J}}
\newcommand{\nn}{\mathcal{N}}

\newcommand{\qq}{\mathcal{Q}}

\newcommand{\FF}{\mathbb{F}}
\newcommand{\CC}{\mathbb{C}}

\newcommand{\NN}{\mathbb{N}}
\newcommand{\PP}{\mathbb{P}}
\newcommand{\QQ}{\mathbb{Q}}
\newcommand{\RR}{\mathbb{R}}

\newcommand{\ZZ}{\mathbb{Z}}

\newcommand{\bfM}{\mathbf{M}}

\newcommand{\CG}{\mathbf{C}}
\newcommand{\DG}{\mathbf{D}}
\newcommand{\UG}{\mathbf{U}}
\newcommand{\VG}{\mathbf{V}}
\newcommand{\SG}{\mathbf{S}}

\def\ie{\textit{i.e. }}

\title[Hyperelliptic curves and their invariants]{Hyperelliptic curves and their invariants: geometric, arithmetic and algorithmic aspects}
\author{Reynald Lercier}

\address{%
  \textsc{DGA MI}, %
  La Roche Marguerite, %
  35174 Bruz, %
  France. %
}%
\address{%
  IRMAR, %
  Universit\'e de Rennes 1, %
  Campus de Beaulieu, %
  35042 Rennes, %
  France. %
}%
\email{reynald.lercier@m4x.org}

\author{Christophe Ritzenthaler}
\address{Institut de Math{\'e}matiques de Luminy, %
         UMR 6206 du CNRS, %
         Luminy, Case 907, 13288 Marseille, France.}
\email{ritzenth@iml.univ-mrs.fr}

\thanks{Both authors acknowledge support by grant ANR-09-BLAN-0020-01.}
\date{\today}

\begin{document}

\begin{abstract}
  We apply classical invariant theory of binary forms to explicitly
  characterize isomorphism classes of hyperelliptic curves of small genus and,
  conversely, propose algorithms for reconstructing hyperelliptic models from
  given invariants. We focus on genus 3 hyperelliptic curves. Both geometric
  and arithmetic aspects are considered.
\end{abstract}

\maketitle

\tableofcontents
\newpage

\section*{Introduction}
\label{sec:introduction}

Invariant theory played a central role in 19th century algebra and geometry,
yet many of its techniques and algorithms were practically forgotten by the
middle of the 20th century and replaced by the abstract and powerful machinery
of modern algebraic geometry. However, motivated by computational applications
to cryptography, robotics, coding theory, \textit{etc.}\,, the classical
invariant theory has come to a renaissance. Among classical groups, the
natural action of $\SL_2$ on binary forms has received most attention. One
reason is the remarkable formalism developed by Gordan in 1868 to compute a
finite set of generators of invariants.  The other reason is the application
of hyperelliptic curves in cryptography, especially for the so-called CM
methods (see~\cite{cohen-handbook}).\smallskip

When $K$ is an algebraically closed field of characteristic $0$ or a prime $p
\ne 2$, hyperelliptic curves of genus $g>1$ are indeed naturally related to
the action of $\SL_2(K)$ on binary forms of even degree $n=2g+2$.  Thus,
isomorphisms between curves with equation $y^2=f(x)$ and $f \in K[x]$ of
degree $n$ are globally determined by a $\SL_2(K)$ action on $f$ (see
Section~\ref{subsec:reconstr}). Thus the $K$-isomorphism classes of
hyperelliptic curves of genus $g$ can be represented by the values of a finite
set of invariants for $SL_2(K)$. In this way, exploring properties of the
invariants lead to effective results on the geometry and arithmetic of the
hyperelliptic moduli space.\medskip

On the geometric level, to make explicit the representation of the classes by
invariants, we have to tackle a double task: compute a set of generators
$\{I_i\}$ of invariants and reciprocally construct a curve from given values
of these invariants. We call the latter the \emph{reconstruction} phase.  As a
by-product, we also want to be able to read some geometric information, in
particular the automorphism group of a curve, from the invariants.\smallskip

The first issue can be addressed thanks to Gordan's method which is based on a
differential operator called \emph{transvectant}, see
Section~\ref{sec:rappel}. Contrary to the genus $2$ case which is described by
$3$ algebraically independent invariants, Shioda \cite{shioda67} gave in genus
$3$ a basis of $9$ invariants (which we call \emph{Shioda invariants}), the
first $6$ being algebraically independent and the last three related to the
others by $5$ explicit relations. Besides this, unlike the genus $2$ case and
the classical Igusa invariants, the discriminant is not an element of this
basis. We therefore decided not to use the usual representation based on
`absolute invariants' and to switch to a weighted projective space of
invariants for which we rely on some specific algorithms to test equality or
to create points. Note that although we restrict to genus $3$ in the present
article, these algorithms apply to any weighted projective space and are
therefore useful for hyperelliptic curves of higher genus too.\smallskip

The main algorithm for the second issue relies on Mestre's method which he
exposes for $g=2$ in~\cite{mestre}. It is based on computations going back
to~\cite[§103]{clebsch} (see Section~\ref{sec:reconstruction}) and uses a
generalization of invariants called \emph{covariants} (see Definition
\ref{def:covariants}).  Roughly speaking, starting from three order $2$
covariants, one constructs a plane conic $\qq$ and a plane degree $g+1$ curve
$\hh$ whose coefficients are invariants, hence expressible in terms of the
generators $\{I_i\}$. After specialization at given values of invariants, if
$\qq$ is not singular, the degree $2$-cover of $\qq$ ramified at the
$2(g+1)=n$ points of intersection of $\qq$ and $\hh$ is a hyperelliptic curve
with invariants equal to the initial values. In order to make this practical,
one of the main computational difficulty is to find the formal expressions of
the coefficients in terms of $\{I_i\}$ as the degree and number of variables
are for $g=3$ already quite large. We by-passed the difficulty using an
evaluation-interpolation strategy.\smallskip

For genus $2$, this algorithm enabled Mestre to reconstruct hyperelliptic
curves $C$ with no \emph{extra-automorphism}, \ie $\Aut(C)$ is generated by
the hyperelliptic involution $\iota$. Indeed, it was proved in \cite{bolza}
that the `classical' $\qq$ is always non-singular in this case. It is however
always singular when $C$ has extra-automorphisms and it can not be used
anymore. Cardona and Quer in~\cite{CAQU} completed the picture by a different
choice of order $2$ covariants which lead to another conic $\qq$ non-singular
in the case $\Aut(C) \simeq (\ZZ/2\ZZ)^2$ (for bigger automorphism groups,
explicit parametrizations were already known). The two authors of the present
article have implemented these constructions with the computational algebra
software \textsc{magma} and generalized them to the remaining fields of small
characteristics ($3$ and the most difficult case $5$)~\cite{LR08}.  Genus $2$
can be considered as solved.\smallskip

It is thus natural to turn to genus $3$. Here also, the issues that we tackle
are naturally stratified by the automorphism group of the curve (see
Table~\ref{tab:auto} and Figure~\ref{fig:lattice} for the lattice of
automorphism groups in genus $3$). In Section~\ref{sec:ident-loci-reconstr},
we find equations for all the strata and show how to reconstruct a curve with
given invariants. These two questions are indeed intertwined.\smallskip

Starting from normal models as in Table~\ref{tab:auto}, we get necessary
algebraic conditions for the invariants to define a curve with a certain
automorphism group. We determine these equations by evaluations and
interpolations, the same approach as the one we followed for finding the
expressions of the coefficients of $\qq$ and $\hh$. Since we only use Gröbner
techniques to reduce bases, we were able to obtain equations for all the
strata, even for the dimension 3 stratum $(\ZZ/2\ZZ)^2$.  In order to check
that these conditions are also sufficient, we reconstruct a curve from given
invariants and check that its Shioda invariants are equal to the original
ones. This last step involves calculations in the quotient ring of $\QQ[J_2,
..., J_{10}]$ by the ideal defined by the stratum equations. Even if modern
computational algebra software can handle them, keeping polynomials reduced to
normal forms modulo the ideal of relations yields several hours of
computations on a powerful computer, at least for the strata of dimension 2 or
3.  A reasonably optimized program written in \textsc{magma} that implements
the corresponding computations is available on the web page of the authors for
independent checks.\smallskip

The reconstruction step for most strata is carried out by `inverting' the
expressions of the invariants in terms of the parameters of normal models
modulo the stratum equations. For strata of dimension less or equal to $1$, we
give models whose coefficients are algebraic expressions in terms of the
invariants. For the dimension $2$ stratum with automorphism group
$(\ZZ/2\ZZ)^3$, we can still work out these computations at the price of a
`cubic extension' (see Lemma \ref{lemma:C2p3}). For the dimension $2$ stratum
$\ZZ/4\ZZ$, we exhibit $5$ conics among which at least one is always
non-singular and use Mestre's method for this stratum.  The dimension $3$
stratum $(\ZZ/2\ZZ)^2$ is more challenging.  One can show
(Lemma~\ref{lemma:reconstr-D4}) that \emph{any} choice of 3 covariants in the
set of the 14 fundamental covariants of order $2$ ($364$ possibilities) leads
to a singular conic, hence Cardona and Quer's patch is not possible for curves
in this stratum. Our computational approach yields, in addition to a set of
$24$ necessary equations for the invariants to define a curve with
automorphism group $(\ZZ/2\ZZ)^2$, an explicit reconstruction at the price
here of a `degree 8 extension' (Lemma \ref{lemma:D4}).  We noticed furthermore
that the singularity of the $364$ conics is equivalent to the nullity of only
$19$ determinants.  Obviously, the locus where these determinants
simultaneously vanish contains the stratum $(\ZZ/2\ZZ)^2$. We show that it is
actually \emph{equal} to the $24$ stratum equations for $(\ZZ/2\ZZ)^2$.  As a
by-product, the reconstruction of curves with no extra-automorphisms can
therefore be achieved thanks to Mestre's method by picking one of the
non-singular conics $\qq$ among the $19$ fixed ones.\medskip

So far we have avoided arithmetic issues by working over an algebraically
closed field but new challenging and deep issues arrive when one consider over
which `minimal field' these constructions can be achieved. Assume for
simplicity that $k$ is a field of characteristic $0$ and let $C$ be a curve
defined over $k$ (see a more general framework in
Section~\ref{subsec:def}). One can consider the intersection of all the
subfields $k'$ of $K=\bar{k}$ over which there exists a curve $K$-isomorphic
to $C$ ($k'$ is called a \emph{field of definition}). This field is called the
\emph{field of moduli} and denoted $\bfM_C$. If it is a field of definition,
it is the minimal field of definition of $C$. Unfortunately, this is not
automatically the case if for instance $C$ has a non-trivial automorphism
group. Therefore hyperelliptic curves are highly concerned with this issue:
when is $\bfM_C$ a field of definition~? On top of this, hyperelliptic curves
lead to a refined question that we want to address too. Let remember that, in
full generality, a hyperelliptic curve $C/k$ is a curve with a degree $2$
morphism from $C$ to a non-singular plane conic $Q$. If $Q$ has a point and,
assuming the characteristic of $k$ not $2$, one can write $C/k : y^2=f(x)$
with $f \in k[x]$. We say that $C$ can be \emph{hyperelliptically defined over
  $k$} or that $C$ admits a \emph{hyperelliptic equation over $k$}. This is
obviously the case if $k$ is algebraically closed or finite (and therefore
people often use this property as a definition) but when $k$ is arbitrary,
there might be again an obstruction.\smallskip

When $g=2$ (or more generally even), Mestre showed that the two questions are
equivalent. Moreover in the case where $C$ has no extra-automorphism, he
showed that $\bfM_C$ is a field of definition if and only if the conic $\qq$
constructed from the covariants has a rational point. When $C$ has
extra-automorphisms, Cardona and Quer were able to exhibit hyperelliptic
equations over $\bfM_C$. We implemented and completed their results for fields
of characteristic $2$, $3$ and $5$ \cite{LR08}.\smallskip

For $g=3$, we address both questions, from a theoretical and computational
point of view. Let us give some of the results that we have obtained according
to the stratification of the moduli space by the automorphism groups (see
Section~\ref{subsec:descentg3}).\smallskip
\begin{itemize}
\item For dimension $0$ and dimension $1$ strata, since $\Aut(C)/\langle \iota
  \rangle$ is not cyclic, Huggins' results in~\cite{hugginsphd} show that
  there is a hyperelliptic equation over $\bfM_C$. This is confirmed and made
  explicit by the computations we performed for the reconstruction. In theses
  cases we can even exhibit parametrized models over $\bfM_C$.\smallskip
\item The dimension $2$ case $\Aut(C) \simeq (\ZZ/2\ZZ)^3$ is theoretically
  covered by~\cite{hugginsphd} and therefore there exists a hyperelliptic
  equation over $\bfM_C$ (see Remark \ref{remark:fuertes} for the controversy
  on this subject).  In the present article, however we only reconstruct the
  curve hyperelliptically over at most a cubic extension of $\bfM_C$. There is
  a real difficulty to perform an explicit descent over $\bfM_C$ as geometric
  isomorphisms between our curve and a model over $\bfM_C$ might be defined
  over a degree $24$ of $\bfM_C$.  In a forthcoming article, we work out an
  algorithm based on covariants to obtain a model over $\bfM_C$ in an
  efficient way.
\item The dimension $2$ case $\Aut(C) \simeq \ZZ/4\ZZ$ is not covered by the
  general result of~\cite{hugginsphd}. However using the special form of the
  ramification signature, we can show that $C$ can always be hyperelliptically
  defined over $\bfM_C$. To make this result explicit, we use the fact that
  there is always a non-singular conic $\qq$. Although this conic has not
  necessarily a rational point, we can make use of the special shapes of $\qq$
  and $\hh$ to perform an explicit hyperelliptic descent.\smallskip
\item For the dimension $3$ case $\Aut(C) \simeq (\ZZ/2\ZZ)^2$, Huggins
  constructed examples of genus $3$ curves over $\bar{\QQ}$ which cannot be
  defined over their field of moduli. Hence, we did not explore this case
  further and our reconstruction takes place over a degree $8$ extension at
  most.\smallskip
\item Finally, for the dimension $5$ case $\Aut(C) \simeq \ZZ/2\ZZ$, we show that
  $\bfM_C$ is always a field of definition. This is more generally true for
  hyperelliptic curves with no extra-automorphisms of odd genus. However, now
  $C$ has not automatically a hyperelliptic model over $\bfM_C$. In genus $3$,
  we show that $C$ has a hyperelliptic equation over $\bfM_C$ if and only if
  one (equivalently, all of the) non-singular conic $\qq$ has a rational point
  over $\bfM_C$. As we proved that such a conic always exists, the method is explicit as well.\smallskip
\end{itemize} 

With a view to applications over finite fields, we want hyperelliptic
equations over $\bfM_C$ in all cases. On one hand, the task is made easier as
there is never any obstruction for the curve to have a hyperelliptic equation
over its field of moduli and we propose an algorithm which makes the
reconstruction over $\bfM_C$ effective. Moreover, starting from the rational
parametrizations that we exhibit for most of the strata, we can state the
exact number of isomorphism classes of rational curves with a given
automorphism group, except for the dimension $3$ stratum $\Aut(C) \simeq
(\ZZ/2\ZZ)^2$. On the other hand, as already illustrated by the genus $2$
case, strange phenomenon can happen when the characteristic $p$ is too
small. This is not so surprising as the stratification itself may be different
when $p \leq 2g+1$. We took special care of denominators when computing
invariants and covariants and our results are then naturally valid for $p >
7=2 \cdot 3+1$. We wonder whether the natural bound $p>2g+1$ may be reached in
this way for general $g$.  \medskip

Finally, this article is only the emerged part of the iceberg. A
\textsc{magma} code\footnote{published under the GNU Lesser General Public
  License} containing the various algorithms to compute invariants and
reconstruct the curves is available on the web page of the authors. We tested
it over $\FF_p$ for $11 \leq p \leq 47$ and checked that we actually obtain
the $p^5$ non-isomorphic curves predicted by the theory. We plan to include
soon the computation of twists to obtain a representative of each
$\FF_q$-isomorphism class as we did for genus $2$.\bigskip

\noindent
{\bf Notation and convention.}
In the following, the integer $p$ is a prime or $0$, $k$ denotes any field of
characteristic $p$ and $K$ an algebraically closed field of characteristic
$p$. For a integer $g$, a \emph{hyperelliptic curve} or a \emph{genus $g$
  curve}  $C$ over a field $k$ is always assumed to be smooth, projective and
absolutely irreducible (if it is given by a singular equation, typically the
case for hyperelliptic curves, $C$ is the smooth model associated to this
equation). However a `curve' (like a conic) may be singular. When we speak of
a morphism from $C$ to $C'$ we always mean  a morphism \emph{defined over
  $k$}. However the notation $\Aut(C)$  stands for $\Aut_{\bar{k}}(C)$.\\

\noindent {\bf Acknowledgments.} This paper has benefit from plenty of
discussions with several mathematicians. During the Huggins /
Fuertes-Gonz\'alez-Diez controversy, the authors wish to thank the help and
patient listening of Nils Bruin, Pierre Dèbes, Yolanda Fuertes, Gabino
Gonz\'alez-Diez, Ruben Hidalgo, Everett Howe, Bonnie Huggins, David Kohel,
Stéphane Louboutin, Enric Nart, Xavier Xarles and finally Alexey Zykin who had
the sharpest eye.  Proposition~\ref{prop:godd} has been suggested to us by
Nils Bruin and we thank Jean-Marc Couveignes for indicating to us the
reference \cite{couveignes-belyi} which lead to
Proposition~\ref{prop:signature} and Jean-François Mestre for the precise
reference to \cite{clebsch}.

\section{Invariants and symbolic computations} \label{sec:rappel}
\subsection{Terminology}

We let $k$ be an infinite field of characteristic $p$ so we can identify the
graded algebra of polynomial functions over a finite $k$-vector space $E$ with
the graded symmetric algebra $S(E^*)$.
If $E=\oplus_{i=1}^m E_i$, a function $\phi$ is homogeneous of degree
$(n_1,\ldots,n_m)$ if it is homogeneous of degree $n_i$ on  $E_i$ for $1 \leq
i \leq m$, \ie if $\phi(f_1,\ldots,\lambda f_i,\ldots,f_m)=\lambda^{n_i} \cdot
\phi(f_1,\ldots,f_m)$ for all $\lambda \in K$ and all $(f_1,\ldots,f_m) \in
E$.\\
 
Let $V=k^2$ be seen as the space of vectors with coordinates $(x,z)$ and let
$E=S^n(V^*)$ be the $(n+1)$-dimension vector space of homogeneous polynomials
functions of degree $n$ in $X=(1,0)^*$ and $Z=(0,1)^*$. We write a generic
element of $E$ as $\sum_{i=0}^{n} a_i X^i Z^{n-i}$ with $a_i \in k$.  We let a
group $G \subset \GL_2(k)$ acts on $V$ in the natural left action, \ie
if $$M=\begin{bmatrix} a & b \\ c & d \end{bmatrix} \in G$$ then
$$M.{ (x,z)}={ (a x + b z,c x + d z)}.$$ Let $\rho$ be the contragredient representation on $V^*$. 
We let the group $G$ act with the tensorial representation $\otimes^n \rho$ on
$E$. More concretely, if $f \in E$ then $M.f$ is defined by $(M.f)
(x,z)=f(M^{-1}.(x,z))$. This can be extended to $E=\oplus_{i=1}^m
S^{n_i}(V^*)$ with $m \geq 1$ and integers $n_i>0$ using the representation
$\bigoplus_i \otimes^{n_i} \rho$.  In this setting, we recover the classical
definition of invariants and covariants.
\begin{definition}
  A homogeneous polynomial function $I : E \to k$ is a \emph{(relative)
    invariant} if there exists $\omega \in \ZZ$ such that for all $M \in G$
  and all $(f_1,\ldots,f_m) \in E$, we have $$I
  (M.f_1,\ldots,M.f_m)=\det(M)^{-\omega} \cdot I(f_1,\ldots,f_m)\,.$$
\end{definition}
\noindent
The space of invariants is classically embedded into a broader framework.
\begin{definition} \label{def:covariants} A homogeneous polynomial function $C
  : E \oplus V \to k$ of degree $(d_1,\ldots,d_m,r)$ is a \emph{(simultaneous)
    covariant} if there exists $\omega \in \ZZ$ such that for all $M \in G$,
  all $(f_1,\ldots,f_m) \in E$ and all $(x,z) \in V$, we have
$$C(M. f_1,\ldots,M. f_m,M.(x,z))= \det(M)^{-\omega} \cdot C(f_1,\ldots,f_m,(x,z))\,.$$
\end{definition}
\noindent
Hence, covariants of order $0$ are invariants.  The integer $\omega$ is called
the \emph{weight} (or index), the $m$-uple $(d_1,\ldots,d_m)$, the
\emph{degree} and the integer $r$ is called the \emph{order} of a covariant.
\begin{example} \label{ex:fundamental} When $E=S^n(V^*)$, the polynomial
  function $(\sum a_i X^i Z^{n-i},(x,z)) \mapsto \sum_{i=0}^n a_i x^i z^{n-i}$
  is a covariant of degree $1$ and of order $n$. We will denote it by
  $\mathfrak{f}$ in the sequel.
\end{example}

The major operation to generate new covariants from given ones is the
transvectant.  For $i,j \in \NN$ two distinct integers, let $(x_i,z_i)$ and
$(x_j,z_j)$ be bases of two copies $V_i$ and $V_j$ of $V$. Let $\Omega_{ij}$
be the differential operator
$$\Omega_{ij}=\frac{\partial}{\partial x_i} \frac{\partial}{\partial z_j} - \frac{\partial}{\partial z_i} \frac{\partial}{\partial x_j}.$$
Let $r_i, r_j>0$ be integers, $f_i \in S^{r_i}(V)=S^{r_i}(V_i)$ and $f_j \in
S^{r_j}(V)=S^{r_j}(V_j)$.  We define the following differential operators
(where composition is denoted multiplicatively).
\begin{definition}
  We define the \emph{$h$-th transvectant} $$(\phantom{f_i},\phantom{f_j})_h :
  S^{r_i}(V) \times S^{r_j}(V) \to S^{r_i+r_j-2 h}(V)$$ as
 $$(f_i,f_j)_h= \frac{(r_i-h)! \, (r_j-h)!}{r_i! \, r_j!} \cdot \left(\Omega_{ij}^h(f_i(x_i,z_i) \cdot f_j(x_j,z_j))\right)_{(x_i,z_i)=(x_j,z_j)=(x,z)}.$$
\end{definition}

In the sequel for all $i \in \NN$ we denote $\ell_i=\ell_i(x,z)=(\alpha_i x+
\beta_i z)$ some linear forms with coefficients $\alpha_i,\beta_i \in k$ and
$[\ell_i,\ell_j]=\alpha_i \beta_j-\alpha_j \beta_i$.  Using the Clebsch-Gordan
formula as in \cite[p.565]{procesi} or by direct computation, one gets the
following lemma which will turn to be central for proving the formula of
Section \ref{sec:reconstruction}.
\begin{lemma}
For any $r_i,r_j$ positive integers,
$$(\ell_i^{r_i},\ell_j^{r_j})_h=[\ell_i,\ell_j]^h \cdot \ell_i^{r_i-h} \ell_j^{r_j-h}.$$
\end{lemma}

\noindent If $C$ is a homogeneous covariant of order $r$, we can see $C$ as an
homogeneous polynomial of degree $r$ in the variables $x, z$, \ie as an
element of $S^{r}(V)$. In particular, $\mathfrak{f}=\sum a_i x^i z^{n-i}$.
Since the $r$-th powers of linear forms span the vector-space $S^r(V)$ and the
$h$-th transvectant is bilinear, we see that the $h$-th transvectant of two
covariants of degree $d_1,d_2$ and of order $r_1,r_2$ is a covariant of degree
$d_1+d_2$ and of order $r_1+r_2-2h$ (see also \cite[Chap.15]{procesi} for a
conceptual explanation).

\subsection{Algebra of invariants and covariants}
\label{sec:algebra-invar-covar}

We focus on invariants and covariants for $G=\GL_2(K)$ and $G=\SL_2(K)$.  In
both groups, using the multiplication by $-1$ scalar matrix, we first notice
that $n \sum d_i -r$ needs to be even in order for a covariant to be defined.
Moreover, if $C$ is a covariant for $\GL_2(k)$, then $C$ is a covariant for
$\SL_2(K)$. Conversely, let $C$ be a covariant for $\SL_2(K)$ and let $g \in
\GL_2(K)$. Since $g'=g/\sqrt{\det(g)} \in \SL_2(K)$, we have that
$$C(g'.f_1,\ldots,g'.f_m,g'.(x,z))=C(f_1,\ldots,f_m,(x,z)).$$
Since $C$ is homogeneous of degree $(d_1,\ldots,d_m)$ and of order $r$, for
any $\lambda \in k$, we have that $t=\lambda \cdot \Id$ acts as
$$ C(t. f_1, \ldots,t.f_m,t.(x,z))=\lambda^{n \sum d_i-r} \cdot C(f_1,\ldots,f_m,(x,z)),$$
and thus
$$C(g.f_1,\ldots,g.f_m,g.(x,z))=\det(g)^{-(n \sum d_i-r)/2} \cdot C(f_1,\ldots,f_m,x,z).$$
This shows that covariants for $\SL_2(k)$ are equal to the covariants for
$\GL_2(k)$, with weight $\omega$ equal to $(n \sum d_i-r)/2$. \medskip

Classically many results are known on the algebra $\cc_n$ generated by the
covariants of all degree and order under the action of $G=\SL_2(\CC)$ on
$E=S^n(V^*)$ and on its sub-algebra $\ii_n$ generated by invariants. Since
Gordan \cite{gordan68}, it is known that $\ii_n$ and $\cc_n$ are finitely
generated. Thanks to the so-called Clebsch-Gordan formula, one can even prove
that $\cc_n$ is generated by a finite number of iterations of transvectants
starting from the single covariant $\mathfrak{f}$ (see
\cite{grace-young}). Thus we can even assume that these generators are defined
over $\QQ$ and this allows to perform exact arithmetic operations with them.
Whereas today this particular case is part of a more general result on
reductive groups, Gordan's method is effective and generators for these
algebra were given for $n \leq 6$ in the nineteenth century (see
\cite{dixmier2}). Generators for $\ii_7$ were determined by von Gall
\cite{vongall} and Dixmier and Lazard \cite{dixmier-lazard} and for $\cc_7$ by
Bedratyuk \cite{bedradyuk}. The case of octics, $n=8$, which is the case we
are interested in this article will be reviewed and developed in
Section~\ref{sec:fund-invar-covar}. Finally, for $n=9$ and $n=10$, only
generators for $\ii_n$ have been determined (see
\cite{cronih,brouwer1},\cite{brouwer2}).

We summarize in Tab.~\ref{tab:In} for each degree $2\leq n \leq 10$, the
number of generators for $\ii_n$ and their degrees. We see that the number and
the degree of the generators increase strongly with $n$, but remains
affordable for some small degree $n$, especially for $n=8$.  \medskip

\begin{table}
  \centering
  \begin{tabular}{r|c|l}
    $n$ & \# Gen. of $\ii_n$ & Generator degrees {\footnotesize\it (exponents
      mean repetition)}\\\hline
    2 & 1 & 2 \\
    3 & 1 & 4 \\
    4 & 2 & 2, 3\\
    5 & 4 & 4, 8, 12, 18\\
    6 & 5 & 2, 4, 6, 10, 15\\
    7 & 30 & 4, 8$^3$, 12$^6$, 14$^4$, 16$^2$, 18$^9$, 20, 22$^2$, 26, 30 \\
    8 & 9 & 2, 3, 4, 5, 6, 7, 8, 9, 10 \\
    9 & 92 & 4$^2$, 8$^5$, 10$^5$, 12$^{14}$, 14$^{17}$, 16$^{21}$, 18$^{25}$,
    20$^2$, 22\\
    10 & 106 & 2, 4, 6$^4$, 8$^5$, 9$^5$, 10$^8$, 11$^8$, 12$^{12}$, 13$^{15}$,
    14$^{13}$, 15$^{19}$, 16$^5$, 17$^5$, 18, 19$^2$, 21$^2$ \\\hline
  \end{tabular}\medskip

  \caption{Number and degrees of fundamental invariants for binary forms of degree $n$}
  \label{tab:In}
\end{table}

The algebra $\ii_n$ is a sufficient tool to grasp the orbits of a generic
binary form.
\begin{theorem}[{\cite[p.~78]{mumford-fogarty},\cite[p.~47]{dixmier2}}]
  \label{th:equivforms}
  Let $f$, $f'$ be binary forms of even degree $n$ greater than or equal to
  $3$ over an algebraically closed field $K$ with no roots of multiplicity
  greater than or equal to $n/2$. Let $\{I_i\}$ be a finite set of homogeneous
  generators of degree $d_i$ for $\ii_n$.  Then $f$ and $f'$ are in the same
  orbit under the action of $\GL_2(K)$ (resp. $\SL_2(K)$) if and only if there
  exists $\lambda \in K$ such that for all $i$, $I_i(f)=\lambda^{d_i} \cdot
  I_i(f')$ (resp. $I_i(f)=I_i(f')$).
\end{theorem}
\begin{remark}
  In the theorem, the exponents are equal to the degrees $d_i$ when one could
  have expected weights $\omega_i$. Actually the two possibilities are
  equivalent since the weight is $n/2$ times the degree.
\end{remark}

Note that the multiplicity condition in Theorem~\ref{th:equivforms} is as good
as possible already for $n=2$ since the only generator of $\ii_2$,
$a_1^2-a_0a_2$, is $0$ for quadratic forms in the two distinct orbits $f=0$
and $f=X^2$. Since we will only deal with polynomials defining hyperelliptic
curves, and hence having simple roots, this is not a limitation.

\subsection{Hyperelliptic curves and invariants of binary
  forms} \label{subsec:reconstr} Let $k$ be a field of characteristic $p$ and
$K=\bar{k}$. A curve $C$ of genus $g \geq 2$ defined over $k$ is said
\emph{hyperelliptic} if $C/K$ has a separable degree $2$ map to
$\PP^1_{K}$. Since the extension $K(C)/K(\PP^1)$ is Galois, the curve $C/K$
has an involution that we denote $\iota$.
\begin{lemma}[{\cite[Prop.IV.5.3]{Hart}}]
  The automorphism $\iota$ is the unique involution on $C/K$ such that
  $C/\langle \iota \rangle$ is of genus $0$. It commutes with all
  automorphisms of $C/K$ and is called the \emph{hyperelliptic involution}.
\end{lemma}
By unicity, $\iota$ is defined over $k$ and induces a morphism $\rho : C \to
Q=C/\langle \iota \rangle$ where $Q/k$ is a genus $0$ curve, not necessarily
isomorphic to $\PP^1$. However, if the curve $Q$ has a rational point then it
is isomorphic to $\PP^1$ and $k(C)$ is a degree $2$ separable cover of
$k(\PP^1) \simeq k(x)$. In that case, $C$ is birationally equivalent to an
affine curve of the form
 $$y^2+h(x) y= f(x)$$
 where $\deg f \leq 2g+2$ and $\deg h \leq g$. We say that $C$ has  a \emph{hyperelliptic equation} if (a curve in the isomorphism class of $C$ (over $k$)) can be written in the form above. A hyperelliptic curve has automatically a hyperelliptic equation when $k$ is algebraically closed or a finite field and we recover the `classical' definition of a hyperelliptic curve. However, as we shall see in Section \ref{subsec:def}  for more general fields and  odd genus, it is not necessarily the case.\\

 Let us assume from now on that $p \ne 2$.  One can assume that $h=0$ and the
 polynomial $f$ has then simple roots. We say that $f$ is a
 \emph{hyperelliptic polynomial}.  By homogenizing, $Y^2 Z^{2g} = f(X,Z)$, we
 have $\deg f(X,Z)=2g+2$, taking into account a `root' at infinity when $\deg
 f=2g+1$. With this convention, the roots of $f$ are the ramification points
 of the cover $C/Q$ and equivalently the Weierstrass points $W$ of $C$. We
 will use these conventions for the roots and degree in the sequel when we
 speak about a hyperelliptic polynomial.
 \begin{remark}
   One cannot always assume that $\deg f=2g+2$ (as a polynomial) if $k$ is not
   infinite. For instance the genus $2$ curve $y^2=x(x-1)(x-2)(x-3)(x-4)$ over
   $\FF_5$ has no model with $\deg f=6$ since all its affine Weierstrass
   points are rational.
 \end{remark}

 The following proposition will give us the link between hyperelliptic curves
 and invariants for the action of $\GL_2$.
 \begin{proposition} \label{prop:iso} Let $C : y^2=f(x)$ and $C' : y^2=f'(x)$
   two hyperelliptic curves of genus $g$ over a field $k$. Every isomorphism
   $\phi : C \to C'$ is given by an expression of the form
$$(x,y) \mapsto \left(\frac{ax+b}{cx+d}, \frac{ey}{(cx+d)^{g+1}}\right),$$
for some $M =\begin{bmatrix} a & b \\ c & d \end{bmatrix} \in \GL_2(k)$ and $e
\in k^*$. The pair $(M,e)$ is unique up to replacement by $(\lambda
M,\lambda^{g+1} e)$ for $\lambda \in k^*$. The composition of isomorphisms
$(M,e)$ and $(M',e')$ is $(M'M,e'e)$.
\end{proposition}
\begin{remark}
  We refer to \cite{CANAPU} for the characteristic $2$ case.
\end{remark}
Let assume now that $k=K$ is algebraically closed. The ramification divisor
$W$ of $\rho : C \to Q$ is a positive divisor of degree $2g+2$.  From
Proposition \ref{prop:iso}, one easily gets
\begin{corollary} \label{cor:W} Let $C,C'$ be two hyperelliptic curves of
  genus $g$ with respective ramification divisor $W$ and $W'$. Let $D=\rho(W)
  \in \Div(Q)$ and $D'=\rho'(W') \in \Div(Q')$. Then $C$ and $C'$ are
  isomorphic if and only if there exists an isomorphism $\varphi : Q \to Q'$
  such that $\varphi(W)=W'$.
\end{corollary}
Checking an isomorphism between $C$ and $C'$ in this way is often difficult
since the points of the ramification divisor can be defined over large degree
extensions. This is where invariants enter the game. Indeed, as in the
proposition let $C : y^2=f(x)$ and $C' : y^2=f'(x)$. Since $f$ and $f'$ have
only simple roots, we can use invariants to decide if the homogeneous binary
forms $f(X,Z)$ and $f'(X,Z)$ of degree $2g+2$ are in the same orbits under the
action of $\GL_2(K)$. If so, one has $f'(X,Z)=f(aX+bZ,cX+dZ)$ for some
$M=\begin{bmatrix} a & b \\ c & d \end{bmatrix} \in \GL_2(K)$ and so $(M,1)$
defines an isomorphism between $C$ and $C'$. Conversely, let $(M,e)$ define an
isomorphism between $C$ and $C'$. we can adjust $\lambda$ in Proposition
\ref{prop:iso} such that $e=1$ and hence $f'$ and $f$ are in the same orbit
under $\GL_2(K)$.  Applying Theorem.~\ref{th:equivforms} to the situation
yields thus the following proposition.
\begin{proposition}\label{prop:isocurves}
  Let $\{I_i\}$ be a set of homogeneous generators of degree $d_i$ for
  $\ii_{2g+2}$. Two hyperelliptic curves $C: y^2=f(x)$ and $C': y^2=f'(x)$ of
  genus $g$ are isomorphic if and only if there exists $\lambda \in K^*$ such
  that $I_i(f')= \lambda^{d_i} \cdot I_i(f)$ for all $i$.
\end{proposition}
Hence, the possible values of a set of generators for $\ii_{2g+2}$ up to this
specific equivalence are in bijection with the points of the coarse moduli
space $\Hm_g$ of hyperelliptic curves of genus $g$. We therefore need
algorithms to handle such `weighted sets'.

\subsection{Algorithms in weighted projective spaces}
\label{sec:algorithms}

In the context of curves of genus $1$ or $2$, one usually prefers to handle
absolute invariants instead of homogeneous ones, by calculating ratios of
homogeneous invariants of the same degree.  In this way, one gets rid of the
constant $\lambda$ in Proposition~\ref{prop:isocurves}.  It becomes then easy
to span the coarse moduli space or to check that two curves are isomorphic.

But, care has to be taken to ensure that the denominators of absolute
invariants do not vanish for some hyperelliptic orbits and a common approach
is to choose as denominator some powers of the discriminant invariant.  We
give up this strategy for higher genus curves, because the degree of the
discriminant is too large (this degree is already equal to 14 for octics) and
selecting invariants of smaller degree as denominators yields too many
technical cases to consider.  We choose instead to work with a (kind of)
weighted projective space, given by $m$-uples $(I_1: I_2: \ldots: I_{m})$ of
weights $d_1$, $d_2$, \ldots, $d_m$. More precisely
\begin{definition}
  Let $k$ be a field and $d_1,\ldots,d_m$ positive integers wit $m \geq 2$. We
  call $\mathbb W$ a \emph{weighted projective space} over $k$ of dimension
  $m-1$ and weights $(d_1,\ldots,d_m)$ the set of elements denoted $(\iota_1 :
  \ldots : \iota_m)$ which are equivalence classes of $m$-uples
  $(\iota_1,\ldots,\iota_m) \in k^m \setminus (0,\ldots,0)$ for the relation
$$(\iota_1, \ldots, \iota_m) \sim (\iota_1',\ldots, \iota_m') \iff \; \exists \, \lambda \in \bar{k}^* \; \textrm{such that} \; \iota_j = \lambda^{d_j} \cdot \iota'_j \; \forall \, 1 \leq j \leq m.$$
\end{definition}

As a first tool, we need an algorithm for testing the equality of two points
in a weighted projective space.
\begin{proposition}\label{prop:WPSequal}
  Let $k$ be a field and $\mathbb W$ a $k$-weighted projective space of
  dimension $m-1$ and weights $(d_1,d_2,\ldots, d_m)$, then
  Algorithm~\ref{WPSEqual} tests if two elements of $k^m$ are in the same
  class of $\mathbb W$.

  If $k$ is a field which admits operations with quasi-linear complexity in
  time and space (multiplications, inverses, tests), then
  Algorithm~\ref{WPSEqual} has quasi linear complexity in time and space.
\end{proposition}
\begin{center}
  \parbox{0.8\linewidth}{%
    \begin{footnotesize}\SetAlFnt{\small\sf}%
      \begin{algorithm}[H]%
        \label{WPSEqual}%
        \caption{Equality in a weighted projective space.}  %
        \SetKwInOut{Input}{Input} \SetKwInOut{Output}{Output} %
        \Input{Two elements $(U_1, U_2, \ldots, U_{m})$ and $(V_1, V_2,
          \ldots, V_{m})$ in $k^m$ and $\mathbb W$ a $k$-weighted projective
          space of dimension $m$ with weights $(d_1, d_2, \ldots, d_m)$.}  %
        \Output{The boolean ``\texttt{true}'' if $(U_1: U_2: \ldots:
          U_{m})=(V_1: V_2: \ldots: V_{m})$, ``\texttt{false}'' otherwise.}
        \BlankLine %
        ${\mathcal S}_U \leftarrow \{i \in \{1,\ldots, m\} \ \mathtt{s.t.}\
        U_i \neq 0\}$;\ ${\mathcal S}_V \leftarrow \{i \in \{1,\ldots, m\} \
        \mathtt{s.t.}\ V_i \neq 0\}$\; \If{${\mathcal S}_U \neq {\mathcal
            S}_V$}{\KwRet{$\mathtt{false}$}} %
        $d,\ (c_i : i \in {\mathcal S}_U) \leftarrow \mathtt{ExtentedGCD}(d_i
        : i \in {\mathcal S}_U)$\; %
        $\Lambda \leftarrow \prod_{i\in {\mathcal S}_U} (V_i/U_i)^{c_i} $\; %
        \KwRet{$ \mathtt{true}\text{ if } V_i/U_i = \Lambda^{d_i/d}\text{ for
            all }i \in {\mathcal S}_U,\ \mathtt{false}\text{ otherwise}.$}
   \end{algorithm}
 \end{footnotesize}
}
\end{center}
\begin{proof}
  If $(U_1, U_2, \ldots, U_{m})$ and $(V_1, V_2, \ldots, V_{m})$ are in the
  same class of $\mathbb W$, then there exists some $\lambda\in\bar{k}^*$ such
  that $V_i = \lambda ^{d_i}\cdot U_i$ for all $i$. So, $\Lambda$ is equal to
  $\prod (V_i/U_i)^{c_i} = \lambda^{\sum c_i d_i } = \lambda ^d$, and thus we
  have $V_i/U_i = \Lambda^{d_i/d} = \lambda ^{d_i} $ for $i\in {\mathcal
    S}_U$. Conversely, if ${\mathcal S}_U={\mathcal S}_V$ and $V_i/U_i =
  \Lambda^{d_i/d}$ for $i\in {\mathcal S}_U$, then let $\lambda$ be a $d$-th
  root of $\Lambda$ in $\bar{k}^*$ and we can easily check that $V_i = \lambda
  ^{d_i} \cdot U_i$.
\end{proof}

We deduce from Algorithm~\ref{WPSEqual} how to associate to a class $(U_1:
U_2: \ldots: U_{m})$ a unique vector of $k^m$ representing the class. Let
${\mathcal S}_U$ be the support of $U$, \textit{i.e.} ${\mathcal S}_U = \{i
\in \{1,\ldots, m\} \ \mathrm{s.t.}\ U_i \neq 0\}$, let $c_i$ be integers such
that $\sum_{i \in {\mathcal S}_U} c_i d_i = d$ with $d = \gcd(d_i : i \in
{\mathcal S}_U)$ and let $\Lambda = \prod_{i\in {\mathcal S}_U} U_i^{c_i}
$. Set $u_i = U_i/\Lambda ^{d_i/d}$ for $i \in {\mathcal S}_U$ and $u_i = 0$
otherwise. Then $(u_1, u_2, \ldots, u_m)$ is the unique representative of
$(U_1: U_2: \ldots: U_{m})$ such that $\prod_{i\in {\mathcal S}_U} u_i^{c_i} =
1$.\medskip

When $k$ is a finite field, enumerating elements in $\mathbb W$ is another
feature that is needed to span $\Hm_g(k)$. In weighted projective spaces, this
can be easily done by enumerating representatives. Typically, for each support
${\mathcal S}_U$, considered in turn, fix, for once, integers $c_i$ such that
$\sum_{i \in {\mathcal S}_U} c_i d_i = \gcd(d_i : i \in {\mathcal S}_U)$. Then
enumerate all the vectors $(u_1, u_2, \ldots, u_m)$ of $k^m$ with support
${\mathcal S}_U$ such that $\prod_{i\in {\mathcal S}_U} u_i^{c_i} = 1$.

\begin{example}
  Let $k=\FF_7$ and $\mathbb W$ be a $k$-weighted projective space of
  dimension $2$ and weights $(5,\ 7)$. To enumerate elements in $\mathbb W$,
  we consider supports in turn.  The supports ${\mathcal S}=\{1\}$ and
  ${\mathcal S}=\{2\}$, yields the classes $(1:0)$ and $(0:1)$. Let us now
  consider the support ${\mathcal S}=\{1, 2\}$. We fix $c_1=3$, $c_2=-2$, so
  that $5 c_1 +7 c_2 = \gcd(5, 7)$. This yields $6$ representatives
  $(U_1:U_2)$ such that $\prod_{i=1,2} U_i^{c_i} = 1$, that is $( 1: 1 )$, $(
  1: 6 )$, $( 2: 1 )$, $( 2: 6 )$, $( 4: 1 )$, $( 4: 6 )$. Another choice of
  $c_1$ and $c_2$ would lead to different representatives: for instance
  $c_1=24$ and $c_2=-17$ gives $(1:i)$ for $1 \leq i \leq 6$.
\end{example}

Finally, enumerating points on a coarse moduli space can be translated into
enumerating points on a variety defined inside a weighted projective
space. This is much more intricate than enumerating the full space, at least
for curves of genus larger than 2. A naive strategy consists in enumerating
all the points in the ambient space and for each point check if it is defined
on the projective moduli variety. It is often inefficient. More sophisticated
methods may be possible, based on a nice description of the variety,
especially when the variety is rational and one has an explicit
parametrization for it.  In this direction, we give in
Section~\ref{sec:fund-invar-covar} an efficient method for the moduli space of
genus 3 hyperelliptic curves.

\subsection{Fundamental invariants and covariants for the binary octics}
\label{sec:fund-invar-covar}

According to the syzygy theorem of Hilbert, the graded ring $\ii_8$ of
invariants of binary octics fits into a finite exact sequence of
$\CC[\Xm]$-module where $\CC[\Xm]:=\CC[X_2,\ldots,X_{10}]$. In what Mumford
called a `tour de force', Shioda managed to find an explicit description.
\begin{theorem}[\cite{shioda67}, Th.~3, p.~1042]\label{th:syzygies}
  The graded ring $\ii_8$ of invariants of binary octics is generated by 9
  elements $J_2$, $J_3$, $\ldots$, $J_{10}$ of degree $2$, $3$, $\ldots$,
  $10$. There exists 5 generating relations, $\mathfrak{R}_i(J)$, of degree
  $15+i$, ($i=1, \ldots, 5$), which, in turn, are connected by 5 fundamental
  first syzygies $\mathfrak{T}_i(R)$ of degree $24+i$, ($i=1, \ldots, 5$). The
  second syzygy $\mathfrak{F}$ is unique up to constants and of degree
  $45$. The syzygy sequence of $\ii_8$ (as $\CC[\Xm]$-module) is given by
  \begin{displaymath}
    0 \rightarrow \CC[\Xm] \mathfrak{F} \rightarrow \sum_{i=1}^{5}\CC[\Xm]
    \mathfrak{T}_i \rightarrow \sum_{i=1}^{5}\CC[\Xm] \mathfrak{R}_i \rightarrow
    \CC[\Xm] \rightarrow \ii_8 \rightarrow 0\,.
  \end{displaymath}
\end{theorem}

Using transvectants and the covariant $\mathfrak{f}$ from Example
\ref{ex:fundamental}, Shioda defined the invariants $J_i$'s of
Theorem~\ref{th:syzygies} in the following way,
\begin{footnotesize}
  \begin{displaymath}
    \begin{array}{c}
      J_2 = (\mathfrak{f},\mathfrak{f})_4,\ J_3 = (f,\mathfrak{g})_8,\ J_4 =
      (\mathfrak{k},\mathfrak{k})_4,\ J_5 = (\mathfrak{m},\mathfrak{k})_4,\ J_6
      = (\mathfrak{k},\mathfrak{h})_4,\\ J_7 = (\mathfrak{m},\mathfrak{h})_4,
      J_8 = (\mathfrak{p},\mathfrak{h})_4,\ J_9 =
      (\mathfrak{n},\mathfrak{h})_4,\ J_{10} = (\mathfrak{q},\mathfrak{h})_4
  \end{array}
\end{displaymath}
\end{footnotesize}
where
\begin{footnotesize}
  \begin{displaymath}
    \mathfrak{g} = (\mathfrak{f},\mathfrak{f})_4,\ \mathfrak{k} = (\mathfrak{f},\mathfrak{f})_4,\ \mathfrak{h} =
    (\mathfrak{k}, \mathfrak{k})_2,\ \mathfrak{m} = (\mathfrak{f}, \mathfrak{k})_4,\
    \mathfrak{n} = (\mathfrak{f}, \mathfrak{h})_4,\ \mathfrak{p} = (\mathfrak{g},
    \mathfrak{k})_4,\
    \mathfrak{q} = (\mathfrak{g}, \mathfrak{h})_4\,.
  \end{displaymath}
\end{footnotesize}
One may remark that none of these transvectant calculations involves
covariants of order greater than 10. These formulae are thus valid, in
addition to fields of characteristic 0, over any field of or characteristic
greater than or equal to 11. We found for a generic form $f = a_8 x^8 + a_7x^7
+ \ldots + a_0$,
\begin{center}
  \parbox{0.85\linewidth}{%
    \begin{dgroup*}[style={\footnotesize}] %
      \begin{dmath*} J_2=1/140\,(280\,{ a_0}\,{ a_8}-35\,{ a_1}\,{ a_7}+10\,{
          a_2}\,{ a_6} -5\,{ a_3}\,{ a_5}+2\,{{ a_4}}^{2})\,, \end{dmath*} %
      \begin{dmath*} J_3 =1/137200\,(11760\,{ a_0}\,{ a_4}\,{ a_8}-7350\,{
          a_0}\,{ a_5}\,{ a_7 }+3150\,{ a_0}\,{{ a_6}}^{2}-7350\,{ a_1}\,{
          a_3}\,{ a_8}+ 2205\,{ a_1}\,{ a_4}\,{ a_7}-525\,{ a_1}\,{ a_5}\,{
          a_6}+ 3150\,{{ a_2}}^{2}{ a_8}-525\,{ a_2}\,{ a_3}\,{ a_7}-330\,{
          a_2}\,{ a_4}\,{ a_6}+225\,{ a_2}\,{{ a_5}}^{2}+225\,{{ a_3}}^{2}{
          a_6}-135\,{ a_3}\,{ a_4}\,{ a_5}+36\,{{ a_4}}^{3})\,, \end{dmath*} %
      \begin{dmath*} J_4=1/3687936\,(8064\,{ a_0}\,{ a_4}\,{{
            a_6}}^{2}-614656\,{ a_0}\,{ a_1}\, { a_7}\,{ a_8}-12544\,{ a_0}\,{
          a_2}\,{ a_6}\,{ a_8}+53312 \,{ a_0}\,{ a_3}\,{ a_5}\,{ a_8}-35280\,{
          a_0}\,{ a_3}\,{ a_6}\,{ a_7}+4704\,{ a_0}\,{ a_4}\,{ a_5}\,{
          a_7}-35280\, { a_1}\,{ a_2}\,{ a_5}\,{ a_8}-13132\,{ a_1}\,{ a_2}\,{
          a_6}\,{ a_7}+4704\,{ a_1}\,{ a_3}\,{ a_4}\,{ a_8}+3626\,{ a_1}\,{
          a_3}\,{ a_5}\,{ a_7}-3864\,{ a_1}\,{ a_4}\,{ a_5}\,{ a_6}-3864\,{
          a_2}\,{ a_3}\,{ a_4}\,{ a_7}-1516\,{ a_2}\,{ a_3}\,{ a_5}\,{
          a_6}-320\,{ a_3}\,{{ a_4}}^{2}{ a_5 }+82320\,{ a_0}\,{ a_2}\,{{
            a_7}}^{2}-3360\,{ a_0}\,{{ a_5}} ^{2}{ a_6}+2458624\,{{
            a_0}}^{2}{{ a_8}}^{2}+2401\,{{ a_1}}^{2 }{{ a_7}}^{2}+1260\,{
          a_1}\,{{ a_5}}^{3}+256\,{{ a_2}}^{2}{{ a_6}}^{2}+1260\,{{ a_3}}^{3}{
          a_7}+589\,{{ a_3}}^{2}{{ a_5} }^{2}-25088\,{ a_0}\,{{ a_4}}^{2}{
          a_8}+82320\,{{ a_1}}^{2}{ a_6}\,{ a_8}+3780\,{ a_1}\,{ a_3}\,{{
            a_6}}^{2}+784\,{ a_1}\,{{ a_4}}^{2}{ a_7}+8064\,{{ a_2}}^{2}{
          a_4}\,{ a_8}+ 3780\,{{ a_2}}^{2}{ a_5}\,{ a_7}-3360\,{ a_2}\,{{
            a_3}}^{2}{ a_8}+1984\,{ a_2}\,{{ a_4}}^{2}{ a_6}-504\,{ a_2}\,{
          a_4} \,{{ a_5}}^{2}-504\,{{ a_3}}^{2}{ a_4}\,{ a_6}+64\,{{ a_4}}^
        {4})\,,
    \end{dmath*}
    \begin{dmath*}
      etc.
    \end{dmath*}
  \end{dgroup*}
}
\end{center}

The invariants $J_2$, \ldots, $J_7$ are algebraically independent (see
\cite[Lemma~4,~p.~1037]{shioda67}). We have moreover 5 relations between
$J_{8}$, $J_{9}$ and $J_{10}$, of the form
\begin{footnotesize}
  \begin{equation}
    \begin{array}{rcl}
      \mathfrak{R}_1 &:& J_{8}^2                 +  A_{6}\,J_{10} +  A_{7}\,J_{9} +  A_{8}\,J_{8} + A_{16} = 0\\
      \mathfrak{R}_2 &:& J_{8}\,J_{9}               +  B_{7}\,J_{10} +  B_{8}\,J_{9} +  B_{9}\,J_{8} + B_{17} = 0\\
      \mathfrak{R}_3 &:& J_{8}\,J_{10} + C_0\,J_{9}^2    +  C_{8}\,J_{10} +  C_{9}\,J_{9} + C_{10}\,J_{8} + C_{18} = 0\\
      \mathfrak{R}_4 &:& J_{9}\,J_{10}             +  D_{9}\,J_{10} + D_{10}\,J_{9} + D_{11}\,J_{8} + D_{19} = 0\\
      \mathfrak{R}_5 &:& J_{10}^2  + E_0\,{J}_{2}\,J_{9}^2 + E_{10}\,J_{10} + E_{11}\,J_{9} + E_{12}\,J_{8} + E_{20} = 0
    \end{array}\label{eq:Syzygies}
\end{equation}
\end{footnotesize}
where the $A_i$'s, $B_i$'s, $C_i$'s and $D_i$'s are invariants of degree $i$,
that can be explicitly determine as functions of $J_2$, $J_3$, \ldots, $J_7$
(see \cite[Th. 1,~p.~1030]{shioda67}).  Moreover, following
Shioda~\cite[p.~1043]{shioda67}), we may consider $\mathfrak{R}_1$,
$\mathfrak{R}_2$, $\mathfrak{R}_3$ and $(J_9-B_9)\mathfrak{R}_2 -
B_7\mathfrak{R}_4$ as linear equations in 1, $J_9$, $J_9^2$ and $J_{10}$ and
we obtain that $J_8$ always satisfies an equation of degree $5$,

  {\rm\setcounter{myequation}{\theequation} 
    \begin{dgroup*}[style={\tiny},spread={-4pt}]
      \begin{dmath}[style={\refstepcounter{myequation}},number={\themyequation},label={eq:J8}]
        {{J_{8}}}^{5}+ ( {A_{8}}+2\,{B_{8}}+{C_{8}} ) {{J_{8}}}^{4}+ (
        -{A_{6}}\,{C_{10}}-{A_{7}}\,{B_{9}}+2\,{A_{8}}\,{B_{8}}+{A_{8}}\,{C_{8}}+{A_{16}}+{B_{7}}\,{B_{9}}\,C_0+
        {B_{7}}\,C_0{D_{9}}-{B_{7}}\,{C_{9}}+{{B_{8}}}^{2}+2\,{B_{8}}\,{C_{8}}
        ) {{J_{8}}}^{3}+ (
        -{A_{6}}\,{B_{7}}\,C_0{D_{11}}-2\,{A_{6}}\,{B_{8}}\,{C_{10}}-{A_{6}}\,{{B_{9}}}^{2}C_0+{A_{6}}\,{B_{9}}\,{C_{9}}-{A_{6}}\,{C_{18}}+{A_{7}}\,{B_{7}}\,{C_{10}}-{A_{7}}\,{B_{8}}\,{B_{9}}-{A_{7}}\,{B_{9}}\,{C_{8}}
        -{A_{7}}\,{B_{17}}+{A_{8}}\,{B_{7}}\,{B_{9}}\,C_0+{A_{8}}\,{B_{7}}\,C_0{D_{9}}-{A_{8}}\,{B_{7}}\,{C_{9}}+{A_{8}}\,{{B_{8}}}^{2}+2\,{A_{8}}\,{B_{8}}\,{C_{8}}+2\,{A_{16}}\,{B_{8}}+{A_{16}}\,{C_{8}}-{{B_{7}}}^{2}C_0{D_{10}}+{B_{7}}\,{B_{8}}\,C_0{D_{9}}-{B_{7}}\,{B_{8}}\,{C_{9}}+{B_{7}}\,{B_{17}}\,C_0+{{B_{8}}}^{2}{C_{8}}
        ) {{J_{8}}}^{2}+ (
        -{A_{6}}\,{B_{7}}\,C_0{D_{19}}+{A_{6}}\,{B_{8}}\,{B_{9}}\,{C_{9}}+{A_{7}}\,{{B_{7}}}^{2}C_0{D_{11}}+{A_{7}}\,{B_{7}}\,{B_{8}}\,{C_{10}}-{A_{7}}\,{B_{8}}\,{B_{9}}\,{C_{8}}
        -{A_{8}}\,{{B_{7}}}^{2}C_0{D_{10}}-{A_{8}}\,{B_{7}}\,{B_{8}}\,{C_{9}}+{A_{8}}\,{B_{7}}\,{B_{17}}\,C_0+{A_{16}}\,{B_{7}}\,{B_{9}}\,C_0+{A_{16}}\,{B_{7}}\,C_0{D_{9}}+{A_{16}}\,{{B_{8}}}^{2}
        -{A_{6}}\,{{B_{8}}}^{2}{C_{10}}+{A_{6}}\,{B_{17}}\,{C_{9}}+{A_{7}}\,{B_{7}}\,{C_{18}}-{A_{7}}\,{B_{8}}\,{B_{17}}-{A_{7}}\,{B_{17}}\,{C_{8}}+{A_{8}}\,{{B_{8}}}^{2}{C_{8}}-{A_{16}}\,{B_{7}}\,{C_{9}}-{A_{6}}\,{B_{7}}\,{B_{8}}\,C_0{D_{11}}+{A_{6}}\,{B_{7}}\,{B_{9}}\,C_0{D_{10}}-{A_{7}}\,{B_{7}}\,{B_{9}}\,C_0{D_{9}}+{A_{8}}\,{B_{7}}\,{B_{8}}\,C_0{D_{9}}-2\,{A_{6}}\,{B_{9}}\,{B_{17}}\,C_0-2\,{A_{6}}\,{B_{8}}\,{C_{18}}+2\,{A_{16}}\,{B_{8}}\,{C_{8}}
        ) {J_{8}}
        -{A_{6}}\,{B_{7}}\,{B_{8}}\,C_0{D_{19}}+{A_{6}}\,{B_{7}}\,{B_{17}}\,C_0{D_{10}}-{A_{6}}\,{{B_{8}}}^{2}{C_{18}}+{A_{6}}\,{B_{8}}\,{B_{17}}\,{C_{9}}-{A_{6}}\,{{B_{17}}}^{2}C_0+{A_{7}}\,{{B_{7}}}^{2}C_0{D_{19}}
        +{A_{7}}\,{B_{7}}\,{B_{8}}\,{C_{18}}-{A_{7}}\,{B_{7}}\,{B_{17}}\,C_0{D_{9}}-{A_{7}}\,{B_{8}}\,{B_{17}}\,{C_{8}}-{A_{16}}\,{{B_{7}}}^{2}C_0{D_{10}}+{A_{16}}\,{B_{7}}\,{B_{8}}\,C_0{D_{9}}-
        {A_{16}}\,{B_{7}}\,{B_{8}}\,{C_{9}}+{A_{16}}\,{B_{7}}\,{B_{17}}\,C_0+{A_{16}}\,{{B_{8}}}^{2}{C_{8}}\,.
      \end{dmath}
  \end{dgroup*}
  \setcounter{equation}{\themyequation} } \medskip

Generically, \textit{i.e.}  when $\delta = { A_6}\,{ J_8}+{ A_6}\,{ B_8}-{
  A_7}\,{ B_7} \neq 0$, $\mathfrak{R}_1$ and $\mathfrak{R}_2$ yields
\begin{tiny}
  \begin{eqnarray*}
    J_9 &=& ( {{ B_7}\,{{ J_8}}^{2}-{ A_6}\,{ B_9}\,{ J_8}-{ A_6}\,{ B_{17}}+{
        B_7}\,{ A_8}\,{ J_8}+{ A_{16}}\,{ B_7}})/\delta\,,\\
    J_{10} &=& -({{{ J_8}}^{3}+{{ J_8}}^{2}{ B_8}-{ A_7}\,{ B_9}\,{ J_8}-{ A_7}\,{ B_{17}}+{ A_8}\,{{ J_8}}^{2}+{ A_8}\,{ J_8}\,{ B_8}+{ A_{16}}\,{ J_8}+{ A_{16}}\,{ B_8}})/\delta\,.
  \end{eqnarray*}
\end{tiny}
Otherwise, \textit{i.e.} when $\delta = 0$, one might have no rational
solution for the system $(\mathfrak{R_i})_{i=1\ldots5}$ in $J_9$ and $J_{10}$
(for instance $J_2=1$, $J_3=0$, $J_4=0$, $J_5=0$, $J_6=8$, $J_7=0$, $J_8=7$
modulo 11) or several (for instance $J_2=9$, $J_3=0$, $J_4=0$, $J_5=0$,
$J_6=2$, $J_7=0$, $J_8=0$ which has 3 solutions modulo 11 for $(J_9, J_{10})$:
$(0, 4)$, $(2, 9)$ or $(9, 9)$).

\begin{remark}
  Instead of $\ii_8$, one may look at $\textrm{Frac}(\ii_8)$ and consider
  absolute invariants. It is known that $\textrm{Frac}(\ii_n)$ is a rational
  function field (see ~\cite{bogomolov86}) and in the case $n=8$, Maeda worked
  out $6$ algebraically independent absolute
  invariants~\cite[Th.~B,~p.~631]{maeda90} which generate
  $\textrm{Frac}(\ii_8)$. Unfortunately their degrees are too large for
  practical computations.
\end{remark}

Now, following Section~\ref{sec:algorithms}, we represent the coarse moduli
space of genus 3 hyperelliptic curves by the projective variety given by
Eq.~\eqref{eq:Syzygies} defined in a weighted projective space of dimension 9
with weights $2$, $3$,\ldots, $10$ the points of which are of the form
$(J_2:J_3:~\ldots:J_{10})$\,. In this setting, enumerating points on the
moduli space over a finite field can be done as follows.
\begin{enumerate}
\item Enumerate representatives for all the points of the weighted projective
  subspace of dimension $6$ defined by the algebraically independent
  invariants $J_2$, $J_3$, \ldots $J_7$ (\textit{cf.}
  Section~\ref{sec:algorithms}).
\item For each such representative, denoted $(j_2,j_3~,~\ldots~,j_7)$, compute
  the gcd $\delta$ of its support $\{d~: j_d \ne 0\}$.
\item For each representative $\pi$ of the quotient $k^*/(k^*)^\delta$ (we
  choose $\pi=1$ when $\delta=1$), compute the roots in $J_8$ of
  Eq.~\eqref{eq:J8} specialized at $J_d = \pi^{(d/\delta)}\,j_d$ for $2\leq d
  \leq 7$.
\item For each root $j_8$, solve Eq.~\eqref{eq:Syzygies} in $j_9$ and
  $j_{10}$.
\item Return representatives for
  $\{(\pi^{(2/\delta)}\,j_2:\pi^{(3/\delta)}\,j_3:~\ldots:\pi^{(7/\delta)}\,j_7:j_8:j_9:j_{10})\}$.
\end{enumerate}
\medskip

We consider in this algorithm several representatives
$(\pi^{(2/\delta)}\,j_2:\pi^{(3/\delta)}\,j_3:~\ldots:\pi^{(7/\delta)}\,j_7)$
starting from the same class $(j_2:j_3:~\ldots:j_7)$ when $\delta \ne 1$, even
if we may, so, encounter several times the same class at the very end of the
enumeration.  We do so because such collisions are straightforward to filter,
and otherwise we may miss points on the coarse moduli space.

For instance, modulo $11$, the representative $j_2=-1$,
$j_3=j_4=j_5=j_6=j_7=0$ yields only one point in the moduli space,
\textit{i.e} $(-1:0:0:0:0:0:0:0:0)$, while another choice of representative
for the class $(-1:0:\ldots:0)$, as for instance $j_2=1$,
$j_3=j_4=j_5=j_6=j_7=0$, yields two points in the moduli space, among which
the (new) point $(1: 0: 0: 0: 0: 0: 8: 2: 7)$. This is mostly due to the fact
that there exists equivalent representatives, here $(-1,0,0,0,0,0)$ and
$(1,0,0,0,0,0)$, for some classes $(j_2:j_3:\ldots:j_7)$, here
$(-1:0:0:0:0:0)$, linked by a projective constant $\lambda$ which may only be
defined in a non-trivial extension of $\FF_{11}$ of degree $\delta$ (in the
example $\delta=2$). In this situation, specializing Eq.~\eqref{eq:Syzygies}
with such representatives yield solutions which are no more necessarily
equivalent.\medskip

Some of the classes encountered while enumerating may be specializations of
Shioda invariants at forms $f$ which do not have simple roots. These classes
are not points of the coarse moduli space. To discriminate them, we may check
that the discriminant is non-zero for forms in the class
$(j_2:j_3:~\ldots:j_{10})$. The discriminant is an invariant $\Delta$, of
degree $14$, and the test can thus be easily done once one has expressed
$\Delta$ as a polynomial in the $J_i$s which is
{\rm\setcounter{myequation}{\theequation}
  \begin{dgroup*}[style={\footnotesize},spread={-5pt}]
    \begin{dmath}[style={\refstepcounter{myequation}},number={\themyequation},label={eq:disc}]
      351157512403353600\,{\it J_{2}}\,{\it J_{3}}\,{\it J_{4}}\,{\it J_{5}}-
      106232524760678400\,{{\it J_{7}}}^{2}+369994358063104\,{{\it
          J_{2}}}^{7}- 66189456176578560\,{{\it J_{2}}}^{3}{\it J_{3}}\,{\it
        J_{5}}+33474799101542400 \,{{\it J_{2}}}^{2}{{\it J_{3}}}^{2}{\it
        J_{4}}-1292776756346880000\,{{\it J_{2}}} ^{2}{\it J_{3}}\,{\it
        J_{7}}+1267160326927810560\,{{\it J_{2}}}^{2}{\it J_{4}}\,{ \it
        J_{6}}-1143872189796188160\,{\it J_{2}}\,{{\it J_{3}}}^{2}{\it J_{6}}-
      16397790463918080\,{{\it J_{2}}}^{5}{\it J_{4}}+10245839137013760\,{{\it
          J_{2}} }^{4}{{\it J_{3}}}^{2}-93857730923593728\,{{\it
          J_{2}}}^{4}{\it J_{6}}+ 313290675380551680\,{{\it J_{2}}}^{3}{{\it
          J_{4}}}^{2}+2345327219866337280\, {{\it J_{2}}}^{3}{\it
        J_{8}}-895451656628551680\,{{\it J_{2}}}^{2}{{\it J_{5}}}^{2
      }+18346188403113984000\,{{\it J_{2}}}^{2}{\it
        J_{10}}+31476303632793600\,{ \it J_{2}}\,{{\it
          J_{3}}}^{4}+968390964336918528\,{\it J_{2}}\,{{\it J_{4}}}^{3}+
      4879526536497070080\,{\it J_{2}}\,{{\it
          J_{6}}}^{2}+212465049521356800\,{{ \it J_{3}}}^{3}{\it
        J_{5}}+169419330714009600\,{{\it J_{3}}}^{2}{{\it J_{4}}}^{2}+
      2549580594256281600\,{{\it J_{3}}}^{2}{\it
        J_{8}}+2751208802098348032\,{{ \it J_{4}}}^{2}{\it
        J_{6}}+29008561427982581760\,{\it J_{4}}\,{\it J_{10}}-
      21756421070986936320\,{\it J_{5}}\,{\it J_{9}}-3021725148748185600\,{\it
        J_{6}} \,{\it J_{8}}-7554312871870464000\,{\it J_{2}}\,{\it
        J_{3}}\,{\it J_{9}}+ 31803657190574653440\,{\it J_{2}}\,{\it
        J_{4}}\,{\it J_{8}}-8849169312564510720 \,{\it J_{2}}\,{\it
        J_{5}}\,{\it J_{7}}-5139293475644375040\,{\it J_{3}}\,{\it J_{4}}
      \,{\it J_{7}}+4582949808934748160\,{\it J_{3}}\,{\it J_{5}}\,{\it
        J_{6}}\,.
    \end{dmath}
  \end{dgroup*}
  \setcounter{equation}{\themyequation} }

\section{Reconstruction of binary forms} \label{sec:reconstruction}

Let $n>2$ be an even positive integer, $k$ be a field of characteristic $0$ or
greater than $n$ and $K=\bar{k}$.

\subsection{Clebsch's identities}
\label{sec:some-clebschs-ident}

In the following section we give proofs of the results in the Example 3 of
\cite{mestre}. They are the fundamental tools to reconstruct a generic binary
form of even degree $n$ from its invariants as it is explained in
\textit{loc. cit.} The classical reference for this section is
\cite[§103]{clebsch} but we wish to give full proofs using the material of
Section \ref{sec:rappel}. Note that all formulae but the last one could also
be checked by a computer algebra software.\smallskip

\subsubsection{Useful identities for three quadratic forms}
Let $q_1,q_2,q_3$ be three elements of $S^2(V)$, \ie binary homogeneous
quadratic forms in $x,z$ over $k$. We denote
$$q_1^*=(q_2,q_3)_1, \quad q_2^*=(q_3,q_1)_1, \quad q_3^*=(q_1,q_2)_1.$$

\begin{lemma} \label{lemma:ortho} We have the relation $$q_1 q_1^* +q_2 q_2^*
  + q_3 q_3^* =0.$$
\end{lemma}
\begin{proof}
  Let us consider the multilinear map $S^2(V^*)^3 \to S^4(V^*)$,
  $(q_1,q_2,q_3) \mapsto \sum_{i=1}^3 q_i q_i^*$. By multilinearity, to check
  the equality to zero, it is enough to verify it on linear powers
  $q_i=\ell_i^2$. Since $(\ell_j^2,\ell_k^2)_1=[\ell_j,\ell_k] \ell_j \ell_k$,
  we get
$$\sum q_i q_i^* = \ell_1 \ell_2 \ell_3 \cdot  ([\ell_2,\ell_3] \ell_1 + [\ell_3,\ell_1] \ell_2 + [\ell_1,\ell_2] \ell_3).$$
Of course $w:={^t
  \left([\ell_2,\ell_3],[\ell_3,\ell_1],[\ell_1,\ell_2]\right)}={^t
  \left(\alpha_1,\alpha_2,\alpha_3\right)} \wedge {^t
  \left(\beta_1,\beta_2,\beta_3\right)}$, so $$w \cdot {^t \left(\alpha_1 x +
    \beta_1 z,\alpha_2 x+ \beta_2 z,\alpha_3 x + \beta_3 z\right)}=0.$$ This
proves the equality to $0$ (this is called the first fundamental identity in
\cite[p.19]{grace-young}).
\end{proof}
\medskip

\begin{lemma} \label{lemma:detform} Let $R(q_1,q_2,q_3)$ be the determinant of
  $q_1,q_2,q_3$ in the basis $x^2, xz,z^2$. Then $$-4 \cdot
  R(q_1,q_2,q_3)=(\Omega_{12} \cdot \Omega_{23} \cdot \Omega_{31})
  (q_1(x_1,z_1)\cdot q_2(x_2,z_2) \cdot q_3(x_3,z_3)).$$
\end{lemma}
\begin{proof}
  The operator $\Omega_{12} \cdot \Omega_{23} \cdot \Omega_{31} : S^2(V^*)^3
  \to k$ is multi-linear and anti-symmetric. Hence it is equal to
  $R(q_1,q_2,q_3)$ up to a constant.
\end{proof}
\medskip

\begin{lemma} \label{lemma:detconic} Let $A_{ij}=(q_i,q_j)_2$ for $i,j \in
  \{1,2,3\}$. Then $$2 \cdot \det ((A_{ij})_{i,j=1,2,3})= R^2(q_1,q_2,q_3).$$
\end{lemma}
\begin{proof}
  We follow \cite[p.79]{grace-young}. Let us consider the multilinear form
  $\phi : S^2(V^*)^6 \to k$ given by
$$(q_1,q_2,q_3,q_4,q_5,q_6) \mapsto \det \begin{bmatrix}  (q_1,q_4)_2 & (q_1, q_5)_2 & (q_1,q_6)_2 \\
  (q_2,q_4)_2 & (q_2, q_5)_2 & (q_2,q_6)_2 \\
  (q_3,q_4)_2 & (q_3, q_5)_2 & (q_3,q_6)_2 \\
\end{bmatrix}.$$ If we evaluate $\phi$ on the elements of the basis
$q_i=\ell_i^2$, since $(q_i,q_j)_2=[\ell_i, \ell_j]^2=(\alpha_i \beta_j
-\alpha_j \beta_i)^2$, we get
\begin{eqnarray*}
  \phi(\ell_1^2,\ldots,\ell_6^2)&=&
  \det \begin{bmatrix} \alpha_1^2 & \alpha_1 \beta_1 & \beta_1^2 \\ \alpha_2^2 & \alpha_2 \beta_2 & \beta_2^2 \\
    \alpha_3^2 & \alpha_3 \beta_3 & \beta_3^2 \end{bmatrix} \cdot \det \begin{bmatrix} \beta_4^2 & \beta_5^2 & \beta_6^2 \\ -2 \alpha_4 \beta_4 & -2 \alpha_5 \beta_5 & -2 \alpha_6 \beta_6  \\ \alpha_4^2 & \alpha_5^2 &\alpha_6^2 \end{bmatrix} \\
  &=&  \frac{1}{2} \cdot  R(q_1,q_2,q_3) \cdot R(q_4,q_5,q_6)
\end{eqnarray*} 
We conclude by letting $q_4=q_1,q_5=q_2$ and $q_6=q_3$.
\end{proof}
\medskip

\begin{lemma}
Let $f$ be any quadratic form. Then 
$$\frac{1}{2} \cdot R(q_1,q_2,q_3) \cdot f = (f,q_1^*)_2 \cdot q_1+ (f,q_2^*)_2 \cdot q_2+ (f,q_3^*)_2 \cdot q_3\,.$$
\end{lemma}
\begin{proof}
  Let us consider the multi-linear map $\psi : S^2(V^*)^4 \to S^2(V^*)$ given
  by the difference of the left and right member of the equality above. We
  want to prove that $\psi$ is the zero map and again it is enough to do this
  on the basis $q_i=\ell_i^2$ and $f=\ell^2=(\mu x+ \nu z)^2$. We proceed like
  in \cite[p.414]{clebsch}.  Obviously
$$\det \begin{bmatrix} \ell_1^2 & \alpha_1^2 & \alpha_1 \beta_1 & \beta_1^2 \\  \ell_2^2 & \alpha_2^2 & \alpha_2 \beta_2 & \beta_2^2 \\  \ell_3^2 & \alpha_3^2 & \alpha_3 \beta_3 & \beta_3^2 \\  f & \mu^2 & \mu \nu & \nu^2 \end{bmatrix}=0.$$
Expanding the determinant with respect to the first column, we get
$$-f R(q_1,q_2,q_3) + q_3 R(q_1,q_2,f) - q_2 R(q_1,q_3,f) + q_1 R(q_2,q_3,f)=0.$$
To get our result, we only need to prove that $R(q_i,q_j,f)=2 \cdot
(f,(q_i,q_j)_1)_2$.  By a direct computation, $R(q_i,q_j,f)=2 [ \ell_i,
\ell_j] \cdot [\ell,\ell_j] \cdot [\ell, \ell_i]$. On the other hand,
\begin{eqnarray*}
  (f,(q_i,q_j)_1)_2&=& (f,[\ell_i,\ell_j] \cdot \ell_i  \ell_j)_2=[\ell_i,\ell_j] \cdot (\ell^2,\ell_i \ell_j)_2\,, \\
  &=& \frac{1}{4} [\ell_i,\ell_j] \cdot \Omega_{12}^2 \left[\ell(x_1,z_1)^2 \ell_i(x_2,z_2) \ell_j(x_2,z_2) \right]\,, \\
  &=& [ \ell_i, \ell_j] \cdot [\ell, \ell_j] \cdot [\ell, \ell_i]\,.
\end{eqnarray*}
\end{proof}
Mestre's formula is somehow dual to this result and can be checked using the same arguments.
\begin{lemma} \label{lemma:decomp} Let $f$ be any quadratic form. Then
$$R(q_1,q_2,q_3) \cdot f = 2 \left((f,q_1)_2 \cdot q_1^*+ (f,q_2)_2 \cdot q_2^*+ (f,q_3)_2 \cdot q_3^*\right).$$
\end{lemma}

\subsubsection{Relations among a binary form and three quadratic forms}

We now come to the two fundamental results.
\begin{proposition}\label{prop:conic}
We have $\sum_{i,j} A_{ij} \cdot  q_i^* \cdot q_j^*=0$.
\end{proposition}
\begin{proof}
  For $i \in \{1,2,3\}$ using Lemma \ref{lemma:decomp}, we get
$$R(q_1,q_2,q_3) q_i = (q_i,q_1)_2 q_1^* + (q_i,q_2)_2 q_2^* + (q_i,q_3)_2 q_3^*.$$ 
We multiply the previous equality by $q_i^*$,
$$R(q_1,q_2,q_3) q_i q_i^* = A_{i1}  q_i^* q_1^* + A_{i2} q_i^* q_2^*+ A_{i3} q_i^* q_3^*,$$
and we sum on $1 \leq i \leq 3$,
$$R(q_1,q_2,q_3) \cdot \sum_{i=1}^3 q_i q_i^* = \sum_{i=1}^3  A_{i1}  q_i^* q_1^* + A_{i2} q_i^* q_2^*+ A_{i3} q_i^* q_3^*=\sum_{1 \leq i,j \leq 3} A_{ij} q_i^* q_j^*,$$
which is zero by Lemma \ref{lemma:ortho}.
\end{proof}
\medskip

\begin{proposition}\label{prop:multic}
  Let $f$ be a binary form of even degree $n$. Then
$$R(q_1,q_2,q_3)^{n/2} \cdot f(x,z) = \frac{1}{n!} \cdot \left(\sum_{i=1}^3  q_i^*(x,z) \delta_i \right)^{n/2} f(x_1,z_1)$$
where $\delta_i$ is the differential operator $\phi(x_1,z_1) \mapsto
\Omega_{12}^2 (\phi(x_1,z_1) \cdot q_i(x_2,z_2))$.
\end{proposition}
\begin{proof}
  As the right member is linear in $f$ we can assume that $f=\ell^n=(\mu
  x_1+\nu z_1)^n$. To make it linear in the $q_i$ we consider $n/2$ triplets
  of quadratic forms $(q_{1j},q_{2j},q_{3j})$ and define
  accordingly $$q_{1j}^*=(q_{2j},q_{3j})_1, \quad q_{2j}^*=(q_{3j},q_{1j})_1,
  \quad q_{3j}^*=(q_{1j},q_{2j})_1.$$ We replace the previous operator
  by $$\Psi:=\frac{1}{n!} \cdot \prod_{j=1}^{n/2} \sum_{i=1}^3 q_{ij}^*(x,z)
  \, \delta_{ij}$$ where $$\delta_{ij}(\phi(x_1,z_1))=
  \Omega_{12}^2(\phi(x_1,z_1) q_{ij}(x_2,z_2)).$$ As $\Psi$ is linear in each
  $q_{ij}$ we can assume that $$q_{ij}=\ell_{ij}^2=(\alpha_{ij} x+\beta_{ij}
  z)^2.$$ For $1 \leq j \leq n/2$, let $i_j \in \{1,2,3\}$.  Observe that for
  all $m$
$$\delta_{i_j j}(\ell^m)=\Omega_{12}^2(\ell^m(x_1,z_1)q_{i_j j}(x_2,z_2))=\frac{2! \, m! }{(m-2)!} \cdot [\ell,\ell_{i_j j}] \cdot \ell^{m-2}.$$ Hence
$$
\frac{1}{n!} \left(\prod_{j=1}^{n/2} \delta_{i_j j}\right)(f)=
\prod_{j=1}^{n/2} 2 [\ell,\ell_{i_j j}]^2 = \prod_{j=1}^{n/2} 2
(\ell^2,q_{i_j,j})_2.$$ Therefore, if we develop the expression of $\Psi$, we
can replace each product of the $\delta_{i_j j}$ operators on $f$ by the right
expression. Re-factoring the new expression, we get that
$$\Psi(f)= \prod_{j=1}^{n/2} 2 \sum_{i=1}^3 q_{ij}^* \cdot (\ell^2,q_{ij})_2.$$
Using Lemma \ref{lemma:decomp},
$$\Psi(f)=\prod_{j=1}^{n/2} \left( R(q_{1j},q_{2j},q_{3j}) \cdot \ell^2 \right) = \left(\prod_{j=1}^{n/2} R(q_{1j},q_{2j},q_{3j})\right) \cdot \ell^n.$$ 
To conclude, we let $q_{ij}=q_i$ for all $j$.

\end{proof}

\subsection{A generic reconstruction algorithm}
\label{sec:general-strategy}

Starting from a weighted projective point $(\iota_1 : \iota_2 : \ldots)$ of
values in $k$ for a finite set of generators $\{I_i\}$ of $\ii_{2g+2}$, we aim
at recovering a hyperelliptic curve $C/K~: y^2=f(x)$ such that $(\iota_1 :
\iota_2 : \ldots)=(I_1(f) : I_2(f) : \ldots)$. Inverting the polynomial system
giving the invariants in terms of a generic polynomial $f$ can be efficiently
done only in very specific cases (see Section
\ref{sec:ident-loci-reconstr}). However, Mestre explained in his Example 3 and
Remark p.~321 \cite{mestre} that for a generic hyperelliptic curves, one can
use the background developed in Section \ref{sec:some-clebschs-ident} in the
following way.\medskip

Let $f(X,Z)$ be a binary form of even degree $2g+2$ defined over $K$.  Using
transvectants one can compute three covariants $q_i(x,z)=q_i(f,(x,z))$ of
order $2$ defined over $K$. With the notation of
Section~\ref{sec:some-clebschs-ident}, we denote
$$q_1^*=(q_2,q_3)_1, \quad q_2^*=(q_3,q_1)_1, \quad q_3^*=(q_1,q_2)_1.$$
From Propositions~\ref{prop:conic} and~\ref{prop:multic}, we see that the
point $P=(x_1^*:x_2^*:x_3^*)=(q_1^*(x,z) : q_2^*(x,z) : q_3^*(x,z))$ is
solution of $$
\begin{cases}
  \sum A_{ij} x_i x_j=0, \\
  \sum_I h_I x_I - R(q_1,q_2,q_3)^{n/2} f(x,z)=0,
\end{cases}$$ where $I$ is a multi-index of cardinality $n/2$.

Clearly the $A_{ij}=(q_i,q_j)_2$ are invariants. In Proposition
 \ref{prop:multic}, each time the operator $\Omega_{ij}^2$ is applied on a
 covariant of order $r$ and on a covariant $q_i$ of order $2$, we get a
 covariant of order $r+2-4=r-2$. As we started with $f$ which is a covariant
 of order $n$ and we applied the operators $n/2$ times, the coefficients $h_I$
 are invariants. Hence, both $A_{ij}$ and $h_I$ being invariants, they can be
 express as polynomials in the $I_i$s with coefficients in the prime field of
 $K$. Thus if the invariants of a given $f(X,Z)$ lives inside a field $k
 \subset K$, the specializations of the $A_{ij}$ and $h_I$ belongs to $k$ as
 well. We refer to Section~\ref{sec:conic-quartic} for a strategy to find the
 formal expressions when $n=8$.

 Now, assume that $R(q_1,q_2,q_3) \ne 0$. Note that this is an invariant as
 well, and so it can be computed as a polynomial in the $I_i$s. Then $P$ lies
 on the curve $\hh/k : \sum_I h_I x_I=0$ if and only if $(x : z)$ is a root of
 $f(x,z)$. Moreover, by Lemma \ref{lemma:detconic}, the conic $\qq/k : \sum
 A_{ij} x_i x_j=0$ is non singular, hence the morphism $(x:z) \mapsto
 (x_1^*:x_2^* : x_3^*) \in \qq$ is a parametrization. As there a quadratic
 extension of $k$, we get a parametrization $\phi : (x,z) \mapsto
 (\chi_1:\chi_2:\chi_3)$.  All parametrizations being $\GL_2$-equivalent, the
 binary polynomial $\sum_I h_I \chi_I$ is $\GL_2$-equivalent to $f$. Hence we
 get the following proposition.
 \begin{proposition} \label{prop:generic-reconstructions} Let $(q_1,q_2,q_3)$
   be three covariants of order $2$ of a binary form $f$ of even degree $n$
   defined over $K$. If $R(q_1,q_2,q_3) \ne 0$, there exists a non-singular
   conic $\qq$ and a plane curve $\hh$ of degree $n/2$ defined over a subfield
   $k$ of $K$ such that there is a $K$-isomorphism $\qq \to \PP^1$ mapping the
   intersections points of $\qq \cap \hh$ on the roots of $f(X,Z)$. Moreover,
   this isomorphism is defined at most over a quadratic extension of $k$ and
   is defined over $k$ as soon as $\qq$ has a $k$-rational point.
 \end{proposition}
 This is Mestre's main result in \cite[Sec.2.4]{mestre}, that for a generic
 hyperelliptic curve, the obstruction to reconstruct the curve from its
 invariant over $k$ is equivalent to $\qq(k)=\emptyset$. We will come back on
 these questions in Section \ref{subsec:def} and relate $k'$ to the field of
 moduli. \medskip

 To conclude, given values $(\iota_1 : \iota_2 : \ldots)$ defined in $k$, the
 general method is as follows.
 \begin{enumerate}
 \item Find a triple $(q_1, q_2, q_3)$ of covariants of order 2 such that that
   the expression of $R(q_1, q_2, q_3)$ as a polynomial in the $I_i$s
   evaluated at the $\iota_is$ is non-zero;
 \item Find a point on the conic $\qq$ defined by the expression of the
   coefficients $A_{ij}(q_1, q_2, q_3)$ as polynomials in the $I_i$s evaluated
   at the $\iota_is$, and get a parametrization for the conic.
 \item Let $\hh$ be the curve defined by the expression of the coefficients
   $h_{I}(q_1, q_2, q_3)$ as polynomials in the $I_i$ evaluated at the
   $\iota_i$s. After a parametrization of the conic $\qq$ return the
   intersection divisor $\qq \cdot \hh$ as a polynomial $f(X, Z)$ over at most
   a quadratic extension of $k$.
 \end{enumerate}

 It might happen that all possible $R(q_1,q_2,q_3)$ evaluated at the
 $\iota_is$ are zero. Even if it is not the case and for a general field $K$,
 $\qq$ might have no $k$-rational point although we know, for theoretical
 reasons, we can reconstruct the curve over $k$.  Section
 \ref{sec:ident-loci-reconstr} and Section \ref{sec:field-definition} will
 deal with these issues.

\subsection{Reconstruction in the hyperelliptic genus 3 case}
\label{sec:conic-quartic}

\begin{table}[htbp]
  \centering
  \begin{sideways}
    \begin{footnotesize}
      \begin{math}\setlength{\arraycolsep}{2pt}
        {\begin{array}{c||c|cccccccc|c|c}
            \begin{array}{c}
              \mathrm{\ \ \ Ord.}\\
              \mathrm{Deg.\ \ \ }
            \end{array}
            &   0 &   2 &   4 &   6 &   8 &  10 &  12 &  14 &  18 & \mathrm{Tot} \\\hline\hline
            &&&&&&&&&&\\ 
            1 &   - &   - &   - &   - &\mathfrak{f}&   - &   - &   - &   - &   1 \\[0.2cm]

            2 &(\,\mathfrak{f},\,\mathfrak{f})_8
            &   - &(\,\mathfrak{f},\,\mathfrak{f})_6
            &   - &(\,\mathfrak{f},\,\mathfrak{f})_4
            &   - &(\,\mathfrak{f},\,\mathfrak{f})_2
            &   - &   - &   4 \\[0.2cm]

            3 &(\,C_{2,8},\,\mathfrak{f})_8
            &   - &(\,C_{2,8},\,\mathfrak{f})_6
            &(\,C_{2,8},\,\mathfrak{f})_5
            &(\,C_{2,8},\,\mathfrak{f})_4
            &(\,C_{2,8},\,\mathfrak{f})_3
            &(\,C_{2,8},\,\mathfrak{f})_2
            &(\,C_{2,8},\,\mathfrak{f})_1
            &(\,C_{2,12},\,\mathfrak{f})_1
            &   8 \\[0.2cm]

            4 &(\,C_{3,8},\,\mathfrak{f})_8
            &   - &\left|\begin{array}{l}
                (\,C_{3,4},\,\mathfrak{f})_4\\
                (\,C_{3,8},\,\mathfrak{f})_6\\
              \end{array}\right.
            &(\,C_{3,4},\,\mathfrak{f})_3
            &(\,C_{3,4},\,\mathfrak{f})_2
            &\left|\begin{array}{l}
                (\,C_{3,4},\,\mathfrak{f})_1\\
                (\,C_{3,8},\,\mathfrak{f})_3\\
              \end{array}\right.
            &(\,C_{3,8},\,\mathfrak{f})_2
            &(\,C_{3,8},\,\mathfrak{f})_1
            &(\,C_{3,12},\,\mathfrak{f})_1
            &  10 \\[0.2cm]

            5 &(\,C_{4,8},\,\mathfrak{f})_8
            &(\,C_{4,10},\,\mathfrak{f})_8
            &\left|\begin{array}{l}
                (\,C_{4,10},\,\mathfrak{f})_7\\
                (\,C_{4,8},\,\mathfrak{f})_6\\
              \end{array}\right.
            &\left|\begin{array}{l}
                (\,C_{4,10},\,\mathfrak{f})_6\\
                (\,C_{4,8},\,\mathfrak{f})_5\\
              \end{array}\right.
            &(\,C_{4,10},\,\mathfrak{f})_5
            &\left|\begin{array}{l}
                (\,C_{4,8},\,\mathfrak{f})_3\\
                (\,C_{4,10},\,\mathfrak{f})_4\\
                (\,C_{4,10}',\,\mathfrak{f})_4\\
              \end{array}\right.
            &   - &(\,C_{4,10},\,\mathfrak{f})_2
            &   - &  11 \\[0.2cm]

            6 &(\,C_{3,4}\,C_{2,4},\,\mathfrak{f})_8
            &(\,C_{5,8},\,\mathfrak{f})_7
            &\left|\begin{array}{l}
                (\,C_{5,8},\,\mathfrak{f})_6\\
                (\,C_{5,4}',\,\mathfrak{f})_4\\
              \end{array}\right.
            &\left|\begin{array}{l}
                (\,C_{5,8},\,\mathfrak{f})_5\\
                (\,C_{5,4}',\,\mathfrak{f})_3\\
                (\,C_{5,10}',\,\mathfrak{f})_6\\
              \end{array}\right.
            &(\,C_{5,4}',\,\mathfrak{f})_2
            &(\,C_{5,4}',\,\mathfrak{f})_1
            &   - &   - &   - &   9 \\[0.2cm]

            7 &(\,C_{2,4}\,C_{4,4}',\,\mathfrak{f})_8
            &\left|\begin{array}{l}
                (\,C_{2,4}\,C_{4,6},\,\mathfrak{f})_8\\
                (\,C_{6,6}'',\,\mathfrak{f})_6\\
              \end{array}\right.
            &\left|\begin{array}{l}
                (\,C_{2,4}\,C_{4,6},\,\mathfrak{f})_7\\
                (\,C_{6,6}'',\,\mathfrak{f})_5\\
              \end{array}\right.
            &\left|\begin{array}{l}
                (\,C_{6,6}'',\,\mathfrak{f})_4\\
                (\,C_{6,2},\,\mathfrak{f})_2\\
                (\,C_{2,4}\,C_{4,6},\,\mathfrak{f})_6\\
              \end{array}\right.
            &   - &   - &   - &   - &   - &   8 \\[0.2cm]

            8 &(\,C_{3,4}\,C_{4,4},\,\mathfrak{f})_8
            &\left|\begin{array}{l}
                (\,C_{2,8}\,C_{5,2},\,\mathfrak{f})_8\\
                (\,C_{3,6}\,C_{4,4},\,\mathfrak{f})_8\\
              \end{array}\right.
            &\left|\begin{array}{l}
                (\,C_{3,6}\,C_{4,4},\,\mathfrak{f})_7\\
                (\,C_{3,4}\,C_{4,6},\,\mathfrak{f})_7\\
              \end{array}\right.
            &\left|\begin{array}{l}
                (\,C_{3,6}\,C_{4,4},\,\mathfrak{f})_6\\
                (\,C_{3,4}\,C_{4,6},\,\mathfrak{f})_6\\
              \end{array}\right.
            &   - &   - &   - &   - &   - &   7 \\[0.2cm]

            9 &(\,C_{2,4}\,C_{6,4},\,\mathfrak{f})_8
            &\left|\begin{array}{l}
                (\,C_{4,6}\,C_{4,4}',\,\mathfrak{f})_8\\
                (\,C_{2,4}\,C_{6,4},\,\mathfrak{f})_7\\
                (\,C_{2,4}\,C_{6,6}',\,\mathfrak{f})_8\\
              \end{array}\right.
            &(\,C_{2,4}\,C_{6,4},\,\mathfrak{f})_6
            &   - &   - &   - &   - &   - &   - &   5 \\

            10 &(\,C_{4,4}\,C_{5,4}',\,\mathfrak{f})_8
            &\left|\begin{array}{l}
                (\,C_{7,2}'\,C_{2,4},\,\mathfrak{f})_6\\
                (\,C_{4,6}\,C_{5,4},\,\mathfrak{f})_8\\
              \end{array}\right.
            &   - &   - &   - &   - &   - &   - &   - &   3 \\[0.2cm]

            11 &   - &\left|\begin{array}{l}
                (\,C_{8,4}'\,C_{2,4},\,\mathfrak{f})_7\\
                (\,C_{5,6}'\,C_{5,4}',\,\mathfrak{f})_8\\
              \end{array}\right.
            &   - &   - &   - &   - &   - &   - &   - &   2 \\[0.2cm]

            12 &   - &(\,C_{6,6}'\,C_{5,4}',\,\mathfrak{f})_8
            &   - &   - &   - &   - &   - &   - &   - &   1 \\[0.2cm]\hline
            &&&&&&&&&&\\ 
            \textrm{Tot} &   9 &  14 &  13 &  12 &   6 &   7 &   3 &   3 &   2 &  69 \\
          \end{array}
        }
      \end{math}
    \end{footnotesize}
  \end{sideways}
  \bigskip
  \caption{Fundamental invariants and covariants for a binary form  of degree $8$}
  \label{tab:fundamentals}
\end{table}

We have computed a system of fundamental generators for invariants and
covariants of octics, using Gordan's algorithm.  Results are given in
Table~\ref{tab:fundamentals}. In this table, generators are all
defined by the mean of transvectants of the form $(\prod_{d,r} C_{d,
  r}, \mathfrak{f})_h$ where we denote recursively by $C_{d, r}$ generators of
degree $d$ and order $r$ given at the intersection of the row $d$ and the
column $r$ of the table. When we have two or three generators of degree $d$ and $r$, we denote them by $C_{d, r}'$ and $C_{d, r}''$. For
instance, $C_{5,10} = (\,C_{4,8}\,,\,\mathfrak{f})_3$, $C_{5,10}' =
(\,C_{4,10}\,,\,\mathfrak{f})_4$ and $C_{5,10}'' = (\,C_{4,10}'\,,\,\mathfrak{f})_4$.
Covariants $C_{2, 0}$, $C_{3, 0}$, $\ldots$, $C_{10, 0}$ given in column
0 are invariants and so, they are related to Shioda invariants
$J_2$, $J_3$, $\ldots$, $J_{10}$ by
\begin{dgroup*}[compact,style={\tiny},spread={0pt}]
  \begin{dmath*}
    C_{2,0} = {J_2}\,,\ C_{3,0} = {J_3}\,,\ C_{4,0} = \frac
    {1}{30}\,{{J_2}}^{2}-{\frac {4}{35}}\,{J_4}\,,\ C_{5,0} = \frac
    {3}{14}\,{J_5}\,,\ 
  \end{dmath*}
  \begin{dmath*}
    C_{6,0} = -{\frac {1}{1890}}\,{{J_2}}^{3}+\frac
    {1}{15}\,{J_2}\,{J_4}+{\frac {1}{63}}\,{{J_3}}^{2}+{\frac
      {8}{21}}\,{J_6}\,,
  \end{dmath*}
  \begin{dmath*}
    C_{7,0} = \frac{1}{42}\,{J_2}\,{J_5}-{\frac
      {3}{245}}\,{J_3}\,{J_4}+\frac{1}{49}\,{J_7}\,,
  \end{dmath*}
  \begin{dmath*}
    C_{8,0} = -{\frac {11}{264600}}\,{{J_2}}^{4}+{\frac
      {59}{14700}}\,{{J_2}}^{2}{J_4}+{\frac
      {11}{8820}}\,{J_2}\,{{J_3}}^{2}+{\frac {22}{735}}\,{J_2}\,{J_6}+{\frac
      {3}{980}}\,{J_3}\,{J_5}+{\frac {198}{8575}}\,{{J_4}}^{2}+{\frac
      {9}{49}}\,{J_8},
  \end{dmath*}
  \begin{dmath*}
    C_{9,0} = -{\frac {1}{45360}}\,{{J_2}}^{3}{J_3}+{\frac
      {1}{1260}}\,{{J_2}}^{2}{J_5}+{\frac
      {61}{17640}}\,{J_2}\,{J_3}\,{J_4}-{\frac
      {17}{4704}}\,{J_2}\,{J_7}+{\frac {1}{1512}}\,{{J_3}}^{3}+{\frac
      {1}{144}}\,{J_3}\,{J_6}-{\frac
      {1}{112}}\,{J_4}\,{J_5}-\frac{1}{28}\,{J_9}\,,
  \end{dmath*}
  \begin{dmath*}
    C_{10,0} = -{\frac {1}{595350}}\,{{J_2}}^{5}+{\frac
      {107}{411600}}\,{{J_2}}^{3}{J_4}+{\frac
      {1}{19845}}\,{{J_2}}^{2}{{J_3}}^{2}+{\frac
      {317}{105840}}\,{{J_2}}^{2}{J_6}+{\frac
      {5}{2352}}\,{J_2}\,{J_3}\,{J_5}+{\frac
      {293}{180075}}\,{J_2}\,{{J_4}}^{2}
    + {\frac {25}{5488}}\,{J_2}\,{J_8}+{\frac
      {5}{6174}}\,{{J_3}}^{2}{J_4}+{\frac {17}{2744}}\,{J_3}\,{J_7}+{\frac
      {19}{5145} }\,{J_4}\,{J_6}-{\frac {153}{27440}}\,{{J_5}}^{2}+{\frac
      {3}{49}}\,{J_{10}}\,.
  \end{dmath*}
\end{dgroup*}

Our main motivation for computing Table~\ref{tab:fundamentals} is to
determine fundamental covariants of order $2$ which can be used as quadratic
forms in Proposition~\ref{prop:generic-reconstructions} to reconstruct a
generic binary octic from its Shioda invariants.  As with Shioda invariants,
we paid attention that none of the transvectant computations in
Table~\ref{tab:fundamentals} involves covariants of order greater than $10$
(except the 2 covariants of order 18, but it does not matter since these 2
covariants are not useful to compute other invariants or covariants of order $2$). These formulas are
thus also valid, in addition to fields of characteristic 0, over any field of or
characteristic greater or equal than 11.\\

The main computational difficulty in the reconstruction method of
Section~\ref{sec:general-strategy} is to write the invariants $A_{ij}$ and $h_I$  as
 polynomials in the $J_i$'s, since their degree may be large (close to $40$ in
our cases).  Writing them as a polynomial with $9$ unknowns
$a_i$'s for a generic form $f = a_8 X^8 + a_7X^7 + \ldots + a_0$ is hopeless.
We follow instead a ``black-box'' approach. 

\begin{figure}[htbp]
  \parbox{0.95\linewidth}{%
    \begin{footnotesize}\SetAlFnt{\small\sf}%
      \begin{algorithm}[H]%
        \label{RecognizeJInvariants}%
        \caption{Write an invariant as a polynomial in the $J_i$'s.}  %
        \SetKwInOut{Input}{Input} \SetKwInOut{Output}{Output} %
        \Input{An invariant $I$ of degree $d$ (given as an evaluation
          program).} %
        \Output{A polynomial $P$ in $\QQ[J_2, \ldots, J_{10}]$
          such that $I(f) = P(J_2(f), J_3(f), \ldots, J_{10}(f))$\,. }
        \BlankLine
        \tcp{A basis for the polynomials of degree $d$ in
          the weighted graded algebra $\QQ[J_2, J_3, \ldots, J_{10}]$}
        ${\mathcal B} \leftarrow \left[\,\prod_{w} J_w ^{e_w}\ \mathtt{s.t.}\ \sum_w
        w\,e_w = d\,\right]$\;
        \BlankLine
        \tcp{Choose at random $\#{\mathcal B}+O(1)$ octics over $\QQ$}
        ${\mathcal F} \leftarrow \left[\, a_8X^8+\ldots +a_1X\,Z^7+a_0Z^8\
          \mathtt{for}\ \#{\mathcal B}+O(1)\ 9\mathtt{-uples}\ (a_0,
          a_1, \ldots, a_8)\mathtt{\ randomly\ chosen\ in\ }\QQ^9\, \right]$\;
        \BlankLine
        \tcp{Evaluate the invariant $I$ and the basis ${\mathcal B}$ at each
          form of ${\mathcal F}$}
        \For{$i = 1$ \KwTo $\#{\mathcal F}$}{
          $V_i \leftarrow I({\mathcal F}_i)$\;
          \For{$j = 1$ \KwTo $\#{\mathcal B}$}{
            $M_{i,j} \leftarrow {\mathcal B}_j({\mathcal F}_i)$
          }
        }
        \BlankLine
        \tcp{Invert the linear system defined by the matrix $M$ and the vector
        $V$}
        Find the vectors $U$ such that $M \times U = V$\;
        \BlankLine
        \KwRet{$\sum_i U_i\,{\mathcal B}_i$} for each $U$
   \end{algorithm}
 \end{footnotesize}
}
\end{figure}

More precisely, since invariants are computed through sequences of covariants
which are the result of transvectant or differential operations, we represent
a covariant no more by a formal expression in the $a_i$'s but as an algorithm
which performs the corresponding sequence of operations.  If one inputs a
generic form $f \in \QQ[a_0, \ldots, a_8][X, Z]$, such an algorithm returns
the formal expression of the covariant as a multivariate homogeneous
polynomial in the $a_i$s, $X$ and $Z$. But, if one inputs a form $f \in \QQ[X,
Z]$, the algorithm returns the covariant as a homogeneous polynomial in
$\QQ[x,z]$, without the use of the formal expression in the $a_i$s.  In other
words, we consider that a covariant is given by an evaluation program. Note
that it is immediate to determine the degree and the order of a covariant from
the sequence of operations which compose its evaluation program.

Now coming back to the question of writing a homogeneous invariant as a
polynomial in the $J_i$s, we propose to construct a basis ${\mathcal B} = \{\
\prod_{w} J_w ^{e_w}~:\ \sum_w w\,e_w = d\ \}$ for the polynomials of degree
$d$ in the weighted graded algebra $\QQ[J_2, J_3, \ldots, J_{10}]$ and we
evaluate this basis (with some given evaluation programs for the $J_i$'s), and
$I$ (given as an evaluation program too), at $\#{\mathcal B}+O(1)$ octics
chosen at random over $\QQ$.  It remains to invert the corresponding linear
system to find $I$ as a polynomial in the $J_i$'s.
Algorithm~\ref{RecognizeJInvariants} summarizes this method.

\begin{remark}
  Working with $\QQ[J_2, J_3, \ldots, J_{10}]$ instead of $\QQ[a_0, a_1,
  \ldots, a_{8}]$ in Algorithm~\ref{RecognizeJInvariants} makes a real
  difference in practice. For instance, an invariant of degree $20$ yields a
  basis ${\mathcal B}$ with only $107$ monomials. This must be compared to the
  $61731$ monomials that otherwise we would have to deal with in the weighted
  projective points $(a_0 : a_1 : \ldots, a_{8})$ where $a_0$, $a_1$, \ldots,
  $a_{8}$ are of weight $0$, $1$, \ldots, $8$. Note that an invariant of
  degree $d$ is of weight $4d$ in the last algebra. For instance, $J_2=2\,{
    a_0}\,{ a_8}- { a_1}\,{ a_7}/4+{ a_2}\,{ a_6}/14 - { a_3}\,{ a_5}/28+{{
      a_4}}^{2}/70$ is of degree 2 and weight 8.
\end{remark}

\begin{example}
  Consider the covariants of order 2 of smallest degree in
  Table~\ref{tab:fundamentals}, that is $q_1=C_{5, 2}$, $q_2=C_{6, 2}$ and
  $q_3=C_{7,2}$. A first call to Algorithm~\ref{RecognizeJInvariants} yields
  that $R=R(q_1,q_2,q_3)$ is equal, up to a constant, to
  \begin{dgroup*}[compact,style={\tiny},spread={-6pt}]
    \begin{dmath*}
      R=-4937630140800\,{{J_9}}^{2}+6172588800000\,{J_8}\,{J_{10}}+1016336160000\,{{J_6}}^{3}-1646487542700\,{J_5}\,{J_6}\,{J_7}
      +475344450\,{{J_5}}^{2}{J_8}-13778100\,{J_4}\,{{J_7}}^{2}+6154254741600\,{J_4}\,{J_6}\,{J_8}+2469123699840\,{J_4}\,{J_5}\,{J_9}
      -3175414824960\,{{J_4}}^{2}{J_{10}}-1028718873000\,{J_3}\,{J_7}\,{J_8}-1555231104000\,{J_3}\,{J_6}\,{J_9}+514676332800\,{J_3}\,{J_5}\,{J_{10}}
      -579162433500\,{J_2}\,{{J_8}}^{2}+231655788000\,{J_2}\,{J_7}\,{J_9}+47632860\,{{J_4}}^{2}{{J_5}}^{2}
      -201602675520\,{{J_4}}^{3}{J_6}-264617457390\,{J_3}\,{J_4}\,{J_5}\,{J_6}+529262244990\,{J_3}\,{{J_4}}^{2}{J_7}
      +4618063800\,{{J_3}}^{2}{{J_6}}^{2}-35210700\,{{J_3}}^{2}{J_5}\,{J_7}-228766979700\,{{J_3}}^{2}{J_4}\,{J_8}
      +38124172800\,{{J_3}}^{3}{J_9}+77149935135\,{J_2}\,{{J_5}}^{2}{J_6}-40603006080\,{J_2}\,{J_4}\,{{J_6}}^{2}-115812049185\,{J_2}\,{J_4}\,{J_5}\,{J_7}-330859026540\,{J_2}\,{{J_4}}^{2}{J_8}+145802916000\,{J_2}\,{J_3}\,{J_6}\,{J_7}-15715198800\,{J_2}\,{J_3}\,{J_4}\,{J_9}+42877447200\,{J_2}\,{{J_3}}^{2}{J_{10}}+53596043550\,{{J_2}}^{2}{{J_7}}^{2}-145802916000\,{{J_2}}^{2}{J_6}\,{J_8}-53606606760\,{{J_2}}^{2}{J_5}\,{J_9}
      +137217628800\,{{J_2}}^{2}{J_4}\,{J_{10}}-36737464140\,{{J_3}}^{2}{{J_4}}^{3}-7824600\,{{J_3}}^{3}{J_4}\,{J_5}
      +11300902200\,{{J_3}}^{4}{J_6}-47249726760\,{J_2}\,{{J_4}}^{4}-12161979900\,{J_2}\,{J_3}\,{{J_4}}^{2}{J_5}+33446455740\,{J_2}\,{{J_3}}^{2}{J_4}\,{J_6}+1760535\,{{J_2}}^{2}{J_4}\,{{J_5}}^{2}+25514097660\,{{J_2}}^{2}{{J_4}}^{2}{J_6}-153935460\,{{J_2}}^{3}{{J_6}}^{2}
      +1173690\,{{J_2}}^{3}{J_5}\,{J_7}+7625565990\,{{J_2}}^{3}{J_4}\,{J_8}-1270805760\,{{J_2}}^{3}{J_3}\,{J_9}
      -1429248240\,{{J_2}}^{4}{J_{10}}
      +289800\,{J_2}\,{{J_3}}^{4}{J_4}+900887400\,{{J_2}}^{2}{{J_3}}^{2}{{J_4}}^{2}+2575261188\,{{J_2}}^{3}{{J_4}}^{3}+260820\,{{J_2}}^{3}{J_3}\,{J_4}\,{J_5}
      -753393480\,{{J_2}}^{3}{{J_3}}^{2}{J_6}-1114881858\,{{J_2}}^{4}{J_4}\,{J_6}-19320\,{{J_2}}^{4}{{J_3}}^{2}{J_4}
      -30029580\,{{J_2}}^{5}{{J_4}}^{2}+12556558\,{{J_2}}^{6}{J_6}+322\,{{J_2}}^{7}{J_4}\,.
    \end{dmath*}
  \end{dgroup*}
  Similarly, we find that the equation of the conic, $\sum_{i,j} A_{ij} \, x_i
  x_j=0$, is equal to
  \begin{dgroup*}[compact,style={\tiny},spread={-6pt}]
    \begin{dmath*}
      ( 9217732608000\,{ J_{10}}-1422489600\,{{ J_2}}^{3}{ J_4}+
      1814283878400\,{ J_4}\,{ J_6}-384072192000\,{ J_3}\,{ J_7}\\+
      42674688000\,{{ J_3}}^{2}{ J_4}-1152216576000\,{{ J_5}}^{2}+
      212154163200\,{ J_2}\,{{ J_4}}^{2}+384072192000\,{ J_2}\,{ J_8 } )\ {{
          x_1}}^{2}\,+
    \end{dmath*}
    \begin{dmath*}
      ( -80015040000\,{{ J_3}}^{2}{ J_5}+ 2667168000\,{{ J_2}}^{3}{
        J_5}-12002256000\,{{ J_2}}^{2}{ J_7} +288054144000\,{ J_2}\,{
        J_9}+216040608000\,{ J_4}\,{ J_7}\\- 102019176000\,{ J_2}\,{ J_4}\,{
        J_5}+138883248000\,{ J_3}\,{{ J_4}}^{2}-48009024000\,{ J_5}\,{
        J_6}+360067680000\,{ J_3}\,{ J_8} )\ { x_1}\,{ x_2}\,+
    \end{dmath*}
    \begin{dmath*}
      ( -12040358400\,{{ J_2}}^{2 }{ J_8}-902039040\,{{ J_2}}^{3}{
        J_6}-24768737280\,{{ J_4}}^{3 }+27061171200\,{{ J_3}}^{2}{
        J_6}+18627840\,{{ J_2}}^{4}{ J_4} -424308326400\,{ J_2}\,{
        J_{10}}\\-5482391040\,{{ J_2}}^{2}{{ J_4} }^{2} -43481733120\,{
        J_2}\,{ J_4}\,{ J_6}+216726451200\,{ J_5} \,{ J_7}+12040358400\,{
        J_2}\,{ J_3}\,{ J_7}-762657638400\,{ J_4}\,{ J_8}\\+36121075200\,{
        J_2}\,{{ J_5}}^{2}+135339724800\, { J_3}\,{ J_9}-162570240000\,{{
          J_6}}^{2}-10516262400\,{ J_3} \,{ J_4}\,{ J_5}-558835200\,{ J_2}\,{{
          J_3}}^{2}{ J_4} )\ { x_1}\,{ x_3}\,+
    \end{dmath*}
    \begin{dmath*}
      ( 135025380000\,{ J_3}\,{ J_9}+ 55566000\,{{ J_2}}^{3}{
        J_6}-15788682000\,{ J_2}\,{ J_4}\,{ J_6}+2813028750\,{{ J_2}}^{2}{
        J_8}-2813028750\,{ J_2}\,{ J_3}\,{ J_7}+149333625\,{{ J_2}}^{4}{
        J_4}\\+8439086250\,{ J_3}\, { J_4}\,{ J_5}-151903552500\,{ J_4}\,{
        J_8}-2509400250\,{{ J_2}}^{2}{{ J_4}}^{2}+75951776250\,{ J_5}\,{
        J_7}-1666980000\,{{ J_3}}^{2}{ J_6}-4480008750\,{ J_2}\,{{ J_3}}^{2}{
        J_4}\\+ 92610000\,{{ J_2}}^{3}{{ J_3}}^{2}-1543500\,{{ J_2}}^{6}-
      1389150000\,{{ J_3}}^{4}-2893401000\,{{ J_4}}^{3}-50009400000\,{{
          J_6}}^{2}-67512690000\,{ J_2}\,{ J_{10}} )\ {{ x_2}}^{2}+
    \end{dmath*}
    \begin{dmath*}
      ( 1434793500\,{{ J_2}}^{2}{ J_4}\,{ J_5}-1629217800\,{ J_2}\,{ J_3}\,{{
          J_4}}^{2}+6460738200\,{ J_2}\,{ J_5}\,{ J_6} +365148000\,{{
          J_3}}^{3}{ J_4}-41806800\,{{ J_2}}^{4}{ J_5}- 12748654800\,{ J_3}\,{
        J_4}\,{ J_6}\\+1254204000\,{ J_2}\,{{ J_3}}^{2}{ J_5}+914457600\,{
        J_2}\,{ J_4}\,{ J_7}-12171600\,{{ J_2}}^{3}{ J_3}\,{
        J_4}-172254600\,{{ J_2}}^{3}{ J_7}- 2400451200\,{{ J_2}}^{2}{
        J_9}-714420000\,{ J_6}\,{ J_7}\\- 5643918000\,{ J_2}\,{ J_3}\,{
        J_8}-4445733600\,{{ J_4}}^{2}{ J_5}+14402707200\,{ J_4}\,{
        J_9}+10811556000\,{{ J_3}}^{2}{ J_7}-63440496000\,{ J_3}\,{
        J_{10}}-44365482000\,{ J_5}\,{ J_8} )\ { x_2}\,{ x_3}+
    \end{dmath*}
    \begin{dmath*}
      ( 94363920\,{{ J_2}}^{3}{ J_8 }+2592705024\,{{ J_4}}^{2}{
        J_6}-32568480\,{{ J_3}}^{2}{{ J_4} }^{2}+57512\,{{ J_2}}^{5}{
        J_4}-283091760\,{{ J_2}}^{2}{{ J_5} }^{2}+4386130560\,{{ J_2}}^{2}{
        J_{10}}-1905120000\,{ J_6}\,{ J_8}\\-40824000\,{{
          J_7}}^{2}+34895088\,{{ J_2}}^{3}{{ J_4}}^{2}+ 1886976000\,{ J_2}\,{{
          J_6}}^{2}-10150479360\,{ J_5}\,{ J_9}- 109801152\,{{ J_2}}^{2}{
        J_4}\,{ J_6}+23227223040\,{ J_4}\,{ J_{10}}\\+21819168\,{{ J_2}}^{4}{
        J_6}+3110425920\,{ J_3}\,{ J_5}\,{ J_6}+15630965280\,{ J_2}\,{ J_4}\,{
        J_8}+164838240\,{ J_2}\,{ J_3}\,{ J_4}\,{ J_5}+635065920\,{ J_2}\,{{
          J_4}}^ {3}-3676609440\,{ J_3}\,{ J_4}\,{ J_7}\\-1725360\,{{
          J_2}}^{2}{{ J_3}}^{2}{ J_4}-3397101120\,{ J_2}\,{ J_5}\,{ J_7}-
      2121396480\,{ J_2}\,{ J_3}\,{ J_9}-94363920\,{{ J_2}}^{2}{ J_3}\,{
        J_7}-654575040\,{ J_2}\,{{ J_3}}^{2}{ J_6})\ {{ x_3}}^{2} = 0\,.
    \end{dmath*}
  \end{dgroup*}
  The beginning of the quartic
  \begin{math}
    \sum_{i,j,k,l} h_{i,j,k,l} x_i x_j x_k x_l
  \end{math}
  (its coefficients are too large to be all written here) is then
  \begin{dgroup*}[compact,style={\tiny},spread={-6pt}]
    \begin{dmath*}
      (20832487200\,{{ J_7}}^{3}-98761420800\,{ J_6}\,{ J_7}\,{ J_8}-
      14814213120\,{{ J_6}}^{2}{ J_9}+140619288600\,{ J_5}\,{{ J_8}}
      ^{2}+21526903440\,{ J_5}\,{ J_7}\,{ J_9}+192584770560\,{ J_4} \,{
        J_8}\,{ J_9}-29628426240\,{ J_3}\,{{ J_9}}^{2}+ 6351593875200000\,{
        J_2}\,{ J_9}\,{ J_{10}}-231472080\,{{ J_5}}^ {3}{ J_6}+17310682368\,{
        J_4}\,{ J_5}\,{{ J_6}}^{2}- 24651776520\,{ J_4}\,{{ J_5}}^{2}{ J_7}
    \end{dmath*}
    \begin{dmath*}
      \cdots
    \end{dmath*}
    \begin{dmath*}
      +653457959280\,{{ J_2}}^{5}{ J_5}\,{ J_6}-653460684660\,{{ J_2}}^{5}{
        J_4}\,{ J_7}+108909756900\,{{ J_2}}^{6}{ J_9}+47040\,{{ J_2}}^{4}{{
          J_3 }}^{3}{ J_4}+56723695560\,{{ J_2}}^{5}{ J_3}\,{{ J_4}}^{2}+
      141120\,{{ J_2}}^{5}{{ J_3}}^{2}{ J_5}+222264\,{{ J_2}}^{6}{ J_4}\,{
        J_5}+117600\,{{ J_2}}^{6}{ J_3}\,{ J_6}+7056\,{{ J_2}}^{7}{
        J_7}-784\,{{ J_2}}^{7}{ J_3}\,{ J_4}-2352\,{{ J_2}}^{8}{ J_5}) \ {{
          x_1}}^{4}\,+\cdots
    \end{dmath*}
  \end{dgroup*}

  These precomputations done, let us look now for an octic $f$ defined over
  $\FF_{11}$ such that for instance $J_2=1$, $J_3=J_4=J_5=J_6=J_7=0$, $J_8=8$,
  $J_9=2$ and $J_{10}=7$. We first check that $R\neq 0$ and that the conic
  equation is equal to
  \begin{displaymath}
    x_1\,x_2 + 3\,x_1\,x_3 + 6\,x_2^{2} + x_2\,x_3 + 8\,x_3^{2} = 0\,.
  \end{displaymath}
  Then, since the point $(1:0:1)$ is on this conic, we have the
  parametrization
  \begin{displaymath}
    (X, Z) \mapsto ( 8\,X^2 + 10\,X\,Z + 6\,Z^2: 8\,X\,Z + 9\,Z^2: X\,Z +
    6\,Z^2)\,.
  \end{displaymath}
  In this case, the quartic equation of $\hh$ is equal to
  \begin{multline*}
    6\,x_1^{4} + 5\,x_1^{3}\,x_2 + 9\,x_1^{3}\,x_3 + 5\,x_1^{*2}\,x_2^{*2} +
    x_1^{2}\,x_2\,x_3 + 7\,x_1^{2}\,x_3^{2} + 8\,x_1\,x_2^{3}\\ +
    10\,x_1\,x_2^{2}\,x_3 + 3\,x_1\,x_2\,x_3^{2} + 3\,x_1\,x_3^{3} +
    7\,x_2^{4} + 7\,x_2^{3}\,x_3 + 9\,x_2\,x_3^{3} + 5\,x_3^{4}
  \end{multline*}
  and we finally find that, up to a constant,
  \begin{displaymath}
    f(X,Z) = 2\,X^8 + 7\,X^7\,Z^1 + 9\,X^6\,Z^2 + 9\,X^5\,Z^3 + 8\,X^4\,Z^4 + 3\,X^3\,Z^5
    + 2\,X^2\,Z^6 + 4\,X\,Z^7 + 8 \,Z^8\,.
  \end{displaymath}
\end{example}

\section{Automorphisms and strata of hyperelliptic genus $3$ curves} \label{sec:field}
In the sequel, $K$ is an algebraically closed  field of characteristic $p$ where $p$ is a prime or $0$.
\subsection{Review on automorphism groups}
Let $C$ be a  hyperelliptic curve of genus $g \geq 2$ over  $K$ and $\iota$ be its hyperelliptic involution. We denote $\Aut(C)$ the automorphism group of $C$ (by convention if $C$ is defined over an arbitrary field $k$ then $\Aut(C)$ is $\Aut(C/\bar{k})$). Since $\iota$ is in the center of $\Aut(C)$ we can defined the \emph{reduced automorphism group} $\overline{\Aut}(C)=\Aut(C)/\langle \iota \rangle$. We say that $C$ has \emph{extra-automorphism} when $\# \overline{\Aut(C)}>1$. The reduced automorphism group  acts on $C/\langle \iota \rangle \simeq \PP^1$.
The list of possible finite groups $\overline{\ggm}$ acting on $\PP^1$ was given in \cite[§ 71-74]{weberbook}. In his PhD thesis, Brandt \cite{brandt} gave the full list of polynomial orbits under any $\overline{\ggm}$ which in turn gives the \emph{normal models} of hyperelliptic curves which automorphism group contains a group $\ggm$ such that $\ggm/\langle \iota \rangle=\overline{\ggm}$. The structure of $\ggm$ itself then depends on the behavior of the exact sequence
$$1 \to \langle \iota \rangle \to \ggm \to \overline{\ggm} \to 1$$
which is measured by $H^2( \overline{\ggm},\mathbb{Z}/2\mathbb{Z})$.  When $p=0$,
the structure of $\ggm$, depending on its signature, has then been worked out in
\cite{brandt-stich}. Finally among the groups $\overline{\ggm}$, one has to
determine the ones which appear as full automorphism group of $C$. When $p=0$,
this can be done using Fuchsian groups \cite{singerman}, Teichm\"uller theory
\cite{ries} or Hurwitz spaces \cite{MSSV}.\medskip

For $g=2$, a complete list of automorphism groups and models can be found in \cite{cardona} for $p
\ne 2$ and in \cite{CANAPU} for $p=2$.\medskip

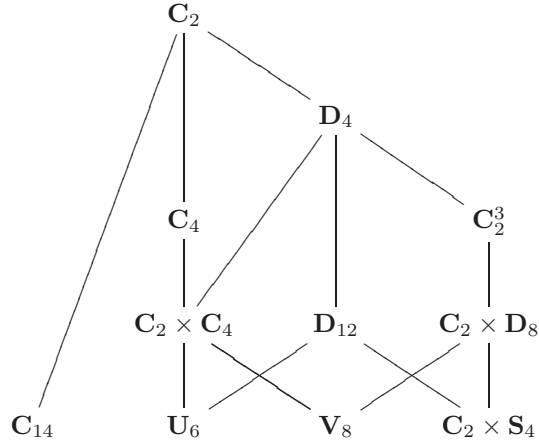
\begin{figure}[htbp]
\begin{center}
\leavevmode
\xymatrix{ & \ar@{-}[ddddl] \ar@{-}[dd] \CG_2 \ar@{-}[rd]  &              & \\
                                                           &              & \DG_4  \ar@{-}[dd] \ar@{-}[ddl]  \ar@{-}[dr] & \\
                                                           &  \CG_4 \ar@{-}[d] & &     \CG_2^3 \ar@{-}[d] \\
                                                   & \CG_2 \times \CG_4 \ar@{-}[d]  \ar@{-}[dr] \ar@{-}[dr]  &  \DG_{12} \ar@{-}[dl]  \ar@{-}[dr] &  \CG_2 \times \DG_8   \ar@{-}[dl] \ar@{-}[d] \\ 
\CG_{14}   &    \UG_6               &   \VG_8       & \CG_2 \times \SG_4 \\}
\end{center}
\caption{Organization of the strata according to inclusion and dimension (see Remark \ref{rem:error})} \label{fig:lattice}
\end{figure}

 For $g=3$ and $p=0$, this work has been achieved explicitly  in several
 papers. We refer to  \cite{MSSV} (or \cite{GSS}) for a historical viewpoint
 and we gather their results  in Table~\ref{tab:auto} (see also Remark
 \ref{rem:error}). Using the signature \cite[Sec.4]{MSSV} (or  the shape of
 the equations), we can deduce the relations between the stratum in the moduli
 space (\textit{cf.} Figure~\ref{fig:lattice}).\medskip

For sake of completeness, we show how these results extend to all $p$. By \cite{roquette}, we know that when $p>3+1=4$ and $p \ne 2 \cdot 3+1=7$, then $C \mapsto C/\Aut(C)$ is tamely ramified and so by \cite[XIII.2.12]{grothendieck}, it can be lifted to characteristic $0$. Hence, Table~\ref{tab:auto} is also valid in these characteristics.
When $p=2$, the possible automorphism groups and models are in
\cite{NS04}. When $p=7$, by \cite{roquette} the curve  $C : y^2=x^7-x$ which
has a group of order $2^5 \cdot 3 \cdot 7$ is the only exceptional
case. Finally for $p=3$, going through the list of \cite[Satz.2.3]{brandt}, it
seems that there is no new automorphism group and moreover the cases for which $3$
divides $\#\Aut(C)$ in Table~\ref{tab:auto} do not appear anymore. 

\begin{remark} \label{rem:error}
Some remarks on the notation and convention for Table~\ref{tab:auto} and Figure~\ref{fig:lattice}.
\begin{itemize}
\item We have the following notation for the groups:
  \begin{itemize}
  \item $\CG_n=\ZZ/n\ZZ$;
  \item $\DG_n$ is the dihedral group with $n$ elements;
  \item $\UG_6$ is a group with $24$ elements defined by $\langle
    S,T\rangle$ with $S^{12}=T^2=1$ and $TST=S^5$;
  \item $\VG_8$ is a group with $32$ elements defined by $V^8=\langle
    S,T\rangle$ with $S^4=T^8=(ST)^2=(S^{-1}T)^2=1$;
  \item $\SG_n$ is the symmetric group over $n$ letters.
  \end{itemize}
\item 
  An exponent in the signature must be understood as repetition: for instance
  $(2^6)=(2,2,2,2,2,2)$.  The dimension $\delta$ of the stratum in $\Hm_3$ is easily
  computed as $-3+\# S$ where $S$ is the signature.
  The column Id. refers to the  GAP or MAGMA library of small groups.
\item 
  The order chosen for Table \ref{tab:auto} (resp. for Figure \ref{fig:lattice}) is by decreasing dimension $\delta$  (resp.  from top to bottom) then by increasing order of the automorphism group (resp. from left to right).
\item 
  The equation of the normal model is valid for  the stratum: for some special values of the parameters, the curve can have more automorphisms.
\item 
  In both references  \cite{MSSV} and \cite{GSS}, Case 11 of Table \ref{tab:auto} is wrongly written as $x^8+14x^2+1$.
\item 
  In \cite[Fig.1]{GSS} and in \cite[Tab.3]{babu}, the organization of the
  strata is wrong.
\end{itemize}
\end{remark}

\begin{table}[htbp]
{\small
\begin{tabular}{|c|c|c|c|c|c|c|}
\hline
\# & $\Aut(C)$ & $\overline{\Aut}(C)$ & signature & $\delta$ & equation $y^2=f(x)$ & Id. \\
\hline
\hline
1 & $\CG_2$ & $\{1\}$ & $(2^8)$ & $5$ & $x(x-1)(x^5+ax^4+bx^3+cx^2+dx+e)$ & $(2,1)$ \\
 2 & $\DG_4$ & $\CG_2$ &  $(2^6)$ &$3$ & $\begin{cases} x^8+a x^6+ b x^4+ c x^2+1 & \textrm{or} \\ (x^2-1)(x^6+ax^4+bx^2+c) \end{cases}$ & $(4,2)$ \\
3 & $\CG_4$ & $\CG_2$ &   $(2^3,4^2)$&   $2$ & $x(x^2-1)(x^4+ax^2+b)$ & $(4,1)$ \\
4 & $\CG_2^3$ & $\DG_4$ & $(2^5)$ & $2$ & $(x^4+ax^2+1)(x^4+bx^2+1)$ & $(8,5)$ \\
5 & $\CG_2 \times \CG_4$ & $\DG_4$ & $(2^2,4^2)$  & $1$  & $\begin{cases} (x^4-1)(x^4+ax^2+1) & \textrm{or} \\ x(x^2-1)(x^4+ax^2+1) \end{cases}$ &  $(8,2)$ \\
6 & $\DG_{12}$  & $\DG_6$ & $(2^3,6)$ & $1$ & $x(x^6+ax^3+1)$& $(12,4)$ \\
7 & $\CG_2 \times \DG_8$ & $\DG_8$ &  $(2^3,4)$ &   $1$  & $x^8+ax^4+1$ & $(16,11)$ \\
8 & $\CG_{14}$ & $\CG_7$ & $(2,7,14)$ & $0$ & $x^7-1$ & $(14,2)$ \\
9 & $\UG_6$ & $\DG_{12}$ & $(2,4,12)$&  $0$ & $x(x^6-1)$ & $(24,5)$ \\
10 & $\VG_8$ & $\DG_{16}$ & $(2,4,8)$& $0$ & $x^8-1$ & $(32,9)$ \\
11 & $\CG_2 \times \SG_4$ &  $\SG_4$ &  $(2,4,6)$ & $0$ & $x^8+14x^4+1$ & $(48,48)$ \\
\hline
\end{tabular}\medskip
}

\caption{Automorphism groups for genus $3$ hyperelliptic curves in characteristic $0$ (see Remark \ref{rem:error})}
\label{tab:auto}
\end{table}

\subsection{Identification of loci and reconstruction with non trivial
  automorphism groups}
\label{sec:ident-loci-reconstr}

Let $C$ be a hyperelliptic curve of genus $3$ defined over a field $k$ of characteristic $p \ne 2,3,5,7$ and let $K=\bar{k}$.
Let $(j_2:j_3:\ldots:j_{10})$ be the Shioda invariants of $C$, our aim is to
retrieve from them a hyperelliptic model $y^2 = f(x)$.
 
The conic and quartic method presented in Section~\ref{sec:conic-quartic},
based on Prop.~\ref{prop:generic-reconstructions} is our main
tool for solving this problem, but unfortunately Lemma~\ref{lemma:reconstr-D4}
shows that this method can not work for most of the curves with a non trivial
automorphism group.
\begin{lemma}\label{lemma:reconstr-D4}
  When the automorphism group of an hyperelliptic curve $C:\ y^2 = f(x)$ contains $\DG_4$, \ie in cases 2,4,5,6,7,9,10,11 of
  Table~\ref{tab:auto}, then for any choice of $3$ quadratic covariants
  $q_i=q_i(f,(x,z))$ of $f=x^8 +a x^6+b x^4+ c x^2+1$, we have that
  $R(q_1,q_2,q_3)=0$.
\end{lemma}
\begin{proof}
  Let $d_i$ be the degree of the covariants $q_i$. Their weight is then $(8
  d_i-2)/2=4d_i-1$ which is always odd. Hence, for $g \in \Aut(C)$ 
  acting as $g.(x,z)=(-x,z)$,  
  $$q_i(g.f,g.(x,z))=-q_i(f,(x,z))= q_i(f,(-x,z)).$$
  Hence for all $i$, $q_i$ has only an $xz$ term. In particular the
  determinant of the $q_i$s in the bases $x^2$, $xz$, $z^2$ is zero.
\end{proof} 

We thus have to develop more specific methods for reconstructing models for
curves in this case.
Given invariants $(j_2:j_3:\ldots:j_{10})$, a prerequisite is to determine
what is the automorphism group of the corresponding curve so that we can
select the right normal model to reconstruct. To this aim, we determine
for each automorphism group equations for the corresponding stratum in
$\Hm_3$. There are very few cases where determining some of these equations is
straightforward, for instance automorphism groups which contains $\CG_4$.
\begin{lemma}\label{lemma:loci-C4}
  When the automorphism group of an hyperelliptic curve $C:\ y^2 = f(x)$
  contains $\CG_4$, \ie in cases 3,5,9 and 10, then $J_i(f)=0$ for all odd index
  $i$.
\end{lemma}
\begin{proof}
  The weight of an invariant of degree $d$ is $8d/2=4d$. If $\Aut(C)$ contains
  $\CG_4$, the curve $C$ admits a model of the form $y^2=f(x)$ with
  $f(x)=x(x^2-1)(x^4+ax^2+1)$. For $g \in \Aut(C)$ acting as $g.(x,z)=(-x,z)$, we get that
  $$J_i(g.f)=\det(g)^{4i} J_i(f)=J_i(f).$$ 
  On the other hand, $g.f=-f$, so if the degree of $J_i$ is odd, then $J_i(g.f)=J_i(-f)=- J_i(f)$.
  Hence we get that $J_i(f)=-J_i(f)$.
\end{proof}

In order to exhibit necessary conditions for all the strata, we have applied
Algorithm~\ref{RecognizeJInvariants}. We choose for $I$ the constant invariant
$0$, seen as an invariant of increasing positive degree $d$, and choose for
$\mathcal F$ (at line $2$) a set of random octics over $\QQ$ of the form the
normal model for the automorphism group that we consider. Increasing one by
one $d$ from $0$ to $30$ yields generators for the ideal of relations which
defines the stratum. The shape of these generators is very close to the one of
a Gröbner basis for the graded reverse lexicographical (or `grevlex') order
$J_{2}< J_3< \ldots <J_{10}$ with weights 2, 3, \ldots, 10. For this order, it
is thus possible to deduce a reduced Gröbner basis. In the easier cases, note
that it is feasible to apply a `change of order algorithm' and to deduce a
reduced Gröbner basis for the lexical order $J_2<J_3<\ldots<J_{10}$ too.

Now, in order to exhibit a model from given invariants, we proceed as above, except
that we slightly modify  Algorithm~\ref{RecognizeJInvariants} to add in the
basis $\mathcal B$ (at line 1) the coefficients $a$, $b$, \textit{etc.} of the
normal model in Tab.~\ref{tab:auto} for the considered stratum. We still
choose for $\mathcal F$ (at line $2$) a set of random octics in the shape of the
normal model. The lowest degree equations found in this way are then enough to
reconstruct a model.\smallskip

All in all, we give below reconstruction lemmas, one for each automorphism
group in Table~\ref{tab:auto}. We precisely state the stratum equations
that we have obtained and a model in terms of the $j_i$s for the curve $C$ which may be defined over a non trivial extension of $k$. 
In order to make this extension as small as possible, we introduce models which have more non-zero coefficients than the ones in Table \ref{tab:auto}. We refer the reader to Section~\ref{sec:field-definition} for the existence of a model over the field of moduli.

The proofs of these lemmas do not depend on the way we have obtained the
equations. They follow essentially all the same principle.
\begin{itemize}
\item We check that the normal models of Tab.~\ref{tab:auto} have Shioda
  invariants that satisfy the stratum equations, so that we are convinced that
  these models are a subset of all the models which satisfy the stratum
  equations.
\item Conversely, we check that the reconstructed model  has
  Shioda invariants in the same weighted projective class as the $9$-uple
  $(j_2 : j_3 : \ldots : j_{10})$ provided in input. Since these models are of
  normal form, we have proved that only normal models satisfy the stratum
  equations. If we can perform this step, we then have checked that the
  equations we found describe the stratum. 
\end{itemize}

Most of these proofs need heavy computations, far too complex to be written down
here, and so we must skip them. But a program written in the \textsc{magma}
computational algebra system that implements the corresponding computations is
available on the web page of the authors for independent checks.\medskip

Incidentally, we succeed in parameterizing the projective variety defined by
the stratum equations of dimension $\leq 2$ (it is often a pencil of rational curves). As a
first consequence, we give parametrized models over the field of moduli for
all these strata, except the cases $\CG_2^3$ where we have to deal
with algebraic extensions. As a second consequence, we give the exact number
of isomorphic classes of curves  over a finite
field $\FF_q$ in these strata.\medskip

\begin{figure}[htbp]
  \begin{center}
  \parbox{0.95\linewidth}{%
    \begin{footnotesize}\SetAlFnt{\small\sf}%
      \begin{algorithm}[H]%
        \label{algo:reconstruct}%
        \caption{Reconstruct a hyperelliptic polynomial from its Shioda invariants}  %
        \SetKwInOut{Input}{Input} \SetKwInOut{Output}{Output} %
        \Input{Shioda invariants $(j_2:j_3:\ldots:j_{10})$.}  %
        \Output{A hyperelliptic polynomial $f$}  \BlankLine
        %
        %
        If $(j_2:j_3:\ldots:j_{10})$ satisfies Eq.~\eqref{eq:2}, then
        reconstruct $f$  with Lemma~\ref{lemma:C2S4}\hfill (case $\CG_2\times \SG_4$)\;
        If $(j_2:j_3:\ldots:j_{10})$ satisfies Eq.~\eqref{eq:3}, then
        reconstruct $f$  with Lemma~\ref{lemma:V8}\hfill (case $\VG_8$)\;
        If $(j_2:j_3:\ldots:j_{10})$ satisfies Eq.~\eqref{eq:4}, then
        reconstruct $f$  with Lemma~\ref{lemma:U6}\hfill (case $\UG_6$)\;
        If $(j_2:j_3:\ldots:j_{10})$ satisfies Eq.~\eqref{eq:5}, then
        reconstruct $f$  with Lemma~\ref{lemma:C14}\hfill (case $\CG_{14}$)\;
        If $(j_2:j_3:\ldots:j_{10})$ satisfies Eq.~\eqref{eq:6}, then
        reconstruct $f$  with Lemma~\ref{lemma:C2D8}\hfill (case $\CG_2\times \DG_8$)\;
        If $(j_2:j_3:\ldots:j_{10})$ satisfies Eq.~\eqref{eq:7}, then
        reconstruct $f$  with Lemma~\ref{lemma:D12}\hfill (case $\DG_{12}$)\;
        If $(j_2:j_3:\ldots:j_{10})$ satisfies Eq.~\eqref{eq:8}, then
        reconstruct $f$  with Lemma~\ref{lemma:C2C4}\hfill (case $\CG_2\times \CG_4$)\;
        If $(j_2:j_3:\ldots:j_{10})$ satisfies Eq.~\eqref{eq:9}, then
        reconstruct $f$  with Lemma~\ref{lemma:C2p3}\hfill (case $\CG_2^3$)\;
        If $(j_2:j_3:\ldots:j_{10})$ satisfies Eq.~\eqref{eq:10}, then
        reconstruct $f$  with Lemma~\ref{lemma:C4}\hfill (case $\CG_4$)\;
        If $(j_2:j_3:\ldots:j_{10})$ satisfies Eq.~\eqref{eq:11}, then
        reconstruct $f$  with Lemma~\ref{lemma:D4}\hfill (case $\DG_4$)\;
        Otherwise, reconstruct $f$ with Lemma~\ref{lemma:C2}\hfill (case
        $\CG_2$)\;
        \KwRet{$f$};
      \end{algorithm}
    \end{footnotesize}
  }
\end{center}
\end{figure}

Algorithm~\ref{algo:reconstruct} summarizes how one can apply the lemmas below
to reconstruct a model in any case. Actually, this algorithm returns more generally a
binary form $f$ of degree 8, even in the cases where $f$ does not have simple
roots. But, when the order of one of these roots is greater than 4, only one
among all the possible orbits, is returned. Typically only $f=x^8$ is returned
for $j_2=j_3=\ldots=j_{10}=0$. Moreover the automorphism group of such a form may in general be bigger than the one indicated in the algorithm.

We state now, by increasing dimension, the reconstruction lemmas.

\subsubsection{Strata of dimension 0}
\label{sec:strata-dimension-0}

\stratum{ $\Aut(C) = \CG_2 \times \SG_4$}

\begin{lemma}\label{lemma:C2S4}
  Let $C$ be a hyperelliptic curve of genus 3 over $k$, let $(j_2 : j_3:\ldots:
  j_{10})$ be its Shioda invariants. Then the automorphism group of $C$ is
  $\CG_2 \times \SG_4$ if and only if
  {\rm\setcounter{myequation}{\theequation}
    \begin{dgroup*}[style={\small}]
      \begin{dmath}[style={\refstepcounter{myequation}},number={\themyequation},label={eq:2}]
        0 = 30\,{{ j_3}}^{2}-{{ j_2}}^{3},
      \end{dmath}
      \begin{dmath*}
        0 = j_4 \hiderel{=} j_5 \hiderel{=} \cdots \hiderel{=} j_{10}\,.
      \end{dmath*}
    \end{dgroup*}
    \setcounter{equation}{\themyequation}
  }
  \noindent
  Furthermore, a curve $C$ with automorphism group $\CG_2 \times \SG_4$ is
  $K$-isomorphic to the curve $y^2 = x^8 + 14\,x^4 + 1$\,.
\end{lemma}

\stratum{ $\Aut(C) = \VG_8$}

\begin{lemma}\label{lemma:V8}
  Let $C$ be a hyperelliptic curve of genus 3 over $k$, let $(j_2: j_3:\ldots:
  j_{10})$ be its Shioda invariants. Then the automorphism group of $C$ is $\VG_8$
  if and only
  {\rm\setcounter{myequation}{\theequation}
    \begin{dgroup*}[style={\small}]
      \begin{dmath*}
        0 = 6\,{ j_4}-{{ j_2}}^{2}\,,
      \end{dmath*}
      \begin{dmath*}
        0 = 36\,{ j_6}+{{ j_2}}^{3}\,,
      \end{dmath*}
      \begin{dmath}[style={\refstepcounter{myequation}},number={\themyequation},label={eq:3}]
        0 = 420\,{ j_8}+{{ j_2}}^{4}\,,
      \end{dmath}
      \begin{dmath*}
        0 = 2520\,{ j_{10}}-{{ j_2}}^{5}\,,
      \end{dmath*}
      \begin{dmath*}
        0 = { j_3} \hiderel{=} { j_5} \hiderel{=} { j_7} \hiderel{=} { j_9}\,.
      \end{dmath*}
    \end{dgroup*}
    \setcounter{equation}{\themyequation}
  }
  \noindent
  Furthermore, a curve $C$ with automorphism group $\VG_8$ is $K$-isomorphic
  to the curve $y^2 = x^8 - 1$\,.
\end{lemma}

\stratum{ $\Aut(C) = \UG_6$}

\begin{lemma}\label{lemma:U6}
  Let $C$ be a hyperelliptic curve of genus 3 over $k$, let $(j_2: j_3:\ldots:
  j_{10})$ be its Shioda invariants. Then the automorphism group of $C$ is
  $\UG_6$ if and only if
  {\rm\setcounter{myequation}{\theequation} 
    \begin{dgroup*}[style={\small}]
      \begin{dmath*}
        0 = 96\,{ j_4}-{{ j_2}}^{2}\,,
      \end{dmath*}
      \begin{dmath*}
        0 = 2304\,{ j_6}+{{ j_2}}^{3}\,,
      \end{dmath*}
      \begin{dmath}[style={\refstepcounter{myequation}},number={\themyequation},label={eq:4}]
        0 = 17920\,{ j_8}-{{ j_2}}^{4}\,,
      \end{dmath}
      \begin{dmath*}
        0 = 430080\,{ j_{10}}+{{ j_2}}^{5}\,,
      \end{dmath*}
      \begin{dmath*}
        0 = { j_3} \hiderel{=} { j_5} \hiderel{=} { j_7} \hiderel{=} { j_9}\,.
      \end{dmath*}
    \end{dgroup*}
    \setcounter{equation}{\themyequation}
  }
  \noindent
  Furthermore, a curve $C$ with automorphism group $\UG_6$ is $K$-isomorphic
  to the curve $y^2 = x^7 - x$\,.
\end{lemma}

\stratum{ $\Aut(C) = \CG_{14}$}

\begin{lemma}\label{lemma:C14}
  Let $C$ be a hyperelliptic curve of genus 3 over $k$, let $(j_2: j_3:\ldots:
  j_{10})$ be its Shioda invariants. Then the automorphism group of $C$ is
  $\CG_{14}$ if and only if
  \begin{equation}
    \label{eq:5}
    \begin{array}{rcl}
      j_2\ =\ j_3\ =\ j_4\ =\ j_5\ =\ j_6\ =\ j_8\ =\ j_9\ =\ j_{10} &=&
      0\hspace*{1cm}(j_7 \neq 0)\,.
    \end{array}
  \end{equation}
  Furthermore, a curve $C$ with automorphism group $\CG_{14}$ is
  $K$-isomorphic to the curve $y^2 = x^7 - 1$\,.
\end{lemma}

\subsubsection{Strata of dimension 1}
\label{sec:strata-dimension-1}

\stratum{ $\Aut(C) = \CG_{2} \times \DG_8$}

\begin{lemma}\label{lemma:C2D8}
  Let $C$ be a hyperelliptic curve of genus 3 over $k$, let $(j_2 : j_3 :\ldots :
  j_{10})$ be its Shioda invariants. Then the automorphism group of $C$ contains
  $\CG_{2} \times \DG_8$ if and only if
  {\rm\setcounter{myequation}{\theequation}
    \begin{dgroup*}[style={\small}]
      \begin{dmath*}
        0=54\,j_4^3 - 81\,j_4^2\,j_2^2 - 1080\,j_4\,j_3^2\,j_2 + 36\,j_4\,j_2^4
        - 3600\,j_3^4 + 240\,j_3^2\,j_2^3 - 4\,j_2^6\,,
      \end{dmath*}
      \begin{dmath*}
        0=60\,j_5\,j_3 + 18\,j_4^2 - 15\,j_4\,j_2^2 - 60\,j_3^2\,j_2 + 2\,j_2^4\,,
      \end{dmath*}
      \begin{dmath*}
        0=3\,j_5\,j_4 - 2\,j_5\,j_2^2 + 15\,j_4\,j_3\,j_2 + 60\,j_3^3 -
        2\,j_3\,j_2^3\,,
      \end{dmath*}
      \begin{dmath*}
        0=5\,j_5^2 - 6\,j_4^2\,j_2 - 30\,j_4\,j_3^2 + j_4\,j_2^3\,,
      \end{dmath*}
      \begin{dmath}[style={\refstepcounter{myequation}},number={\themyequation},label={eq:6}]
        0=18\,j_6 - 9\,j_4\,j_2 - 60\,j_3^2 + 2\,j_2^3\,,
      \end{dmath}
      \begin{dmath*}
        0=3\,j_7 - j_5\,j_2 + 3\,j_4\,j_3\,,
      \end{dmath*}
      \begin{dmath*}
        0=1260\,j_8 - 18\,j_4^2 - 63\,j_4\,j_2^2 - 420\,j_3^2\,j_2 + 14\,j_2^4\,,
      \end{dmath*}
      \begin{dmath*}
        0=18\,j_9 - 2\,j_5\,j_2^2 + 15\,j_4\,j_3\,j_2 + 60\,j_3^3 -
        2\,j_3\,j_2^3\,,
      \end{dmath*}
      \begin{dmath*}
        0=630\,j_{10} - 15\,j_4^2\,j_2 - 30\,j_4\,j_3^2 + j_4\,j_2^3\,.
      \end{dmath*}
    \end{dgroup*}
    \setcounter{equation}{\themyequation}
  }
  \noindent
  Furthermore, a curve $C$ with automorphism group $\CG_{2} \times \DG_8$ is
  $K$-isomorphic to the curve $y^2 = x^8 + a_4\,x^4 + a_0$ where $a_0 =
  -a_4^2/140+j_2/2$ and
  \begin{small}
    \begin{displaymath}
      a_4 = %
      \left\{\begin{array}{ll}
          {{35({ j_5}\,{ j_2}+6\,{ j_4}\,{ j_3})}/({-60
              \,{{ j_3}}^{2}+{{2 j_2}}^{3}})} & \text{ if }{-30\,{{ j_3}}^{2}+{{
                j_2}}^{3}} \neq 0\,,\\
          35 j_5\, /\, 3 j_4 & \text{ otherwise}\,.\\
        \end{array}\right.
    \end{displaymath}
  \end{small}
\end{lemma}
\medskip

\begin{remark}
  When $-30\,{{ j_3}}^{2}+{{ j_2}}^{3}$ and $j_4$ are both equal to 0,
  Eq.~\eqref{eq:6} can be reduced to Eq.~\eqref{eq:2}. This means that $C$ has
  a larger automorphism group, specifically $\CG_{2} \times \SG_4$. So
  Lemma~\ref{lemma:C2S4} can be used, instead, for reconstructing a model.
\end{remark}
\medskip

Actually the quotient field of $\ii_8$ modulo the ideal defined by
Eq.~\eqref{eq:6} is obtained by adjoining $j_4$ to $k[j_2, j_3]$ and $j_4$
satisfies an irreducible monic equation of degree $3$. 
The invariants $j_5$, $j_6$, \ldots, $j_{10}$ are then rational in $j_2$,
$j_3$ and $j_4$.

From a more geometric viewpoint, the projective variety defined by
Eq.~\eqref{eq:6} has two singularities, the $\CG_2\times \SG_4$ and $\VG_8$
points. Moreover, it is birationally equivalent to the conic
\begin{math}
  {{X_2}}^{2}-400\,{X_1}\,{X_3}+({{4000}/{3}})\,{{X_3}}^{2}
\end{math}
of discriminant
\begin{math}
  -2^8\cdot 5^4\,.
\end{math}
We can classically parameterized it by the slope of lines which
all intersect at the $\VG_8$ point. This yields
\begin{footnotesize}
  \begin{displaymath}
    t \mapsto \left(%
      \frac{30\,t+1}{30}:%
      \frac{45\,t+1}{900}:%
      \frac{t^{2}}{6}:%
      {\frac {t^{2}}{60}}:%
      \frac{-t^3}{36}:%
      {\frac {-t^3}{360}}:%
      {\frac {-{t}^{3} \left( 180\,t+7 \right)}{2^{4}\cdot3^{3}\cdot5^{2}\cdot7
}}\,:%
      {\frac {{t}^{4}}{2^{4}\cdot3^{3}\cdot5}}\,:%
      {\frac {{t}^{4} \left( 180\,t+7%
      \right)}{2^{5}\cdot3^{4}\cdot5^{2}\cdot7}}\,
    \right)\,.
  \end{displaymath}
\end{footnotesize}

The corresponding curves are of the form
\begin{displaymath}
y^2 = {x}^{8}+7/6\,{x}^{4}+1/2\,t+{{1}/{144}}\,.
\end{displaymath}
The values $t=-1/72$ and $t=2/3$ give polynomials with multiple roots, the
curve with automorphism group $\CG_2\times \SG_4$ is obtained with $t=0$,
otherwise the automorphism group is always $\CG_{2} \times \DG_8$ since the point $\VG_8$
is not reachable with this conic parametrization we chose.
Incidentally, this shows that over a finite field $\FF_{q}$ of characteristic at
least 11 there are $q-3$ isomorphic classes of curves with automorphism group
$\CG_2\times \SG_4$.

\stratum{ $\Aut(C) = \DG_{12}$}

\begin{lemma}\label{lemma:D12}
  Let $C$ be a hyperelliptic curve of genus 3 over $k$, let $(j_2: j_3:\ldots:
  j_{10})$ be its Shioda invariants. Then the automorphism group of $C$
  contains $\DG_{12}$ if and only if
  {\rm\setcounter{myequation}{\theequation} 
    \begin{dgroup*}[style={\small},spread={-7pt}]
      \begin{dmath*}
        0 = 55296\,j_4^3 - 5184\,j_4^2\,j_2^2 - 4320\,j_4\,j_3^2\,j_2 +
        144\,j_4\,j_2^4 - 900\,j_3^4 + 60\,j_3^2\,j_2^3 - j_2^6\,,
      \end{dmath*}
      \begin{dmath*}
        0 = 720\,j_5\,j_3 - 2304\,j_4^2 + 120\,j_4\,j_2^2 + 30\,j_3^2\,j_2 -
        j_2^4\,,
      \end{dmath*}
      \begin{dmath*}
        0 = 576\,j_5\,j_4 - 24\,j_5\,j_2^2 - 120\,j_4\,j_3\,j_2 - 30\,j_3^3 +
        j_3\,j_2^3\,,
      \end{dmath*}
      \begin{dmath*}
        0 = 180\,j_5^2 - 96\,j_4^2\,j_2 - 30\,j_4\,j_3^2 + j_4\,j_2^3\,,
      \end{dmath*}
      \begin{dmath}[style={\refstepcounter{myequation}},number={\themyequation},label={eq:7}]
        0 = 576\,j_6 - 72\,j_4\,j_2 - 30\,j_3^2 + j_2^3\,,
      \end{dmath}
      \begin{dmath*}
        0 = 12\,j_7 - j_5\,j_2 - 2\,j_4\,j_3\,,
      \end{dmath*}
      \begin{dmath*}
        0 = 60480\,j_8 - 14976\,j_4^2 + 504\,j_4\,j_2^2 + 210\,j_3^2\,j_2 -
        7\,j_2^4\,,
      \end{dmath*}
      \begin{dmath*}
        0 = 3456\,j_9 - 24\,j_5\,j_2^2 - 120\,j_4\,j_3\,j_2 - 30\,j_3^3 +
        j_3\,j_2^3\,,
      \end{dmath*}
      \begin{dmath*}
        0 = 30240\,j_{10} - 600\,j_4^2\,j_2 - 390\,j_4\,j_3^2 + 13\,j_4\,j_2^3\,.
      \end{dmath*}
    \end{dgroup*}
    \setcounter{equation}{\themyequation}
  }
  \noindent
  Furthermore, a curve $C$ with automorphism group $\DG_{12}$ is
  $K$-isomorphic to the curve $y^2 = x \, (x^6 + a_4\,x^3 +
  a_1)$ where $a_1 = {2}\,{{a_4}}^{2}\,/\,35-4\,{j_2}$ and
  \begin{small}
    \begin{displaymath}
      a_4 = %
      \left\{\begin{array}{ll}
          280\,({{-{j_5}\,{j_2}+4\,{j_4}\,{j_3}}})/({30\,{{j_3}}
            ^{2}-{{j_2}}^{3}})
          & \text{ if }{-30\,{{ j_3}}^{2}+{{
                j_2}}^{3}} \neq 0\,,\\
          35 j_5\, /\, 3 j_4 & \text{ otherwise}\,.\\
        \end{array}\right.
    \end{displaymath}
  \end{small}
\end{lemma}
\medskip

\begin{remark}
  When $-30\,{{ j_3}}^{2}+{{ j_2}}^{3}$ and $j_4$ are both equal to 0,
  Eq.~\eqref{eq:7} can be reduced to Eq.~\eqref{eq:2}. Again, this means that
  $C$ has a larger automorphism group, specifically $\CG_{2} \times \SG_4$.  So,
  Lemma~\ref{lemma:C2S4} can be used, instead, for reconstructing a model.
\end{remark}
\medskip

Similarly to $\CG_{2} \times \DG_8$, the quotient field of $\ii_8$ modulo the
ideal defined by Eq.~\eqref{eq:7} is obtained by adjoining $j_4$ to $k[j_2,
j_3]$ and $j_4$ satisfies an irreducible monic equation of degree $3$.  The
invariants $j_5$, $j_6$, \ldots, $j_{10}$ are then rational in $j_2$, $j_3$
and $j_4$.

More geometrically, the projective variety defined by Eq.~\eqref{eq:7} has two
singularities, the $\CG_2\times \SG_4$ and $\UG_6$ points. Moreover, it is
birationally equivalent to the conic
\begin{math}
  {{X_2}}^{2}-({{25}/{16}})\,{X_1}\,{X_3}+({{125}/{6144}})\,{{X_3}}^{2}
\end{math}
of discriminant
\begin{math}
  -5^4/ 2^8\,.
\end{math}
We can classically parameterized it by the slope of lines which all intersect
at the $\UG_6$ point, this yields
\begin{footnotesize}
  \begin{displaymath}
    t \mapsto \left(%
      \frac{360\,t+1}{30}:%
      \frac{540\,t+1}{900}:%
      \frac{3\,{t}^{2}}{2}:%
      \frac{-{t}^{2}}{10}:%
      \frac{-3\,{t}^{3}}{4}:%
      \frac{{t}^{3}}{20}:%
      {\frac {{t}^{3} \left( 4860\,t+7 \right)}{
          2^{3}\cdot3\cdot5^{2}\cdot7}}:%
      \frac{-{t}^{4}}{40}:%
      {\frac {-{t}^{4} \left( 4860\,t+7 \right)}{
          2^{4}\cdot3\cdot5^{2}\cdot7
        }}%
    \right)\,.
  \end{displaymath}
\end{footnotesize}

The corresponding curves are of the form
\begin{displaymath}
y^2 = x\,({x}^{6}-{ {7}/{9}}\,{x}^{3} - 48\,t-{{8}/{81}})\,.
\end{displaymath}
The values $t=-1/(2\cdot3^{5})$ and
$t=-1/(2^{6}\cdot3)$ give polynomials with multiple roots, the curve with automorphism group $\CG_2\times \SG_4$ is
obtained with $t=0$, otherwise the automorphism group is always $\DG_{12}$ since
the point $\UG_6$ is not reachable with this conic parametrization we chose.
Incidentally, this shows that over a finite field $\FF_{q}$ of characteristic at
least 11 there are $q-3$ isomorphic classes of curves with automorphism group
$\DG_{12}$.

\stratum{ $\Aut(C) = \CG_{2} \times \CG_4$}

\begin{lemma}\label{lemma:C2C4}
  Let $C$ be a hyperelliptic curve of genus 3 over $k$, let $(j_2: j_3:\ldots:
  j_{10})$ be its Shioda invariants. Then the automorphism group of $C$ contains
  $\CG_{2} \times \CG_4$ if and only if
  {\rm\setcounter{myequation}{\theequation}
    \begin{dgroup*}[style={\small}]
      \begin{dmath*}
        0=648000\,j_6^2 - 356400\,j_6\,j_4\,j_2 + 4275\,j_6\,j_2^3 -
        839808\,j_4^3 + 75249\,j_4^2\,j_2^2 - 1449\,j_4\,j_2^4 + 8\,j_2^6\,,
      \end{dmath*}
      \begin{dmath}[style={\refstepcounter{myequation}},number={\themyequation},label=eq:8]
        0=85050\,j_8 + 13860\,j_6\,j_2 + 37908\,j_4^2 - 2961\,j_4\,j_2^2 +
        28\,j_2^4\,,
      \end{dmath}
      \begin{dmath*}
        0=2041200\,j_{10} + 560520\,j_6\,j_4 + 2205\,j_6\,j_2^2 +
        98091\,j_4^2\,j_2 - 5607\,j_4\,j_2^3 + 56\,j_2^5\,,
      \end{dmath*}
      \begin{dmath*}
        0={ j_3} \hiderel{=} { j_5} \hiderel{=} { j_7} \hiderel{=} { j_9}\,.
      \end{dmath*}
    \end{dgroup*}
    \setcounter{equation}{\themyequation}
  }
  \noindent
  Furthermore, a curve $C$ with automorphism group $\CG_{2} \times \CG_4$ is
  $K$-isomorphic to the curve $y^2 =
  {a}^{2}{x}^{8}+2\,{a}^{2}{x}^{6}+8\,a{x}^{2}-16$ with
  \begin{small}
    \begin{equation}\label{eq:modelC2C4}
      a = %
      \left\{\begin{array}{ll}
          196/3
          & \text{ if }6\,{j_4}-{{j_2}}^{2} = 0\,,\\
          -84
          & \text{ if }147\,{j_4}-2\,{{j_2}}^{2} = 0\,,\\
          \displaystyle{ {98}\,{\frac {36288\,{{j_4}}^{2}-3906\,{j_4}\,{{
                    j_2}}^{2}+14400\,{j_6}\,{j_2}+43\,{{j_2}}^{4}}{ 9 ( 96\,{
                  j_4}-{{j_2}}^{2} )  ( 147\,{j_4}-2\,{{j_2}}^{2}
                ) }}} & \text{ otherwise}\,.\\
        \end{array}\right.
    \end{equation}
  \end{small}
\end{lemma}
\medskip

\begin{remark}
  There is no difficulty to check that curves that annihilate denominators in
  Eq.~\eqref{eq:modelC2C4} have a larger automorphism group than $\CG_{2} \times
  \CG_4$.\medskip
  \begin{itemize}
  \item When $6\,{j_4}-{{j_2}}^{2}$ and $36\,{j_6}+{{j_2}}^{3}$ are both equal
    to 0, Eq.~\eqref{eq:8} can be reduced to Eq.~\eqref{eq:3} and
    Lemma~\ref{lemma:V8} can be used, instead, for reconstructing a
    model. Otherwise, \textit{i.e.} when $6\,{j_4}-{{j_2}}^{2} =0 $ but
    $36\,{j_6}+{{j_2}}^{3} \neq 0$ , then now $576\,{j_6}-65\,{{j_2}}^{3} = 0$
    but this lemma can still be applied.\medskip
  \item When $96\,{j_4}-{{j_2}}^{2}=0$, Eq.~\eqref{eq:8} can be reduced to
    Eq.~\eqref{eq:4} and Lemma~\ref{lemma:U6} can now be used for
    reconstructing a model.\medskip
  \item When $147\,{j_4}-2\,{{j_2}}^{2} = 0$ and
    $3087\,{j_6}-2\,{{j_2}}^{3}=0$, then we can choose for $f(x)$ the
    polynomial $x ( x-1 ) ( x+1 ) ( {x}^{2}+1 ) ^{2}$. Its discriminant is
    obviously equal to zero and no hyperelliptic curve exists with such
    invariants. Otherwise, \textit{i.e.}  when $147\,{j_4}-2\,{{j_2}}^{2} = 0$
    but $3087\,{j_6}-2\,{{j_2}}^{3} \neq 0$, then now
    $197568\,{j_6}-47\,{{j_2}}^{3} = 0$ but this lemma can still be applied.    %
  \end{itemize}
\end{remark}
\medskip

Here, the quotient field of $\ii_8$ modulo the ideal defined by
Eq.~\eqref{eq:8} is obtained by adjoining $j_6$ to $k[j_2, j_4]$ and $j_6$
satisfies an irreducible monic equation of degree $2$.  The invariants $j_8$
and $j_{10}$ are then rational in $j_2$, $j_4$ and $j_6$ (the invariants
$j_3$, $j_5$, $j_7$ and $j_9$ are trivial).

Geometrically, the projective variety defined by Eq.~\eqref{eq:8} has two
singularities, the $\VG_8$ and $\UG_6$ points. Moreover, it is birationally
equivalent to the conic
\begin{math}
{{X_1}}^{2}-{{11}/{20}}\,\,{X_1}\,{X_3}-{{162}/{125}}
\,\,{X_2}\,{X_3}+{{713}/{8000}}\,{{X_3}}^{2}
\end{math}
of discriminant
\begin{math}
  -2^2\cdot 3^8/ 5^6\,.
\end{math}
We can classically parameterized it by the slope of lines which all intersect
at the $\VG_8$ point, this yields
\begin{footnotesize}
  \begin{multline*}
    t \mapsto \left(%
      t~:%
      0~:%
      \frac{{t}^{2}}{6}-\frac{3\,t}{2}+{\frac {18}{5}}~:%
      0~:%
      -{\frac {\, \left( 5\,t-27 \right)  \left(
        25\,{t}^{2}-405\,t+1296 \right)}{4500}}~:%
      0~:\right.\\%
      \left.-{\frac {{t}^{4}}{420}}+{\frac {51\,{t}^{3}}{700}}-{\frac
        {684\,{t}^{2}}{875}}+{\frac {15516\,t}{4375}}-{\frac {25272}{4375}}~:%
      0~:%
      {\frac {{t}^{5}}{2520}}-{\frac {11\,{t}^{4}}{560}}+\frac{3\,{t}^{3}}{10}-{\frac {813\,{t}^{2}}{400}}+{\frac {28026\,t}{4375}}-{\frac {168156}{21875}}%
    \right)\,.
  \end{multline*}
\end{footnotesize}

The corresponding curves are of the form
\begin{displaymath}
  y^2 =  ( {x}^{4}+2\,{x}^{2}+{{5}/{147}}\,t-1/7 )  ( {x}^{4}-{{5}/{147}}\,t+1/7 )\,.
\end{displaymath}
The values $t=21/5$ and $t=168/5$ give polynomials with multiple roots, the
curve with automorphism group $\UG_6$ is obtained with $t=24/5$, otherwise the
automorphism group is always $\CG_{2} \times \CG_4$ since the point $\VG_8$ is
not reachable with this conic parametrization.
Incidentally, this shows again that over a finite field $\FF_{q}$ of
characteristic at least 11 there are $q-3$ isomorphic classes of curves with
automorphism group $\CG_{2} \times \CG_4$.

\subsubsection{Strata of dimension 2}
\label{sec:strata-dimension-2}

\stratum{ $\Aut(C) = \CG_{2}^3$}

\begin{lemma}\label{lemma:C2p3}
  Let $C$ be a hyperelliptic curve of genus 3 over $k$, let $(j_2: j_3:\ldots:
  j_{10})$ be its Shioda invariants. Then the automorphism group of $C$
  contains $\CG_{2}^3$ if and only if
  {\rm\setcounter{myequation}{\theequation}
    \begin{dgroup*}[style={\footnotesize},spread={-3pt}]
      \begin{dmath*}
        0=97200\,j_5^4 - 656100\,j_5^3\,j_3\,j_2 - 437400\,j_5^2\,j_4^2\,j_2 +
        2916000\,j_5^2\,j_4\,j_3^2 + 212625\,j_5^2\,j_4\,j_2^3 +
        753300\,j_5^2\,j_3^2\,j_2^2 - 25110\,j_5^2\,j_2^5 +
        2494800\,j_5\,j_4^3\,j_3 - 1336500\,j_5\,j_4^2\,j_3\,j_2^2 -
        5427000\,j_5\,j_4\,j_3^3\,j_2 + 180900\,j_5\,j_4\,j_3\,j_2^4 -
        3024000\,j_5\,j_3^5 + 201600\,j_5\,j_3^3\,j_2^3 - 3360\,j_5\,j_3\,j_2^6
        + 381024\,j_4^5 - 259200\,j_4^4\,j_2^2 - 186300\,j_4^3\,j_3^2\,j_2 +
        66960\,j_4^3\,j_2^4 + 2592000\,j_4^2\,j_3^4 +
        151200\,j_4^2\,j_3^2\,j_2^3 - 7920\,j_4^2\,j_2^6 +
        351000\,j_4\,j_3^4\,j_2^2 - 23400\,j_4\,j_3^2\,j_2^5 + 390\,j_4\,j_2^8 +
        108000\,j_3^6\,j_2 - 10800\,j_3^4\,j_2^4 + 360\,j_3^2\,j_2^7 -
        4\,j_2^{10}\,,
      \end{dmath*}
      \begin{dmath*}
        0=7620480\,j_6\,j_4^3 - 4626720\,j_6\,j_4^2\,j_2^2 -
        6609600\,j_6\,j_4\,j_3^2\,j_2 + 929070\,j_6\,j_4\,j_2^4 +
        10368000\,j_6\,j_3^4 + 1502550\,j_6\,j_3^2\,j_2^3 - 61605\,j_6\,j_2^6 +
        874800\,j_5^3\,j_3 - 2449440\,j_5^2\,j_4^2 + 947700\,j_5^2\,j_4\,j_2^2 -
        3207600\,j_5^2\,j_3^2\,j_2 - 89910\,j_5^2\,j_2^4 +
        3440880\,j_5\,j_4^2\,j_3\,j_2 + 16912800\,j_5\,j_4\,j_3^3 -
        982935\,j_5\,j_4\,j_3\,j_2^3 - 1555200\,j_5\,j_3^3\,j_2^2 +
        51840\,j_5\,j_3\,j_2^5 + 843696\,j_4^4\,j_2 - 4581360\,j_4^3\,j_3^2 -
        559278\,j_4^3\,j_2^3 - 168480\,j_4^2\,j_3^2\,j_2^2 +
        125091\,j_4^2\,j_2^5 - 64800\,j_4\,j_3^4\,j_2 +
        305370\,j_4\,j_3^2\,j_2^4 - 10107\,j_4\,j_2^7 + 432000\,j_3^6 +
        104400\,j_3^4\,j_2^3 - 8400\,j_3^2\,j_2^6 + 148\,j_2^9\,,
      \end{dmath*}
      \begin{dmath*}
        0=16200\,j_6\,j_5\,j_3 + 45360\,j_6\,j_4^2 -
        17550\,j_6\,j_4\,j_2^2 - 49950\,j_6\,j_3^2\,j_2 + 1665\,j_6\,j_2^4
        - 14580\,j_5^2\,j_4 + 2430\,j_5^2\,j_2^2 +
        22275\,j_5\,j_4\,j_3\,j_2 + 43200\,j_5\,j_3^3 -
        1440\,j_5\,j_3\,j_2^3 + 5022\,j_4^3\,j_2 - 34560\,j_4^2\,j_3^2 -
        2223\,j_4^2\,j_2^3 - 7650\,j_4\,j_3^2\,j_2^2 + 255\,j_4\,j_2^5 -
        3600\,j_3^4\,j_2 + 240\,j_3^2\,j_2^4 - 4\,j_2^7\,,
      \end{dmath*}
      \begin{dmath*}
        0=22680\,j_6\,j_5\,j_4 - 4995\,j_6\,j_5\,j_2^2 -
        14850\,j_6\,j_4\,j_3\,j_2 - 86400\,j_6\,j_3^3 +
        2880\,j_6\,j_3\,j_2^3 - 7290\,j_5^3 + 26730\,j_5^2\,j_3\,j_2 +
        2511\,j_5\,j_4^2\,j_2 - 80460\,j_5\,j_4\,j_3^2 -
        693\,j_5\,j_4\,j_2^3 - 360\,j_5\,j_3^2\,j_2^2 + 12\,j_5\,j_2^5 -
        10206\,j_4^3\,j_3 + 1350\,j_4^2\,j_3\,j_2^2 -
        4500\,j_4\,j_3^3\,j_2 + 150\,j_4\,j_3\,j_2^4 - 3600\,j_3^5 +
        240\,j_3^3\,j_2^3 - 4\,j_3\,j_2^6\,,
      \end{dmath*}
      \begin{dmath*}
        0=194400\,j_6\,j_5^2 + 2034720\,j_6\,j_4^2\,j_2 +
        2073600\,j_6\,j_4\,j_3^2 - 777870\,j_6\,j_4\,j_2^3 -
        2193750\,j_6\,j_3^2\,j_2^2 + 73125\,j_6\,j_2^5 -
        913680\,j_5^2\,j_4\,j_2 + 518400\,j_5^2\,j_3^2 +
        179550\,j_5^2\,j_2^3 + 1516320\,j_5\,j_4^2\,j_3 +
        419175\,j_5\,j_4\,j_3\,j_2^2 + 1555200\,j_5\,j_3^3\,j_2 -
        51840\,j_5\,j_3\,j_2^4 + 244944\,j_4^4 + 112590\,j_4^3\,j_2^2 -
        1170720\,j_4^2\,j_3^2\,j_2 - 80451\,j_4^2\,j_2^4 +
        86400\,j_4\,j_3^4 - 306810\,j_4\,j_3^2\,j_2^3 + 10131\,j_4\,j_2^6
        - 147600\,j_3^4\,j_2^2 + 9840\,j_3^2\,j_2^5 - 164\,j_2^8\,,
      \end{dmath*}
      \begin{dmath}[style={\refstepcounter{myequation}},number={\themyequation},label=eq:9]
        0=32400\,j_6^2 + 28350\,j_6\,j_4\,j_2 + 87750\,j_6\,j_3^2 -
        2925\,j_6\,j_2^3 - 21870\,j_5^2\,j_2 + 76545\,j_5\,j_4\,j_3 +
        13122\,j_4^3 - 405\,j_4^2\,j_2^2 + 5130\,j_4\,j_3^2\,j_2 -
        171\,j_4\,j_2^4 + 3600\,j_3^4 - 240\,j_3^2\,j_2^3 + 4\,j_2^6\,,
      \end{dmath}
      \begin{dmath*}
        0=12150\,j_7\,j_3 + 7560\,j_6\,j_4 - 1665\,j_6\,j_2^2 -
        2430\,j_5^2 + 1620\,j_5\,j_3\,j_2 + 837\,j_4^2\,j_2 +
        1530\,j_4\,j_3^2 - 231\,j_4\,j_2^3 - 120\,j_3^2\,j_2^2 +
        4\,j_2^5\,,
      \end{dmath*}
      \begin{dmath*}
        0=4860\,j_7\,j_4 - 810\,j_7\,j_2^2 - 1080\,j_6\,j_5 +
        3330\,j_6\,j_3\,j_2 - 837\,j_5\,j_4\,j_2 - 2880\,j_5\,j_3^2
        -12\,j_5\,j_2^3 + 2916\,j_4^2\,j_3 + 360\,j_4\,j_3\,j_2^2 +
        240\,j_3^3\,j_2 - 8\,j_3\,j_2^4\,,
      \end{dmath*}
      \begin{dmath*}
        0=18225\,j_7\,j_5 + 7425\,j_6\,j_4\,j_2 + 43200\,j_6\,j_3^2 -
        1440\,j_6\,j_2^3 - 10935\,j_5^2\,j_2 + 42525\,j_5\,j_4\,j_3 +
        5103\,j_4^3 - 675\,j_4^2\,j_2^2 + 2250\,j_4\,j_3^2\,j_2 -
        75\,j_4\,j_2^4 + 1800\,j_3^4 - 120\,j_3^2\,j_2^3 + 2\,j_2^6\,,
      \end{dmath*}
      \begin{dmath*}
        0=145800\,j_7\,j_6 + 4050\,j_7\,j_2^3 - 64800\,j_6\,j_5\,j_2 +
        83700\,j_6\,j_4\,j_3 - 16650\,j_6\,j_3\,j_2^2 + 77760\,j_5^2\,j_3
        + 13122\,j_5\,j_4^2 - 5751\,j_5\,j_4\,j_2^2 +
        8640\,j_5\,j_3^2\,j_2 + 252\,j_5\,j_2^4 - 10260\,j_4^2\,j_3\,j_2 +
        2880\,j_4\,j_3^3 - 1896\,j_4\,j_3\,j_2^3 - 1200\,j_3^3\,j_2^2 +
        40\,j_3\,j_2^5\,,
      \end{dmath*}
      \begin{dmath*}
        0=546750\,j_7^2 - 963900\,j_6\,j_4^2 + 453600\,j_6\,j_4\,j_2^2 +
        1404000\,j_6\,j_3^2\,j_2 - 46800\,j_6\,j_2^4 + 320760\,j_5^2\,j_4
        - 196830\,j_5^2\,j_2^2 + 204120\,j_5\,j_4\,j_3\,j_2 -
        680400\,j_5\,j_3^3 + 22680\,j_5\,j_3\,j_2^3 - 12312\,j_4^3\,j_2 +
        537570\,j_4^2\,j_3^2 + 32301\,j_4^2\,j_2^3 +
        157680\,j_4\,j_3^2\,j_2^2 - 5256\,j_4\,j_2^5 + 82800\,j_3^4\,j_2 -
        5520\,j_3^2\,j_2^4 + 92\,j_2^7\,,
      \end{dmath*}
      \begin{dmath*}
        0=85050\,j_8 - 11655\,j_6\,j_2 + 11340\,j_5\,j_3 + 2187\,j_4^2 -
        1260\,j_4\,j_2^2 - 840\,j_3^2\,j_2 + 28\,j_2^4\,,
      \end{dmath*}
      \begin{dmath*}
        0=9720\,j_9 + 810\,j_7\,j_2 - 3330\,j_6\,j_3 - 2187\,j_5\,j_4 +
        108\,j_5\,j_2^2 - 360\,j_4\,j_3\,j_2 - 240\,j_3^3 +
        8\,j_3\,j_2^3\,,
      \end{dmath*}
      \begin{dmath*}
        0=170100\,j_{10} + 37890\,j_6\,j_4 - 22680\,j_5^2 +
        4221\,j_4^2\,j_2 + 1680\,j_4\,j_3^2 - 56\,j_4\,j_2^3\,.
      \end{dmath*}
    \end{dgroup*}
    \setcounter{equation}{\themyequation} 
  }
  \noindent Furthermore, a curve $C$ with automorphism group $\CG_{2}^3$ is
  $K$-isomorphic to the curve $y^2 = {a_8}\,{x}^{8} + {a_6}\,{x}^{6} +
  {a_4}\,{x}^{4} + \lambda{a_6}\,{x}^{2} + {\lambda}^{2}{a_8}$ where $a_4$ is
  any root of the degree 3 equation {\rm\setcounter{myequation}{\theequation}
    \begin{dgroup*}[style={\footnotesize}]
      \begin{dmath}[style={\refstepcounter{myequation}},number={\themyequation},label={eq:13}]
        0 = 192\, ( -60\,{{j_3}}^{2}+2\,{{j_2}}^{3}+18\,{j_6}-9\,{j_4}\,{j_2} )\,
        {x}^{3}-90720\, ( 3\,{j_4}\,{j_3}+3\,{j_7}-{j_5}\,{j_2} )\,
        {x}^{2}+294\, (
        765\,{j_4}\,{{j_2}}^{2}+1440\,{j_6}\,{j_2}-5940\,{j_5}\,{j_3}-1782\,{{j_4}}^{2}+1140\,{{j_3}}^{2}{j_2}-38\,{{j_2}}^{4}
        )\,
        x+6637050\,{j_6}\,{j_3}+1250235\,{j_5}\,{j_4}-833490\,{j_5}\,{{j_2}}^{2}+2932650\,{j_4}\,{j_3}\,{j_2}+2881200\,{{j_3}}^{3}-96040\,{j_3}\,{{j_2}}^{3}
      \end{dmath}
    \end{dgroup*}
    \setcounter{equation}{\themyequation}
  }
  and
  \begin{math}
    a_6 = -28\,{{ \nu}}^{2}-{{a_4}}^{2}/5+14\,{j_2},\ a_8 = \nu\,
    a_6\text{ and }  \lambda = 1/a_6 
  \end{math}
  with
  \begin{footnotesize}
    \begin{displaymath}
      \nu = {\frac {18\,{j_6}\,{a_4}-9\,{a_4}\,{j_4}\,{j_2}-
          60\,{a_4}\,{{j_3}}^{2}+2\,{a_4}\,{{j_2}}^{3}-810\,{j_7}
          +270\,{j_5}\,{j_2}-810\,{j_4}\,{j_3}}{10\,(-18\,{j_6}+9\,{
            j_4}\,{j_2}+60\,{{j_3}}^{2}-2\,{{j_2}}^{3})}}\,.
    \end{displaymath}
  \end{footnotesize}
\end{lemma}

Eq.~\eqref{eq:13} may not have a root in the base field $k$. In
this case, Lemma~\ref{lemma:C2p3} yields a model over an extension $k'$ of
degree 3 (see
Section~\ref{sec:field} for the rationality issue).

\begin{remark}
  When $-18\,{j_6}+9\,{ j_4}\,{j_2}+60\,{{j_3}}^{2}-2\,{{j_2}}^{3} = 0$,
  Eq.~\eqref{eq:9} can be reduced to Eq.~\eqref{eq:6} and
  Lemma~\ref{lemma:C2D8} can be used, instead, for reconstructing a
  model.

\end{remark}
\medskip

The quotient field of $\ii_8$ modulo the ideal defined by Eq.~\eqref{eq:9} is
obtained by adjoining $j_5$ to $k[j_2, j_3, j_4]$ and $j_5$ satisfies an
irreducible monic equation of degree $4$.  The invariants $j_6$, $j_7$~,
\ldots $j_{10}$ are then rational in $j_2$, $j_3$, $j_4$ and $j_5$.

More geometrically, singularities of the projective variety defined by
Eq.~\eqref{eq:9} are on the curve $\CG_2\times \DG_8$.  Moreover, for non-zero
$J_2$, $J_3$, we can set $J_3=J_2=t$ for some parameter $t$ and the
intersection of the stratum with this hyperplane is a pencil of rational
curves that we can parametrize by an other parameter $u$.  So, we find for $t
\neq 0$ the parametrization
\begin{math}
  (t,u) \mapsto (j_2(t,u): j_3(t,u): \cdots : j_{10}(t,u))
\end{math}
where
{\rm\setcounter{myequation}{\theequation} 
  \begin{dgroup*}[style={\tiny},spread={-1pt},compact]
  \begin{dmath*}
    j_2(t,u) = t\,,\ j_3(t,u) = t\,, j_4(t,u) =
    -9800\,{u}^{4}+85\,t{u}^{2}-{\frac {35}{3}}\,tu+{\frac
      {2}{147}}\,{t}^{2}\,,
  \end{dmath*}
  \begin{dmath*}
    j_5(t,u) = -137200\,{u}^{5}+700\,t{u}^{3}+{\frac {350}{3}}\,t{u}^{2}-{\frac
      {25}{28}}\,{t}^{2}u+\frac{1}{28}\,{t}^{2},
  \end{dmath*}
  \begin{dmath*}
    j_6(t,u) = -617400\,{u}^{6}+2030\,t{u}^{4}+1470\,t{u}^{3}-{\frac
      {689}{56}}\,{t}^{2}{u}^{2}+{\frac
      {11}{12}}\,{t}^{2}u-\frac{1}{24}\,{t}^{2}+{\frac {2}{3087}}\,{t}^{3},
  \end{dmath*}
  \begin{dmath*}
    j_7(t,u) = -1920800\,{u}^{7}-{\frac {72520}{3}}\,t{u}^{5}+{\frac
      {43120}{3}}\,t{u}^{4}+{\frac {377}{6}}\,{t}^{2}{u}^{3}-{\frac
      {226}{9}}\,{t}^{2}{u}^{2}+{\frac {1}{588}}\,{t}^{2} ( 17\,t+686
    ) u-{\frac {1}{588}}\,{t}^{3},
  \end{dmath*}
  \begin{dmath}[style={\refstepcounter{myequation}},number={\themyequation},label={eqparaC2p3}]
    j_8(t,u) = -2469600\,{u}^{8}-{\frac {125300}{3}}\,t{u}^{6}+{\frac
      {37240}{3}}\,t{u}^{5}-{\frac {643}{14}}\,{t}^{2}{u}^{4}+{\frac
      {1432}{9}}\,{t}^{2}{u}^{3}-{\frac {1}{246960}}\,{t}^{2} (
      4705960+120087\,t ) {u}^{2}+{\frac {1411}{17640}}\,{t}^{3}u-{\frac
      {11}{252105}}\,{t}^{4}-{\frac {1}{1680}}\,{t}^{3},
  \end{dmath}
  \begin{dmath*}
    j_9(t,u) = 302526000\,{u}^{9}-{\frac {12022150}{3}}\,t{u}^{7}-{\frac
      {325850}{3}}\,t{u}^{6}+{\frac {665105}{36}}\,{t}^{2}{u}^{5}-{\frac
      {9905}{18}}\,{t}^{2}{u}^{4}-{\frac {1}{2016}}\,{t}^{2} (
      -397880+56341\,t ) {u}^{3}+{\frac
      {18829}{6048}}\,{t}^{3}{u}^{2}+{\frac {1}{98784}}\,{t}^{3} (
      472\,t-30527 ) u+{\frac {1}{96}}\,{t}^{3}-{\frac
      {1}{4116}}\,{t}^{4},
  \end{dmath*}
  \begin{dmath*}
    j_{10}(t,u) = 1162084000\,{u}^{10}-11872700\,t{u}^{8}-{\frac
      {7991900}{3}}\,t{u}^{7}+{\frac {1366535}{18}}\,{t}^{2}{u}^{6}-{\frac
      {51835}{9}}\,{t}^{2}{u}^{5}-{\frac {1}{2352}}\,{t}^{2} (
      -13267240+268943\,t ) {u}^{4}-{\frac
      {77657}{3024}}\,{t}^{3}{u}^{3}+{\frac {1}{1037232}}\,{t}^{3} (
      67571+105696\,t ) {u}^{2}-{\frac {1}{148176}}\,{t}^{3} (
      824\,t-1029 ) u-{\frac {11}{5294205}}\,{t}^{5}+{\frac
      {1}{6174}}\,{t}^{4}\,.
  \end{dmath*}
\end{dgroup*}
\setcounter{equation}{\themyequation}
}
\noindent 
The lines
\begin{small}
  \begin{math}
   280\,{u}^{2}-t = 0
  \end{math}
\end{small}
and
\begin{small}
  \begin{math}
    14280\,{u}^{3}-33\,tu-7\,t = 0
  \end{math}
\end{small}
are the loci of forms with multiple roots (discriminant equal to 0).
The rational curve
\begin{footnotesize}
  \begin{math}
    -564715200\,{u}^{6}+6338640\,{u}^{4}t+1344560\,{u}^{3}t-50127\,{u}^{2}
    {t}^{2}+6174\,u{t}^{2}+96\,{t}^{3}-3087\,{t}^{2} =0
  \end{math}
\end{footnotesize}
is the locus of forms with automorphism groups which contain $\CG_2\times
\DG_8$\,.    \medskip

In the case $J_2=0$ and non-zero $J_3$, $J_4$, setting $J_4=J_3$ yields a rational
curve that we parametrize as 
\begin{math}
  t \mapsto (j_2(t) : j_3(t): \ldots : j_{10}(t))\,,
\end{math}
with
\begin{math}
  t \neq -105/1260\,,
\end{math}
where
\begin{dgroup*}[style={\tiny},compact]
  \begin{dmath*}
    j_2(t) = 0\,,\ 
    j_{3}(t) = 1260\,t+105\,,\ 
    j_{4}(t) = 7350\,t\,,\ 
    j_{5}(t) = 36750\,t+7350\,,\ 
    j_{6}(t) = -66150\,{t}^{2}-242550\,t-29400\,,
  \end{dmath*}
  \begin{dmath*}
    j_{7}(t) = -926100\,{t}^{2}+977550\,t+102900\,,\ 
    j_{8}(t) = -7563150\,{t}^{2}-1749300\,t-102900\,,\ 
  \end{dmath*}
  \begin{dmath*}
    j_{9}(t) = 20837250\,{t}^{3}-33957000\,{t}^{2}-8232000\,t-1029000\,,\ 
    j_{10}(t) = -6945750\,{t}^{3}+557975250\,{t}^{2}+119364000\,t+7203000\,.
  \end{dmath*}
\end{dgroup*}
The point $t = -2/7$ yields a form with multiple roots, the roots of
$81\,{t}^{2}+17\,t+1=0$ yield forms with
automorphism groups which contain $\CG_2\times \DG_8$\,.
\medskip

Similarly, in the case $J_3=0$ and non-zero $J_2$, $J_5$, setting $J_5=J_2^2$
yields the rational curve
\begin{dgroup*}[style={\tiny},compact]
  \begin{dmath*}
    j_{2}(t) = {\frac {1010}{3}}\, ( 52\,t+59 ) ( t+1 ) \,,\
    j_{3}(t) = 0\,,\
    j_{4}(t) = { \frac {50}{27}}\, ( 52\,t+59 ) ^{2} ( 11272\,{t}^{2}+
      23147\,t+11767 ) \,,\
  \end{dmath*}
  \begin{dmath*}
    j_{5}(t) = {\frac {25}{9}}\, ( 52\,t+59 ) ^{3} ( 17\,t+29 )
    ^{2}\,,\
    j_{6}(t) = -{\frac {25}{486}}\, ( 52\,t+59 ) ^{3} (
      27670948\,{t}^{3}+88650255\,{t}^{2}+94796346\,t+ 33852031 ) \,,\
  \end{dmath*}
  \begin{dmath*}
    j_{7}(t) = {\frac {50}{81}}\, ( 52\,t+59 ) ^{4} ( 17\,t+29 )
    ( 38117\,{t}^{2}+81691\,t+43682 ) \,,\
  \end{dmath*}
  \begin{dmath*}
    j_{8}(t) = -{\frac {25}{3402}}\, ( 52\,t+59 ) ^{4} (
      6266813468\,{t}^{4}+26795826605\,{t}^{3}+42971724291\,{t}^{2}+
      30632960351\,t+8190389165 ) \,,\
  \end{dmath*}
  \begin{dmath*}
    j_{9}(t) = -{\frac {125}{2916}}\, ( 52\,t +59 ) ^{5} ( 17\,t+29 )
    ( 11612756\,{t}^{3}+ 34093287\,{t}^{2}+33246138\,t+10800599 ) \,,\
  \end{dmath*}
  \begin{dmath*}
    j_{10}(t) = {\frac {125}{30618}}\, ( 52\,t+59 ) ^{5} (
      794827443664\,{t}^{5}+
      4302413461544\,{t}^{4}+9310923963145\,{t}^{3}+10069939679311\,{t}^{2}+
      5442845374787\,t+1176260576869 )\,.
  \end{dmath*}
\end{dgroup*}
\noindent
The point
\begin{tiny}
  \begin{math}
\left({\frac {1010}{39}}:0:{\frac {563600}{4563}}:{\frac {14450}{1521}}:-{
\frac {345886850}{533871}}:{\frac {64798900}{177957}}:-{\frac {
78335168350}{48582261}}:-{\frac {12338553250}{20820969}}:{\frac {
49676715229000}{5684124537}}\right),
  \end{math}
\end{tiny}
when defined, is the only point which is not reached by the parametrization.  Moreover, the
points $t=-4259/2831$, $t=-59/52$ and $t=1/83$ (when defined) yield forms with
multiple roots. The points $t = -59/52$, $t = -29/17$ (when defined) and the
roots of $9844\,t^2 + 19487\,t + 9679 = 0$ yield forms with automorphism
groups which contain $\CG_2\times \DG_8$\,.  \medskip

Modulo 101, the latter degenerates to only one projective point. In this case,
we thus consider instead the rational curve
\begin{dgroup*}[style={\tiny},compact]
  \begin{dmath*}
    j_{2}(t) = 32\, ( t+44 ) ( t+12 ) \,,\
    j_{3}(t) = 0\,,\
    j_{4}(t) = 57\, ( {t}^{2}+ 84\,t+59 ) ( t+12 ) ^{2}\,,\
    j_{5}(t) = 14\, ( t+12 ) ^{ 3} ( t+84 ) ^{2}\,,\
  \end{dmath*}
  \begin{dmath*}
    j_{6}(t) = 22\,{\frac { ( t+62 ) ( t +12 ) ^{3} ( t+22
        ) ( {t}^{2}+73\,t+23 ) ( t+54 )
      }{2\,{t}^{2}+45\,t+46}}\,,\
    j_{7}(t) = {\frac {97}{51}}\, ( t+12 ) ^{4} (
      {t}^{2}+7\,t+86 ) ( t+84 ) \,,\
  \end{dmath*}
  \begin{dmath*}
    j_{8}(t) = 81\,{\frac { ( t+12 ) ^{4} ( {t}^{2}+51\,t+ 54 )
        ( {t}^{2}+73\,t+23 ) ( {t}^{2}+64\,t+62 )
      }{2\,{t}^{2}+45\,t+46}}\,,\
    j_{9}(t) = {\frac {2}{51}}\, ( t+12 ) ^{5} (
      {t}^{2}+3\,t+19 ) ( t+73 ) ( t+84 ) \,,\
  \end{dmath*}
  \begin{dmath*}
    j_{10}(t) = 3\,{\frac { ( {t}^{2}+3\,t+75 ) ( {t}^{2}+26\,t+79 )
        ( t+12 ) ^{5} ( {t}^ {2}+73\,t+23 ) ( t+68
        ) }{2\,{t}^{2}+45\,t+46}}\,.
  \end{dmath*}
\end{dgroup*}
\noindent
The point
\begin{small}
  \begin{math}
(32: 0: 57: 14: 11: 93: 91: 4: 52)
  \end{math}
\end{small}
is the only point which is not reached by the parametrization.  Moreover, the
points $t=20$, $t=56$ and $t=89$ yield forms with
multiple roots. The points $t = 17$, $t = 89$  and the
roots of $t^2 + 14\,t + 87=0$ yield forms with automorphism
groups which contain $\CG_2\times \DG_8$\,.  \medskip

We can now easily check that the rational projective points defined on this
variety are none other than the ones given by these
parametrizations. Incidentally, this shows that over a finite field $\FF_{q}$
of characteristic at least 11 there are $q^2-2\,q+2$ isomorphic classes of
curves with automorphism group $\CG_2^3$.

\stratum{ $\Aut(C) = \CG_{4}$}
\ \\

Unlike the previous strata, curves with automorphism group $\CG_4$ can be
reconstructed with the conic and quartic method of
Section~\ref{sec:conic-quartic}.

\begin{lemma}\label{lemma:C4}
  Let $C$ be a hyperelliptic curve of genus 3 over $k$, let $(j_2: j_3:\ldots:
  j_{10})$ be its Shioda invariants. Then the automorphism group of $C$
  contains $\CG_{4}$ if and only if
  {\rm\setcounter{myequation}{\theequation} 
    \begin{dgroup*}[style={\footnotesize},spread={-3pt}]
      \begin{dmath*}
        0 = 1093705578000 \,j_8^3 + 291654820800 \,j_8^2 \,j_6 \,j_2 -
        1104121821600 \,j_8^2 \,j_4^2 + 6076142100 \,j_8^2 \,j_4 \,j_2^2 -
        1469076134400 \,j_8 \,j_6^2 \,j_4 + 36726903360 \,j_8 \,j_6^2 \,j_2^2 +
        231472080000 \,j_8 \,j_6 \,j_4^2 \,j_2 - 6631246440 \,j_8 \,j_6 \,j_4 \,j_2^3
        + 30005640 \,j_8 \,j_6 \,j_2^5 + 267846264000 \,j_8 \,j_4^4 -
        23803045560 \,j_8 \,j_4^3 \,j_2^2 + 540101520 \,j_8 \,j_4^2 \,j_2^4 -
        3333960 \,j_8 \,j_4 \,j_2^6 - 256048128000 \,j_6^4 +
        42247941120 \,j_6^3 \,j_4 \,j_2 - 405409536 \,j_6^3 \,j_2^3 +
        8888527872 \,j_6^2 \,j_4^3 - 5893679232 \,j_6^2 \,j_4^2 \,j_2^2 +
        274718304 \,j_6^2 \,j_4 \,j_2^4 - 3093174 \,j_6^2 \,j_2^6 -
        9523422720 \,j_6 \,j_4^4 \,j_2 + 1237946976 \,j_6 \,j_4^3 \,j_2^3 -
        60328800 \,j_6 \,j_4^2 \,j_2^5 + 1037232 \,j_6 \,j_4 \,j_2^7 -
        5488 \,j_6 \,j_2^9 + 5509980288 \,j_4^6 - 1031704128 \,j_4^5 \,j_2^2 +
        58796766 \,j_4^4 \,j_2^4 - 1091475 \,j_4^3 \,j_2^6 +
        6174 \,j_4^2 \,j_2^8,
      \end{dmath*}
      \begin{dmath}[style={\refstepcounter{myequation}},number={\themyequation},label={eq:10}]
        0 = 2540160 \,j_{10} \,j_6 - 317520 \,j_{10} \,j_4 \,j_2 + 4410 \,j_{10} \,j_2^3
        + 1786050 \,j_8^2 + 238140 \,j_8 \,j_6 \,j_2 - 884520 \,j_8 \,j_4^2 +
        133056 \,j_6^2 \,j_4 + 47880 \,j_6 \,j_4^2 \,j_2 - 1092 \,j_6 \,j_4 \,j_2^3 -
        17496 \,j_4^4 - 2205 \,j_4^3 \,j_2^2 + 49 \,j_4^2 \,j_2^4,
      \end{dmath}
      \begin{dmath*}
        0 = 1543147200 \,j_{10} \,j_8 - 793618560 \,j_{10} \,j_4^2 +
        34292160 \,j_{10} \,j_4 \,j_2^2 - 357210 \,j_{10} \,j_2^4 -
        144670050 \,j_8^2 \,j_2 + 1538248320 \,j_8 \,j_6 \,j_4 -
        36435420 \,j_8 \,j_6 \,j_2^2 - 82668600 \,j_8 \,j_4^2 \,j_2 +
        1905120 \,j_8 \,j_4 \,j_2^3 + 254016000 \,j_6^3 -
        10160640 \,j_6^2 \,j_4 \,j_2 - 38808 \,j_6^2 \,j_2^3 - 50388480 \,j_6 \,j_4^3
        + 6373080 \,j_6 \,j_4^2 \,j_2^2 - 278460 \,j_6 \,j_4 \,j_2^4 +
        3136 \,j_6 \,j_2^6 - 11809800 \,j_4^4 \,j_2 + 643545 \,j_4^3 \,j_2^3 -
        7497 \,j_4^2 \,j_2^5,
      \end{dmath*}
      \begin{dmath*}
        0=4800902400 \,j_{10}^2 + 262906560 \,j_{10} \,j_4^2 \,j_2 -
        9485910 \,j_{10} \,j_4 \,j_2^3 + 92610 \,j_{10} \,j_2^5 -
        3541737150 \,j_8^2 \,j_4 + 37507050 \,j_8^2 \,j_2^2 -
        555660000 \,j_8 \,j_6^2 - 94382820 \,j_8 \,j_6 \,j_4 \,j_2 +
        1049580 \,j_8 \,j_6 \,j_2^3 + 133698600 \,j_8 \,j_4^3 -
        793800 \,j_8 \,j_4^2 \,j_2^2 - 23756544 \,j_6^2 \,j_4^2 +
        3087000 \,j_6^2 \,j_4 \,j_2^2 + 98280 \,j_6 \,j_4^3 \,j_2 +
        47460 \,j_6 \,j_4^2 \,j_2^3 - 980 \,j_6 \,j_4 \,j_2^5 + 3831624 \,j_4^5 +
        2010015 \,j_4^4 \,j_2^2 - 102312 \,j_4^3 \,j_2^4 + 1029 \,j_4^2 \,j_2^6,
      \end{dmath*}
      \begin{dmath*}
        0 = { j_3} \hiderel{=}{ j_5} \hiderel{=}{ j_7} \hiderel{=}{ j_9}.
      \end{dmath*}
    \end{dgroup*}
    \setcounter{equation}{\themyequation}
  }
  
  \noindent
  Furthermore, at least one among the five determinants
  \begin{footnotesize}
    \begin{equation}\label{eq:14}
      R(C_{5,2}, C_{6,2}, C_{7,2}),\  
      R(C_{5,2}, C_{6,2}, C_{7,2}'),\ 
      R(C_{5,2},C_{7,2},C_{8,2}'),\   
      R(C_{5,2}, C_{8,2}, C_{9,2}) \text{ and } 
      R(C_{5,2}, C_{6,2}, C_{9,2}'')            
    \end{equation}
  \end{footnotesize}
  is non-zero for a curve $C$ with automorphism group $\CG_4$ and if the
  corresponding conic has a $k$-rational point, then $C$ is $K$-isomorphic to
  a curve with a $k$-rational model constructed with the conic and quartic
  method of Prop.~\ref{prop:multic}.
\end{lemma}

\begin{remark}
  To check that at least one among the five determinants of Eq.~\eqref{eq:14}
  is non-zero for a curve with automorphism group $\CG_4$, we solve in $a$ and
  $b$ (over the integers to ensure that the result is still true modulo any
  positive prime $p$) the system obtained by evaluating these determinants at
  normal forms $x\,(x^2-1)\,(x^4+a\,x^2+b)$.  We found a finite number of
  irreducible components, possibly defined in a $k$-extension, but all with an
  automorphism group larger than $\CG_4$ (\textit{cf.}  Tab.~\ref{tab:c4}).
\end{remark}

\begin{table}
  \centering
  \begin{small}
    \begin{tabular}{l|l|c}
      \multicolumn{2}{c|}{$f(x)=x\,(x^2-1)\,(x^4+a\,x^2+b)$} & $\Aut(C)$ \\\hline\hline
      $a=0$ & $b=0$ & $\UG_6, \VG_8, \CG_2\times \SG_4\subset \Aut(C)$\\\hline
      $any$ $a$ & $b=1$ & $\CG_2\times \CG_4\subset \Aut(C)$\\
      $any$ $a$ & $b^2+(-3\,a+1)\,b+a^3=0$ & $\CG_2\times \CG_4\subset \Aut(C)$\\
      \hline
      $a=-6$ & $b=-27$ & $\UG_6\subset \Aut(C)$\\
      $a=10/3$ & $b=1$ & $\UG_6\subset \Aut(C)$\\
      $a=2/9$ & $b=-1/27$ & $\UG_6\subset \Aut(C)$\\
      \hline
      $a=-6$ & $b=1$ & $\VG_8\subset \Aut(C)$\\
      ${{a}}^{2}+40\,{a}+8 = 0$ & $b=-5\,a-1$
      & $\VG_8\subset \Aut(C)$ \\
      \hline
    \end{tabular}\medskip
  \end{small}
  \caption{$\CG_4$-hyperelliptic polynomials that cancel all the determinants of Eq.~\eqref{eq:14}}
\label{tab:c4}
\end{table}

\begin{remark}
  In fields of characteristic $0$, we can conclude with only the 4
  determinants $R(C_{5,2}, C_{6,2}, C_{7,2})$,\ 
  $R(C_{5,2},C_{7,2},C_{8,2}')$,\ 
  $R(C_{5,2}, C_{8,2}, C_{9,2})$ and 
  $R(C_{5,2}, C_{6,2}, C_{9,2}'')$.  
  These four determinants are enough in positive characteristic $p$ too,
  except for a finite number of primes the smallest of which is $p=47$ (modulo
  $47$, the point $(1~: 0~: 1: 0: 3: 0: 43: 0: 18)$ cancel all the
  determinants of Eq.~\eqref{eq:14} except $R(C_{5,2}, C_{6,2}, C_{7,2}')$).
\end{remark}

The quotient field of $\ii_8$ modulo the ideal defined by Eq.~\eqref{eq:10} is
obtained by adjoining $j_8$ to $k[j_2, j_4, j_6]$ and $j_8$ satisfies an
irreducible monic equation of degree $3$.  The invariant $j_{10}$ is then
rational in $j_2$, $j_4$, $j_6$ and $j_8$ (the invariants $j_3$, $j_5$, $j_7$
and $j_9$ are trivial).

Geometrically, singularities of the projective variety defined by
Eq.~\eqref{eq:10} are on the curve $\CG_2\times \CG_4$. Moreover, for non-zero
$J_2$, we can set $J_2=1$, $J_4=t$ for some parameter $t$ and the intersection
of the stratum with this hyperplane is a pencil a rational curves that can we
parametrize with an other parameter $u$. So, we find for $t\neq 0$ the
parametrization
\begin{math}
  (t,u) \mapsto (j_2(t,u): 0: j_4(t,u): 0: j_6(t,u): 0: j_8(t,u): 0: j_{10}(t,u))\,,
\end{math}
where
{\rm\setcounter{myequation}{\theequation} 
  \begin{dgroup*}[style={\tiny},spread={-1pt},compact]
  \begin{dmath*}
    j_2(t,u) = t\,,\ j_4(t,u) = {\frac {1}{96}}\, ( t+1 ) {t}^{2}\,,
  \end{dmath*}
  \begin{dmath*}
    j_6(t,u) = -{\frac {1}{124416000}}\, ( 720\,u+63\,t-25 )  ( -518400\,{u}^{2}+ ( -90720\,t+36000 ) u+2700\,{t}^{3}-1269\,{t}^{2}+3150\,t-
625 )
\,,
  \end{dmath*}
  \begin{dmath}[style={\refstepcounter{myequation}},number={\themyequation},label={eqparaC4}]
    j_8(t,u) = \frac{8}{3}\,{u}^{4}+ ( 2/3\,t-{\frac {10}{27}} ) {u}^{3}+ ( -{\frac {1}{360}}\,{t}^{3}+{\frac {139}{3600}}\,{t}^{2}-{\frac {5}{72}}\,t+{\frac {25
}{1296}} ) {u}^{2}+ ( -{\frac {17}{14400}}\,{t}^{4}+{\frac {241}{432000}}\,{t}^{3}-{\frac {139}{51840}}\,{t}^{2}+{\frac {25}{10368}}\,t-{\frac 
{125}{279936}} ) u-{\frac {1}{483840}}\,{t}^{6}-{\frac {9707}{145152000}}\,{t}^{5}+{\frac {82091}{967680000}}\,{t}^{4}-{\frac {299}{18662400}}\,{t}^
{3}+{\frac {139}{2985984}}\,{t}^{2}-{\frac {125}{4478976}}\,t+{\frac
  {625}{161243136}}\,,
  \end{dmath}
  \begin{dmath*}
    j_{10}(t,u) =-\frac{5}{3}\,{u}^{5}+ ( -{\frac {31}{48}}\,t+{\frac {125}{432}} ) {u}^{4}+ ( -{\frac {23}{3360}}\,{t}^{3}-{\frac {3349}{40320}}\,{t}^{2}+{\frac {
155}{1728}}\,t-{\frac {625}{31104}} ) {u}^{3}+ ( -{\frac {1}{345600}}\,{t}^{4}-{\frac {94723}{29030400}}\,{t}^{3}+{\frac {3349}{387072}}\,{t}^{
2}-{\frac {775}{165888}}\,t+{\frac {3125}{4478976}} ) {u}^{2}+ ( {\frac {61}{3870720}}\,{t}^{6}+{\frac {86893}{580608000}}\,{t}^{5}-{\frac {
735179}{13934592000}}\,{t}^{4}+{\frac {105073}{418037760}}\,{t}^{3}-{\frac {16745}{55738368}}\,{t}^{2}+{\frac {3875}{35831808}}\,t-{\frac {15625}{
1289945088}} ) u+{\frac {1}{44236800}}\,{t}^{7}+{\frac {734581}{139345920000}}\,{t}^{6}-{\frac {23681183}{3344302080000}}\,{t}^{5}+{\frac {736579}{
401316249600}}\,{t}^{4}-{\frac {108523}{24078974976}}\,{t}^{3}+{\frac {83725}{24078974976}}\,{t}^{2}-{\frac {19375}{20639121408}}\,t+{\frac {15625}{
185752092672}}\,.
  \end{dmath*}
\end{dgroup*}
\setcounter{equation}{\themyequation}
}
\noindent The intersections of the stratum with the $J_4=J_2^2/96$ and with
the $J_2=0$, $J_6=J_4$ hyperplanes are two curves which are not reachable by
this parametrization.
Furthermore, the rational curves 
\begin{small}
  \begin{math}
    5040\,u-579\,t-175=0
  \end{math}
\end{small}
and
\begin{small}
  \begin{math}
    25401600\,{u}^{2}+ \left( -12186720\,t-1764000 \right) u+185220\,{t}^{
      3}-2167119\,{t}^{2}+423150\,t+30625=0
  \end{math}
\end{small}
are the loci of forms with multiple roots (discriminant equal to 0), the rational curve
\begin{tiny}
  \begin{multline*}
    268738560000\,{u}^{4}+1492992000\, ( 78\,t-25 )\,
    {u}^{3}-3110400\, ( 90\,{t}^{3}-5409\,{t}^{2}+3900\,t-625 )\,
    {u}^{2}\\-2880\, (
      130410\,{t}^{4}-323352\,{t}^{3}+405675\,{t}^{2}-146250\,t+15625 )\,
    u\\+4665600\,{t}^{6}-55992060\,{t}^{5}+29516481\,{t}^{4}-31997700\,{t}^{3}+20283750
    \,{t}^{2}-4875000\,t+390625 = 0
  \end{multline*}
\end{tiny}
is the locus of forms with automorphism groups which contain $\CG_2\times
\CG_4$\,.    \medskip

A parametrization for the intersection with $J_4=J_2^2/96$ is of the form
\begin{math}
  t \mapsto (\,j_2(t): 0: j_4(t): 0: j_6(t): 0: j_8(t): 0: j_{10}(t)\,)
\end{math}
where
\begin{dgroup*}[style={\tiny},spread={-1pt},compact]
  \begin{dmath*}
    j_2(t) = t\,,\ 
    j_4(t) = {\frac {1}{96}}\,{t}^{2}\,,\ 
    j_6(t) = {\frac {1}{7776000}}\, ( 71\,t+2 ) ( 2183\,{t}^{2}+142\,t+2 )\,,
  \end{dmath*}
  \begin{dmath*} 
    j_{10}(t) = {\frac {16012139}{8817984000}}\,{t }^{4}+{\frac
      {51401}{196830000}}\,{t}^{3}+{\frac {7301}{524880000}}\,{t}^{2}+{\frac
      {31}{98415000}}\,t+{\frac {1}{393660000}},0,-{\frac {905961157}{
        3174474240000}}\,{t}^{5}-{\frac {672811}{14171760000}}\,{t}^{4}-{\frac
      {2449033}{793618560000}}\,{t}^{3}-{\frac
      {77179}{793618560000}}\,{t}^{2}-{\frac { 337}{226748160000}}\,t-{\frac
      {1}{113374080000}}
  \end{dmath*}
\end{dgroup*}
with
\begin{math}
  t \neq 0\,.
\end{math} 
The point
\begin{tiny}
  \begin{math}
    \left(1:0:{\frac {1}{96}}:0:{\frac {154993}{7776000}}:0:{\frac {16012139}{8817984000}}:0:-{\frac {905961157}{3174474240000}}\right)
  \end{math}
\end{tiny}
is the only point on this line which is not reached by the parametrization.
The points $13\,t - 14 = 0$ and $199904\,{t}^{2}+2296\,t-49=0$, when defined,
yield forms with multiple roots. The values $101\,t+2=0$ and $28\,t+1=0$, when
defined too, yield forms with automorphism groups which contain $\CG_2\times \CG_4$\,.
\medskip

Similarly, a parametrization for the intersection with $J_2=0$ and $J_6=J_4$ is
\begin{dgroup*}[style={\tiny},spread={-1pt},compact]
  \begin{dmath*}
    j_2(t) = 0\,,\ j_4(t) = {\frac {4}{3969}}\, ( t-9 ) (
      11\,t+6 )\,,\ j_6(t) = {\frac {4}{3969}}\, ( t+1 )
    ( 11\,t+6 ) ^{2}\,,
  \end{dmath*}
  \begin{dmath*}
    j_8(t) = {\frac {32}{36756909}}\, ( 373\,{t}^{2}+440\,t+127 ) 
    ( 11\,t+6 ) ^{2}\,,\ 
    j_{10}(t) = -{\frac {8}{992436543}}\, ( 835\,{
     t}^{2}+646\,t+51 )  ( 11\,t+6 ) ^{3}
  \end{dmath*}
\end{dgroup*}
with
\begin{math}
  11\,t+6 \neq 0\,.
\end{math}
The point
\begin{tiny}
  \begin{math}
    \left(0:0:{\frac {4}{43659}}:0:{\frac {4}{43659}}:0:{\frac {11936}{4447585989}}:0:-{\frac {6680}{120084821703}}\right),
  \end{math}
\end{tiny}
when defined, is the only point on this line which is not reached by the
parametrization.  Moreover, the point
\begin{small}
  \begin{math}
    t=1/39,
  \end{math}
\end{small}
when defined, yields a form with multiple roots. The
roots of
\begin{small}
  \begin{math}
    667\,{t}^{2}+769\,t+302 =0,
  \end{math}
\end{small}
when defined too, yield forms with automorphism groups which contain $\CG_2\times
\CG_4$\,.  \medskip

 We can now easily check that the rational projective points defined on this
 variety are none other than the ones given by these
 parametrizations. Incidentally, this shows that over a finite field $\FF_{q}$
 of characteristic at least 11 there are $q^2-2\,q+2$ isomorphic classes of
 curves with automorphism group $\CG_4$.

\begin{remark}
  We are also able to deduce from these parametrizations explicit models over the field of moduli (see 
   Section~\ref{subsec:descentg3}), but such expressions are far too
  large to be written down here.
\end{remark}

\subsubsection{Stratum of dimension 3}
\label{sec:strata-dimension-3}

The stratum of dimension 3 corresponds to curves with automorphism group
$\DG_4$.  Contrary to the other strata, we only give generators for the
grevlex order $J_{2}< J_3< \ldots <J_{10}$ with weights 2, 3, \ldots, 10
(since generators for the lexical order are far too huge to be explicitly
written here).

\begin{lemma}\label{lemma:D4}
  Let $C$ be a hyperelliptic curve of genus 3 over $k$, let $(j_2, j_3,\ldots,
  j_{10})$ be its Shioda invariants. If the automorphism group of $C$ contains
  $\DG_{4}$, then $j_2$, $j_3$, \ldots, $j_{10}$ satisfy the following 24
  equations, from degree 16 up to degree 24 (\textit{cf.}
  Appendix~\ref{sec:strat-equat-autm} for the complete set),
  {\rm\setcounter{myequation}{\theequation} 
    \begin{dgroup*}[style={\footnotesize},spread={-4pt}]
      \begin{dmath*}
        0 = 49\,j_2^{4}\,j_4^{2} - 1470\,j_2\,j_3^{2}\,j_4^{2} - 2205\,j_2^{2}\,j_4^{3} - 17496\,j_4^{4} +
        19845\,j_3\,j_4^{2}\,j_5 - 1092\,j_2^{3}\,j_4\,j_6 + 32760\,j_3^{2}\,j_4\,j_6 +
        47880\,j_2\,j_4^{2}\,j_6 - 238140\,j_5^{2}\,j_6 + 133056\,j_4\,j_6^{2} + 357210\,j_4\,j_5\,j_7 -
        238140\,j_3\,j_6\,j_7 - 59535\,j_2\,j_7^{2} - 884520\,j_4^{2}\,j_8 + 238140\,j_2\,j_6\,j_8 +
        1786050\,j_8^{2} + 119070\,j_3\,j_4\,j_9 + 59535\,j_2\,j_5\,j_9 - 714420\,j_7\,j_9 +
        4410\,j_2^{3}\,j_{10} - 132300\,j_3^{2}\,j_{10} - 317520\,j_2\,j_4\,j_{10} + 2540160\,j_6\,j_{10},
      \end{dmath*}
      \begin{dmath}[style={\refstepcounter{myequation}},number={\themyequation},label={eq:11}]
        0=    -14\,j_2^{3}\,j_3\,j_4^{2} + 420\,j_3^{3}\,j_4^{2} + 630\,j_2\,j_3\,j_4^{3} - 189\,j_2^{2}\,j_4^{2}\,j_5 + 
        360\,j_3\,j_4^{2}\,j_6 - 168\,j_2^{3}\,j_5\,j_6 + 5040\,j_3^{2}\,j_5\,j_6 + 3348\,j_2\,j_4\,j_5\,j_6 +
        12096\,j_5\,j_6^{2} + 252\,j_2^{3}\,j_4\,j_7 - 7560\,j_3^{2}\,j_4\,j_7 - 2916\,j_2\,j_4^{2}\,j_7 +
        2268\,j_2^{2}\,j_6\,j_7 - 22032\,j_4\,j_6\,j_7 + 34020\,j_3\,j_7^{2} - 68040\,j_3\,j_6\,j_8 -
        34020\,j_2\,j_7\,j_8 - 42\,j_2^{4}\,j_9 + 1260\,j_2\,j_3^{2}\,j_9 + 756\,j_2^{2}\,j_4\,j_9 +
        7776\,j_4^{2}\,j_9 - 34020\,j_3\,j_5\,j_9 - 24192\,j_2\,j_6\,j_9 + 408240\,j_8\,j_9 +
        34020\,j_3\,j_4\,j_{10} + 17010\,j_2\,j_5\,j_{10} - 204120\,j_7\,j_{10},
      \end{dmath}
      \begin{dmath*}
        0=    14\,j_2^{7}\,j_4 - 840\,j_2^{4}\,j_3^{2}\,j_4 + 12600\,j_2\,j_3^{4}\,j_4 - 1260\,j_2^{5}\,j_4^{2} + 
        37800\,j_2^{2}\,j_3^{2}\,j_4^{2} + 25191\,j_2^{3}\,j_4^{3} + 94770\,j_3^{2}\,j_4^{3} -
        743580\,j_2\,j_4^{4} + 11340\,j_2^{3}\,j_3\,j_4\,j_5 - 340200\,j_3^{3}\,j_4\,j_5 -
        510300\,j_2\,j_3\,j_4^{2}\,j_5 + 76545\,j_2^{2}\,j_4\,j_5^{2} + 1771470\,j_4^{2}\,j_5^{2} +
        546\,j_2^{6}\,j_6 - 32760\,j_2^{3}\,j_3^{2}\,j_6 + 491400\,j_3^{4}\,j_6 - 58446\,j_2^{4}\,j_4\,j_6 +
        1753380\,j_2\,j_3^{2}\,j_4\,j_6 + 1524420\,j_2^{2}\,j_4^{2}\,j_6 - 2449440\,j_4^{3}\,j_6 -
        3768930\,j_3\,j_4\,j_5\,j_6 - 2985255\,j_2\,j_5^{2}\,j_6 + 47250\,j_2^{3}\,j_6^{2} -
        1417500\,j_3^{2}\,j_6^{2} + 3418200\,j_2\,j_4\,j_6^{2} + 13608000\,j_6^{3} +
        8726130\,j_3\,j_4^{2}\,j_7 + 51030\,j_2^{3}\,j_5\,j_7 - 1530900\,j_3^{2}\,j_5\,j_7 +
        688905\,j_2\,j_4\,j_5\,j_7 - 22963500\,j_5\,j_6\,j_7 + 221130\,j_2^{3}\,j_4\,j_8 -
        6633900\,j_3^{2}\,j_4\,j_8 - 18676980\,j_2\,j_4^{2}\,j_8 + 20667150\,j_5^{2}\,j_8 +
        61236000\,j_4\,j_6\,j_8 - 34020\,j_2^{3}\,j_3\,j_9 + 1020600\,j_3^{3}\,j_9 +
        1530900\,j_2\,j_3\,j_4\,j_9 - 459270\,j_2^{2}\,j_5\,j_9 + 11022480\,j_4\,j_5\,j_9 -
        30618000\,j_3\,j_6\,j_9 + 34020\,j_2^{4}\,j_{10} - 1020600\,j_2\,j_3^{2}\,j_{10} -
        1530900\,j_2^{2}\,j_4\,j_{10} - 44089920\,j_4^{2}\,j_{10} + 13778100\,j_3\,j_5\,j_{10} +
        30618000\,j_2\,j_6\,j_{10},
      \end{dmath*}
      \begin{dmath*}
        0 = \cdots
      \end{dmath*}
    \end{dgroup*}
    \setcounter{equation}{\themyequation}
  }\medskip

  \noindent Furthermore, a curve $C$ with automorphism group $\DG_4$ is
  $K$-isomorphic to the curve $y^2 = {a_0}\,{x}^{8} + {a_6}\,{x}^{6} +
  {a_4}\,{x}^{4} + {a_2}\,{x}^{2} + {a_0}$ where:
  \begin{itemize}
  \item $a_4$ is the solution of the linear equation
    {\rm\setcounter{myequation}{\theequation}
    \begin{dgroup*}[style={\tiny},spread={-7pt}]
      \begin{dmath}[style={\refstepcounter{myequation}},number={\themyequation},label={eqd4a40}]
        0 = ( -11022480{ j_5}{ j_4}{ j_3}+2933280{ j_4}{{
            j_3}}^{2}{ j_2}+784{{ j_2}}^{6}+55339200{{ j_6}}^{2}+
        4408992{{ j_4}}^{3}+705600{{ j_3}}^{4}-36741600{ j_9}{
          j_3}-4286520{{ j_5}}^{2}{ j_2}+32315760{ j_6}{ j_4}
        { j_2}+36741600{ j_{10}}{ j_2}-64297800{ j_7}{ j_5}+
        10054800{ j_6}{{ j_3}}^{2}-335160{ j_6}{{ j_2}}^{3}-
        47040{{ j_3}}^{2}{{ j_2}}^{3}+2812320{{ j_4}}^{2}{{ j_2}}^
        {2}-97776{ j_4}{{ j_2}}^{4} ) X-410791500{ j_7}{
          j_6}-10716300{ j_9}{{ j_2}}^{2}+4365900{ j_7}{{ j_2
          }}^{3}+150793650{ j_6}{ j_4}{ j_3}+39590775{ j_6}{
          j_5}{ j_2}-196465500{ j_7}{ j_4}{ j_2}+10716300{{
            j_4}}^{2}{ j_3}{ j_2}+1205583750{ j_8}{ j_5}-
        1157360400{ j_9}{ j_4}+254435580{ j_5}{{ j_4}}^{2}+
        321489000{ j_{10}}{ j_3}-75014100{{ j_5}}^{2}{ j_3}-
        130977000{ j_7}{{ j_3}}^{2}+7144200{ j_4}{{ j_3}}^{3}+
        5120010{ j_5}{ j_4}{{ j_2}}^{2}-238140{ j_4}{ j_3}
        {{ j_2}}^{3}-185220{ j_5}{{ j_2}}^{4}+5556600{ j_5}{
          { j_3}}^{2}{ j_2}\condition[]{if this equation is non trivial,}
      \end{dmath}
      \begin{dmath}[style={\refstepcounter{myequation}},number={\themyequation},label={eqd4a41}]
        0 = ( -483840{ j_4}{{ j_3}}^{2}{ j_2}+2449440{ j_5}{
          j_4}{ j_3}-112{{ j_2}}^{6}-9072000{{ j_6}}^{2}-629856
        {{ j_4}}^{3}-100800{{ j_3}}^{4}+612360{{ j_5}}^{2}{ j_2}-
        6713280{ j_6}{ j_4}{ j_2}+9185400{ j_7}{ j_5}-
        2797200{ j_6}{{ j_3}}^{2}+93240{ j_6}{{ j_2}}^{3}+6720
        {{ j_3}}^{2}{{ j_2}}^{3}-498960{{ j_4}}^{2}{{ j_2}}^{2}+
        16128{ j_4}{{ j_2}}^{4}+18370800{ j_8}{ j_4} ) X
        +89302500{ j_7}{ j_6}+26460{ j_5}{{ j_2}}^{4}-2381400
        {{ j_4}}^{2}{ j_3}{ j_2}-22027950{ j_6}{ j_4}{ 
          j_3}-11013975{ j_6}{ j_5}{ j_2}+41079150{ j_7}{ j_4
        }{ j_2}+10716300{{ j_5}}^{2}{ j_3}-241116750{ j_8}{
          j_5}+257191200{ j_9}{ j_4}-56030940{ j_5}{{ j_4}}^{
          2}+23814000{ j_7}{{ j_3}}^{2}-1587600{ j_4}{{ j_3}}^{3
        }-793800{ j_7}{{ j_2}}^{3}-793800{ j_5}{{ j_3}}^{2}{
          j_2}-476280{ j_5}{ j_4}{{ j_2}}^{2}+52920{ j_4}{j_3}{{ j_2}}^{3}
        \condition[]{otherwise;}
    \end{dmath}
    \end{dgroup*}
    \setcounter{equation}{\themyequation}
  }
  \item  $a_0$ is any root of the quadratic equation
    {\rm\setcounter{myequation}{\theequation}
    \begin{dgroup*}[style={\tiny},spread={-7pt}]
      \begin{dmath}[style={\refstepcounter{myequation}},number={\themyequation},label={eqd4a00}]
        0= ( -529079040{ j_5}{ j_4}{ j_3}+1551156480{ j_6}
        { j_4}{ j_2}+140797440{ j_4}{{ j_3}}^{2}{ j_2}-
        2257920{{ j_3}}^{2}{{ j_2}}^{3}-4693248{ j_4}{{ j_2}}^{4
        }-1763596800{ j_9}{ j_3}+1763596800{ j_{10}}{ j_2}+
        2656281600{{ j_6}}^{2}+33868800{{ j_3}}^{4}+37632{{ j_2}}^
        {6}+211631616{{ j_4}}^{3}+134991360{{ j_4}}^{2}{{ j_2}}^{2}-
        16087680{ j_6}{{ j_2}}^{3}-205752960{{ j_5}}^{2}{ j_2}+
        482630400{ j_6}{{ j_3}}^{2}-3086294400{ j_7}{ j_5}
        ) {X}^{2}-14647500{ j_4}{{ j_3}}^{2}{{ j_2}}^{2}-
        307355310{ j_6}{ j_4}{{ j_2}}^{2}-140277690{ j_6}{{
            j_3}}^{2}{ j_2}-382571910{ j_7}{ j_4}{ j_3}+
        1130212440{ j_6}{ j_5}{ j_3}-105019740{ j_8}{ j_4}
        { j_2}+555660{ j_5}{ j_3}{{ j_2}}^{3}+206671500{ 
          j_9}{ j_3}{ j_2}+36996750{ j_5}{ j_4}{ j_3}{ 
          j_2}+3472081200{{ j_7}}^{2}-3724{{ j_2}}^{7}+507952620{ j_7
        }{ j_5}{ j_2}+8024365440{ j_{10}}{ j_4}+530456850{{
            j_5}}^{2}{ j_4}-472252032{ j_6}{{ j_4}}^{2}-4089340080{
          j_9}{ j_5}+39293100{{ j_4}}^{2}{{ j_3}}^{2}-16669800{
          j_5}{{ j_3}}^{3}-400029840{{ j_6}}^{2}{ j_2}-114388848{
          { j_4}}^{3}{ j_2}-3351600{{ j_3}}^{4}{ j_2}-206671500{ 
          j_{10}}{{ j_2}}^{2}+27862380{{ j_5}}^{2}{{ j_2}}^{2}+10835370
        { j_8}{{ j_2}}^{3}-15739920{{ j_4}}^{2}{{ j_2}}^{3}+
        4675923{ j_6}{{ j_2}}^{4}+223440{{ j_3}}^{2}{{ j_2}}^{4}
        +488250{ j_4}{{ j_2}}^{5}-325061100{ j_8}{{ j_3}}^{2}-
        13974055200{ j_8}{ j_6}\condition[]{if this equation is non trivial,}
      \end{dmath}
      \begin{dmath}[style={\refstepcounter{myequation}},number={\themyequation},label={eqd4a01}]
        0= ( 117573120{ j_5}{ j_4}{ j_3}-322237440{ j_6}{
          j_4}{ j_2}-23224320{ j_4}{{ j_3}}^{2}{ j_2}+322560{
          { j_3}}^{2}{{ j_2}}^{3}+774144{ j_4}{{ j_2}}^{4}+881798400
        { j_8}{ j_4}-435456000{{ j_6}}^{2}-4838400{{ j_3}}^{4}
        -5376{{ j_2}}^{6}-30233088{{ j_4}}^{3}-23950080{{ j_4}}^{2
        }{{ j_2}}^{2}+4475520{ j_6}{{ j_2}}^{3}+29393280{{ j_5}}
        ^{2}{ j_2}-134265600{ j_6}{{ j_3}}^{2}+440899200{ j_7}
        { j_5} ) {X}^{2}+151099830{ j_7}{ j_4}{ j_3}+
        2457000{ j_4}{{ j_3}}^{2}{{ j_2}}^{2}+34838370{ j_6}{{
            j_3}}^{2}{ j_2}-257905620{ j_6}{ j_5}{ j_3}-184779630
        { j_8}{ j_4}{ j_2}-79380{ j_5}{ j_3}{{ j_2}}^{
          3}-72564660{ j_7}{ j_5}{ j_2}-10206000{ j_5}{ j_4}
        { j_3}{ j_2}-868020300{{ j_7}}^{2}+532{{ j_2}}^{7}-
        2469035520{ j_{10}}{ j_4}-96446700{{ j_5}}^{2}{ j_4}+
        62215776{ j_6}{{ j_4}}^{2}+1080203040{ j_9}{ j_5}+
        66418380{ j_6}{ j_4}{{ j_2}}^{2}+2381400{ j_5}{{ 
            j_3}}^{3}+20230560{{ j_6}}^{2}{ j_2}+15422724{{ j_4}}^{3}{
          j_2}+478800{{ j_3}}^{4}{ j_2}-3980340{{ j_5}}^{2}{{ j_2
          }}^{2}-1547910{ j_8}{{ j_2}}^{3}+2999430{{ j_4}}^{2}{{ 
            j_2}}^{3}-1161279{ j_6}{{ j_2}}^{4}-31920{{ j_3}}^{2}{{ 
            j_2}}^{4}-81900{ j_4}{{ j_2}}^{5}+46437300{ j_8}{{ j_3}
        }^{2}+3815002800{ j_8}{ j_6}-11736900{{ j_4}}^{2}{{ j_3}
        }^{2}\condition[]{otherwise;}
    \end{dmath}
    \end{dgroup*}
    \setcounter{equation}{\themyequation}
  }
  \item $a_2$ is any root of
    \begin{dgroup*}[style={\tiny},spread={-6pt}]
      \begin{dmath*}
        0 = 15750\,a_0\,X^4 + (105000\,a_0^2\,a_4 + 510\,a_4^3 -
        23100\,a_4\,j_2 - 686000\,j_3)\,X^2 + 705600\,a_0^3\,a_4^2 -
        10804500\,a_0^3\,j_2 + 3024\,a_0\,a_4^4 - 244755\,a_0\,a_4^2\,j_2 -
        1440600\,a_0\,a_4\,j_3 + 15126300\,a_0\,j_4 +
        2881200\,a_0\,j_2^2\condition[]{if $a_0\neq 0$,}
      \end{dmath*}
      \begin{dmath*}
        0 = X-1\condition[]{otherwise;}
      \end{dmath*}
    \end{dgroup*}
  \item and $a_6$ is any root of
    \begin{dgroup*}[style={\tiny},spread={-6pt}]
      \begin{dmath*}
        0 = 5\,a_2\,X + 140\,a_0^2 + a_4^2 - 70\,j_2 \condition[]{if $a_2\neq
          0$,}
      \end{dmath*}
      \begin{dmath*}
        0 = 1575\,a_0\,X^2 - 24\,a_4^3 + 2940\,a_4\,j_2 - 68600\,j_3\condition[]{otherwise.}
      \end{dmath*}
    \end{dgroup*}
  \end{itemize}
\end{lemma}

Similarly to the automorphism group $\CG_2^3$, Lemma~\ref{lemma:C2p3} yields
most of the time a model over an extension of $k$ (of degree at most 8, see
Section~\ref{sec:field-definition} for the rationality issue).
 
\begin{remark}
  In case the $j_i$s correspond to polynomials with multiple roots, it may
  happen that some of equations in Lemma~\ref{lemma:D4} for $a_4$, $a_0$,
  $a_2$ or $a_6$ are inconsistent (typically `$0=1$'). It means that none of the 
  polynomials ${a_0}\,{x}^{8} + {a_6}\,{x}^{6} + {a_4}\,{x}^{4} +
  {a_2}\,{x}^{2} + {a_0}$ have such invariants.  Instead, we construct
  polynomials of the form
  \begin{math}
    f(x) = {a_8}\,{x}^{8} + {a_6}\,{x}^{6} + {a_4}\,{x}^{4} + {a_2}\,{x}^{2}
  \end{math}
  where $a_4$ is the root of
    \begin{dgroup*}[style={\tiny},spread={-4pt}]
      \begin{dmath*}
        0 = ( 11741206932\,{ j_7}+7172848242\,{ j_4}\,{ j_3}-407219400
          \,{j_3}\,{{j_2}}^{2}-340313967\,{j_5}\,{j_2} ) X -
        770798216160\,{j_8}+48652305408\,{j_6}\,{j_2}-21652069824\,{{
            j_4}}^{2}-15830753400\,{j_5}\,{j_3}+1196327160\,{j_4}\,{{
            j_2}}^{2}-77421296\,{{j_2}}^{4}-6442857120\,{{j_3}}^{2}{
          j_2}\condition[]{if this equation is non trivial,}
      \end{dmath*}
      \begin{dmath*}
        0 = ( 488750830848\,{{ j_4}}^{2}+19912489511040\,{ j_8}+
        119825959260\,{{ j_3}}^{2}{ j_2}+272672590848\,{ j_4}\,{{ j_2}
        }^{2}-11523989242\,{{ j_2}}^{4} ) X+663046656007680\,{ j_9}+
        29540232087660\,{ j_7}\,{ j_2}-147027580680960\,{ j_5}\,{ j_4}
        -238603106730240\,{ j_6}\,{ j_3}-32453930490450\,{ j_4}\,{ j_3 }\,{
          j_2}+5621022644895\,{ j_5}\,{{ j_2}}^{2}+590440221000\,{ j_3}\,{{
            j_2}}^{3}j_2\condition[]{otherwise, if this equation is non
          trivial,}
    \end{dmath*}
    \begin{dmath*}
      0 = ( 7023805824\,{ j_6}\,{ j_2}-83616540\,{{ j_3}}^{2}{ j_2
        }+1090412288\,{ j_4}\,{{ j_2}}^{2}+1936166400\,{{ j_4}}^{2}-
        40257382\,{{ j_2}}^{4} ) X+231759118080\,{ j_9}-96001584000
      \,{ j_5}\,{ j_4}-23998355150\,{ j_4}\,{ j_3}\,{ j_2}-
      191099623680\,{ j_6}\,{ j_3}-32526010860\,{ j_7}\,{ j_2}+
      310611000\,{ j_3}\,{{ j_2}}^{3}-617600655\,{ j_5}\,{{ j_2}}^{2}
      \condition[]{otherwise, if this equation is non trivial,}
    \end{dmath*}
    \begin{dmath*}
      0 = X\condition[]{otherwise;}
    \end{dmath*}
    \end{dgroup*}
  and
  \begin{tiny}
    \begin{displaymath}
      a_2 = 1\,,\ 
      a_6 = -\,{( {{{ a_4}}^{2}-70\,{ j_2}})/{{ 5 a_2}}}\,,\ 
      a_8 = -{ {17}}\,{{ a_4}}^{3}/525+{ {22}}\,{ a_4}\,{ j_2
      }/15+{{392}}\,{ j_3}/9\,.
    \end{displaymath}
  \end{tiny}
\end{remark}

\begin{remark} It may happen that a curve $C$ annihilates both
  Eq.~\eqref{eqd4a40} and Eq.~\eqref{eqd4a41} too, so that we can not
  determine $a_4$. But, when it happens, the automorphism group is larger than
  $\DG_4$ (see Table~\ref{tab:d4a4}) and so, lemmas for larger automorphism
  group apply.
  \begin{table}[htbp]
    \centering
    \begin{footnotesize}
      \begin{tabular}{c|c|c}
        Curves $C: y^2={a_0}\,{x}^{8} +
        {a_6}\,{x}^{6} + {a_4}\,{x}^{4} + {a_2}\,{x}^{2} + {a_0}$
        & Singular & $\Aut C$ \\\hline\hline
        $a_0 = 0 $ & Yes & $ \CG_2^3 \subset \Aut C$ \\\hline
        $a_6 + a_2 = 0$ & No & $ \CG_2^3 \subset \Aut C$ \\\hline
        $a_6 - a_2 = 0$ & No & $ \CG_2^3 \subset \Aut C$ \\\hline
        $\left\{\begin{array}{l} %
            308{a_0}{a_2}{a_6}{a_4}-224{a_0}{{a_4}}
            ^{3}+3{{a_2}}^{3}{a_6}+15{{a_2}}^{2}{{a_4}}^{2}+3{
              a_2}{{a_6}}^{3}+15{{a_6}}^{2}{{a_4}}^{2}=0\\
            21{a_0} {{a_2}}^{2}+21{a_0}{{a_6}}^{2}+11{a_2}{a_6}
            {a_4}-8{{a_4}}^{3}=0\\
            196{{a_0}}^{2}-{a_2}{a_6}-5{{a_4}}^{2}=0
          \end{array}\right.$ &
        No & $ \DG_{12} \subset \Aut C$ \\\hline
      \end{tabular}
        \end{footnotesize}\smallskip
    \caption{Automorphism groups of curves $C$ which cancel
      both Eq.~\eqref{eqd4a40} and Eq.~\eqref{eqd4a41}}
    \label{tab:d4a4}
  \end{table}
\end{remark}

\begin{remark} It may happen that a curve $C$ cancels both Eq.~\eqref{eqd4a00}
  and Eq.~\eqref{eqd4a01} too, so that we can not determine $a_0$. But again,
  when it happens, the automorphism group is larger than $\DG_4$ (see
  Table~\ref{tab:d4a0}) and so, lemmas for larger automorphism group apply.
  \begin{table}[htbp]
    \centering
    \begin{footnotesize}
      \begin{tabular}{c|c|c}
        Curves $C: y^2={a_0}\,{x}^{8} +
        {a_6}\,{x}^{6} + {a_4}\,{x}^{4} + {a_2}\,{x}^{2} + {a_0}$
        & Singular & $ \CG_2^3 \subset \Aut C$ \\\hline\hline
        $a_0 = 0 $ & Yes & $ \CG_2^3 \subset \Aut C$ \\\hline
        $a_6 + a_2 = 0$ & No & $ \CG_2^3 \subset \Aut C$ \\\hline
        $a_6 - a_2 = 0$ & No & $ \CG_2^3 \subset \Aut C$ \\\hline
        $\left\{\begin{array}{l} %
            308{a_0}{a_2}{a_6}{a_4}-224{a_0}{{a_4}}
            ^{3}+3{{a_2}}^{3}{a_6}+15{{a_2}}^{2}{{a_4}}^{2}+3{
              a_2}{{a_6}}^{3}+15{{a_6}}^{2}{{a_4}}^{2}=0\\
            21{a_0} {{a_2}}^{2}+21{a_0}{{a_6}}^{2}+11{a_2}{a_6}
            {a_4}-8{{a_4}}^{3}=0\\
            196{{a_0}}^{2}-{a_2}{a_6}-5{{a_4}}^{2}=0
          \end{array}\right.$ &
        No & $ \DG_{12} \subset \Aut C$ \\\hline
        $\left\{\begin{array}{l} %
        14{a_0}+2{a_6}-{a_4} = 0\\
        14{a_0}+2{a_2}-{a_4} = 0
          \end{array}\right.$&
        No& $ \CG_2 \times \DG_8 \subset \Aut C$ \\\hline
        $\left\{\begin{array}{l} %
        14{a_0}-2{a_6}-{a_4} = 0\\
        14{a_0}-2{a_2}-{a_4} = 0
          \end{array}\right.$&
        No& $\CG_2 \times \DG_8 \subset \Aut C$ \\\hline
      \end{tabular}
        \end{footnotesize}\smallskip
    \caption{Automorphism groups of curves $C$ which cancel
      both Eq.~\eqref{eqd4a00} and Eq.~\eqref{eqd4a01}}
    \label{tab:d4a0}
  \end{table}
\end{remark}

Here, the quotient field of $\ii_8$ modulo the ideal defined by
Eq.~\eqref{eq:11} is obtained by adjoining $j_6$ to $k[j_2, j_3, j_4, j_5]$
and $j_6$ satisfies an irreducible monic equation of degree $10$.  The
invariants $j_7$, $j_8$, $j_9$ and $j_{10}$ are then rational in $j_2$, $j_3$,
$j_4$ and $j_5$. 

Geometrically and contrary to the other strata, we did not find any rational
intersection with trivial hyperplane of the projective variety defined by
Eq.~\eqref{eq:11}. For instance, intersecting with $j_2=j_3=0$ is geometrically
birationally equivalent to the elliptic curve $y^2=x^3-x$, of $j$-invariant
1728 (the curve with complex multiplication by $\sqrt{-1}$). Nevertheless, it
is straightforward to adapt the enumeration algorithm at the end of
Section~\ref{sec:fund-invar-covar} to the case of a field of degree 10 over
$k[j_2, j_3, j_4, j_5]$ (instead of a field of degree 5 over $k[j_2, j_3, j_4,
j_5, j_6, j_7]$ for the full $\ii_8$), in order to find non-trivial curves
with automorphism group contained in $\DG_4$.

At last, we emphasize that we have encountered unexpected computational
difficulties to check with \textsc{magma} that a reconstructed model for
$\DG_4$ has the same weighted projective class as the $9$-uple $(j_2 : j_3 :
\ldots : j_{10})$ provided in input. Our first implementations took about one week only to  define the degree 8 extensions needed for the
reconstructed model (due to some unexpected internal Gröbner basis
computations) that we had to design a more optimized version.  We perform by
`hand', \ie directly over the quotient field of $\ii_8$ modulo the ideal
defined by Eq.~\eqref{eq:11}, the calculation of the Shioda invariants of the
reconstructed model. This requires explicit expressions for each
coordinate of the Shioda invariants in the basis $1$, $a_0$, $a_2$,
$a_0\,a_2$, $a_2^2$,\ldots, $a_0a_2^3$.  Moreover, we only check the equality
of $j_2$, \ldots, $j_{6}$ since $j_7$, $j_8$, $j_9$ and $j_{10}$ are rational
in $j_2$, $j_3$, $j_4$, $j_5$ and $j_6$. The corresponding routine took about
half a day on a computer to check finally that everything was fine.

\subsubsection{Generic stratum}
\label{sec:generic-stratum}

Our main tool for reconstructing curves $C$ with automorphism group $\CG_2$ is
the conic and quartic method developed in
Proposition.~\ref{prop:generic-reconstructions}. It requires that $C$ does not cancel
all the 364 fundamental conic determinants $R(q_1, q_2,q_3)$ (or equivalently
Eq.~\eqref{eq:11}), otherwise Algorithm~\ref{algo:reconstruct} fails.  To
check that this cannot happen, we have computed a Gröbner basis (for the
grevlex order $J_{2}< J_3< \ldots <J_{10}$ with weights 2, 3, \ldots, 10) of a
polynomial system which concatenates Relations~\eqref{eq:Syzygies} and
the \emph{fundamental determinants} $R(q_1, q_2,q_3)$ with $\{q_1, q_2,q_3\}$ any subset of size 3 of
the set of 14 covariants of order 2 given in
Tab.~\ref{tab:fundamentals}. From Lemma~\ref{lemma:reconstr-D4}, these $364$
fundamental determinants all vanish for curves with automorphism group
contained in $\DG_4$. Surprisingly, the Gröbner basis that we obtain is equal
to the system of 24 polynomials~\eqref{eq:11} of Lemma~\ref{lemma:D4}.

After some further Gröbner computations, still for the `grevlex' order
$J_{2}< J_3< \ldots <J_{10}$ with weights 2, 3, \ldots, 10, we observe that we
can reduce the set of 364 fundamental conic determinants $R(q_1, q_2,q_3)$ to
a subset of only 19 elements, since the determinants of these 19 conics
generate the same ideal as the one defined by Eq.~\eqref{eq:11}.

\begin{lemma}\label{lemma:C2}
  Let $C$ be a hyperelliptic curve of genus 3 over $k$, let $(j_2: j_3:\ldots:
  j_{10})$ be its Shioda invariants.  If $j_2$, $j_3$, $\ldots$ $j_{10}$ do
  not cancel Eq.~\eqref{eq:11}, nor Eq.~\eqref{eq:10} and Eq.~\eqref{eq:5},
  then its automorphism group is $\CG_2$. Some of the $364$ fundamental
  determinants are non-zero among which at least one of the following $19$
  determinants,
  \begin{displaymath}
    \begin{small}
      \begin{array}{lllll}
        R( C_{5, 2}, C_{6, 2}, C_{7, 2}  ), & 
        R( C_{5, 2}, C_{8,2}, C_{9, 2}), &    
        R( C_{5, 2}, C_{7, 2}, C'_{8,2} ), &  
        R( C_{5, 2}, C_{6, 2}, C''_{9, 2}), & 
        R( C_{5, 2}, C'_{7, 2}, C'_{9, 2}), \\ 
        R( C_{5, 2}, C_{6, 2}, C'_{7, 2} ), & 
        R( C_{5, 2}, C'_{8,2}, C_{11, 2}), & 
        R( C_{5, 2}, C_{7, 2}, C_{11, 2}), & 
        R( C_{5, 2}, C_{7, 2}, C'_{11, 2}), & 
        R( C_{5, 2}, C'_{7, 2}, C'_{10, 2}), \\ 
        R( C_{6, 2}, C'_{7, 2}, C_{9, 2}), & 
        R( C_{5, 2}, C_{7, 2}, C_{10, 2}), & 
        R( C_{5, 2}, C_{7, 2}, C_{9, 2}), & 
        R( C_{5, 2}, C_{6, 2}, C'_{10, 2}), & 
        R( C_{5, 2}, C_{6, 2}, C_{10, 2}), \\ 
        R( C_{5, 2}, C'_{7, 2}, C'_{8,2} ), & 
        R( C_{5, 2}, C_{6, 2}, C'_{9, 2}), & 
        R( C_{5, 2}, C_{6, 2}, C'_{8,2} ), & 
        R( C_{5, 2}, C_{6, 2}, C_{8,2} ).& 
      \end{array}
    \end{small}
  \end{displaymath}
  With one of the corresponding non-singular conics, we can use Proposition.~\ref{prop:generic-reconstructions} to reconstruct the curve over $k$.
\end{lemma}

There are examples of curves with automorphism group $\CG_2$ for which many
determinants $R$ are zero.
\begin{itemize}
\item The curve defined over $\QQ$ by
  \begin{footnotesize}
    \begin{multline*}
      y^2 =
      2581381040\,{x}^{8}+7704021083264\,{x}^{7}+101567018399840\,{x}^{6}
      -96172789044745280\,{x}^{5}-1962596803409291000\,{x}^{4}\\
      -15122894980514300000\,{x}^{3}+266225999081701799000\,{x}^{2}
      +3116782046067938990000\,x+191614097302577354466875
    \end{multline*}
  \end{footnotesize}
  \noindent
  has Shioda invariants in the class $(0:1:0:1:0:0:0:0:0)$ for which the only
  non-zero determinants in the list of Lemma~\ref{lemma:C2} are $R( C_{5, 2},
  C_{6, 2}, C'_{9, 2})$ and $R( C_{5, 2}, C_{6, 2}, C''_{9, 2})$. These two
  discriminants are zero modulo $19$ and $113$ too, but then the automorphism group
  is larger than $\CG_2$.\smallskip
\item Similarly, the singular form $8\,{x}^{5}z^3-125z^8$ has Shioda
  invariants in the class $( 0: 0: 0: 1: 0: 0: 0: 0: 0 )$
  (resp. $16\,{x}^{7}z-35\,{x}^{2}z^6$, in the class $(0: 0: 0: 1: 0: 0: 0: 0:
  3/32)$) for which $18$ determinants in the list are zero in fields of
  characteristic $17$ (resp. characteristic $13$). Note that over $\QQ$, these
  two forms already have $17$ determinants in the list which are zero.
\end{itemize}

\section{Field of definition}
\label{sec:field-definition}

\subsection{Field of definition and field of moduli: general facts} \label{subsec:def}
Let $K$ be an algebraically closed field and let $k \subset K$ be a subfield. Let $C$ be a curve defined over $K$ of genus $g \geq 2$.
\begin{definition}
 We say that $k$ is a \emph{field of definition} of $C$ if there exists a curve $\cc/k$ such that $\cc$ is $K$-isomorphic to $C$. The curve $\cc/k$ is a model of $C$ over $k$.\\
The \emph{field of moduli} of $C$, denoted $\bfM_C$, is the intersection of all fields of definition of $C$.
\end{definition}
This is the definition given by \cite{koizumi}. Usually, one defines the field of moduli of a curve $C$ as the subfield $K^H$ of $K$ fixed by 
$$H=\{ \sigma \in \Aut(K), \; C \simeq {^\sigma C} \}.$$
As  shown in \cite{koizumi} (or \cite[Th.1.5.8]{hugginsphd}), $K^H$ is a purely inseparable extension of $\bfM_C$. Hence, when $\bfM_C$ is a perfect field, both notions coincide.\\
Another common definition is field of moduli relative to a Galois extension $K/k$, as used for instance in \cite{debes-ensalem}. In \cite{hammer}, it is proved that there is always a field of definition for $C$ which is a finite extension of $\bfM_C$ and  \cite[Cor.1.5.9]{hugginsphd} shows that this extension can be chosen separable. Hence, $\bfM_C$ is the field of moduli relative to the Galois extension $\bfM_C^{\textrm{sep}}/\bfM_C$.
\begin{remark}
As a consequence of Torelli theorem (see for instance the formulation of  \cite{matsusaka}), $C$ can be defined over $k$ if and only if the principally polarized abelian variety $\Jac C$ can be. Hence the field of moduli of $C$ coincides with the field of moduli of $\Jac C$. Note however that, when $C$ is non hyperelliptic, a model $A$  of $\Jac C$  over $k$ is not necessarily the Jacobian of a curve  over $k$ (see \cite{LRZ}).
\end{remark}
A classical problem is to know whether $\bfM_C$ is a field of definition or if there is an obstruction. One finds several sufficient conditions in the literature. For instance, the moduli field $\bfM_C$ is a field of definition
\begin{itemize}
\item when $C$ has no automorphism (or more generally when $\Aut(C)$ has no center and has a complement in the automorphism group of $\Aut(C)$ \cite[Cor.4.3]{debes-ensalem}); 
\item when $K$ is the algebraic closure of a finite field  \cite[Cor.2.11]{huggins};
\item when $K \subset \CC$ and $\Jac(C)$ has complex multiplication \cite{milne-moduli,milne-corr}.
\end{itemize}
The difference between a field of moduli and a field of definition boils down to Weil's cocycle relations.
\begin{theorem}[{\cite[proof of Th.1]{wei56}}]
Assume that $K/k$ is a Galois extension. Then $k$ is a field of definition of $C$ if and only if there exists a finite extension $k'$ of $k$ such that for all $\sigma \in \Gal(K/k)$, there exists a $k'$-isomorphism $f_{\sigma} : C \to {^{\sigma} C}$ such that for all $\sigma,\tau \in \Gal(K/k)$, $$f_{\tau}^{\sigma} \circ f_{\sigma}=f_{\tau \sigma}.$$
\end{theorem}
Assume now that $K/\bfM_C$ is a Galois extension and let $\Gamma:=\Gal(K/\bfM_C)$.  By definition of the field of moduli, for all $\sigma \in \Gamma$, there exists,  a $K$-isomorphism $F_{\sigma} : C \to {^\sigma C}$. But only when  $C$ has a  trivial automorphism group, this forces Weil's cocycle relations. Following this remark, we consider the curve $B=C/\Aut(C)$. The isomorphism $F_{\sigma}$ induces an isomorphism $f^{\sigma} : B \to {^\sigma C}/\Aut({^\sigma C})={^{\sigma} B}$ and the diagram
\begin{center}
\leavevmode
\xymatrix{ C \ar[r]^{F_{\sigma}} \ar[d] &  {^\sigma C} \ar[d]  \\ B \ar[r]^{f_{\sigma}}  &  {^\sigma B}
}
\end{center}
is commutative.
The curve $B$ admits a model over $\bfM_C$. Indeed we have to check the cocycle relations  $$\forall \sigma,\tau \in G, \; f_{\tau}^{\sigma} \circ f_{\sigma}=f_{\tau \sigma}.$$ But 
$$F_{\tau \sigma}^{-1} \circ (F_{\tau}^{\sigma} \circ F_{\sigma}) \in \Aut(C)$$
so $f_{\tau \sigma}^{-1} \circ (f_{\tau}^{\sigma} \circ f_{\sigma})=\Id$. Hence there exists a model $\bb/\bfM_C$ and a $K$-isomorphism $\phi : B \to \bb$ such that for all $\sigma \in \Gamma$, $f_{\sigma}= (\phi^{-1})^{\sigma} \circ \phi$. One of the main results of  \cite{debes-ensalem} is the following.
\begin{theorem}[{\cite[Cor.4.3(c)]{debes-ensalem}}] \label{th:debes}
If $\# \Aut(C)$ is coprime to $p$ and $\bb(\bfM_C) \ne \emptyset$ then $\bfM_C$ is a field of definition of $C$.
\end{theorem}
If $\bb$ is a genus $0$ curve, its canonical divisor is of degree $-2$. Hence to get a rational point on $\bb$, it is enough that $\bb$ has a rational divisor of odd degree \cite[Lem.4.0.4]{hugginsphd} (this is the case for instance if $\bb(k') \ne \emptyset$ for $k'/\bfM_C$ an extension of odd degree \cite[Lem.4.0.5]{hugginsphd}). \\

Finally, the field of moduli is related to the moduli space as follows. Assume that the genus $g$ of $C$ is greater than or equal to $2$ and let $\Mm_g$ be the coarse moduli space of curves of genus $g$ viewed as a scheme over the prime field of $K$. The curve $C$ gives a morphism $\Spec K \to \Mm_g$ whose image $[C]$ is a closed point of $\Mm_g$. Let $K([C])$ be the residue field at $[C]$. If $K([C])$ is perfect then  we have the equality $\bfM_C=K([C])$. Note that counterexamples are known if $K([C])$ is not perfect (in that case $\bfM_C$ is a purely inseparable extension of $K([C])$, see \cite{seki}).

\subsection{The hyperelliptic case} \label{subsec:hyp} The problem to know if
the field of moduli is a field of definition has been first addressed by
Mestre in \cite{mestre} under the conditions that the the genus is even. Note
that in this case, he showed \cite[p.322]{mestre}, that if a hyperelliptic
curve $C$ is defined over a field $k$, it has a hyperelliptic equation (in the
sense of Section \ref{subsec:reconstr}) over $k$. However, in general, this is
not true and we will have to distinguish two problems: is the field of moduli
a field of definition ? And if so, has the model over the field of moduli a
hyperelliptic equation ?  This motivates the following terminology.
 
\begin{definition}
  An hyperelliptic curve $C$ over a field $k$ can be \emph{hyperelliptically
    defined} over a subfield $k \subset K$ if there exists a model of $C$ over
  $k$ which has a hyperelliptic equation.
\end{definition}

From now on, let $K$ be an algebraically closed field of characteristic $p \ne
2$. Let $g \geq 2$ be an integer, $n=2g+2$ and let $f \in K[x]$ be a
hyperelliptic polynomial of degree $n$. Let $C$ be the hyperelliptic curve of
genus $g$ defined by the equation $y^2=f(x)$. In the sequel we assume that
$K$ is a Galois extension of $\bfM_C$ and that $K([C])$ (and $\bfM_C$) are
perfect fields (we are typically interested in the cases $p=0$ or $K$ is the
Galois closure of a finite field). We denote  $\Gamma:=\Gal(K/\bfM_C)$. \\

From the form of the isomorphisms between hyperelliptic equations and Weil
cocycle relations, we get the following result.
\begin{lemma}
  $C$ has a model over $k$ if and only if there exists a finite Galois
  extension $k'/k$ such that for all $\sigma \in \Gal(k'/k)$, there exist
  $(M_{\sigma},e_{\sigma}) \in (\GL_2(k'),{k'}^*)$ such that
  \begin{itemize} \label{lem:descent}
  \item $(M_{\sigma},e_{\sigma})$ is an isomorphism between $C$ and ${^{\sigma} C}$; 
  \item for all $\sigma, \tau \in \Gal(k'/k)$, there exists $\lambda \in {k'}^*$ such that 
    $$M_{\tau \sigma}= \lambda M_{\tau}^{\sigma} M_{\sigma}, \quad e_{\tau \sigma}=\lambda^{g+1} e_{\tau}^{\sigma} e_{\sigma}.$$ 
  \end{itemize}
\end{lemma}
If we wish to obtain a model with a hyperelliptic equation, we get the
following version of the lemma.
\begin{lemma} \label{lem:hdescent}
  $C$ can be hyperelliptically defined over $k$ if and only if there exists a
  finite Galois extension $k'/k$ such that for all $\sigma \in \Gal(k'/k)$,
  there exist $M_{\sigma} \in \GL_2(k')$ such that
  \begin{itemize}
  \item there exists $e_{\sigma} \in {k'}^*$ such that
    $(M_{\sigma},e_{\sigma})$ is an isomorphism between $C$ and ${^{\sigma}
      C}$;
  \item for all $\sigma, \tau \in \Gal(k'/k)$,  
    $$M_{\tau \sigma}=  M_{\tau}^{\sigma} M_{\sigma}.$$
  \end{itemize}
\end{lemma}
\begin{proof}
  If $C$ can be hyperelliptically defined over $k$, there exists a curve $\cc
  : y^2=\tilde{f}(x)$ with $\tilde{f} \in k[x]$, a finite extension $k'/k$ and
  $(A,a) \in (\GL_2(k'),{k'}^*)$ such that $(A,a)$ is a $k'$-isomorphism
  between $C$ and $\cc$. Defined then $M_{\sigma}=(A^{\sigma})^{-1} A$ and
  $e_{\sigma}=e/e^{\sigma}$. One can check that the $(M_{\sigma},e_{\sigma})$
  satisfy the hypotheses.\\
  Conversely, we can assume that there is a finite Galois extension $k''/k$
  containing $k'$ and a splitting field of $f$. We denote $r_i$ the roots of
  $f$. Following the proof of Hilbert 90 as in
  \cite[p.159,Prop.3]{serre-local}, we can find a matrix $P \in M_2(k'')$ such
  that the matrix
  $$A=\sum_{\tau \in \Gal(k''/k)}  P^{\tau} M_{\tau}.$$
  is invertible. Now for any lift $\bar{\sigma}$ of an element $\sigma$ in
  $\Gal(k''/k)$ to $\Gamma$, we get
  $$A^{\bar{\sigma}}=\sum_{\tau \in \Gal(k''/k)}  P^{\tau \sigma} M_{\tau}^{\sigma}= \left(\sum_{\tau \in \Gal(k''/k)} P^{\tau \sigma} M_{\tau \sigma}\right) M_{\sigma}^{-1}=A M_{\sigma}^{-1}$$
  since $P$ and the $M_{\tau}$ are defined over $k''$.  Let us denote
  $\tilde{r}_i=A.r_i$. For all $\sigma \in \Gal(k''/k)$, letting
  $\sigma(\infty)=\infty$, we get
  $$\{\tilde{r}_i^{\sigma}\}=\{A^{\sigma}.r_i^{\sigma}\}=\{A
  M_{\sigma}^{-1}.r_i^{\sigma}\}=\{A. r_i\}=\{\tilde{r}_i\}.$$
  Hence the polynomial $\tilde{f}=\prod (x-\tilde{r}_i) \in k[x]$ and the
  curve $\cc/k : y^2=\tilde{f}(x)$ is a model of $C$.
\end{proof}
\begin{remark}
  The previous lemmas give an easy proof that for $g$ even, $C$ can be defined
  over $k$ if and only if $C$ can be hyperelliptically defined over
  $k$. Indeed, assuming the cocycle condition of Lemma \ref{lem:descent} we
  can write
  $$\det M_{\tau \sigma} = \lambda^2 (\det M_{\tau})^{\sigma} (\det M_{\sigma}), \quad e_{\tau \sigma}= \lambda (\lambda^2)^{g/2} e_{\tau}^{\sigma} e_{\sigma}.$$
  Hence we can obtain a good representative satisfying the cocycle condition of
  Lemma \ref{lem:hdescent} letting
  $$M'_{\sigma} = e_{\sigma}^{-1} \cdot (\det M_{\sigma})^{g/2} \cdot M_{\sigma}.$$
\end{remark}

Let now $I_1$ and $I_2$ be two homogeneous invariants of the same degree, of
degree $2g+2$ binary forms, defined over the prime field of $K$ and $I_2(f)
\ne 0$. Their quotient $\iota=I_1/I_2$ satisfies $\iota(f)=\iota(M.f)$ for all
$M \in \GL_2(K)$. The function $\iota$ is therefore a function on the moduli
space of hyperelliptic curves of genus $g$. Hence $\iota(f) \in \bfM_C$. More
specifically, we can prove the following result.
\begin{proposition} \label{prop:invariant-kc} Let $\jj=(I_1 : \ldots : I_m)$
  be a $m$-uple of invariants of degree $d_i$, of degree $n$ binary forms, and
  suppose each $I_i$ is defined over the prime field of $K$. Then, there
  exists $\lambda \in K^*$ such that for $1 \leq i \leq m$,
  $I_i(f)/\lambda^{d_i} \in \bfM_C$.
\end{proposition}
\begin{proof}
  Let $d$ be the gcd of the degree $d_i$ of the invariants $I_i$ which value
  at $f$ is not zero. Then there exists $c_i \in \ZZ$ with $c_i=0$ if
  $I_i(f)=0$ such that $\sum c_i d_i=d$. We then define $I=\prod_i
  I_i^{c_i}$. From what we said above, the value $\mu_i=I_i(f)/I(f)^{d_i/d}
  \in \bfM_C$. If we let $\lambda=I(f)^{1/d}$ for any choice of a $d$th root
  of $I(f)$ we get the result.
\end{proof}
Note that the corresponding representative
$(I_1(f)/\lambda^{d_1},\ldots,I_m(f)/\lambda^{d_m})$ for the m-uple of
invariants is the one constructed in Section \ref{sec:algorithms} and used in
our algorithms.
\begin{corollary}
  With the notation of Section \ref{sec:general-strategy}, there exists
  $\lambda \in K^*$ such that the change of variables $x_i = \lambda^{e_i}
  \cdot x_i'$ for some integers $e_1,e_2,e_3$ makes the conic $\qq$ and the
  curve $\hh$ defined over $\bfM_C$.
\end{corollary}
\begin{proof}
  Thanks to Proposition \ref{prop:invariant-kc} we can find $\lambda \in K^*$
  for the uple of all the $A_{ij}$ and $h_I$. Let $d_1,d_2,d_3$ be the degree
  of $q_1, q_2$ and $q_3$.  From the formula for the coefficients $A_{ij}$ and
  $h_I$, one sees easily that
  $\deg(A_{ij})+\deg(x_i^*)+\deg(x_j^*)=2(d_1+d_2+d_3)$ is a constant $c$ and
  $\deg(h_I)+ \sum_{i \in I} \deg(x_i^*)=1+n/2 \cdot d_1 \cdot d_2 \cdot d_3$
  is a constant $c'$. Hence, if we make the change of variables
  $x_i=\lambda^{e_i} \cdot x_i'$ where $e_i=\deg(x_i^*)$, we get that
$$\sum A_{ij} \lambda^{\deg(x_i^*)+\deg(x_j^*)} x_i' x_j'=\lambda^c \cdot \sum \left[A_{ij}/\lambda^{\deg(A_{ij})}\right] x_i x_j$$
and similarly
$$\sum_I h_I x_I = \lambda^{c'} \cdot \sum \left[h_I/\lambda^{\deg(h_I)}\right] x_I'.$$
This gives the result.
\end{proof}
We then find a refinement of Proposition \ref{prop:generic-reconstructions}.
\begin{proposition} \label{prop:generic-reconstruction-kc} Let $(q_1,q_2,q_3)$
  be three covariants of order $2$ for degree $n$ binary forms and assume they
  are defined over the prime field of $K$. If $R(q_1,q_2,q_3)$ evaluated at a
  form $f$ in non-zero, there exists a non-singular conic $\qq$ and a plane
  curve $\hh$ of degree $n/2$ defined over $\bfM_C$ such that there is a
  $K$-isomorphism $\qq \to \PP^1$ mapping the intersections points of $\qq
  \cap \hh$ on the roots of $f(X,Z)$. In particular the field of moduli is a
  field of definition if $\qq$ has a $\bfM_C$-rational point and in this case
  $C$ can be hyperelliptically defined over $\bfM_C$.
\end{proposition}

\subsection{The hyperelliptic case with no
  extra-automorphism} \label{subsec:hypnoauto}
Besides the hypotheses of Section \ref{subsec:hyp}, we assume that $\Aut(C)=\langle \iota \rangle$. In particular $C/\langle \iota \rangle=B$ is the genus $0$ curve $Q$ defined in Section \ref{subsec:reconstr}.\\
Let $D$ be the image of the ramification divisor $W$ of $C$ by the map $\phi
\circ \rho : C \to B \to \bb$. The divisor $D$ is defined over
$\bfM_C$. Indeed for any $\sigma \in \Gamma$
$${^\sigma D}=\phi^{\sigma} \circ \rho^{\sigma}({^{\sigma} W})=\phi^{\sigma} \circ \rho^{\sigma}({F_{\sigma}}(W))=\phi^{\sigma} \circ f_{\sigma}  \circ \rho(W)=\phi \circ \rho(W)=D.$$ 
Assume now that there exists $\Phi : C \to C'$ an isomorphism onto a model
$C'/\bfM_C$ of $C$. Let $Q'$ be the quotient of $C'$ by the hyperelliptic
involution and $D'$ the image of the ramification divisor $W'$ on $C'$. The
induced $K$-isomorphism $\tilde{\Phi} : \bb \to Q'$ maps the divisor $D$ onto
$D'$. Now for any $\sigma \in \Gamma$, we have $\tilde{\Phi}^{-1} \circ
\tilde{\Phi}^{\sigma}(D)=D$. Since $C$ has no extra-automorphism, this means
that $\tilde{\Phi}^{-1} \circ \tilde{\Phi}^{\sigma}=\Id$ hence $\tilde{\Phi}$
is defined over $\bfM_C$. Hence we get

\begin{proposition} \label{prop:noautomestre} The curve $\bb$ has a rational
  point if and only if $C$ can be hyperelliptically defined over $\bfM_C$.
\end{proposition}

With the notation and hypotheses of Proposition
\ref{prop:generic-reconstruction-kc}, using $Q'=\qq$, we get a sufficient and
necessary condition.
\begin{corollary} \label{cor:noautomestre} Let $(q_1,q_2,q_3)$ be three
  covariants of order $2$ for degree $n$ binary forms and assume they are
  defined over the prime field of $K$. If $R(q_1,q_2,q_3)$ evaluated at a form
  $f$ in non-zero, there exists a non-singular conic $\qq$ and a plane curve
  $\hh$ of degree $n/2$ defined over $\bfM_C$ such that there is a
  $K$-isomorphism $\qq \to \PP^1$ mapping the intersections points of $\qq
  \cap \hh$ on the roots of $f(X,Z)$. In particular the field of moduli is a
  field of definition if and only if $\qq$ has a $\bfM_C$-rational point and
  in this case $C$ can be hyperelliptically defined over $\bfM_C$.
\end{corollary}
Note that the equivalence does not hold if $C$ has extra-automorphism (see Remark \ref{remark:noarithinfo}).\\

If `to be defined' and `to be hyperelliptically defined' over $\bfM_C$ are the
same problem when $g$ is even, it is not the case for odd genus. Indeed,
\begin{proposition} \label{prop:godd} Let $g$ be odd. Then $\bfM_C$ is a field
  of definition for $C$.
\end{proposition}
\begin{proof}
  Consider the rational canonical divisor $\kappa$ on $\bb$ of degree $-2$. It
  is the negative of the intersection of a line with a place model of
  $\bb$. If $\bfM_C$ is finite, we get the result directly from Section
  \ref{subsec:def}. Hence we assume that $\bfM_C$ is infinite. We can find a
  line such that $\textrm{Supp}(-\kappa) \cap \textrm{Supp}(D) = \emptyset$.
  Since $\deg D=2g+2$, $D-2 \cdot \frac{g+1}{2} (-\kappa)$ is a divisor of
  degree $0$, it is the divisor of a function $u \in \bfM_C(\tilde{Q})$. If we
  consider the degree $2$ extension of the form $z^2=u$, this defines a
  hyperelliptic curve over $\bfM_C$ with the same ramification as $C$, hence
  $C$ admits a model over $\bfM_C$ (note that the crucial fact is that each
  point in the divisor $2 \cdot \frac{g+1}{2} (-\kappa)$ has even multiplicity
  and hence does not contribute to the ramification).
\end{proof}

On the other hand, it is easy to exhibit examples for which $\QQ$ is a field
of moduli of a genus $3$ curve $C$ with no extra automorphism but $C$ cannot
be hyperelliptically defined over $\QQ$. Indeed, consider the hyperelliptic
genus $3$ curve over $\QQ$ which is a degree $2$ cover of the conic $Q :
x_1^2+x_2^2+x_3^2=0$ ramified over the intersection points of $Q$ and $$H :
x_2^4+x_1 x_2^3-(x_1^4+x_1^3 x_3+x_1^2 x_3^2+x_1 x_3^3+2 x_3^4)=0.$$ The curve
$C$ has a hyperelliptic equation $$y^2=-x^8 + (2i + 2) x^7 - 8 x^6 + (-2 i -
6) x^5 - 14 x^4 + (-2 i + 6) x^3 - 8 x^2 + (2 i - 2)x - 1$$ over $\QQ(i)$ and
we can check using our programs that $C$ has no extra-automorphism. Since $Q
\simeq \bb$ and $Q(\QQ) =\emptyset$, Proposition \ref{prop:noautomestre} shows
that $C$ cannot be hyperelliptically defined over $\QQ$.

\subsection{The hyperelliptic case with
  extra-automorphisms} \label{subsec:hypwithauto} When $C$ has
extra-automorphisms, the issues are less clear. In \cite{huggins} and
\cite{hugginsphd}, the following results are proved.
\begin{proposition}[{\cite[Th.5.4]{huggins},
    \cite[Prop.4.2.2]{hugginsphd}}] \label{prop:non-cyclic} If
  $\overline{\Aut}(C)$ is not cyclic or is cyclic and of order divisible by
  $p$ then $C$ can be hyperelliptically defined over its field of moduli.
\end{proposition}
\begin{remark} \label{remark:fuertes} In \cite{fuertes-gonzalez} and
  \cite{fuertes}, one can find a family of hyperelliptic curves, among which
  for instance the genus $5$ curve
\begin{small}
  \begin{multline*}
    \hspace*{2cm}C: y^2 = (x^4-2\, (1-2\, \frac{r_3-r_1}{r_3-r_2}\,
    \frac{q_4-r_2}{q_4-r_1})\, x^2+1) \cdot\\ (x^4-2\, (1-2\, \frac{r_3-r_1}{
      r_3-r_2}\, \frac{q_5-r_2}{q_5-r_1})\, x^2+1) \cdot (x^4-2\, (1-2\,
    \frac{r_3-r_1}{r_3-r_2}\, \frac{q_6-r_2}{q_6-r_1})\, x^2+1)\,,
    \end{multline*}
  \end{small}
  with $q_4=3,q_5=-1,q_6=7$ and $r_1,r_2,r_3$ the conjugate roots of $X^3-3
  X+1=0$ (these values of $q_i$ are only to fix ideas). This curve was claimed
  to be a curve with automorphism group $(\ZZ/2\ZZ)^3$, field of moduli $\QQ$
  but with no hyperelliptic equation over $\QQ$, contradicting Huggins'
  general result without the respective authors being aware of it. It took us
  a long month of discussions to realize that there was a subtle gap at page
  406 point d) of \cite{fuertes-gonzalez}: it is claimed there that a certain
  extension is Galois which is not always the case. It took us even longer to
  find an explicit hyperelliptic model of the curve over $\QQ$ because,
  although the curve is defined over a cubic Galois extension of $\QQ$, the
  $\bar{\QQ}$-isomorphism we found between $C$ and a model over $\QQ$ is
  defined over a $12$ degree extension of $\QQ$. A faster and more systematic
  way to derive an hyperelliptic equation over the field of moduli is worked
  out in a work in progress by the authors of the present article.
\end{remark}
We can add the following result which is implicitly contained in
\cite{couveignes-belyi}.
\begin{proposition} \label{prop:signature} If the signature of $C$ has at
  least one odd exponent (with the convention of Table \ref{tab:auto}), then
  the field of moduli is a field of definition.
\end{proposition}
\begin{proof}
  Let $s$ be one of the ramification index with odd signature exponent
  $e$. Consider the corresponding part $Z$ of the ramification divisor
  $W$. Since Galois action respect the ramification index, the image of $Z$ in
  the curve $\bb$ (with the notation of Section \ref{subsec:def}) is a
  rational divisor of odd degree $e$. Hence by the remark following Theorem
  \ref{th:debes}, we get the result.
\end{proof}

\subsection{The hyperelliptic case of genus $3$} \label{subsec:descentg3} Let
now $C$ be a genus $3$ hyperelliptic curve. In this section, we want to
address several issues which will be cut out according to the automorphism
groups reviewed in Table~\ref{tab:auto}.
\begin{description}
\item[Issue \textsf{I}] Can we compute the automorphism group from the invariants~?
\item[Issue \textsf{II}] Is the field of moduli automatically a field of
  definition~?
\item[Issue \textsf{III}] Can the curve be hyperelliptically defined over its
  field of moduli~?
\item[Issue \textsf{IV}] Can we (hyperelliptically) reconstruct the curve from
  its invariants over $K$~?
\item[Issue \textsf{V}] Can we reconstruct the curve over its field of
  moduli~?
\item[Issue \textsf{VI}] Can we hyperelliptically reconstruct the curve over
  its field of moduli~?
\item[Issue \textsf{VII}] What is the number of $\bar{F}_q$-isomorphism
  classes of hyperelliptic genus $3$ curves with given automorphism group over
  a finite field $\FF_{q}$~?
\end{description}
\medskip

The following table gather our state of knowledge (we \textbf{emphasize} what
is proved in the present article). Note that the practical results are valid
for $p=0$ or $p \geq 11$.\medskip
\begin{table}[htbp]
  \centering
  {\small
    \begin{tabular}{|c|c|c|c|c|c|c|}
      \hline
      \# & Dim. $0$ & Dim. $1$ & $\CG_2^3$ & $\CG_4$ & $\DG_4$ & $\CG_2$ \\
      \hline \hline
      \textsf{I} & yes & \textbf{yes} & \textbf{yes} & \textbf{yes} & \textbf{yes} & \textbf{yes} \\
      \hline 
      \textsf{II} & yes & yes & yes & \textbf{yes} &
      \begin{tabular}{c} no\\ {\tiny (counterex. exist)} \end{tabular}
      & \textbf{yes} \\
      \hline
      \textsf{III} & yes & yes & yes & \textbf{yes} & no & \begin{tabular}{c}
        \textbf{computable} \\ {\tiny (possible theoretical obstruction)}  \end{tabular} \\
      \hline
      \textsf{IV} & yes & \textbf{yes} & \textbf{yes} & \textbf{yes} & \textbf{yes} & \textbf{yes} \\
      \hline
      \textsf{V} & yes & \textbf{yes} &
      \begin{tabular}{c}
        no\\
        {\tiny \textbf{(over a cubic ext.)}}
      \end{tabular}
      & \textbf{yes} &
      \begin{tabular}{c}
        no\\
        {\tiny \textbf{(over a deg. $8$ ext.)}}
      \end{tabular}  & \textbf{yes} \\
      \hline
      \textsf{VI} & yes & \textbf{ yes} &
      \begin{tabular}{c}
        no \\
        {\tiny \textbf{(over a cubic ext.)}}
      \end{tabular}&
      \textbf{yes} & 
      \begin{tabular}{c}
        no\\
        {\tiny \textbf{(over a deg. $8$ ext.)}}
      \end{tabular} &\begin{tabular}{c}
        \textbf{yes } \\  {\tiny \textbf{(if no theoretical obstruction)}}
      \end{tabular} \\
      \hline
      \textsf{VII} & {\tiny $1$} & {\tiny $\mathbf{q-3}$} & {\tiny
        $\mathbf{q^2-2\,q+2}$} & {\tiny $\mathbf{q^2-2\,q+2}$} & {\tiny
        \begin{tabular}{c}
          \\
          q$^3-$2\,q$^2+$3\\
          (unproven)
        \end{tabular}
      } & {\tiny
        \begin{tabular}{c}
          \\
          q$^5-$q$^3+$q$-2$\\
          (unproven)
        \end{tabular}
      } \\
      \hline
    \end{tabular} 
  }\medskip
  \caption{Issues addressed in the present paper.\label{tab:issues}
  } 
\end{table}

\begin{remark}
  Similar issues for genus $2$ are completely solved by \cite{CAQU} (see also
  \cite[Th.5]{CANAPU} for $p=2$). In particular if the curve has
  extra-automorphism, then $\bfM_C$ is always a field of definition and
  moreover $C$ can be hyperelliptically defined over $\bfM_C$. Actually, in
  \cite{CAQU} for the case of automorphism group $\DG_4$, the authors assumed
  that $p \ne 2,3,5$. A careful analysis and some improvements made during the
  implementation of their work in MAGMA by the two authors removed this
  restriction (see \cite{LR08}).
\end{remark}
\medskip

\stratum{of dimension $0$} This corresponds to cases 8,9,10 and 11 in Table
\ref{tab:auto}. There is nothing to prove for them.

\stratum{of dimension $1$} By Proposition \ref{prop:non-cyclic}, we have that
cases 5,6 and 7 are hyperelliptically defined over their field of
moduli. Moreover in Section \ref{sec:strata-dimension-1} explicit equations
for their stratum and an explicit hyperelliptic equation over the field of
moduli are given.  \stratum{$\CG_2^3$} By Proposition \ref{prop:non-cyclic},
we know that this case is hyperelliptically defined over its field of
moduli. In Lemma \ref{lemma:C2p3}, explicit equations for the stratum is
given. Moreover an explicit hyperelliptic equation over an extension at most
cubic of the field of moduli is given. To write a hyperelliptic equation over
the field of moduli is a work in progress as mentioned in the introduction.
\stratum{$\CG_4$} In Lemma \ref{lemma:C4} explicit equations for the stratum
are given. Moreover since the signature is $(2^3,4^2)$, the following result
is a consequence of Proposition \ref{prop:signature}.
\begin{proposition} \label{prop:C4descent} Let $C/k$ be a hyperelliptic curve
  of genus $3$ with $\Aut(C) \simeq \CG_4$. Then $C$ can be defined over its
  field of moduli.
\end{proposition}
\begin{remark}
  Proposition \ref{prop:C4descent} has no analogue for even genus greater than
  $3$, since there can be an obstruction as shown in \cite[Chap5]{hugginsphd}.
 \end{remark}

 Thanks to Lemma \ref{lemma:C4}, the procedure is effective when $p \geq 11$
 or $p=0$. Indeed let $\{j_i\}$ be the Shioda invariants of such a curve
 $C/k$. In Lemma \ref{lemma:C4}, it is proved that we can find three order $2$
 covariants $q_i$, such that $R(q_1,q_2,q_3)$ is non-zero. Now, we can use the
 first part of Proposition \ref{prop:generic-reconstruction-kc} to construct a
 conic $\qq$ and an quartic $\hh$ over $\bfM_C$ such that there exists a
 $K$-isomorphism $\qq \to \PP^1$ mapping the intersection divisor of $\qq \cap
 \hh$ on the ramification divisor of $C$. Since $g=3$ is odd, we can moreover
 proceed as in the proof of Proposition \ref{prop:godd} and get a curve over
 $\bfM_C$, $K$-isomorphic to $C$. Hence $C$ can be defined over its field of
 moduli. Moreover we can effectively construct a hyperelliptic equation over
 at most a quadratic extension of the field of moduli. Actually one can do
 better.
 \begin{proposition} \label{prop:C4hypdescent} Let $C/k$ be a hyperelliptic
   curve of genus $3$ with $\Aut(C) \simeq \CG_4$. Then $C$ can be
   hyperelliptically defined over its field of moduli.
\end{proposition}
\begin{proof}
  By Proposition \ref{prop:C4descent}, we know that there is a model $\cc$ of $C$ over $\bfM_C$. Let $\bb=\cc/\langle \iota \rangle$ be seen as a plane non singular conic. The signature of $C$ singles out on $\bb$ a degree $2$ effective rational divisor $Z$. If the points of $Z$ are rational then we are done since then $\bb \simeq \PP^1$. Otherwise, after a quadratic extension $F$ of $\bfM_C$, there is a map $\phi : \bb \to \PP^1$ sending one of this point to $0$ and the other to $\infty$. Since these points were branches of $\cc/\bb$, there exists a $K$-isomorphism $\Phi$ lifting $\phi$ such that the curve $D=\Phi(\cc)$ is of the form $y^2=x g(x)$ with $g \in F[x]$ of degree $6$. Moreover since there is an involution which fixes $0$ and $\infty$ and which is then  $x \mapsto -x$, the polynomial $g$ has only even power monomials. Let $\sigma$ be the Galois involution of $F/\bfM_C$. Since $D$ has a model over $\bfM_C$, there exists a $K$-isomorphism $(M_{\sigma},e_{\sigma}) : D \to {^\sigma D}$ and $M_{\sigma}$ is actually given by $\phi^{\sigma} \circ \phi^{-1}$ and hence is defined over $F$. Since the isomorphism preserves the ramification indexes, $M_{\sigma}$ maps $\{0, \infty\}$ onto $\{0,\infty\}$. It is then easy to see that $M_{\sigma}$ is either the map $x \mapsto ax$ or $x \mapsto a/x$ for some $a \in F$.  By Lemma \ref{lem:descent}, since there exists $\lambda \in F^*$ such that $M_{\sigma}^{\sigma} M_{\sigma}=\lambda \cdot \Id$, we get in the first case $a=\lambda=1$ (and we are done) and $a=a^{\sigma}$ in the later case and so $a \in \bfM_C$.\\
  If $b=\sqrt{a} \notin F$, let $F'=F(b)$ and $\sigma'$ and $\tau$ be the two generators of the $\CG_2 \times \CG_2$ Galois extension $F'/\bfM_C$ defined by $\sigma'_{|F}=\sigma$, $b^{\sigma'}=b$ and $\tau_{|F}=\Id$ and $b^{\tau}=-b$. We define $$M'_{\sigma'}=\frac{1}{b} \cdot M_{\sigma}=\begin{bmatrix} 0 & b \\ \frac{1}{b} & 0 \end{bmatrix} \; \textrm{and} \; M'_{\tau}=\begin{bmatrix} -1 & 0 \\ 0 & 1 \end{bmatrix}.$$ One can check that the conditions of \ref{lem:hdescent} are satisfied.\\
  If $b \in F$, then we define
 $$M'_{\sigma}=\frac{1}{b} \cdot M_{\sigma} \begin{bmatrix} -1 & 0 \\ 0 & 1 \end{bmatrix}=\begin{bmatrix} 0 & b \\ -\frac{1}{b} & 0 \end{bmatrix}.$$
 One can again check that the conditions of Lemma \ref{lem:hdescent} are
 satisfied.
\end{proof}
The previous proof can be turned out into an effective method. Let assume that
among the five determinants of Lemma \ref{lemma:C4},
$R(C_{5,2},C_{6,2},C_{7,2})$ is not zero (the other cases are obtained by
permutation). Then, due to the degree of the covariants and the action of
$\CG_4$, the conic and quartic have the following forms
\begin{eqnarray*}
  \qq & : &  A_{11}  x_1^2 + A_{13}   x_1 x_3 - A_{22}   x_2^2 - A_{33}  x_3^2=0 \\
  \hh & : & x_2 \cdot  (h_{11}   x_1^3 + h_{13}   x_1^2  x_3 + h_{31}  x_3^2  x_1 + h_{33}  x_3^3 + h_{12}  x_1  x_2^2 + h_{32}  x_2^2  x_3)=0.
\end{eqnarray*}
By a linear change of variables in $x_1$ and $x_3$ we can assume that
$A_{13}=0$ and $A_{11}=1$ and we write
$$\qq : x_1^2 - a x_2^2 + c x_3^2=0.$$
The conic and the quartic have the involution $(x_1:x_2:x_3) \mapsto
(x_1:-x_2:x_3)$ which hence stabilizes the ramification divisor of $C$ and is
the involution of $\overline{\Aut}(C)$. The divisor $\{x_2=0\} \cap \qq$ is
the fixed divisor $Z$ of the proof. The points of $Z$ are $(\pm \alpha:0:1)$
where $\alpha$ is a root of $X^2=c$. The isomorphism $\phi^{-1} : \PP^1 \to
\qq$ is given by $$(t:u) \mapsto (\alpha (a t^2 + u^2) : 2 \alpha t u : u^2- a
t^2).$$ From this, we see that $$M_{\sigma}= \phi^{\sigma} \circ \phi^{-1}
= \begin{bmatrix} 0 & a \\ 1 & 0 \end{bmatrix}.$$ We can then apply the
formula in the proof of Lemma \ref{lem:hdescent} to get the result. As
alternative, we have noticed that there is another choice of parametrization
which only requires a quadratic extension and that we implemented. We
normalize
$$\qq : a x_1^2 - x_2^2 + c x_3^2$$
and parametrize the conic $\qq$ with respect to the point $(0:\sqrt{c}:1)$ by
$$ (t:u) \mapsto (2 \sqrt{c} t u : \sqrt{c} (a t^2+u^2) : a t^2-u^2).$$
Plugging these expressions in $\hh$ leads to a polynomial $g \in \bfM_C[x]$
times $\sqrt{c}$, which is what we need.
\begin{remark} \label{remark:noarithinfo} Note that contrary to the case with
  no extra-automorphism, it might happen that none of the conics $\qq$
  contains arithmetic information, \ie that all of them have no points over
  $\bfM_C$. In this way, for instance, we noticed that none of the
  non-degenerate conics among the 19 conics of Lemma~\ref{lemma:C2} have a
  rational point for $(j_2:j_3:\ldots:j_9)$ equal to $(1: 0: 6: 0: 6: 0:
  2963/2835: 0: 2963/2835)$ over $\QQ$.
\end{remark}

\stratum{$\DG_4$} In \cite[Chap.5]{hugginsphd}, Huggins constructs examples of
curves
$$C : y^2=(x^2-a_1)(x^2+\frac{1}{\bar{a_1}})(x^2-a_2)(x^2+\frac{1}{\bar{a_2}})$$
with $a_1,a_2 \in \QQ(i) \setminus (\QQ \cup i \QQ)$, $|a_i| > 1$ and
$|a_1/a_2| \ne 1$, which are defined over $\QQ(i)$ with geometric group of
automorphism $(\ZZ/2\ZZ)^2$. By \cite[Prop.5.0.5]{huggins}, $C$ has field of
moduli $\QQ(i) \cap \RR=\QQ$ but cannot be defined over $\QQ$. Hence the field
of moduli is not always a field of definition and there is not always a
hyperelliptic equation over the field of moduli. However Lemma \ref{lemma:D4}
gives explicit equations for the stratum and shows how to construct a
hyperelliptic equation over an extension of the field of moduli of degree at
most $8$.

\begin{remark}
  In \cite[Cor.2]{GSS}, it is asserted that as soon as $\# \Aut(C)>2$, the
  field of moduli of $C$ is a field of definition. However, the `moduli
  problem' is distinct from the usual one, which explains the contradiction
  with Huggins' result.
\end{remark}

\stratum{$\CG_2$} This is the generic case for which we can apply the results
of Section \ref{subsec:hypnoauto}. Note that as we can find a non singular
conic $\qq$ in this case, we can compute the obstruction using Corollary
\ref{cor:noautomestre}. When the obstruction is trivial, we apply Mestre's
method to obtain a model with a hyperelliptic equation over the field of
moduli.

\subsection{The case of curves over finite fields}
\label{sec:finite-field-case}

\begin{proposition}
  Let $C$ be a (hyperelliptic) curve of genus $g \geq 2$ defined over a finite
  field $k=\FF_q$, such that all the $\bar{k}$-automorphisms of $C$ are
  defined over $k$. Let $k'$ be the extension of $k$ of degree equal to the
  exponent $e$ of $\Aut(C)$. Then there exists a $k'$-isomorphism from $C$ to
  a curve defined over $\bfM_C$.
\end{proposition}
\begin{proof}
  Let $\cc$ be a model of $C$ over $\bfM_C$ and $\Phi : C \to \cc$. We want to
  show that for the Frobenius automorphism $(x \mapsto x^{q^e})=\tau \in
  \Gal(\bar{k'}/k')$ we have ${\Phi^{\tau}}=\Phi$. If $\gamma=x \mapsto x^q$,
  we have that $(\Phi^{\gamma})^{-1} \circ \Phi=f \in \Aut(C)$. Since $\gamma$
  acts trivially on the automorphism group, we get
$$(\Phi^{\tau})^{-1} \circ \Phi=((\Phi^{\gamma^{e}})^{-1} \circ {\Phi^{\gamma^{e-1}}}) \cdots \circ ((\Phi^{\gamma})^{-1} \circ {\Phi^{\gamma}})=f \circ \cdots \circ f=1.$$
\end{proof}

Assume one knows how to compute a $k'$-isomorphism $F_{\sigma} : C \to {^\sigma C}$ for $\sigma$ a generator of $\Gal(k'/\bfM_C)$ and that we know the automorphism of $C$ explicitly. If $C$ is given by a hyperelliptic equation, we can represent all these morphisms by  couple $(M_{\sigma},e_{\sigma})$ and we will actually only care about the matrix $M_{\sigma}$. Two different $F_{\sigma}$ differ by an automorphism of $C$, hence we have a finite explicit set $\ff$ of $k'$-isomorphisms represented by $M_{\sigma}$.\\

Now pick one of the $M_{\sigma}$ in $\ff$. Let $m=\deg(k'/\bfM_C)$ If
$$M_{\sigma}^{\sigma^{m-1}} \cdots M_{\sigma}^{\sigma} \cdot M_{\sigma} \ne \lambda \cdot \Id$$
for $\lambda \in \bfM_C$ then pick another $M_{\sigma} \in \ff$. Since there
exists a $k'$-isomorphism from $C$ to a model over $\bfM_C$, we know that
there is at least one solution to the previous equation and the procedure
produces a good $M_{\sigma}$ and $\lambda$ after a finite number of trials
(less than $\# \overline{\Aut}(C)$). We let $M'_{\sigma}=\frac{1}{a} \cdot
M_{\sigma}$ for $a$ a solution of $\textrm{Norm}_{k'/\bfM_C}(a)=\lambda$
(which always exists in a finite field). We can then apply the formula of
Lemma \ref{lem:hdescent}
$$A = P+ \sum_{i=1}^{m-1} P^{\sigma^i} \cdot M^{\sigma^{i}} \cdots M^{\sigma} \cdot M$$
for a random matrix $P \in M_2(k')$. We repeat this last part until we find $A$ invertible. We can then apply $A$ to the hyperelliptic polynomial of $C$ to get a hyperelliptic equation defined over $\bfM_C$.\\

Practically, in the hyperelliptic case, $M_{\sigma}$ can be computed thanks
for instance to the function \texttt{IsGL2\-Equivalent} of Magma. However,
when $C$ has extra-automorphisms, a better strategy is available. Since the
automorphisms of $C$ are defined over $k'$, we can assume that
$\overline{\Aut}(C)$ is given by one of the canonical expression of
\cite[Lem.2.2.1]{hugginsphd}. Moreover, if $p \nmid \#\Aut(C)$ then by
\cite[Lem.3.1.2]{hugginsphd}, we have that $M_{\sigma}$ belongs to the
normalizer $\nn$ of $\overline{\Aut}(C)$ in $\PGL_2(k')$. Moreover when
$\overline{\Aut}(C)$ is not cyclic, then $\nn$ is finite so we can easily find
the isomorphism.

\begin{remark}
  If $C$ is non hyperelliptic, we can assume that it is canonically
  embedded. Then, isomorphisms and automorphisms are given by $(g+1) \times
  (g+1)$ matrices and the same strategy can be applied.
\end{remark}

   As an application, we use the algorithms developed in this article
   to exhibit a hyperelliptic equation over $\FF_p$, $p=11$, $13$, \ldots,
   $47$, for all the $\bar{\FF}_p$-isomorphism classes of hyperelliptic curves
   of genus $3$. As expected (\textit{cf.} Table~\ref{tab:issues}), we have
   obtained
   \begin{itemize}
   \item one model for each stratum of dimension 0,
   \item $p-3$ models for each stratum of dimension 1,
   \item $p^2-2\,p+2$ models for each stratum of dimension 2,
   \item $p^3-2\,p^2+3$ models for the stratum $\DG_4$,
   \item $p^5-p^3+p-2$ models for the stratum $\CG_2$,
   \end{itemize}
   that is a total number of $p^5$ non-$\bar{\FF}_p$-isomorphic curves as
   predicted in~\cite{BG01}.

\section{Conclusion}
Although a big part of the tasks we addressed in the introduction is now
completed, we wish to pursue our researches in two directions. The first one
is geometric and would be to give a conceptual proof of the existence of a non
singular conic $\qq$ when the curve has no extra-automorphism. We could then
hope to prove such a statement without restricting to the genus $3$ case. The
second one is arithmetic and is related to the explicit hyperelliptic descent
over the field of moduli in the case $\CG_2^3$. Here also, as we saw for
instance with the genus $5$ example of Fuertes and Gonz\'alez-Diez, the
problem is not limited to genus $3$ curves, neither is it to this particular
automorphism group nor to hyperelliptic curves. We will use covariants in a
forthcoming article to get both theoretical and practical results on this
issue in a general context.

\bibliographystyle{alpha}

\newcommand{\etalchar}[1]{$^{#1}$}

\newpage
\appendix

\section{Stratum equations for the automorphism group {\small $\DG_4$}}
\label{sec:strat-equat-autm}

{ \rm\setcounter{myequation}{\theequation}
  \begin{dgroup*}[style={{\tiny}},spread={-7pt}]
    \begin{dmath*}
      0 = 49\,j_2^{4}\,j_4^{2} - 1470\,j_2\,j_3^{2}\,j_4^{2} -
      2205\,j_2^{2}\,j_4^{3} - 17496\,j_4^{4} + 19845\,j_3\,j_4^{2}\,j_5 -
      1092\,j_2^{3}\,j_4\,j_6 + 32760\,j_3^{2}\,j_4\,j_6 +
      47880\,j_2\,j_4^{2}\,j_6 - 238140\,j_5^{2}\,j_6 + 133056\,j_4\,j_6^{2} +
      357210\,j_4\,j_5\,j_7 - 238140\,j_3\,j_6\,j_7 - 59535\,j_2\,j_7^{2} -
      884520\,j_4^{2}\,j_8 + 238140\,j_2\,j_6\,j_8 + 1786050\,j_8^{2} +
      119070\,j_3\,j_4\,j_9 + 59535\,j_2\,j_5\,j_9 - 714420\,j_7\,j_9 +
      4410\,j_2^{3}\,j_{10} - 132300\,j_3^{2}\,j_{10} -
      317520\,j_2\,j_4\,j_{10} + 2540160\,j_6\,j_{10},
    \end{dmath*}
    \begin{dmath*}
      0= -14\,j_2^{3}\,j_3\,j_4^{2} + 420\,j_3^{3}\,j_4^{2} +
      630\,j_2\,j_3\,j_4^{3} - 189\,j_2^{2}\,j_4^{2}\,j_5 +
      360\,j_3\,j_4^{2}\,j_6 - 168\,j_2^{3}\,j_5\,j_6 +
      5040\,j_3^{2}\,j_5\,j_6 + 3348\,j_2\,j_4\,j_5\,j_6 + 12096\,j_5\,j_6^{2}
      + 252\,j_2^{3}\,j_4\,j_7 - 7560\,j_3^{2}\,j_4\,j_7 -
      2916\,j_2\,j_4^{2}\,j_7 + 2268\,j_2^{2}\,j_6\,j_7 - 22032\,j_4\,j_6\,j_7
      + 34020\,j_3\,j_7^{2} - 68040\,j_3\,j_6\,j_8 - 34020\,j_2\,j_7\,j_8 -
      42\,j_2^{4}\,j_9 + 1260\,j_2\,j_3^{2}\,j_9 + 756\,j_2^{2}\,j_4\,j_9 +
      7776\,j_4^{2}\,j_9 - 34020\,j_3\,j_5\,j_9 - 24192\,j_2\,j_6\,j_9 +
      408240\,j_8\,j_9 + 34020\,j_3\,j_4\,j_{10} + 17010\,j_2\,j_5\,j_{10} -
      204120\,j_7\,j_{10},
    \end{dmath*}
    \begin{dmath*}
      0= 14\,j_2^{7}\,j_4 - 840\,j_2^{4}\,j_3^{2}\,j_4 +
      12600\,j_2\,j_3^{4}\,j_4 - 1260\,j_2^{5}\,j_4^{2} +
      37800\,j_2^{2}\,j_3^{2}\,j_4^{2} + 25191\,j_2^{3}\,j_4^{3} +
      94770\,j_3^{2}\,j_4^{3} - 743580\,j_2\,j_4^{4} +
      11340\,j_2^{3}\,j_3\,j_4\,j_5 - 340200\,j_3^{3}\,j_4\,j_5 -
      510300\,j_2\,j_3\,j_4^{2}\,j_5 + 76545\,j_2^{2}\,j_4\,j_5^{2} +
      1771470\,j_4^{2}\,j_5^{2} + 546\,j_2^{6}\,j_6 -
      32760\,j_2^{3}\,j_3^{2}\,j_6 + 491400\,j_3^{4}\,j_6 -
      58446\,j_2^{4}\,j_4\,j_6 + 1753380\,j_2\,j_3^{2}\,j_4\,j_6 +
      1524420\,j_2^{2}\,j_4^{2}\,j_6 - 2449440\,j_4^{3}\,j_6 -
      3768930\,j_3\,j_4\,j_5\,j_6 - 2985255\,j_2\,j_5^{2}\,j_6 +
      47250\,j_2^{3}\,j_6^{2} - 1417500\,j_3^{2}\,j_6^{2} +
      3418200\,j_2\,j_4\,j_6^{2} + 13608000\,j_6^{3} +
      8726130\,j_3\,j_4^{2}\,j_7 + 51030\,j_2^{3}\,j_5\,j_7 -
      1530900\,j_3^{2}\,j_5\,j_7 + 688905\,j_2\,j_4\,j_5\,j_7 -
      22963500\,j_5\,j_6\,j_7 + 221130\,j_2^{3}\,j_4\,j_8 -
      6633900\,j_3^{2}\,j_4\,j_8 - 18676980\,j_2\,j_4^{2}\,j_8 +
      20667150\,j_5^{2}\,j_8 + 61236000\,j_4\,j_6\,j_8 -
      34020\,j_2^{3}\,j_3\,j_9 + 1020600\,j_3^{3}\,j_9 +
      1530900\,j_2\,j_3\,j_4\,j_9 - 459270\,j_2^{2}\,j_5\,j_9 +
      11022480\,j_4\,j_5\,j_9 - 30618000\,j_3\,j_6\,j_9 +
      34020\,j_2^{4}\,j_{10} - 1020600\,j_2\,j_3^{2}\,j_{10} -
      1530900\,j_2^{2}\,j_4\,j_{10} - 44089920\,j_4^{2}\,j_{10} +
      13778100\,j_3\,j_5\,j_{10} + 30618000\,j_2\,j_6\,j_{10},
    \end{dmath*}
    \begin{dmath*}
      0= -1078\,j_2^{7}\,j_4 + 64680\,j_2^{4}\,j_3^{2}\,j_4 -
      970200\,j_2\,j_3^{4}\,j_4 + 83790\,j_2^{5}\,j_4^{2} -
      2513700\,j_2^{2}\,j_3^{2}\,j_4^{2} - 417312\,j_2^{3}\,j_4^{3} -
      35108640\,j_3^{2}\,j_4^{3} - 3061800\,j_2\,j_4^{4} -
      873180\,j_2^{3}\,j_3\,j_4\,j_5 + 26195400\,j_3^{3}\,j_4\,j_5 +
      33934950\,j_2\,j_3\,j_4^{2}\,j_5 - 5893965\,j_2^{2}\,j_4\,j_5^{2} +
      8266860\,j_4^{2}\,j_5^{2} - 18522\,j_2^{6}\,j_6 +
      1111320\,j_2^{3}\,j_3^{2}\,j_6 - 16669800\,j_3^{4}\,j_6 +
      1907262\,j_2^{4}\,j_4\,j_6 - 57217860\,j_2\,j_3^{2}\,j_4\,j_6 -
      47979540\,j_2^{2}\,j_4^{2}\,j_6 - 30093120\,j_4^{3}\,j_6 -
      22095990\,j_3\,j_4\,j_5\,j_6 + 165566835\,j_2\,j_5^{2}\,j_6 +
      476280\,j_2^{3}\,j_6^{2} - 14288400\,j_3^{2}\,j_6^{2} -
      324505440\,j_2\,j_4\,j_6^{2} - 47304810\,j_3\,j_4^{2}\,j_7 -
      3929310\,j_2^{3}\,j_5\,j_7 + 117879300\,j_3^{2}\,j_5\,j_7 +
      43401015\,j_2\,j_4\,j_5\,j_7 + 64297800\,j_2\,j_3\,j_6\,j_7 +
      32148900\,j_5\,j_6\,j_7 + 32148900\,j_2^{2}\,j_7^{2} -
      289340100\,j_4\,j_7^{2} - 2738610\,j_2^{3}\,j_4\,j_8 +
      82158300\,j_3^{2}\,j_4\,j_8 + 409362660\,j_2\,j_4^{2}\,j_8 -
      1591370550\,j_5^{2}\,j_8 - 64297800\,j_2^{2}\,j_6\,j_8 +
      3735396000\,j_4\,j_6\,j_8 - 964467000\,j_3\,j_7\,j_8 +
      1428840\,j_2^{3}\,j_3\,j_9 - 42865200\,j_3^{3}\,j_9 -
      96446700\,j_2\,j_3\,j_4\,j_9 + 3214890\,j_2^{2}\,j_5\,j_9 +
      308629440\,j_4\,j_5\,j_9 + 514382400\,j_3\,j_6\,j_9 -
      192893400\,j_2\,j_7\,j_9 - 2619540\,j_2^{4}\,j_{10} +
      78586200\,j_2\,j_3^{2}\,j_{10} + 150028200\,j_2^{2}\,j_4\,j_{10} -
      1014068160\,j_4^{2}\,j_{10} - 578680200\,j_3\,j_5\,j_{10} -
      1200225600\,j_2\,j_6\,j_{10} + 11573604000\,j_8\,j_{10},
    \end{dmath*}
    \begin{dmath*}
      0= 14\,j_2^{7}\,j_4 - 840\,j_2^{4}\,j_3^{2}\,j_4 +
      12600\,j_2\,j_3^{4}\,j_4 - 1050\,j_2^{5}\,j_4^{2} +
      31500\,j_2^{2}\,j_3^{2}\,j_4^{2} - 5994\,j_2^{3}\,j_4^{3} +
      746820\,j_3^{2}\,j_4^{3} - 126360\,j_2\,j_4^{4} +
      11340\,j_2^{3}\,j_3\,j_4\,j_5 - 340200\,j_3^{3}\,j_4\,j_5 -
      425250\,j_2\,j_3\,j_4^{2}\,j_5 + 76545\,j_2^{2}\,j_4\,j_5^{2} +
      4067820\,j_4^{2}\,j_5^{2} + 546\,j_2^{6}\,j_6 -
      32760\,j_2^{3}\,j_3^{2}\,j_6 + 491400\,j_3^{4}\,j_6 -
      55926\,j_2^{4}\,j_4\,j_6 + 1677780\,j_2\,j_3^{2}\,j_4\,j_6 +
      1405620\,j_2^{2}\,j_4^{2}\,j_6 + 77760\,j_4^{3}\,j_6 -
      852930\,j_3\,j_4\,j_5\,j_6 - 4005855\,j_2\,j_5^{2}\,j_6 +
      117180\,j_2^{3}\,j_6^{2} - 3515400\,j_3^{2}\,j_6^{2} +
      2371680\,j_2\,j_4\,j_6^{2} + 2894130\,j_3\,j_4^{2}\,j_7 +
      51030\,j_2^{3}\,j_5\,j_7 - 1530900\,j_3^{2}\,j_5\,j_7 +
      2219805\,j_2\,j_4\,j_5\,j_7 - 1020600\,j_2\,j_3\,j_6\,j_7 -
      1530900\,j_5\,j_6\,j_7 - 510300\,j_2^{2}\,j_7^{2} -
      4592700\,j_4\,j_7^{2} + 221130\,j_2^{3}\,j_4\,j_8 -
      6633900\,j_3^{2}\,j_4\,j_8 - 16635780\,j_2\,j_4^{2}\,j_8 +
      20667150\,j_5^{2}\,j_8 + 1020600\,j_2^{2}\,j_6\,j_8 +
      12247200\,j_4\,j_6\,j_8 + 15309000\,j_3\,j_7\,j_8 -
      15120\,j_2^{3}\,j_3\,j_9 + 453600\,j_3^{3}\,j_9 +
      1190700\,j_2\,j_3\,j_4\,j_9 + 51030\,j_2^{2}\,j_5\,j_9 -
      44089920\,j_4\,j_5\,j_9 - 13608000\,j_3\,j_6\,j_9 -
      3061800\,j_2\,j_7\,j_9 + 146966400\,j_9^{2} + 34020\,j_2^{4}\,j_{10} -
      1020600\,j_2\,j_3^{2}\,j_{10} - 2041200\,j_2^{2}\,j_4\,j_{10} -
      19595520\,j_4^{2}\,j_{10} + 6123600\,j_3\,j_5\,j_{10} +
      24494400\,j_2\,j_6\,j_{10},
    \end{dmath*}
    \begin{dmath*}
      0= -14\,j_2^{6}\,j_3\,j_4 + 840\,j_2^{3}\,j_3^{3}\,j_4 -
      12600\,j_3^{5}\,j_4 + 1260\,j_2^{4}\,j_3\,j_4^{2} -
      37800\,j_2\,j_3^{3}\,j_4^{2} - 28350\,j_2^{2}\,j_3\,j_4^{3} +
      1180980\,j_3\,j_4^{4} - 378\,j_2^{5}\,j_4\,j_5 +
      11340\,j_2^{2}\,j_3^{2}\,j_4\,j_5 + 41229\,j_2^{3}\,j_4^{2}\,j_5 -
      726570\,j_3^{2}\,j_4^{2}\,j_5 - 1061424\,j_2\,j_4^{3}\,j_5 -
      76545\,j_2\,j_3\,j_4\,j_5^{2} + 2066715\,j_4\,j_5^{3} +
      16866\,j_2^{3}\,j_3\,j_4\,j_6 - 505980\,j_3^{3}\,j_4\,j_6 -
      758970\,j_2\,j_3\,j_4^{2}\,j_6 + 22302\,j_2^{4}\,j_5\,j_6 -
      669060\,j_2\,j_3^{2}\,j_5\,j_6 - 775899\,j_2^{2}\,j_4\,j_5\,j_6 -
      2758536\,j_4^{2}\,j_5\,j_6 + 9032310\,j_3\,j_5^{2}\,j_6 -
      7510320\,j_3\,j_4\,j_6^{2} + 1701000\,j_2\,j_5\,j_6^{2} +
      126\,j_2^{6}\,j_7 - 7560\,j_2^{3}\,j_3^{2}\,j_7 + 113400\,j_3^{4}\,j_7 -
      22491\,j_2^{4}\,j_4\,j_7 + 674730\,j_2\,j_3^{2}\,j_4\,j_7 +
      261954\,j_2^{2}\,j_4^{2}\,j_7 + 629856\,j_4^{3}\,j_7 -
      10563210\,j_3\,j_4\,j_5\,j_7 - 192780\,j_2^{3}\,j_6\,j_7 +
      5783400\,j_3^{2}\,j_6\,j_7 + 3363120\,j_2\,j_4\,j_6\,j_7 +
      6889050\,j_5\,j_7^{2} + 14849730\,j_3\,j_4^{2}\,j_8 +
      51030\,j_2^{3}\,j_5\,j_8 - 1530900\,j_3^{2}\,j_5\,j_8 +
      3750705\,j_2\,j_4\,j_5\,j_8 - 41334300\,j_5\,j_6\,j_8 -
      18370800\,j_4\,j_7\,j_8 + 1134\,j_2^{5}\,j_9 -
      34020\,j_2^{2}\,j_3^{2}\,j_9 - 81648\,j_2^{3}\,j_4\,j_9 +
      918540\,j_3^{2}\,j_4\,j_9 + 4531464\,j_2\,j_4^{2}\,j_9 +
      459270\,j_2\,j_3\,j_5\,j_9 - 19289340\,j_5^{2}\,j_9 +
      1632960\,j_2^{2}\,j_6\,j_9 + 14696640\,j_4\,j_6\,j_9 +
      9185400\,j_3\,j_7\,j_9 - 34020\,j_2^{3}\,j_3\,j_{10} +
      1020600\,j_3^{3}\,j_{10} + 1530900\,j_2\,j_3\,j_4\,j_{10} -
      459270\,j_2^{2}\,j_5\,j_{10} + 29393280\,j_4\,j_5\,j_{10} -
      48988800\,j_3\,j_6\,j_{10} - 9185400\,j_2\,j_7\,j_{10},
    \end{dmath*}
    \begin{dmath*}
      0= -28\,j_2^{6}\,j_3\,j_4 + 1680\,j_2^{3}\,j_3^{3}\,j_4 -
      25200\,j_3^{5}\,j_4 + 2520\,j_2^{4}\,j_3\,j_4^{2} -
      75600\,j_2\,j_3^{3}\,j_4^{2} - 56700\,j_2^{2}\,j_3\,j_4^{3} +
      1771470\,j_3\,j_4^{4} - 756\,j_2^{5}\,j_4\,j_5 +
      22680\,j_2^{2}\,j_3^{2}\,j_4\,j_5 + 76788\,j_2^{3}\,j_4^{2}\,j_5 -
      1283040\,j_3^{2}\,j_4^{2}\,j_5 - 2009853\,j_2\,j_4^{3}\,j_5 -
      153090\,j_2\,j_3\,j_4\,j_5^{2} + 4133430\,j_4\,j_5^{3} +
      18612\,j_2^{3}\,j_3\,j_4\,j_6 - 558360\,j_3^{3}\,j_4\,j_6 -
      837540\,j_2\,j_3\,j_4^{2}\,j_6 + 29484\,j_2^{4}\,j_5\,j_6 -
      884520\,j_2\,j_3^{2}\,j_5\,j_6 - 1075518\,j_2^{2}\,j_4\,j_5\,j_6 -
      6333552\,j_4^{2}\,j_5\,j_6 + 11941020\,j_3\,j_5^{2}\,j_6 -
      11088900\,j_3\,j_4\,j_6^{2} + 1275750\,j_2\,j_5\,j_6^{2} +
      252\,j_2^{6}\,j_7 - 15120\,j_2^{3}\,j_3^{2}\,j_7 + 226800\,j_3^{4}\,j_7
      - 37422\,j_2^{4}\,j_4\,j_7 + 1122660\,j_2\,j_3^{2}\,j_4\,j_7 +
      591948\,j_2^{2}\,j_4^{2}\,j_7 + 1653372\,j_4^{3}\,j_7 -
      11941020\,j_3\,j_4\,j_5\,j_7 - 226800\,j_2^{3}\,j_6\,j_7 +
      6804000\,j_3^{2}\,j_6\,j_7 + 6123600\,j_2\,j_4\,j_6\,j_7 +
      5103000\,j_6^{2}\,j_7 + 17452260\,j_3\,j_4^{2}\,j_8 +
      102060\,j_2^{3}\,j_5\,j_8 - 3061800\,j_3^{2}\,j_5\,j_8 +
      1377810\,j_2\,j_4\,j_5\,j_8 - 45927000\,j_5\,j_6\,j_8 +
      2268\,j_2^{5}\,j_9 - 68040\,j_2^{2}\,j_3^{2}\,j_9 -
      197316\,j_2^{3}\,j_4\,j_9 + 2857680\,j_3^{2}\,j_4\,j_9 +
      7225848\,j_2\,j_4^{2}\,j_9 + 918540\,j_2\,j_3\,j_5\,j_9 -
      24800580\,j_5^{2}\,j_9 + 2041200\,j_2^{2}\,j_6\,j_9 +
      24494400\,j_4\,j_6\,j_9 - 68040\,j_2^{3}\,j_3\,j_{10} +
      2041200\,j_3^{3}\,j_{10} + 3061800\,j_2\,j_3\,j_4\,j_{10} -
      918540\,j_2^{2}\,j_5\,j_{10} + 22044960\,j_4\,j_5\,j_{10} -
      61236000\,j_3\,j_6\,j_{10},
    \end{dmath*}
    \begin{dmath*}
      0= -56\,j_2^{6}\,j_3\,j_4 + 3360\,j_2^{3}\,j_3^{3}\,j_4 -
      50400\,j_3^{5}\,j_4 + 5670\,j_2^{4}\,j_3\,j_4^{2} -
      170100\,j_2\,j_3^{3}\,j_4^{2} - 141750\,j_2^{2}\,j_3\,j_4^{3} +
      1574640\,j_3\,j_4^{4} - 1512\,j_2^{5}\,j_4\,j_5 +
      45360\,j_2^{2}\,j_3^{2}\,j_4\,j_5 + 128061\,j_2^{3}\,j_4^{2}\,j_5 -
      1545480\,j_3^{2}\,j_4^{2}\,j_5 - 2204496\,j_2\,j_4^{3}\,j_5 -
      306180\,j_2\,j_3\,j_4\,j_5^{2} + 8266860\,j_4\,j_5^{3} +
      22104\,j_2^{3}\,j_3\,j_4\,j_6 - 663120\,j_3^{3}\,j_4\,j_6 -
      1010880\,j_2\,j_3\,j_4^{2}\,j_6 + 36288\,j_2^{4}\,j_5\,j_6 -
      1088640\,j_2\,j_3^{2}\,j_5\,j_6 - 1145016\,j_2^{2}\,j_4\,j_5\,j_6 -
      10392624\,j_4^{2}\,j_5\,j_6 + 11634840\,j_3\,j_5^{2}\,j_6 -
      9979200\,j_3\,j_4\,j_6^{2} - 3470040\,j_2\,j_5\,j_6^{2} +
      504\,j_2^{6}\,j_7 - 30240\,j_2^{3}\,j_3^{2}\,j_7 + 453600\,j_3^{4}\,j_7
      - 63504\,j_2^{4}\,j_4\,j_7 + 1905120\,j_2\,j_3^{2}\,j_4\,j_7 +
      702756\,j_2^{2}\,j_4^{2}\,j_7 + 944784\,j_4^{3}\,j_7 -
      5511240\,j_3\,j_4\,j_5\,j_7 - 260820\,j_2^{3}\,j_6\,j_7 +
      4762800\,j_3^{2}\,j_6\,j_7 + 8485560\,j_2\,j_4\,j_6\,j_7 -
      1530900\,j_2\,j_3\,j_7^{2} + 22657320\,j_3\,j_4^{2}\,j_8 +
      204120\,j_2^{3}\,j_5\,j_8 - 6123600\,j_3^{2}\,j_5\,j_8 -
      6429780\,j_2\,j_4\,j_5\,j_8 + 3061800\,j_2\,j_3\,j_6\,j_8 +
      9185400\,j_5\,j_6\,j_8 + 1530900\,j_2^{2}\,j_7\,j_8 -
      27556200\,j_4\,j_7\,j_8 + 6426\,j_2^{5}\,j_9 -
      192780\,j_2^{2}\,j_3^{2}\,j_9 - 428652\,j_2^{3}\,j_4\,j_9 +
      5715360\,j_3^{2}\,j_4\,j_9 + 9552816\,j_2\,j_4^{2}\,j_9 +
      3367980\,j_2\,j_3\,j_5\,j_9 - 49601160\,j_5^{2}\,j_9 +
      3946320\,j_2^{2}\,j_6\,j_9 + 57386880\,j_4\,j_6\,j_9 -
      36741600\,j_3\,j_7\,j_9 - 136080\,j_2^{3}\,j_3\,j_{10} +
      4082400\,j_3^{3}\,j_{10} + 4592700\,j_2\,j_3\,j_4\,j_{10} -
      2602530\,j_2^{2}\,j_5\,j_{10} - 66134880\,j_4\,j_5\,j_{10} -
      85730400\,j_3\,j_6\,j_{10} + 27556200\,j_2\,j_7\,j_{10} +
      440899200\,j_9\,j_{10},
    \end{dmath*}
    \begin{dmath*}
      0= -56\,j_2^{10} + 5040\,j_2^{7}\,j_3^{2} - 151200\,j_2^{4}\,j_3^{4} +
      1512000\,j_2\,j_3^{6} + 7560\,j_2^{8}\,j_4 -
      453600\,j_2^{5}\,j_3^{2}\,j_4 + 6804000\,j_2^{2}\,j_3^{4}\,j_4 -
      138780\,j_2^{6}\,j_4^{2} - 1879200\,j_2^{3}\,j_3^{2}\,j_4^{2} +
      181278000\,j_3^{4}\,j_4^{2} - 12009060\,j_2^{4}\,j_4^{3} +
      513361800\,j_2\,j_3^{2}\,j_4^{3} + 362167200\,j_2^{2}\,j_4^{4} +
      476171136\,j_4^{5} - 68040\,j_2^{6}\,j_3\,j_5 +
      4082400\,j_2^{3}\,j_3^{3}\,j_5 - 61236000\,j_3^{5}\,j_5 +
      6123600\,j_2^{4}\,j_3\,j_4\,j_5 - 183708000\,j_2\,j_3^{3}\,j_4\,j_5 -
      137781000\,j_2^{2}\,j_3\,j_4^{2}\,j_5 + 411374700\,j_3\,j_4^{3}\,j_5 -
      918540\,j_2^{5}\,j_5^{2} + 27556200\,j_2^{2}\,j_3^{2}\,j_5^{2} +
      68890500\,j_2^{3}\,j_4\,j_5^{2} - 826686000\,j_3^{2}\,j_4\,j_5^{2} -
      2341292850\,j_2\,j_4^{2}\,j_5^{2} - 124002900\,j_2\,j_3\,j_5^{3} +
      3348078300\,j_5^{4} - 15750\,j_2^{7}\,j_6 +
      945000\,j_2^{4}\,j_3^{2}\,j_6 - 14175000\,j_2\,j_3^{4}\,j_6 +
      1417500\,j_2^{5}\,j_4\,j_6 - 42525000\,j_2^{2}\,j_3^{2}\,j_4\,j_6 -
      67845600\,j_2^{3}\,j_4^{2}\,j_6 + 1078555500\,j_3^{2}\,j_4^{2}\,j_6 +
      1261461600\,j_2\,j_4^{3}\,j_6 - 12757500\,j_2^{3}\,j_3\,j_5\,j_6 +
      382725000\,j_3^{3}\,j_5\,j_6 + 574087500\,j_2\,j_3\,j_4\,j_5\,j_6 -
      86113125\,j_2^{2}\,j_5^{2}\,j_6 - 4842018000\,j_4\,j_5^{2}\,j_6 -
      14742000\,j_2^{4}\,j_6^{2} + 442260000\,j_2\,j_3^{2}\,j_6^{2} +
      663390000\,j_2^{2}\,j_4\,j_6^{2} + 1539648000\,j_4^{2}\,j_6^{2} -
      5970510000\,j_3\,j_5\,j_6^{2} + 48478500\,j_2^{3}\,j_3\,j_4\,j_7 -
      1454355000\,j_3^{3}\,j_4\,j_7 - 2181532500\,j_2\,j_3\,j_4^{2}\,j_7 +
      22963500\,j_2^{4}\,j_5\,j_7 - 688905000\,j_2\,j_3^{2}\,j_5\,j_7 -
      378897750\,j_2^{2}\,j_4\,j_5\,j_7 - 6200145000\,j_4^{2}\,j_5\,j_7 +
      9300217500\,j_3\,j_5^{2}\,j_7 + 4363065000\,j_3\,j_4\,j_6\,j_7 +
      1033357500\,j_2\,j_5\,j_6\,j_7 + 850500\,j_2^{6}\,j_8 -
      51030000\,j_2^{3}\,j_3^{2}\,j_8 + 765450000\,j_3^{4}\,j_8 -
      125023500\,j_2^{4}\,j_4\,j_8 + 3750705000\,j_2\,j_3^{2}\,j_4\,j_8 +
      3903795000\,j_2^{2}\,j_4^{2}\,j_8 + 23147208000\,j_4^{3}\,j_8 -
      19633792500\,j_3\,j_4\,j_5\,j_8 - 4650108750\,j_2\,j_5^{2}\,j_8 +
      76545000\,j_2^{3}\,j_6\,j_8 - 2296350000\,j_3^{2}\,j_6\,j_8 -
      7807590000\,j_2\,j_4\,j_6\,j_8 - 11573604000\,j_3\,j_4^{2}\,j_9 -
      91854000\,j_2^{3}\,j_5\,j_9 + 2755620000\,j_3^{2}\,j_5\,j_9 +
      826686000\,j_2\,j_4\,j_5\,j_9 + 16533720000\,j_5\,j_6\,j_9 -
      244944000\,j_2^{3}\,j_4\,j_{10} + 7348320000\,j_3^{2}\,j_4\,j_{10} +
      22596084000\,j_2\,j_4^{2}\,j_{10} - 37200870000\,j_5^{2}\,j_{10} -
      22044960000\,j_4\,j_6\,j_{10},
    \end{dmath*}
    \begin{dmath*}
      0= 616\,j_2^{10} - 55440\,j_2^{7}\,j_3^{2} + 1663200\,j_2^{4}\,j_3^{4} -
      16632000\,j_2\,j_3^{6} - 76860\,j_2^{8}\,j_4 +
      4611600\,j_2^{5}\,j_3^{2}\,j_4 - 69174000\,j_2^{2}\,j_3^{4}\,j_4 +
      1205280\,j_2^{6}\,j_4^{2} + 22939200\,j_2^{3}\,j_3^{2}\,j_4^{2} -
      1772928000\,j_3^{4}\,j_4^{2} + 126996660\,j_2^{4}\,j_4^{3} -
      5111164800\,j_2\,j_3^{2}\,j_4^{3} - 4209100200\,j_2^{2}\,j_4^{4} -
      5237882496\,j_4^{5} + 748440\,j_2^{6}\,j_3\,j_5 -
      44906400\,j_2^{3}\,j_3^{3}\,j_5 + 673596000\,j_3^{5}\,j_5 -
      62256600\,j_2^{4}\,j_3\,j_4\,j_5 + 1867698000\,j_2\,j_3^{3}\,j_4\,j_5 +
      1285956000\,j_2^{2}\,j_3\,j_4^{2}\,j_5 - 3186677700\,j_3\,j_4^{3}\,j_5 +
      10103940\,j_2^{5}\,j_5^{2} - 303118200\,j_2^{2}\,j_3^{2}\,j_5^{2} -
      723350250\,j_2^{3}\,j_4\,j_5^{2} + 9093546000\,j_3^{2}\,j_4\,j_5^{2} +
      25724696850\,j_2\,j_4^{2}\,j_5^{2} + 1364031900\,j_2\,j_3\,j_5^{3} -
      36828861300\,j_5^{4} + 570150\,j_2^{7}\,j_6 -
      34209000\,j_2^{4}\,j_3^{2}\,j_6 + 513135000\,j_2\,j_3^{4}\,j_6 -
      58903200\,j_2^{5}\,j_4\,j_6 + 1767096000\,j_2^{2}\,j_3^{2}\,j_4\,j_6 +
      1798151400\,j_2^{3}\,j_4^{2}\,j_6 - 9061834500\,j_3^{2}\,j_4^{2}\,j_6 -
      9782013600\,j_2\,j_4^{3}\,j_6 + 160744500\,j_2^{3}\,j_3\,j_5\,j_6 -
      4822335000\,j_3^{3}\,j_5\,j_6 - 10307331000\,j_2\,j_3\,j_4\,j_5\,j_6 -
      947244375\,j_2^{2}\,j_5^{2}\,j_6 + 18659484000\,j_4\,j_5^{2}\,j_6 +
      151956000\,j_2^{4}\,j_6^{2} - 4558680000\,j_2\,j_3^{2}\,j_6^{2} -
      3458376000\,j_2^{2}\,j_4\,j_6^{2} - 14696640000\,j_4^{2}\,j_6^{2} +
      38578680000\,j_3\,j_5\,j_6^{2} - 543469500\,j_2^{3}\,j_3\,j_4\,j_7 +
      16304085000\,j_3^{3}\,j_4\,j_7 + 31138506000\,j_2\,j_3\,j_4^{2}\,j_7 -
      229635000\,j_2^{4}\,j_5\,j_7 + 6889050000\,j_2\,j_3^{2}\,j_5\,j_7 +
      5717911500\,j_2^{2}\,j_4\,j_5\,j_7 + 103985289000\,j_4^{2}\,j_5\,j_7 -
      102302392500\,j_3\,j_5^{2}\,j_7 - 1082565000\,j_3\,j_4\,j_6\,j_7 -
      27900652500\,j_2\,j_5\,j_6\,j_7 + 688905000\,j_2^{3}\,j_7^{2} -
      20667150000\,j_3^{2}\,j_7^{2} - 22733865000\,j_2\,j_4\,j_7^{2} -
      9355500\,j_2^{6}\,j_8 + 561330000\,j_2^{3}\,j_3^{2}\,j_8 -
      8419950000\,j_3^{4}\,j_8 + 1536003000\,j_2^{4}\,j_4\,j_8 -
      46080090000\,j_2\,j_3^{2}\,j_4\,j_8 - 56857626000\,j_2^{2}\,j_4^{2}\,j_8
      - 300913704000\,j_4^{3}\,j_8 + 199437997500\,j_3\,j_4\,j_5\,j_8 +
      60451413750\,j_2\,j_5^{2}\,j_8 - 841995000\,j_2^{3}\,j_6\,j_8 +
      25259850000\,j_3^{2}\,j_6\,j_8 + 110684070000\,j_2\,j_4\,j_6\,j_8 +
      186004350000\,j_5\,j_7\,j_8 - 15309000\,j_2^{4}\,j_3\,j_9 +
      459270000\,j_2\,j_3^{3}\,j_9 + 688905000\,j_2^{2}\,j_3\,j_4\,j_9 +
      99202320000\,j_3\,j_4^{2}\,j_9 + 482233500\,j_2^{3}\,j_5\,j_9 -
      20667150000\,j_3^{2}\,j_5\,j_9 + 13226976000\,j_2\,j_4\,j_5\,j_9 -
      22044960000\,j_2\,j_3\,j_6\,j_9 - 4133430000\,j_2^{2}\,j_7\,j_9 -
      198404640000\,j_4\,j_7\,j_9 + 15309000\,j_2^{5}\,j_{10} -
      459270000\,j_2^{2}\,j_3^{2}\,j_{10} + 2602530000\,j_2^{3}\,j_4\,j_{10} -
      98743050000\,j_3^{2}\,j_4\,j_{10} - 289891224000\,j_2\,j_4^{2}\,j_{10} +
      6200145000\,j_2\,j_3\,j_5\,j_{10} + 372008700000\,j_5^{2}\,j_{10} +
      22044960000\,j_2^{2}\,j_6\,j_{10} + 308629440000\,j_4\,j_6\,j_{10} +
      124002900000\,j_3\,j_7\,j_{10},
    \end{dmath*}
    \begin{dmath*}
      0= 56\,j_2^{10} - 5040\,j_2^{7}\,j_3^{2} + 151200\,j_2^{4}\,j_3^{4} -
      1512000\,j_2\,j_3^{6} - 4830\,j_2^{8}\,j_4 +
      289800\,j_2^{5}\,j_3^{2}\,j_4 - 4347000\,j_2^{2}\,j_3^{4}\,j_4 -
      91800\,j_2^{6}\,j_4^{2} + 8343000\,j_2^{3}\,j_3^{2}\,j_4^{2} -
      167670000\,j_3^{4}\,j_4^{2} + 14360085\,j_2^{4}\,j_4^{3} -
      418045050\,j_2\,j_3^{2}\,j_4^{3} - 422528400\,j_2^{2}\,j_4^{4} -
      51018336\,j_4^{5} + 68040\,j_2^{6}\,j_3\,j_5 -
      4082400\,j_2^{3}\,j_3^{3}\,j_5 + 61236000\,j_3^{5}\,j_5 -
      3912300\,j_2^{4}\,j_3\,j_4\,j_5 + 117369000\,j_2\,j_3^{3}\,j_4\,j_5 +
      38272500\,j_2^{2}\,j_3\,j_4^{2}\,j_5 - 897544800\,j_3\,j_4^{3}\,j_5 +
      918540\,j_2^{5}\,j_5^{2} - 27556200\,j_2^{2}\,j_3^{2}\,j_5^{2} -
      53964225\,j_2^{3}\,j_4\,j_5^{2} + 826686000\,j_3^{2}\,j_4\,j_5^{2} +
      2604060900\,j_2\,j_4^{2}\,j_5^{2} + 124002900\,j_2\,j_3\,j_5^{3} -
      3348078300\,j_5^{4} + 112140\,j_2^{7}\,j_6 -
      6728400\,j_2^{4}\,j_3^{2}\,j_6 + 100926000\,j_2\,j_3^{4}\,j_6 -
      11907270\,j_2^{5}\,j_4\,j_6 + 357218100\,j_2^{2}\,j_3^{2}\,j_4\,j_6 +
      370978380\,j_2^{3}\,j_4^{2}\,j_6 - 1867041900\,j_3^{2}\,j_4^{2}\,j_6 -
      2435443200\,j_2\,j_4^{3}\,j_6 + 4592700\,j_2^{3}\,j_3\,j_5\,j_6 -
      137781000\,j_3^{3}\,j_5\,j_6 - 941612850\,j_2\,j_3\,j_4\,j_5\,j_6 -
      551124000\,j_2^{2}\,j_5^{2}\,j_6 + 7581891600\,j_4\,j_5^{2}\,j_6 +
      21120750\,j_2^{4}\,j_6^{2} - 633622500\,j_2\,j_3^{2}\,j_6^{2} +
      130734000\,j_2^{2}\,j_4\,j_6^{2} - 4776408000\,j_4^{2}\,j_6^{2} +
      4822335000\,j_3\,j_5\,j_6^{2} - 44396100\,j_2^{3}\,j_3\,j_4\,j_7 +
      1331883000\,j_3^{3}\,j_4\,j_7 + 3699419850\,j_2\,j_3\,j_4^{2}\,j_7 -
      13012650\,j_2^{4}\,j_5\,j_7 + 390379500\,j_2\,j_3^{2}\,j_5\,j_7 +
      568346625\,j_2^{2}\,j_4\,j_5\,j_7 + 744017400\,j_4^{2}\,j_5\,j_7 -
      9300217500\,j_3\,j_5^{2}\,j_7 + 3444525000\,j_3\,j_4\,j_6\,j_7 -
      3444525000\,j_2\,j_5\,j_6\,j_7 - 850500\,j_2^{6}\,j_8 +
      51030000\,j_2^{3}\,j_3^{2}\,j_8 - 765450000\,j_3^{4}\,j_8 +
      164061450\,j_2^{4}\,j_4\,j_8 - 4921843500\,j_2\,j_3^{2}\,j_4\,j_8 -
      7362098100\,j_2^{2}\,j_4^{2}\,j_8 - 5456127600\,j_4^{3}\,j_8 +
      17980420500\,j_3\,j_4\,j_5\,j_8 + 8680203000\,j_2\,j_5^{2}\,j_8 +
      76545000\,j_2^{3}\,j_6\,j_8 - 2296350000\,j_3^{2}\,j_6\,j_8 +
      5051970000\,j_2\,j_4\,j_6\,j_8 + 13778100000\,j_6^{2}\,j_8 -
      6633900\,j_2^{4}\,j_3\,j_9 + 199017000\,j_2\,j_3^{3}\,j_9 +
      298525500\,j_2^{2}\,j_3\,j_4\,j_9 + 7936185600\,j_3\,j_4^{2}\,j_9 -
      6889050\,j_2^{3}\,j_5\,j_9 - 2480058000\,j_3^{2}\,j_5\,j_9 +
      165337200\,j_2\,j_4\,j_5\,j_9 - 5970510000\,j_2\,j_3\,j_6\,j_9 +
      6633900\,j_2^{5}\,j_{10} - 199017000\,j_2^{2}\,j_3^{2}\,j_{10} -
      169929900\,j_2^{3}\,j_4\,j_{10} - 3857868000\,j_3^{2}\,j_4\,j_{10} -
      22320522000\,j_2\,j_4^{2}\,j_{10} + 2686729500\,j_2\,j_3\,j_5\,j_{10} +
      33480783000\,j_5^{2}\,j_{10} + 5970510000\,j_2^{2}\,j_6\,j_{10} -
      33067440000\,j_4\,j_6\,j_{10},
    \end{dmath*}
    \begin{dmath*}
      0= 616\,j_2^{10} - 55440\,j_2^{7}\,j_3^{2} + 1663200\,j_2^{4}\,j_3^{4} -
      16632000\,j_2\,j_3^{6} - 49140\,j_2^{8}\,j_4 +
      2948400\,j_2^{5}\,j_3^{2}\,j_4 - 44226000\,j_2^{2}\,j_3^{4}\,j_4 -
      1546560\,j_2^{6}\,j_4^{2} + 113205600\,j_2^{3}\,j_3^{2}\,j_4^{2} -
      2004264000\,j_3^{4}\,j_4^{2} + 179105580\,j_2^{4}\,j_4^{3} -
      4990442400\,j_2\,j_3^{2}\,j_4^{3} - 4995545400\,j_2^{2}\,j_4^{4} -
      1836660096\,j_4^{5} + 748440\,j_2^{6}\,j_3\,j_5 -
      44906400\,j_2^{3}\,j_3^{3}\,j_5 + 673596000\,j_3^{5}\,j_5 -
      39803400\,j_2^{4}\,j_3\,j_4\,j_5 + 1194102000\,j_2\,j_3^{3}\,j_4\,j_5 +
      275562000\,j_2^{2}\,j_3\,j_4^{2}\,j_5 - 9733243500\,j_3\,j_4^{3}\,j_5 +
      10103940\,j_2^{5}\,j_5^{2} - 303118200\,j_2^{2}\,j_3^{2}\,j_5^{2} -
      571791150\,j_2^{3}\,j_4\,j_5^{2} + 9093546000\,j_3^{2}\,j_4\,j_5^{2} +
      29604216150\,j_2\,j_4^{2}\,j_5^{2} + 1364031900\,j_2\,j_3\,j_5^{3} -
      36828861300\,j_5^{4} + 1318590\,j_2^{7}\,j_6 -
      79115400\,j_2^{4}\,j_3^{2}\,j_6 + 1186731000\,j_2\,j_3^{4}\,j_6 -
      137884680\,j_2^{5}\,j_4\,j_6 + 4136540400\,j_2^{2}\,j_3^{2}\,j_4\,j_6 +
      4083926040\,j_2^{3}\,j_4^{2}\,j_6 - 16477805700\,j_3^{2}\,j_4^{2}\,j_6 -
      20814991200\,j_2\,j_4^{3}\,j_6 + 95426100\,j_2^{3}\,j_3\,j_5\,j_6 -
      2862783000\,j_3^{3}\,j_5\,j_6 - 12074864400\,j_2\,j_3\,j_4\,j_5\,j_6 -
      5921138475\,j_2^{2}\,j_5^{2}\,j_6 + 57112192800\,j_4\,j_5^{2}\,j_6 +
      274428000\,j_2^{4}\,j_6^{2} - 8232840000\,j_2\,j_3^{2}\,j_6^{2} +
      239112000\,j_2^{2}\,j_4\,j_6^{2} - 41290560000\,j_4^{2}\,j_6^{2} +
      78994440000\,j_3\,j_5\,j_6^{2} - 510810300\,j_2^{3}\,j_3\,j_4\,j_7 +
      15324309000\,j_3^{3}\,j_4\,j_7 + 41435339400\,j_2\,j_3\,j_4^{2}\,j_7 -
      128595600\,j_2^{4}\,j_5\,j_7 + 3857868000\,j_2\,j_3^{2}\,j_5\,j_7 +
      4767222600\,j_2^{2}\,j_4\,j_5\,j_7 + 44393038200\,j_4^{2}\,j_5\,j_7 -
      102302392500\,j_3\,j_5^{2}\,j_7 - 21880935000\,j_3\,j_4\,j_6\,j_7 -
      42367657500\,j_2\,j_5\,j_6\,j_7 - 229635000\,j_2^{3}\,j_7^{2} +
      6889050000\,j_3^{2}\,j_7^{2} + 2066715000\,j_2\,j_4\,j_7^{2} +
      41334300000\,j_6\,j_7^{2} - 9355500\,j_2^{6}\,j_8 +
      561330000\,j_2^{3}\,j_3^{2}\,j_8 - 8419950000\,j_3^{4}\,j_8 +
      1839121200\,j_2^{4}\,j_4\,j_8 - 55173636000\,j_2\,j_3^{2}\,j_4\,j_8 -
      82264442400\,j_2^{2}\,j_4^{2}\,j_8 - 159385060800\,j_4^{3}\,j_8 +
      227545321500\,j_3\,j_4\,j_5\,j_8 + 101372370750\,j_2\,j_5^{2}\,j_8 +
      535815000\,j_2^{3}\,j_6\,j_8 - 16074450000\,j_3^{2}\,j_6\,j_8 +
      102417210000\,j_2\,j_4\,j_6\,j_8 - 82668600\,j_2^{4}\,j_3\,j_9 +
      2480058000\,j_2\,j_3^{3}\,j_9 + 3720087000\,j_2^{2}\,j_3\,j_4\,j_9 +
      119704132800\,j_3\,j_4^{2}\,j_9 - 41334300\,j_2^{3}\,j_5\,j_9 -
      32240754000\,j_3^{2}\,j_5\,j_9 + 992023200\,j_2\,j_4\,j_5\,j_9 -
      66134880000\,j_2\,j_3\,j_6\,j_9 + 4133430000\,j_2^{2}\,j_7\,j_9 -
      99202320000\,j_4\,j_7\,j_9 + 82668600\,j_2^{5}\,j_{10} -
      2480058000\,j_2^{2}\,j_3^{2}\,j_{10} - 2431069200\,j_2^{3}\,j_4\,j_{10}
      - 38670534000\,j_3^{2}\,j_4\,j_{10} - 262114574400\,j_2\,j_4^{2}\,j_{10}
      + 33480783000\,j_2\,j_3\,j_5\,j_{10} + 342248004000\,j_5^{2}\,j_{10} +
      66134880000\,j_2^{2}\,j_6\,j_{10} + 44089920000\,j_4\,j_6\,j_{10} -
      124002900000\,j_3\,j_7\,j_{10},
    \end{dmath*}
    \begin{dmath*}
      0= -4312\,j_2^{10} + 388080\,j_2^{7}\,j_3^{2} -
      11642400\,j_2^{4}\,j_3^{4} + 116424000\,j_2\,j_3^{6} +
      861420\,j_2^{8}\,j_4 - 51685200\,j_2^{5}\,j_3^{2}\,j_4 +
      775278000\,j_2^{2}\,j_3^{4}\,j_4 - 39329010\,j_2^{6}\,j_4^{2} +
      849536100\,j_2^{3}\,j_3^{2}\,j_4^{2} + 9910026000\,j_3^{4}\,j_4^{2} -
      199178595\,j_2^{4}\,j_4^{3} + 33391225350\,j_2\,j_3^{2}\,j_4^{3} +
      13993956900\,j_2^{2}\,j_4^{4} - 22856214528\,j_4^{5} -
      5239080\,j_2^{6}\,j_3\,j_5 + 314344800\,j_2^{3}\,j_3^{3}\,j_5 -
      4715172000\,j_3^{5}\,j_5 + 697750200\,j_2^{4}\,j_3\,j_4\,j_5 -
      20932506000\,j_2\,j_3^{3}\,j_4\,j_5 -
      20387760750\,j_2^{2}\,j_3\,j_4^{2}\,j_5 + 57992022900\,j_3\,j_4^{3}\,j_5
      - 70727580\,j_2^{5}\,j_5^{2} + 2121827400\,j_2^{2}\,j_3^{2}\,j_5^{2} +
      6831641250\,j_2^{3}\,j_4\,j_5^{2} - 63654822000\,j_3^{2}\,j_4\,j_5^{2} -
      131195068200\,j_2\,j_4^{2}\,j_5^{2} - 9548223300\,j_2\,j_3\,j_5^{3} +
      257802029100\,j_5^{4} + 4740750\,j_2^{7}\,j_6 -
      284445000\,j_2^{4}\,j_3^{2}\,j_6 + 4266675000\,j_2\,j_3^{4}\,j_6 -
      457531200\,j_2^{5}\,j_4\,j_6 + 13725936000\,j_2^{2}\,j_3^{2}\,j_4\,j_6 +
      5540157000\,j_2^{3}\,j_4^{2}\,j_6 + 162696397500\,j_3^{2}\,j_4^{2}\,j_6
      + 217767463200\,j_2\,j_4^{3}\,j_6 - 696559500\,j_2^{3}\,j_3\,j_5\,j_6 +
      20896785000\,j_3^{3}\,j_5\,j_6 + 27801144000\,j_2\,j_3\,j_4\,j_5\,j_6 -
      40145938875\,j_2^{2}\,j_5^{2}\,j_6 - 939941982000\,j_4\,j_5^{2}\,j_6 -
      358202250\,j_2^{4}\,j_6^{2} + 10746067500\,j_2\,j_3^{2}\,j_6^{2} +
      99499995000\,j_2^{2}\,j_4\,j_6^{2} + 546924960000\,j_4^{2}\,j_6^{2} -
      57868020000\,j_3\,j_5\,j_6^{2} + 3589960500\,j_2^{3}\,j_3\,j_4\,j_7 -
      107698815000\,j_3^{3}\,j_4\,j_7 - 130754169000\,j_2\,j_3\,j_4^{2}\,j_7 +
      2786238000\,j_2^{4}\,j_5\,j_7 - 83587140000\,j_2\,j_3^{2}\,j_5\,j_7 -
      67994923500\,j_2^{2}\,j_4\,j_5\,j_7 + 494771571000\,j_4^{2}\,j_5\,j_7 +
      716116747500\,j_3\,j_5^{2}\,j_7 - 4822335000\,j_2^{2}\,j_3\,j_6\,j_7 -
      1107070335000\,j_3\,j_4\,j_6\,j_7 - 16878172500\,j_2\,j_5\,j_6\,j_7 -
      12055837500\,j_2^{3}\,j_7^{2} + 289340100000\,j_3^{2}\,j_7^{2} +
      65101522500\,j_2\,j_4\,j_7^{2} + 65488500\,j_2^{6}\,j_8 -
      3929310000\,j_2^{3}\,j_3^{2}\,j_8 + 58939650000\,j_3^{4}\,j_8 -
      7215642000\,j_2^{4}\,j_4\,j_8 + 216469260000\,j_2\,j_3^{2}\,j_4\,j_8 +
      143384094000\,j_2^{2}\,j_4^{2}\,j_8 - 1180507608000\,j_4^{3}\,j_8 -
      585913702500\,j_3\,j_4\,j_5\,j_8 + 54251268750\,j_2\,j_5^{2}\,j_8 +
      10716300000\,j_2^{3}\,j_6\,j_8 - 176818950000\,j_3^{2}\,j_6\,j_8 -
      22504230000\,j_2\,j_4\,j_6\,j_8 + 72335025000\,j_2\,j_3\,j_7\,j_8 -
      589396500\,j_2^{4}\,j_3\,j_9 + 17681895000\,j_2\,j_3^{3}\,j_9 +
      28934010000\,j_2^{2}\,j_3\,j_4\,j_9 - 11573604000\,j_3\,j_4^{2}\,j_9 -
      8680203000\,j_2^{3}\,j_5\,j_9 + 57868020000\,j_3^{2}\,j_5\,j_9 -
      23147208000\,j_2\,j_4\,j_5\,j_9 - 183248730000\,j_2\,j_3\,j_6\,j_9 +
      159137055000\,j_2^{2}\,j_7\,j_9 - 1041624360000\,j_4\,j_7\,j_9 +
      678699000\,j_2^{5}\,j_{10} - 20360970000\,j_2^{2}\,j_3^{2}\,j_{10} -
      54170896500\,j_2^{3}\,j_4\,j_{10} + 636548220000\,j_3^{2}\,j_4\,j_{10} +
      1157360400000\,j_2\,j_4^{2}\,j_{10} +
      238705582500\,j_2\,j_3\,j_5\,j_{10} - 4687309620000\,j_5^{2}\,j_{10} +
      234686970000\,j_2^{2}\,j_6\,j_{10} + 9126613440000\,j_4\,j_6\,j_{10} -
      5208121800000\,j_3\,j_7\,j_{10} + 20832487200000\,j_{10}^{2},
    \end{dmath*}
    \begin{dmath*}
      0= 14\,j_2^{7}\,j_3\,j_4 - 840\,j_2^{4}\,j_3^{3}\,j_4 +
      12600\,j_2\,j_3^{5}\,j_4 - 1260\,j_2^{5}\,j_3\,j_4^{2} +
      37800\,j_2^{2}\,j_3^{3}\,j_4^{2} + 22050\,j_2^{3}\,j_3\,j_4^{3} +
      189000\,j_3^{3}\,j_4^{3} - 897480\,j_2\,j_3\,j_4^{4} +
      252\,j_2^{6}\,j_4\,j_5 - 3780\,j_2^{3}\,j_3^{2}\,j_4\,j_5 -
      113400\,j_3^{4}\,j_4\,j_5 - 37125\,j_2^{4}\,j_4^{2}\,j_5 +
      603450\,j_2\,j_3^{2}\,j_4^{2}\,j_5 + 1046844\,j_2^{2}\,j_4^{3}\,j_5 +
      2834352\,j_4^{4}\,j_5 + 76545\,j_2^{2}\,j_3\,j_4\,j_5^{2} -
      2930580\,j_3\,j_4^{2}\,j_5^{2} - 1377810\,j_2\,j_4\,j_5^{3} +
      294\,j_2^{6}\,j_3\,j_6 - 17640\,j_2^{3}\,j_3^{3}\,j_6 +
      264600\,j_3^{5}\,j_6 - 43326\,j_2^{4}\,j_3\,j_4\,j_6 +
      1299780\,j_2\,j_3^{3}\,j_4\,j_6 + 1354320\,j_2^{2}\,j_3\,j_4^{2}\,j_6 +
      162000\,j_3\,j_4^{3}\,j_6 - 14364\,j_2^{5}\,j_5\,j_6 +
      430920\,j_2^{2}\,j_3^{2}\,j_5\,j_6 + 527094\,j_2^{3}\,j_4\,j_5\,j_6 -
      3252150\,j_3^{2}\,j_4\,j_5\,j_6 - 1411344\,j_2\,j_4^{2}\,j_5\,j_6 -
      7424865\,j_2\,j_3\,j_5^{2}\,j_6 + 20207880\,j_5^{3}\,j_6 -
      7560\,j_2^{3}\,j_3\,j_6^{2} + 226800\,j_3^{3}\,j_6^{2} +
      7850520\,j_2\,j_3\,j_4\,j_6^{2} - 1803060\,j_2^{2}\,j_5\,j_6^{2} -
      3110400\,j_4\,j_5\,j_6^{2} - 84\,j_2^{7}\,j_7 +
      5040\,j_2^{4}\,j_3^{2}\,j_7 - 75600\,j_2\,j_3^{4}\,j_7 +
      14742\,j_2^{5}\,j_4\,j_7 - 442260\,j_2^{2}\,j_3^{2}\,j_4\,j_7 +
      648\,j_2^{3}\,j_4^{2}\,j_7 + 31590\,j_3^{2}\,j_4^{2}\,j_7 -
      3254256\,j_2\,j_4^{3}\,j_7 + 17010\,j_2^{3}\,j_3\,j_5\,j_7 -
      510300\,j_3^{3}\,j_5\,j_7 + 8190315\,j_2\,j_3\,j_4\,j_5\,j_7 -
      38578680\,j_4\,j_5^{2}\,j_7 + 128520\,j_2^{4}\,j_6\,j_7 -
      3855600\,j_2\,j_3^{2}\,j_6\,j_7 + 447120\,j_2^{2}\,j_4\,j_6\,j_7 -
      15863040\,j_4^{2}\,j_6\,j_7 + 10716300\,j_3\,j_5\,j_6\,j_7 +
      29087100\,j_3\,j_4\,j_7^{2} + 119070\,j_2^{3}\,j_3\,j_4\,j_8 -
      3572100\,j_3^{3}\,j_4\,j_8 - 20207880\,j_2\,j_3\,j_4^{2}\,j_8 -
      34020\,j_2^{4}\,j_5\,j_8 + 1020600\,j_2\,j_3^{2}\,j_5\,j_8 -
      2908710\,j_2^{2}\,j_4\,j_5\,j_8 + 123084360\,j_4^{2}\,j_5\,j_8 +
      6889050\,j_3\,j_5^{2}\,j_8 - 27556200\,j_3\,j_4\,j_6\,j_8 +
      4592700\,j_2\,j_5\,j_6\,j_8 + 510300\,j_2^{3}\,j_7\,j_8 -
      15309000\,j_3^{2}\,j_7\,j_8 - 33679800\,j_2\,j_4\,j_7\,j_8 -
      378\,j_2^{6}\,j_9 - 11340\,j_2^{3}\,j_3^{2}\,j_9 + 680400\,j_3^{4}\,j_9
      + 6804\,j_2^{4}\,j_4\,j_9 + 1326780\,j_2\,j_3^{2}\,j_4\,j_9 -
      3102624\,j_2^{2}\,j_4^{2}\,j_9 - 1959552\,j_4^{3}\,j_9 -
      459270\,j_2^{2}\,j_3\,j_5\,j_9 - 22657320\,j_3\,j_4\,j_5\,j_9 +
      8266860\,j_2\,j_5^{2}\,j_9 - 1360800\,j_2^{3}\,j_6\,j_9 -
      8164800\,j_3^{2}\,j_6\,j_9 - 12247200\,j_2\,j_4\,j_6\,j_9 -
      9185400\,j_2\,j_3\,j_7\,j_9 + 110224800\,j_5\,j_7\,j_9 +
      34020\,j_2^{4}\,j_3\,j_{10} - 1020600\,j_2\,j_3^{3}\,j_{10} -
      1530900\,j_2^{2}\,j_3\,j_4\,j_{10} + 12247200\,j_3\,j_4^{2}\,j_{10} -
      51030\,j_2^{3}\,j_5\,j_{10} + 15309000\,j_3^{2}\,j_5\,j_{10} +
      13471920\,j_2\,j_4\,j_5\,j_{10} + 48988800\,j_2\,j_3\,j_6\,j_{10} -
      220449600\,j_5\,j_6\,j_{10} + 9185400\,j_2^{2}\,j_7\,j_{10} -
      146966400\,j_4\,j_7\,j_{10},
    \end{dmath*}
    \begin{dmath*}
      0= 350\,j_2^{7}\,j_3\,j_4 - 21000\,j_2^{4}\,j_3^{3}\,j_4 +
      315000\,j_2\,j_3^{5}\,j_4 - 31500\,j_2^{5}\,j_3\,j_4^{2} +
      945000\,j_2^{2}\,j_3^{3}\,j_4^{2} + 587655\,j_2^{3}\,j_3\,j_4^{3} +
      3632850\,j_3^{3}\,j_4^{3} - 16694100\,j_2\,j_3\,j_4^{4} +
      4536\,j_2^{6}\,j_4\,j_5 + 11340\,j_2^{3}\,j_3^{2}\,j_4\,j_5 -
      4422600\,j_3^{4}\,j_4\,j_5 - 487944\,j_2^{4}\,j_4^{2}\,j_5 +
      1880820\,j_2\,j_3^{2}\,j_4^{2}\,j_5 + 12203460\,j_2^{2}\,j_4^{3}\,j_5 +
      11337408\,j_4^{4}\,j_5 + 1913625\,j_2^{2}\,j_3\,j_4\,j_5^{2} +
      12006630\,j_3\,j_4^{2}\,j_5^{2} - 24800580\,j_2\,j_4\,j_5^{3} +
      6426\,j_2^{6}\,j_3\,j_6 - 385560\,j_2^{3}\,j_3^{3}\,j_6 +
      5783400\,j_3^{5}\,j_6 - 810990\,j_2^{4}\,j_3\,j_4\,j_6 +
      24329700\,j_2\,j_3^{3}\,j_4\,j_6 + 23481900\,j_2^{2}\,j_3\,j_4^{2}\,j_6
      + 583200\,j_3\,j_4^{3}\,j_6 - 195048\,j_2^{5}\,j_5\,j_6 +
      5851440\,j_2^{2}\,j_3^{2}\,j_5\,j_6 + 7635060\,j_2^{3}\,j_4\,j_5\,j_6 -
      59960250\,j_3^{2}\,j_4\,j_5\,j_6 - 2099520\,j_2\,j_4^{2}\,j_5\,j_6 -
      114128595\,j_2\,j_3\,j_5^{2}\,j_6 + 99202320\,j_5^{3}\,j_6 +
      425250\,j_2^{3}\,j_3\,j_6^{2} - 12757500\,j_3^{3}\,j_6^{2} +
      119475000\,j_2\,j_3\,j_4\,j_6^{2} - 10206000\,j_2^{2}\,j_5\,j_6^{2} -
      143467200\,j_4\,j_5\,j_6^{2} - 1512\,j_2^{7}\,j_7 +
      90720\,j_2^{4}\,j_3^{2}\,j_7 - 1360800\,j_2\,j_3^{4}\,j_7 +
      233604\,j_2^{5}\,j_4\,j_7 - 7008120\,j_2^{2}\,j_3^{2}\,j_4\,j_7 -
      2980152\,j_2^{3}\,j_4^{2}\,j_7 + 84046410\,j_3^{2}\,j_4^{2}\,j_7 -
      14696640\,j_2\,j_4^{3}\,j_7 + 663390\,j_2^{3}\,j_3\,j_5\,j_7 -
      19901700\,j_3^{3}\,j_5\,j_7 + 84276045\,j_2\,j_3\,j_4\,j_5\,j_7 -
      99202320\,j_4\,j_5^{2}\,j_7 + 1701000\,j_2^{4}\,j_6\,j_7 -
      51030000\,j_2\,j_3^{2}\,j_6\,j_7 - 34700400\,j_2^{2}\,j_4\,j_6\,j_7 -
      22963500\,j_3\,j_5\,j_6\,j_7 + 2602530\,j_2^{3}\,j_3\,j_4\,j_8 -
      78075900\,j_3^{3}\,j_4\,j_8 - 335267100\,j_2\,j_3\,j_4^{2}\,j_8 -
      612360\,j_2^{4}\,j_5\,j_8 + 18370800\,j_2\,j_3^{2}\,j_5\,j_8 -
      11941020\,j_2^{2}\,j_4\,j_5\,j_8 + 363741840\,j_4^{2}\,j_5\,j_8 +
      268672950\,j_3\,j_5^{2}\,j_8 + 275562000\,j_3\,j_4\,j_6\,j_8 +
      137781000\,j_2\,j_5\,j_6\,j_8 - 12096\,j_2^{6}\,j_9 -
      124740\,j_2^{3}\,j_3^{2}\,j_9 + 14628600\,j_3^{4}\,j_9 +
      789264\,j_2^{4}\,j_4\,j_9 + 14594580\,j_2\,j_3^{2}\,j_4\,j_9 -
      36741600\,j_2^{2}\,j_4^{2}\,j_9 - 17635968\,j_4^{3}\,j_9 -
      11481750\,j_2^{2}\,j_3\,j_5\,j_9 - 66134880\,j_3\,j_4\,j_5\,j_9 +
      148803480\,j_2\,j_5^{2}\,j_9 - 17690400\,j_2^{3}\,j_6\,j_9 -
      234738000\,j_3^{2}\,j_6\,j_9 - 24494400\,j_2\,j_4\,j_6\,j_9 +
      489888000\,j_6^{2}\,j_9 + 850500\,j_2^{4}\,j_3\,j_{10} -
      25515000\,j_2\,j_3^{3}\,j_{10} - 38272500\,j_2^{2}\,j_3\,j_4\,j_{10} -
      330674400\,j_3\,j_4^{2}\,j_{10} + 6123600\,j_2^{3}\,j_5\,j_{10} +
      160744500\,j_3^{2}\,j_5\,j_{10} - 55112400\,j_2\,j_4\,j_5\,j_{10} +
      765450000\,j_2\,j_3\,j_6\,j_{10} - 1102248000\,j_5\,j_6\,j_{10},
    \end{dmath*}
    \begin{dmath*}
      0= -2702\,j_2^{7}\,j_3\,j_4 + 162120\,j_2^{4}\,j_3^{3}\,j_4 -
      2431800\,j_2\,j_3^{5}\,j_4 + 243180\,j_2^{5}\,j_3\,j_4^{2} -
      7295400\,j_2^{2}\,j_3^{3}\,j_4^{2} - 5174442\,j_2^{3}\,j_3\,j_4^{3} -
      8913240\,j_3^{3}\,j_4^{3} + 153148320\,j_2\,j_3\,j_4^{4} -
      61236\,j_2^{6}\,j_4\,j_5 + 1485540\,j_2^{3}\,j_3^{2}\,j_4\,j_5 +
      10546200\,j_3^{4}\,j_4\,j_5 + 6241941\,j_2^{4}\,j_4^{2}\,j_5 -
      88770330\,j_2\,j_3^{2}\,j_4^{2}\,j_5 - 162170424\,j_2^{2}\,j_4^{3}\,j_5
      + 8503056\,j_4^{4}\,j_5 - 14773185\,j_2^{2}\,j_3\,j_4\,j_5^{2} -
      28737180\,j_3\,j_4^{2}\,j_5^{2} + 334807830\,j_2\,j_4\,j_5^{3} -
      13230\,j_2^{6}\,j_3\,j_6 + 793800\,j_2^{3}\,j_3^{3}\,j_6 -
      11907000\,j_3^{5}\,j_6 + 2873358\,j_2^{4}\,j_3\,j_4\,j_6 -
      86200740\,j_2\,j_3^{3}\,j_4\,j_6 - 102510360\,j_2^{2}\,j_3\,j_4^{2}\,j_6
      - 7639920\,j_3\,j_4^{3}\,j_6 + 2374596\,j_2^{5}\,j_5\,j_6 -
      71237880\,j_2^{2}\,j_3^{2}\,j_5\,j_6 - 89948394\,j_2^{3}\,j_4\,j_5\,j_6
      + 174223710\,j_3^{2}\,j_4\,j_5\,j_6 - 299461536\,j_2\,j_4^{2}\,j_5\,j_6
      + 1034046405\,j_2\,j_3\,j_5^{2}\,j_6 - 545612760\,j_5^{3}\,j_6 +
      340200\,j_2^{3}\,j_3\,j_6^{2} - 10206000\,j_3^{3}\,j_6^{2} -
      1055899800\,j_2\,j_3\,j_4\,j_6^{2} + 111755700\,j_2^{2}\,j_5\,j_6^{2} +
      545875200\,j_4\,j_5\,j_6^{2} + 20412\,j_2^{7}\,j_7 -
      1224720\,j_2^{4}\,j_3^{2}\,j_7 + 18370800\,j_2\,j_3^{4}\,j_7 -
      3024378\,j_2^{5}\,j_4\,j_7 + 90731340\,j_2^{2}\,j_3^{2}\,j_4\,j_7 +
      48522240\,j_2^{3}\,j_4^{2}\,j_7 - 205074990\,j_3^{2}\,j_4^{2}\,j_7 +
      60466176\,j_2\,j_4^{3}\,j_7 - 1581930\,j_2^{3}\,j_3\,j_5\,j_7 +
      47457900\,j_3^{3}\,j_5\,j_7 - 939896055\,j_2\,j_3\,j_4\,j_5\,j_7 +
      545612760\,j_4\,j_5^{2}\,j_7 - 18370800\,j_2^{4}\,j_6\,j_7 +
      551124000\,j_2\,j_3^{2}\,j_6\,j_7 + 468018000\,j_2^{2}\,j_4\,j_6\,j_7 -
      451396800\,j_4^{2}\,j_6\,j_7 - 390379500\,j_3\,j_5\,j_6\,j_7 +
      344452500\,j_3\,j_4\,j_7^{2} - 5358150\,j_2^{3}\,j_3\,j_4\,j_8 +
      160744500\,j_3^{3}\,j_4\,j_8 + 1833405840\,j_2\,j_3\,j_4^{2}\,j_8 +
      8266860\,j_2^{4}\,j_5\,j_8 - 248005800\,j_2\,j_3^{2}\,j_5\,j_8 +
      85883490\,j_2^{2}\,j_4\,j_5\,j_8 - 1173894120\,j_4^{2}\,j_5\,j_8 -
      640681650\,j_3\,j_5^{2}\,j_8 + 413343000\,j_3\,j_4\,j_6\,j_8 -
      2824510500\,j_2\,j_5\,j_6\,j_8 + 7654500\,j_2^{3}\,j_7\,j_8 -
      229635000\,j_3^{2}\,j_7\,j_8 - 413343000\,j_2\,j_4\,j_7\,j_8 +
      2755620000\,j_6\,j_7\,j_8 + 184842\,j_2^{6}\,j_9 -
      4524660\,j_2^{3}\,j_3^{2}\,j_9 - 30618000\,j_3^{4}\,j_9 -
      15683220\,j_2^{4}\,j_4\,j_9 + 175032900\,j_2\,j_3^{2}\,j_4\,j_9 +
      579537504\,j_2^{2}\,j_4^{2}\,j_9 - 123451776\,j_4^{3}\,j_9 +
      88639110\,j_2^{2}\,j_3\,j_5\,j_9 + 259028280\,j_3\,j_4\,j_5\,j_9 -
      2008846980\,j_2\,j_5^{2}\,j_9 + 175543200\,j_2^{3}\,j_6\,j_9 +
      367416000\,j_3^{2}\,j_6\,j_9 + 1775844000\,j_2\,j_4\,j_6\,j_9 -
      137781000\,j_2\,j_3\,j_7\,j_9 - 6565860\,j_2^{4}\,j_3\,j_{10} +
      196975800\,j_2\,j_3^{3}\,j_{10} + 295463700\,j_2^{2}\,j_3\,j_4\,j_{10} +
      308629440\,j_3\,j_4^{2}\,j_{10} - 81596970\,j_2^{3}\,j_5\,j_{10} -
      211264200\,j_3^{2}\,j_5\,j_{10} + 1499057280\,j_2\,j_4\,j_5\,j_{10} -
      5633712000\,j_2\,j_3\,j_6\,j_{10} + 5511240000\,j_5\,j_6\,j_{10} +
      137781000\,j_2^{2}\,j_7\,j_{10} - 3306744000\,j_4\,j_7\,j_{10},
    \end{dmath*}
    \begin{dmath*}
      0= -3374\,j_2^{7}\,j_3\,j_4 + 202440\,j_2^{4}\,j_3^{3}\,j_4 -
      3036600\,j_2\,j_3^{5}\,j_4 + 303660\,j_2^{5}\,j_3\,j_4^{2} -
      9109800\,j_2^{2}\,j_3^{3}\,j_4^{2} - 6342498\,j_2^{3}\,j_3\,j_4^{3} -
      14695560\,j_3^{3}\,j_4^{3} + 186990120\,j_2\,j_3\,j_4^{4} -
      68796\,j_2^{6}\,j_4\,j_5 + 1394820\,j_2^{3}\,j_3^{2}\,j_4\,j_5 +
      20071800\,j_3^{4}\,j_4\,j_5 + 7003395\,j_2^{4}\,j_4^{2}\,j_5 -
      87119550\,j_2\,j_3^{2}\,j_4^{2}\,j_5 - 180571842\,j_2^{2}\,j_4^{3}\,j_5
      + 30862944\,j_4^{4}\,j_5 - 18447345\,j_2^{2}\,j_3\,j_4\,j_5^{2} -
      80940870\,j_3\,j_4^{2}\,j_5^{2} + 376142130\,j_2\,j_4\,j_5^{3} -
      27510\,j_2^{6}\,j_3\,j_6 + 1650600\,j_2^{3}\,j_3^{3}\,j_6 -
      24759000\,j_3^{5}\,j_6 + 4605246\,j_2^{4}\,j_3\,j_4\,j_6 -
      138157380\,j_2\,j_3^{3}\,j_4\,j_6 -
      151528320\,j_2^{2}\,j_3\,j_4^{2}\,j_6 - 6972480\,j_3\,j_4^{3}\,j_6 +
      2696652\,j_2^{5}\,j_5\,j_6 - 80899560\,j_2^{2}\,j_3^{2}\,j_5\,j_6 -
      102476502\,j_2^{3}\,j_4\,j_5\,j_6 + 296199990\,j_3^{2}\,j_4\,j_5\,j_6 -
      301430808\,j_2\,j_4^{2}\,j_5\,j_6 + 1242554985\,j_2\,j_3\,j_5^{2}\,j_6 -
      464781240\,j_5^{3}\,j_6 - 604800\,j_2^{3}\,j_3\,j_6^{2} +
      18144000\,j_3^{3}\,j_6^{2} - 1279508400\,j_2\,j_3\,j_4\,j_6^{2} +
      129616200\,j_2^{2}\,j_5\,j_6^{2} + 524880000\,j_4\,j_5\,j_6^{2} +
      22932\,j_2^{7}\,j_7 - 1375920\,j_2^{4}\,j_3^{2}\,j_7 +
      20638800\,j_2\,j_3^{4}\,j_7 - 3412206\,j_2^{5}\,j_4\,j_7 +
      102366180\,j_2^{2}\,j_3^{2}\,j_4\,j_7 + 53584416\,j_2^{3}\,j_4^{2}\,j_7
      - 405316710\,j_3^{2}\,j_4^{2}\,j_7 + 83800008\,j_2\,j_4^{3}\,j_7 -
      3010770\,j_2^{3}\,j_3\,j_5\,j_7 + 90323100\,j_3^{3}\,j_5\,j_7 -
      1084106835\,j_2\,j_3\,j_4\,j_5\,j_7 + 349963740\,j_4\,j_5^{2}\,j_7 -
      20865600\,j_2^{4}\,j_6\,j_7 + 625968000\,j_2\,j_3^{2}\,j_6\,j_7 +
      513434700\,j_2^{2}\,j_4\,j_6\,j_7 - 451396800\,j_4^{2}\,j_6\,j_7 +
      68890500\,j_3\,j_5\,j_6\,j_7 + 229635000\,j_3\,j_4\,j_7^{2} +
      688905000\,j_7^{3} - 11141550\,j_2^{3}\,j_3\,j_4\,j_8 +
      334246500\,j_3^{3}\,j_4\,j_8 + 2512513080\,j_2\,j_3\,j_4^{2}\,j_8 +
      9287460\,j_2^{4}\,j_5\,j_8 - 278623800\,j_2\,j_3^{2}\,j_5\,j_8 +
      105172830\,j_2^{2}\,j_4\,j_5\,j_8 - 1260236880\,j_4^{2}\,j_5\,j_8 -
      1219361850\,j_3\,j_5^{2}\,j_8 - 520506000\,j_3\,j_4\,j_6\,j_8 -
      3191926500\,j_2\,j_5\,j_6\,j_8 + 7654500\,j_2^{3}\,j_7\,j_8 -
      229635000\,j_3^{2}\,j_7\,j_8 - 298525500\,j_2\,j_4\,j_7\,j_8 +
      205254\,j_2^{6}\,j_9 - 4116420\,j_2^{3}\,j_3^{2}\,j_9 -
      61236000\,j_3^{4}\,j_9 - 17064432\,j_2^{4}\,j_4\,j_9 +
      142986060\,j_2\,j_3^{2}\,j_4\,j_9 + 641875752\,j_2^{2}\,j_4^{2}\,j_9 -
      219469824\,j_4^{3}\,j_9 + 110684070\,j_2^{2}\,j_3\,j_5\,j_9 +
      264539520\,j_3\,j_4\,j_5\,j_9 - 2256852780\,j_2\,j_5^{2}\,j_9 +
      204120000\,j_2^{3}\,j_6\,j_9 + 979776000\,j_3^{2}\,j_6\,j_9 +
      1837080000\,j_2\,j_4\,j_6\,j_9 - 137781000\,j_2\,j_3\,j_7\,j_9 -
      8198820\,j_2^{4}\,j_3\,j_{10} + 245964600\,j_2\,j_3^{3}\,j_{10} +
      368946900\,j_2^{2}\,j_3\,j_4\,j_{10} + 1179405360\,j_3\,j_4^{2}\,j_{10}
      - 91700910\,j_2^{3}\,j_5\,j_{10} - 569494800\,j_3^{2}\,j_5\,j_{10} +
      1603770840\,j_2\,j_4\,j_5\,j_{10} - 7103376000\,j_2\,j_3\,j_6\,j_{10} +
      4408992000\,j_5\,j_6\,j_{10} + 137781000\,j_2^{2}\,j_7\,j_{10} -
      551124000\,j_4\,j_7\,j_{10},
    \end{dmath*}
    \begin{dmath*}
      0= -392\,j_2^{11} + 35280\,j_2^{8}\,j_3^{2} - 1058400\,j_2^{5}\,j_3^{4}
      + 10584000\,j_2^{2}\,j_3^{6} + 52220\,j_2^{9}\,j_4 -
      3427200\,j_2^{6}\,j_3^{2}\,j_4 + 64638000\,j_2^{3}\,j_3^{4}\,j_4 -
      264600000\,j_3^{6}\,j_4 - 1466060\,j_2^{7}\,j_4^{2} +
      44871600\,j_2^{4}\,j_3^{2}\,j_4^{2} - 26694000\,j_2\,j_3^{4}\,j_4^{2} -
      34824420\,j_2^{5}\,j_4^{3} + 1478487600\,j_2^{2}\,j_3^{2}\,j_4^{3} +
      1906270200\,j_2^{3}\,j_4^{4} + 10926981000\,j_3^{2}\,j_4^{4} -
      6936954048\,j_2\,j_4^{5} - 476280\,j_2^{7}\,j_3\,j_5 +
      28576800\,j_2^{4}\,j_3^{3}\,j_5 - 428652000\,j_2\,j_3^{5}\,j_5 +
      34360200\,j_2^{5}\,j_3\,j_4\,j_5 -
      1030806000\,j_2^{2}\,j_3^{3}\,j_4\,j_5 -
      414963000\,j_2^{3}\,j_3\,j_4^{2}\,j_5 -
      5003370000\,j_3^{3}\,j_4^{2}\,j_5 - 3794372100\,j_2\,j_3\,j_4^{3}\,j_5 -
      6429780\,j_2^{6}\,j_5^{2} + 192893400\,j_2^{3}\,j_3^{2}\,j_5^{2} +
      417170250\,j_2^{4}\,j_4\,j_5^{2} -
      5557167000\,j_2\,j_3^{2}\,j_4\,j_5^{2} -
      15345304200\,j_2^{2}\,j_4^{2}\,j_5^{2} + 20076660000\,j_4^{3}\,j_5^{2} -
      868020300\,j_2^{2}\,j_3\,j_5^{3} + 18600435000\,j_3\,j_4\,j_5^{3} +
      23436548100\,j_2\,j_5^{4} - 330750\,j_2^{8}\,j_6 +
      19845000\,j_2^{5}\,j_3^{2}\,j_6 - 297675000\,j_2^{2}\,j_3^{4}\,j_6 +
      35514000\,j_2^{6}\,j_4\,j_6 - 858330000\,j_2^{3}\,j_3^{2}\,j_4\,j_6 -
      6212700000\,j_3^{4}\,j_4\,j_6 - 1287423000\,j_2^{4}\,j_4^{2}\,j_6 +
      1452802500\,j_2\,j_3^{2}\,j_4^{2}\,j_6 +
      13651027200\,j_2^{2}\,j_4^{3}\,j_6 - 29029363200\,j_4^{4}\,j_6 +
      233887500\,j_2^{4}\,j_3\,j_5\,j_6 - 7016625000\,j_2\,j_3^{3}\,j_5\,j_6 -
      5401890000\,j_2^{2}\,j_3\,j_4\,j_5\,j_6 -
      53512974000\,j_3\,j_4^{2}\,j_5\,j_6 + 960001875\,j_2^{3}\,j_5^{2}\,j_6 +
      120175650000\,j_3^{2}\,j_5^{2}\,j_6 -
      46983321000\,j_2\,j_4\,j_5^{2}\,j_6 - 97524000\,j_2^{5}\,j_6^{2} +
      2925720000\,j_2^{2}\,j_3^{2}\,j_6^{2} +
      3212136000\,j_2^{3}\,j_4\,j_6^{2} - 133688880000\,j_3^{2}\,j_4\,j_6^{2}
      - 5787936000\,j_2\,j_4^{2}\,j_6^{2} - 1224720000\,j_2\,j_3\,j_5\,j_6^{2}
      + 40415760000\,j_5^{2}\,j_6^{2} + 2646000\,j_2^{6}\,j_3\,j_7 -
      158760000\,j_2^{3}\,j_3^{3}\,j_7 + 2381400000\,j_3^{5}\,j_7 -
      60385500\,j_2^{4}\,j_3\,j_4\,j_7 + 1811565000\,j_2\,j_3^{3}\,j_4\,j_7 -
      13778100000\,j_2^{2}\,j_3\,j_4^{2}\,j_7 + 58484754000\,j_3\,j_4^{3}\,j_7
      + 168399000\,j_2^{5}\,j_5\,j_7 - 5051970000\,j_2^{2}\,j_3^{2}\,j_5\,j_7
      - 5370907500\,j_2^{3}\,j_4\,j_5\,j_7 -
      130279590000\,j_3^{2}\,j_4\,j_5\,j_7 +
      46845540000\,j_2\,j_4^{2}\,j_5\,j_7 +
      69234952500\,j_2\,j_3\,j_5^{2}\,j_7 - 111602610000\,j_5^{3}\,j_7 -
      2517480000\,j_2^{3}\,j_3\,j_6\,j_7 + 75524400000\,j_3^{3}\,j_6\,j_7 +
      24308505000\,j_2\,j_3\,j_4\,j_6\,j_7 +
      14581822500\,j_2^{2}\,j_5\,j_6\,j_7 + 48945060000\,j_4\,j_5\,j_6\,j_7 -
      382725000\,j_2^{4}\,j_7^{2} + 11481750000\,j_2\,j_3^{2}\,j_7^{2} +
      5740875000\,j_2^{2}\,j_4\,j_7^{2} - 33067440000\,j_4^{2}\,j_7^{2} +
      5197500\,j_2^{7}\,j_8 - 311850000\,j_2^{4}\,j_3^{2}\,j_8 +
      4677750000\,j_2\,j_3^{4}\,j_8 - 896427000\,j_2^{5}\,j_4\,j_8 +
      26892810000\,j_2^{2}\,j_3^{2}\,j_4\,j_8 +
      38918880000\,j_2^{3}\,j_4^{2}\,j_8 + 60980850000\,j_3^{2}\,j_4^{2}\,j_8
      - 155110788000\,j_2\,j_4^{3}\,j_8 + 153090000\,j_2^{3}\,j_3\,j_5\,j_8 -
      4592700000\,j_3^{3}\,j_5\,j_8 - 44434372500\,j_2\,j_3\,j_4\,j_5\,j_8 -
      41850978750\,j_2^{2}\,j_5^{2}\,j_8 + 301740390000\,j_4\,j_5^{2}\,j_8 +
      161595000\,j_2^{4}\,j_6\,j_8 - 4847850000\,j_2\,j_3^{2}\,j_6\,j_8 -
      37813230000\,j_2^{2}\,j_4\,j_6\,j_8 + 173910240000\,j_4^{2}\,j_6\,j_8 -
      564902100000\,j_3\,j_5\,j_6\,j_8 - 68890500000\,j_3\,j_4\,j_7\,j_8 +
      25515000\,j_2^{5}\,j_3\,j_9 - 765450000\,j_2^{2}\,j_3^{3}\,j_9 -
      2393307000\,j_2^{3}\,j_3\,j_4\,j_9 + 37353960000\,j_3^{3}\,j_4\,j_9 +
      44824752000\,j_2\,j_3\,j_4^{2}\,j_9 - 803722500\,j_2^{4}\,j_5\,j_9 +
      34445250000\,j_2\,j_3^{2}\,j_5\,j_9 +
      10103940000\,j_2^{2}\,j_4\,j_5\,j_9 + 88179840000\,j_4^{2}\,j_5\,j_9 -
      421609860000\,j_3\,j_5^{2}\,j_9 + 36741600000\,j_2^{2}\,j_3\,j_6\,j_9 +
      367416000000\,j_3\,j_4\,j_6\,j_9 + 5051970000\,j_2^{3}\,j_7\,j_9 +
      55112400000\,j_3^{2}\,j_7\,j_9 - 55112400000\,j_2\,j_4\,j_7\,j_9 -
      30051000\,j_2^{6}\,j_{10} + 1037610000\,j_2^{3}\,j_3^{2}\,j_{10} -
      4082400000\,j_3^{4}\,j_{10} + 1316574000\,j_2^{4}\,j_4\,j_{10} -
      5051970000\,j_2\,j_3^{2}\,j_4\,j_{10} +
      83770848000\,j_2^{2}\,j_4^{2}\,j_{10} - 387991296000\,j_4^{3}\,j_{10} -
      10333575000\,j_2^{2}\,j_3\,j_5\,j_{10} +
      655837560000\,j_3\,j_4\,j_5\,j_{10} - 173604060000\,j_2\,j_5^{2}\,j_{10}
      - 18779040000\,j_2^{3}\,j_6\,j_{10} - 538876800000\,j_3^{2}\,j_6\,j_{10}
      + 102876480000\,j_2\,j_4\,j_6\,j_{10} -
      206671500000\,j_2\,j_3\,j_7\,j_{10} + 496011600000\,j_5\,j_7\,j_{10},
    \end{dmath*}
    \begin{dmath*}
      0= -8624\,j_2^{11} + 776160\,j_2^{8}\,j_3^{2} -
      23284800\,j_2^{5}\,j_3^{4} + 232848000\,j_2^{2}\,j_3^{6} +
      836430\,j_2^{9}\,j_4 - 61299000\,j_2^{6}\,j_3^{2}\,j_4 +
      1419579000\,j_2^{3}\,j_3^{4}\,j_4 - 10001880000\,j_3^{6}\,j_4 +
      2396520\,j_2^{7}\,j_4^{2} + 169759800\,j_2^{4}\,j_3^{2}\,j_4^{2} -
      7249662000\,j_2\,j_3^{4}\,j_4^{2} - 1755542565\,j_2^{5}\,j_4^{3} +
      33823449450\,j_2^{2}\,j_3^{2}\,j_4^{3} + 57104013420\,j_2^{3}\,j_4^{4} +
      610877651400\,j_3^{2}\,j_4^{4} - 63190933056\,j_2\,j_4^{5} -
      10478160\,j_2^{7}\,j_3\,j_5 + 628689600\,j_2^{4}\,j_3^{3}\,j_5 -
      9430344000\,j_2\,j_3^{5}\,j_5 + 377451900\,j_2^{5}\,j_3\,j_4\,j_5 -
      11323557000\,j_2^{2}\,j_3^{3}\,j_4\,j_5 +
      17861010300\,j_2^{3}\,j_3\,j_4^{2}\,j_5 -
      408842154000\,j_3^{3}\,j_4^{2}\,j_5 -
      515879620200\,j_2\,j_3\,j_4^{3}\,j_5 - 141455160\,j_2^{6}\,j_5^{2} +
      4243654800\,j_2^{3}\,j_3^{2}\,j_5^{2} +
      8045262225\,j_2^{4}\,j_4\,j_5^{2} -
      164923857000\,j_2\,j_3^{2}\,j_4\,j_5^{2} -
      381143580300\,j_2^{2}\,j_4^{2}\,j_5^{2} + 244888012800\,j_4^{3}\,j_5^{2}
      - 19096446600\,j_2^{2}\,j_3\,j_5^{3} + 1328071059000\,j_3\,j_4\,j_5^{3}
      + 515604058200\,j_2\,j_5^{4} - 15245370\,j_2^{8}\,j_6 +
      914722200\,j_2^{5}\,j_3^{2}\,j_6 - 13720833000\,j_2^{2}\,j_3^{4}\,j_6 +
      1616110650\,j_2^{6}\,j_4\,j_6 - 41880037500\,j_2^{3}\,j_3^{2}\,j_4\,j_6
      - 198098460000\,j_3^{4}\,j_4\,j_6 - 50739668460\,j_2^{4}\,j_4^{2}\,j_6 -
      30550786200\,j_2\,j_3^{2}\,j_4^{2}\,j_6 +
      343850637600\,j_2^{2}\,j_4^{3}\,j_6 - 316389265920\,j_4^{4}\,j_6 +
      9708967800\,j_2^{4}\,j_3\,j_5\,j_6 -
      291269034000\,j_2\,j_3^{3}\,j_5\,j_6 -
      248928167250\,j_2^{2}\,j_3\,j_4\,j_5\,j_6 -
      2166618034800\,j_3\,j_4^{2}\,j_5\,j_6 +
      71491116375\,j_2^{3}\,j_5^{2}\,j_6 +
      4288020282000\,j_3^{2}\,j_5^{2}\,j_6 -
      1124706303000\,j_2\,j_4\,j_5^{2}\,j_6 - 3118641750\,j_2^{5}\,j_6^{2} +
      93559252500\,j_2^{2}\,j_3^{2}\,j_6^{2} +
      6336225000\,j_2^{3}\,j_4\,j_6^{2} - 4275803700000\,j_3^{2}\,j_4\,j_6^{2}
      + 640073664000\,j_2\,j_4^{2}\,j_6^{2} -
      308629440000\,j_2\,j_3\,j_5\,j_6^{2} + 100018800\,j_2^{6}\,j_3\,j_7 -
      6001128000\,j_2^{3}\,j_3^{3}\,j_7 + 90016920000\,j_3^{5}\,j_7 -
      7372814400\,j_2^{4}\,j_3\,j_4\,j_7 +
      221184432000\,j_2\,j_3^{3}\,j_4\,j_7 -
      305977155750\,j_2^{2}\,j_3\,j_4^{2}\,j_7 +
      1008391582800\,j_3\,j_4^{3}\,j_7 + 2470107150\,j_2^{5}\,j_5\,j_7 -
      74103214500\,j_2^{2}\,j_3^{2}\,j_5\,j_7 -
      91538251875\,j_2^{3}\,j_4\,j_5\,j_7 -
      4395489462000\,j_3^{2}\,j_4\,j_5\,j_7 +
      308767221000\,j_2\,j_4^{2}\,j_5\,j_7 +
      1484314713000\,j_2\,j_3\,j_5^{2}\,j_7 - 1406192886000\,j_5^{3}\,j_7 -
      83587140000\,j_2^{3}\,j_3\,j_6\,j_7 + 2507614200000\,j_3^{3}\,j_6\,j_7 +
      1928934000000\,j_2\,j_3\,j_4\,j_6\,j_7 +
      450888322500\,j_2^{2}\,j_5\,j_6\,j_7 + 1512835380000\,j_4\,j_5\,j_6\,j_7
      + 121451400\,j_2^{7}\,j_8 - 7287084000\,j_2^{4}\,j_3^{2}\,j_8 +
      109306260000\,j_2\,j_3^{4}\,j_8 - 23588362350\,j_2^{5}\,j_4\,j_8 +
      707650870500\,j_2^{2}\,j_3^{2}\,j_4\,j_8 +
      1093636687500\,j_2^{3}\,j_4^{2}\,j_8 +
      4713488010000\,j_3^{2}\,j_4^{2}\,j_8 - 1533337192800\,j_2\,j_4^{3}\,j_8
      + 28934010000\,j_2^{3}\,j_3\,j_5\,j_8 - 868020300000\,j_3^{3}\,j_5\,j_8
      - 1961725878000\,j_2\,j_3\,j_4\,j_5\,j_8 -
      1304200500750\,j_2^{2}\,j_5^{2}\,j_8 + 6138143550000\,j_4\,j_5^{2}\,j_8
      - 6429780000\,j_2^{4}\,j_6\,j_8 + 192893400000\,j_2\,j_3^{2}\,j_6\,j_8 -
      932318100000\,j_2^{2}\,j_4\,j_6\,j_8 - 6811892640000\,j_4^{2}\,j_6\,j_8
      - 15624365400000\,j_3\,j_5\,j_6\,j_8 + 1696747500\,j_2^{5}\,j_3\,j_9 -
      50902425000\,j_2^{2}\,j_3^{3}\,j_9 -
      120664007100\,j_2^{3}\,j_3\,j_4\,j_9 + 1329311088000\,j_3^{3}\,j_4\,j_9
      + 1750920948000\,j_2\,j_3\,j_4^{2}\,j_9 - 8342639550\,j_2^{4}\,j_5\,j_9
      + 937461924000\,j_2\,j_3^{2}\,j_5\,j_9 +
      188815082400\,j_2^{2}\,j_4\,j_5\,j_9 + 273798403200\,j_4^{2}\,j_5\,j_9 -
      11874517704000\,j_3\,j_5^{2}\,j_9 +
      1527072750000\,j_2^{2}\,j_3\,j_6\,j_9 +
      12201885360000\,j_3\,j_4\,j_6\,j_9 - 1153788300\,j_2^{6}\,j_{10} +
      18324873000\,j_2^{3}\,j_3^{2}\,j_{10} + 488663280000\,j_3^{4}\,j_{10} +
      51993956700\,j_2^{4}\,j_4\,j_{10} +
      730790424000\,j_2\,j_3^{2}\,j_4\,j_{10} +
      2438723700000\,j_2^{2}\,j_4^{2}\,j_{10} - 2095152998400\,j_4^{3}\,j_{10}
      - 687182737500\,j_2^{2}\,j_3\,j_5\,j_{10} +
      12008440836000\,j_3\,j_4\,j_5\,j_{10} -
      3749847696000\,j_2\,j_5^{2}\,j_{10} - 819796950000\,j_2^{3}\,j_6\,j_{10}
      - 21218274000000\,j_3^{2}\,j_6\,j_{10} +
      1686439440000\,j_2\,j_4\,j_6\,j_{10} + 9258883200000\,j_6^{2}\,j_{10},
    \end{dmath*}
    \begin{dmath*}
      0= -392\,j_2^{11} + 35280\,j_2^{8}\,j_3^{2} - 1058400\,j_2^{5}\,j_3^{4}
      + 10584000\,j_2^{2}\,j_3^{6} + 40740\,j_2^{9}\,j_4 -
      2856000\,j_2^{6}\,j_3^{2}\,j_4 + 61362000\,j_2^{3}\,j_3^{4}\,j_4 -
      370440000\,j_3^{6}\,j_4 - 124140\,j_2^{7}\,j_4^{2} +
      5936400\,j_2^{4}\,j_3^{2}\,j_4^{2} - 66366000\,j_2\,j_3^{4}\,j_4^{2} -
      76735620\,j_2^{5}\,j_4^{3} + 1800273600\,j_2^{2}\,j_3^{2}\,j_4^{3} +
      2536560360\,j_2^{3}\,j_4^{4} + 19717846200\,j_3^{2}\,j_4^{4} -
      4249568448\,j_2\,j_4^{5} - 476280\,j_2^{7}\,j_3\,j_5 +
      28576800\,j_2^{4}\,j_3^{3}\,j_5 - 428652000\,j_2\,j_3^{5}\,j_5 +
      21886200\,j_2^{5}\,j_3\,j_4\,j_5 - 656586000\,j_2^{2}\,j_3^{3}\,j_4\,j_5
      + 417133800\,j_2^{3}\,j_3\,j_4^{2}\,j_5 -
      13126374000\,j_3^{3}\,j_4^{2}\,j_5 - 15611462100\,j_2\,j_3\,j_4^{3}\,j_5
      - 6429780\,j_2^{6}\,j_5^{2} + 192893400\,j_2^{3}\,j_3^{2}\,j_5^{2} +
      378897750\,j_2^{4}\,j_4\,j_5^{2} -
      6934977000\,j_2\,j_3^{2}\,j_4\,j_5^{2} -
      17324320500\,j_2^{2}\,j_4^{2}\,j_5^{2} + 17620221600\,j_4^{3}\,j_5^{2} -
      868020300\,j_2^{2}\,j_3\,j_5^{3} + 45881073000\,j_3\,j_4\,j_5^{3} +
      23436548100\,j_2\,j_5^{4} - 580230\,j_2^{8}\,j_6 +
      34813800\,j_2^{5}\,j_3^{2}\,j_6 - 522207000\,j_2^{2}\,j_3^{4}\,j_6 +
      60975720\,j_2^{6}\,j_4\,j_6 - 1609459200\,j_2^{3}\,j_3^{2}\,j_4\,j_6 -
      6594372000\,j_3^{4}\,j_4\,j_6 - 1956375720\,j_2^{4}\,j_4^{2}\,j_6 +
      1731464100\,j_2\,j_3^{2}\,j_4^{2}\,j_6 +
      14895705600\,j_2^{2}\,j_4^{3}\,j_6 - 17473605120\,j_4^{4}\,j_6 +
      314174700\,j_2^{4}\,j_3\,j_5\,j_6 - 9425241000\,j_2\,j_3^{3}\,j_5\,j_6 -
      7624611000\,j_2^{2}\,j_3\,j_4\,j_5\,j_6 -
      83295831600\,j_3\,j_4^{2}\,j_5\,j_6 + 2483247375\,j_2^{3}\,j_5^{2}\,j_6
      + 147915558000\,j_3^{2}\,j_5^{2}\,j_6 -
      47125038600\,j_2\,j_4\,j_5^{2}\,j_6 - 138348000\,j_2^{5}\,j_6^{2} +
      4150440000\,j_2^{2}\,j_3^{2}\,j_6^{2} +
      1659204000\,j_2^{3}\,j_4\,j_6^{2} - 161876880000\,j_3^{2}\,j_4\,j_6^{2}
      + 10552032000\,j_2\,j_4^{2}\,j_6^{2} -
      25412940000\,j_2\,j_3\,j_5\,j_6^{2} + 40415760000\,j_5^{2}\,j_6^{2} +
      3704400\,j_2^{6}\,j_3\,j_7 - 222264000\,j_2^{3}\,j_3^{3}\,j_7 +
      3333960000\,j_3^{5}\,j_7 - 195785100\,j_2^{4}\,j_3\,j_4\,j_7 +
      5873553000\,j_2\,j_3^{3}\,j_4\,j_7 -
      13915881000\,j_2^{2}\,j_3\,j_4^{2}\,j_7 + 49257363600\,j_3\,j_4^{3}\,j_7
      + 122472000\,j_2^{5}\,j_5\,j_7 - 3674160000\,j_2^{2}\,j_3^{2}\,j_5\,j_7
      - 3784384800\,j_2^{3}\,j_4\,j_5\,j_7 -
      147976794000\,j_3^{2}\,j_4\,j_5\,j_7 +
      13998549600\,j_2\,j_4^{2}\,j_5\,j_7 +
      67581580500\,j_2\,j_3\,j_5^{2}\,j_7 - 66961566000\,j_5^{3}\,j_7 -
      2789640000\,j_2^{3}\,j_3\,j_6\,j_7 + 83689200000\,j_3^{3}\,j_6\,j_7 +
      79552125000\,j_2\,j_3\,j_4\,j_6\,j_7 +
      18255982500\,j_2^{2}\,j_5\,j_6\,j_7 - 41334300000\,j_4\,j_5\,j_6\,j_7 +
      76545000\,j_2^{4}\,j_7^{2} - 2296350000\,j_2\,j_3^{2}\,j_7^{2} -
      1148175000\,j_2^{2}\,j_4\,j_7^{2} - 16533720000\,j_4^{2}\,j_7^{2} +
      5499900\,j_2^{7}\,j_8 - 329994000\,j_2^{4}\,j_3^{2}\,j_8 +
      4949910000\,j_2\,j_3^{4}\,j_8 - 1024682400\,j_2^{5}\,j_4\,j_8 +
      30740472000\,j_2^{2}\,j_3^{2}\,j_4\,j_8 +
      45892299600\,j_2^{3}\,j_4^{2}\,j_8 + 129177342000\,j_3^{2}\,j_4^{2}\,j_8
      - 70948029600\,j_2\,j_4^{3}\,j_8 + 949158000\,j_2^{3}\,j_3\,j_5\,j_8 -
      28474740000\,j_3^{3}\,j_5\,j_8 - 105333574500\,j_2\,j_3\,j_4\,j_5\,j_8 -
      55491297750\,j_2^{2}\,j_5^{2}\,j_8 + 265366206000\,j_4\,j_5^{2}\,j_8 -
      178605000\,j_2^{4}\,j_6\,j_8 + 5358150000\,j_2\,j_3^{2}\,j_6\,j_8 -
      44549190000\,j_2^{2}\,j_4\,j_6\,j_8 + 95528160000\,j_4^{2}\,j_6\,j_8 -
      537345900000\,j_3\,j_5\,j_6\,j_8 + 13778100000\,j_3\,j_4\,j_7\,j_8 +
      62937000\,j_2^{5}\,j_3\,j_9 - 1888110000\,j_2^{2}\,j_3^{3}\,j_9 -
      4599844200\,j_2^{3}\,j_3\,j_4\,j_9 + 53030376000\,j_3^{3}\,j_4\,j_9 +
      46294416000\,j_2\,j_3\,j_4^{2}\,j_9 - 549593100\,j_2^{4}\,j_5\,j_9 +
      41977278000\,j_2\,j_3^{2}\,j_5\,j_9 +
      11867536800\,j_2^{2}\,j_4\,j_5\,j_9 + 31744742400\,j_4^{2}\,j_5\,j_9 -
      461290788000\,j_3\,j_5^{2}\,j_9 + 53887680000\,j_2^{2}\,j_3\,j_6\,j_9 +
      448247520000\,j_3\,j_4\,j_6\,j_9 - 765450000\,j_2^{3}\,j_7\,j_9 -
      18370800000\,j_3^{2}\,j_7\,j_9 + 220449600000\,j_6\,j_7\,j_9 -
      46607400\,j_2^{6}\,j_{10} + 908334000\,j_2^{3}\,j_3^{2}\,j_{10} +
      14696640000\,j_3^{4}\,j_{10} + 2088147600\,j_2^{4}\,j_4\,j_{10} +
      22320522000\,j_2\,j_3^{2}\,j_4\,j_{10} +
      107652888000\,j_2^{2}\,j_4^{2}\,j_{10} - 218686003200\,j_4^{3}\,j_{10} -
      25489485000\,j_2^{2}\,j_3\,j_5\,j_{10} +
      447512688000\,j_3\,j_4\,j_5\,j_{10} - 163683828000\,j_2\,j_5^{2}\,j_{10}
      - 31842720000\,j_2^{3}\,j_6\,j_{10} - 661348800000\,j_3^{2}\,j_6\,j_{10}
      + 88179840000\,j_2\,j_4\,j_6\,j_{10} +
      41334300000\,j_2\,j_3\,j_7\,j_{10},
    \end{dmath*}
    \begin{dmath*}
      0= 5096\,j_2^{11} - 458640\,j_2^{8}\,j_3^{2} +
      13759200\,j_2^{5}\,j_3^{4} - 137592000\,j_2^{2}\,j_3^{6} -
      443520\,j_2^{9}\,j_4 + 28627200\,j_2^{6}\,j_3^{2}\,j_4 -
      520128000\,j_2^{3}\,j_3^{4}\,j_4 + 1814400000\,j_3^{6}\,j_4 -
      7755720\,j_2^{7}\,j_4^{2} + 552661200\,j_2^{4}\,j_3^{2}\,j_4^{2} -
      9599688000\,j_2\,j_3^{4}\,j_4^{2} + 1264264740\,j_2^{5}\,j_4^{3} -
      32927002200\,j_2^{2}\,j_3^{2}\,j_4^{3} - 37554396120\,j_2^{3}\,j_4^{4} -
      90975263400\,j_3^{2}\,j_4^{4} + 26493213024\,j_2\,j_4^{5} +
      6191640\,j_2^{7}\,j_3\,j_5 - 371498400\,j_2^{4}\,j_3^{3}\,j_5 +
      5572476000\,j_2\,j_3^{5}\,j_5 - 304819200\,j_2^{5}\,j_3\,j_4\,j_5 +
      9144576000\,j_2^{2}\,j_3^{3}\,j_4\,j_5 -
      836357400\,j_2^{3}\,j_3\,j_4^{2}\,j_5 +
      60454512000\,j_3^{3}\,j_4^{2}\,j_5 + 15308343900\,j_2\,j_3\,j_4^{3}\,j_5
      + 83587140\,j_2^{6}\,j_5^{2} - 2507614200\,j_2^{3}\,j_3^{2}\,j_5^{2} -
      4675368600\,j_2^{4}\,j_4\,j_5^{2} +
      78535170000\,j_2\,j_3^{2}\,j_4\,j_5^{2} +
      234850338900\,j_2^{2}\,j_4^{2}\,j_5^{2} - 78582409200\,j_4^{3}\,j_5^{2}
      + 11284263900\,j_2^{2}\,j_3\,j_5^{3} - 193444524000\,j_3\,j_4\,j_5^{3} -
      304675125300\,j_2\,j_5^{4} + 10104570\,j_2^{8}\,j_6 -
      606274200\,j_2^{5}\,j_3^{2}\,j_6 + 9094113000\,j_2^{2}\,j_3^{4}\,j_6 -
      1065657060\,j_2^{6}\,j_4\,j_6 + 31054897800\,j_2^{3}\,j_3^{2}\,j_4\,j_6
      + 27444420000\,j_3^{4}\,j_4\,j_6 + 32724232200\,j_2^{4}\,j_4^{2}\,j_6 -
      115775932500\,j_2\,j_3^{2}\,j_4^{2}\,j_6 -
      202754242800\,j_2^{2}\,j_4^{3}\,j_6 + 121545411840\,j_4^{4}\,j_6 -
      1047645900\,j_2^{4}\,j_3\,j_5\,j_6 + 31429377000\,j_2\,j_3^{3}\,j_5\,j_6
      - 28485456300\,j_2^{2}\,j_3\,j_4\,j_5\,j_6 +
      438804491400\,j_3\,j_4^{2}\,j_5\,j_6 -
      47277636525\,j_2^{3}\,j_5^{2}\,j_6 - 663369588000\,j_3^{2}\,j_5^{2}\,j_6
      + 630393345900\,j_2\,j_4\,j_5^{2}\,j_6 + 2134188000\,j_2^{5}\,j_6^{2} -
      64025640000\,j_2^{2}\,j_3^{2}\,j_6^{2} -
      5592888000\,j_2^{3}\,j_4\,j_6^{2} + 829981080000\,j_3^{2}\,j_4\,j_6^{2}
      - 272957040000\,j_2\,j_4^{2}\,j_6^{2} +
      524486340000\,j_2\,j_3\,j_5\,j_6^{2} - 303118200000\,j_5^{2}\,j_6^{2} -
      18144000\,j_2^{6}\,j_3\,j_7 + 1088640000\,j_2^{3}\,j_3^{3}\,j_7 -
      16329600000\,j_3^{5}\,j_7 - 1674294300\,j_2^{4}\,j_3\,j_4\,j_7 +
      50228829000\,j_2\,j_3^{3}\,j_4\,j_7 +
      289794777300\,j_2^{2}\,j_3\,j_4^{2}\,j_7 -
      382003727400\,j_3\,j_4^{3}\,j_7 - 1241559900\,j_2^{5}\,j_5\,j_7 +
      37246797000\,j_2^{2}\,j_3^{2}\,j_5\,j_7 +
      46060953750\,j_2^{3}\,j_4\,j_5\,j_7 +
      700386750000\,j_3^{2}\,j_4\,j_5\,j_7 +
      33356780100\,j_2\,j_4^{2}\,j_5\,j_7 -
      863680198500\,j_2\,j_3\,j_5^{2}\,j_7 + 468730962000\,j_5^{3}\,j_7 +
      14441490000\,j_2^{3}\,j_3\,j_6\,j_7 - 433244700000\,j_3^{3}\,j_6\,j_7 -
      468160155000\,j_2\,j_3\,j_4\,j_6\,j_7 -
      294506887500\,j_2^{2}\,j_5\,j_6\,j_7 + 106091370000\,j_4\,j_5\,j_6\,j_7
      - 1377810000\,j_2^{4}\,j_7^{2} + 41334300000\,j_2\,j_3^{2}\,j_7^{2} +
      2066715000\,j_2^{2}\,j_4\,j_7^{2} + 97430850000\,j_4^{2}\,j_7^{2} -
      74220300\,j_2^{7}\,j_8 + 4453218000\,j_2^{4}\,j_3^{2}\,j_8 -
      66798270000\,j_2\,j_3^{4}\,j_8 + 14603255100\,j_2^{5}\,j_4\,j_8 -
      438097653000\,j_2^{2}\,j_3^{2}\,j_4\,j_8 -
      671522140800\,j_2^{3}\,j_4^{2}\,j_8 -
      391144950000\,j_3^{2}\,j_4^{2}\,j_8 + 122992506000\,j_2\,j_4^{3}\,j_8 -
      3490452000\,j_2^{3}\,j_3\,j_5\,j_8 + 104713560000\,j_3^{3}\,j_5\,j_8 +
      1620511231500\,j_2\,j_3\,j_4\,j_5\,j_8 +
      818729147250\,j_2^{2}\,j_5^{2}\,j_8 - 2372175477000\,j_4\,j_5^{2}\,j_8 +
      4618215000\,j_2^{4}\,j_6\,j_8 - 138546450000\,j_2\,j_3^{2}\,j_6\,j_8 +
      708653610000\,j_2^{2}\,j_4\,j_6\,j_8 - 789944400000\,j_4^{2}\,j_6\,j_8 +
      2335387950000\,j_3\,j_5\,j_6\,j_8 + 310007250000\,j_3\,j_4\,j_7\,j_8 +
      1860043500000\,j_7^{2}\,j_8 - 757285200\,j_2^{5}\,j_3\,j_9 +
      22718556000\,j_2^{2}\,j_3^{3}\,j_9 + 43385706000\,j_2^{3}\,j_3\,j_4\,j_9
      - 279236160000\,j_3^{3}\,j_4\,j_9 + 330453950400\,j_2\,j_3\,j_4^{2}\,j_9
      + 2571912000\,j_2^{4}\,j_5\,j_9 - 383857866000\,j_2\,j_3^{2}\,j_5\,j_9 -
      61009426800\,j_2^{2}\,j_4\,j_5\,j_9 - 513206668800\,j_4^{2}\,j_5\,j_9 +
      2462697594000\,j_3\,j_5^{2}\,j_9 - 631955520000\,j_2^{2}\,j_3\,j_6\,j_9
      - 2006091360000\,j_3\,j_4\,j_6\,j_9 + 24800580000\,j_2^{3}\,j_7\,j_9 -
      347208120000\,j_2\,j_4\,j_7\,j_9 + 713059200\,j_2^{6}\,j_{10} -
      20064996000\,j_2^{3}\,j_3^{2}\,j_{10} - 39803400000\,j_3^{4}\,j_{10} -
      28627830000\,j_2^{4}\,j_4\,j_{10} -
      163500120000\,j_2\,j_3^{2}\,j_4\,j_{10} -
      1750039149600\,j_2^{2}\,j_4^{2}\,j_{10} + 1724797670400\,j_4^{3}\,j_{10}
      + 306700506000\,j_2^{2}\,j_3\,j_5\,j_{10} -
      2787585192000\,j_3\,j_4\,j_5\,j_{10} +
      2403176202000\,j_2\,j_5^{2}\,j_{10} + 531528480000\,j_2^{3}\,j_6\,j_{10}
      + 3012811200000\,j_3^{2}\,j_6\,j_{10} -
      837708480000\,j_2\,j_4\,j_6\,j_{10} -
      744017400000\,j_2\,j_3\,j_7\,j_{10},
    \end{dmath*}
    \begin{dmath*}
      0= 5488\,j_2^{10}\,j_3 - 493920\,j_2^{7}\,j_3^{3} +
      14817600\,j_2^{4}\,j_3^{5} - 148176000\,j_2\,j_3^{7} +
      79870\,j_2^{8}\,j_3\,j_4 - 4792200\,j_2^{5}\,j_3^{3}\,j_4 +
      71883000\,j_2^{2}\,j_3^{5}\,j_4 - 59795820\,j_2^{6}\,j_3\,j_4^{2} +
      2371912200\,j_2^{3}\,j_3^{3}\,j_4^{2} - 17341128000\,j_3^{5}\,j_4^{2} +
      2488273830\,j_2^{4}\,j_3\,j_4^{3} - 39790472400\,j_2\,j_3^{3}\,j_4^{3} -
      70269330600\,j_2^{2}\,j_3\,j_4^{4} - 86770852128\,j_3\,j_4^{5} +
      6667920\,j_2^{6}\,j_3^{2}\,j_5 - 400075200\,j_2^{3}\,j_3^{4}\,j_5 +
      6001128000\,j_3^{6}\,j_5 + 9887220\,j_2^{7}\,j_4\,j_5 -
      528538500\,j_2^{4}\,j_3^{2}\,j_4\,j_5 +
      6957657000\,j_2\,j_3^{4}\,j_4\,j_5 - 993245085\,j_2^{5}\,j_4^{2}\,j_5 +
      13383553050\,j_2^{2}\,j_3^{2}\,j_4^{2}\,j_5 +
      23895147420\,j_2^{3}\,j_4^{3}\,j_5 - 49969888200\,j_3^{2}\,j_4^{3}\,j_5
      + 68630684400\,j_2\,j_4^{4}\,j_5 + 90016920\,j_2^{5}\,j_3\,j_5^{2} -
      2700507600\,j_2^{2}\,j_3^{3}\,j_5^{2} -
      2420990775\,j_2^{3}\,j_3\,j_4\,j_5^{2} +
      85730400000\,j_3^{3}\,j_4\,j_5^{2} +
      286878412800\,j_2\,j_3\,j_4^{2}\,j_5^{2} +
      12152284200\,j_2\,j_3^{2}\,j_5^{3} - 56180202750\,j_2^{2}\,j_4\,j_5^{3}
      - 122762871000\,j_4^{2}\,j_5^{3} - 328111673400\,j_3\,j_5^{4} +
      16321410\,j_2^{7}\,j_3\,j_6 - 979284600\,j_2^{4}\,j_3^{3}\,j_6 +
      14689269000\,j_2\,j_3^{5}\,j_6 - 1954956150\,j_2^{5}\,j_3\,j_4\,j_6 +
      58648684500\,j_2^{2}\,j_3^{3}\,j_4\,j_6 +
      58880448720\,j_2^{3}\,j_3\,j_4^{2}\,j_6 -
      118748316600\,j_3^{3}\,j_4^{2}\,j_6 -
      142003951200\,j_2\,j_3\,j_4^{3}\,j_6 - 398249460\,j_2^{6}\,j_5\,j_6 +
      12973867200\,j_2^{3}\,j_3^{2}\,j_5\,j_6 - 30791502000\,j_3^{4}\,j_5\,j_6
      + 14393005830\,j_2^{4}\,j_4\,j_5\,j_6 -
      137182503150\,j_2\,j_3^{2}\,j_4\,j_5\,j_6 +
      65446908120\,j_2^{2}\,j_4^{2}\,j_5\,j_6 +
      203172649920\,j_4^{3}\,j_5\,j_6 -
      236672164575\,j_2^{2}\,j_3\,j_5^{2}\,j_6 +
      43915397400\,j_3\,j_4\,j_5^{2}\,j_6 + 150746192100\,j_2\,j_5^{3}\,j_6 +
      2635416000\,j_2^{4}\,j_3\,j_6^{2} - 79062480000\,j_2\,j_3^{3}\,j_6^{2} +
      190968435000\,j_2^{2}\,j_3\,j_4\,j_6^{2} +
      78253776000\,j_3\,j_4^{2}\,j_6^{2} - 19735852500\,j_2^{3}\,j_5\,j_6^{2}
      + 788719680000\,j_3^{2}\,j_5\,j_6^{2} -
      238861224000\,j_2\,j_4\,j_5\,j_6^{2} - 3360420\,j_2^{8}\,j_7 +
      201625200\,j_2^{5}\,j_3^{2}\,j_7 - 3024378000\,j_2^{2}\,j_3^{4}\,j_7 +
      489145230\,j_2^{6}\,j_4\,j_7 - 18940805100\,j_2^{3}\,j_3^{2}\,j_4\,j_7 +
      127993446000\,j_3^{4}\,j_4\,j_7 - 6300367920\,j_2^{4}\,j_4^{2}\,j_7 +
      388146682350\,j_2\,j_3^{2}\,j_4^{2}\,j_7 -
      31090741920\,j_2^{2}\,j_4^{3}\,j_7 - 69044814720\,j_4^{4}\,j_7 -
      619759350\,j_2^{4}\,j_3\,j_5\,j_7 + 18592780500\,j_2\,j_3^{3}\,j_5\,j_7
      + 196807528575\,j_2^{2}\,j_3\,j_4\,j_5\,j_7 +
      912744012600\,j_3\,j_4^{2}\,j_5\,j_7 - 707275800\,j_2^{3}\,j_5^{2}\,j_7
      - 890203041000\,j_3^{2}\,j_5^{2}\,j_7 -
      118918781100\,j_2\,j_4\,j_5^{2}\,j_7 + 3083913000\,j_2^{5}\,j_6\,j_7 -
      92517390000\,j_2^{2}\,j_3^{2}\,j_6\,j_7 -
      64430478000\,j_2^{3}\,j_4\,j_6\,j_7 -
      617105790000\,j_3^{2}\,j_4\,j_6\,j_7 -
      98100072000\,j_2\,j_4^{2}\,j_6\,j_7 -
      333544837500\,j_2\,j_3\,j_5\,j_6\,j_7 - 202538070000\,j_5^{2}\,j_6\,j_7
      - 3214890000\,j_2^{3}\,j_3\,j_7^{2} + 96446700000\,j_3^{3}\,j_7^{2} +
      2411167500\,j_2\,j_3\,j_4\,j_7^{2} - 83349000\,j_2^{6}\,j_3\,j_8 +
      5000940000\,j_2^{3}\,j_3^{3}\,j_8 - 75014100000\,j_3^{5}\,j_8 +
      18237356550\,j_2^{4}\,j_3\,j_4\,j_8 -
      547120696500\,j_2\,j_3^{3}\,j_4\,j_8 -
      1115245341000\,j_2^{2}\,j_3\,j_4^{2}\,j_8 -
      2316980408400\,j_3\,j_4^{3}\,j_8 - 1360970100\,j_2^{5}\,j_5\,j_8 +
      40829103000\,j_2^{2}\,j_3^{2}\,j_5\,j_8 -
      19876440150\,j_2^{3}\,j_4\,j_5\,j_8 +
      2254694211000\,j_3^{2}\,j_4\,j_5\,j_8 +
      323344450800\,j_2\,j_4^{2}\,j_5\,j_8 +
      1116129435750\,j_2\,j_3\,j_5^{2}\,j_8 - 286446699000\,j_5^{3}\,j_8 +
      5358150000\,j_2^{3}\,j_3\,j_6\,j_8 - 160744500000\,j_3^{3}\,j_6\,j_8 +
      962859555000\,j_2\,j_3\,j_4\,j_6\,j_8 +
      417131977500\,j_2^{2}\,j_5\,j_6\,j_8 + 2615083380000\,j_4\,j_5\,j_6\,j_8
      - 803722500\,j_2^{4}\,j_7\,j_8 + 24111675000\,j_2\,j_3^{2}\,j_7\,j_8 +
      33756345000\,j_2^{2}\,j_4\,j_7\,j_8 - 347208120000\,j_4^{2}\,j_7\,j_8 -
      29569050\,j_2^{7}\,j_9 - 220279500\,j_2^{4}\,j_3^{2}\,j_9 +
      33220530000\,j_2\,j_3^{4}\,j_9 + 2362518900\,j_2^{5}\,j_4\,j_9 +
      18873445500\,j_2^{2}\,j_3^{2}\,j_4\,j_9 -
      87591974400\,j_2^{3}\,j_4^{2}\,j_9 +
      1071201348000\,j_3^{2}\,j_4^{2}\,j_9 - 237424219200\,j_2\,j_4^{3}\,j_9 -
      15994077750\,j_2^{3}\,j_3\,j_5\,j_9 - 327918780000\,j_3^{3}\,j_5\,j_9 -
      137367657000\,j_2\,j_3\,j_4\,j_5\,j_9 +
      337081216500\,j_2^{2}\,j_5^{2}\,j_9 + 957302388000\,j_4\,j_5^{2}\,j_9 -
      30148524000\,j_2^{4}\,j_6\,j_9 - 745854480000\,j_2\,j_3^{2}\,j_6\,j_9 -
      354985092000\,j_2^{2}\,j_4\,j_6\,j_9 - 890616384000\,j_4^{2}\,j_6\,j_9 +
      72335025000\,j_2^{2}\,j_3\,j_7\,j_9 - 810152280000\,j_3\,j_4\,j_7\,j_9 +
      1994422500\,j_2^{5}\,j_3\,j_{10} - 59832675000\,j_2^{2}\,j_3^{3}\,j_{10}
      - 77458947300\,j_2^{3}\,j_3\,j_4\,j_{10} -
      368701956000\,j_3^{3}\,j_4\,j_{10} -
      2576504700000\,j_2\,j_3\,j_4^{2}\,j_{10} +
      12854201850\,j_2^{4}\,j_5\,j_{10} +
      422115057000\,j_2\,j_3^{2}\,j_5\,j_{10} -
      151393762800\,j_2^{2}\,j_4\,j_5\,j_{10} -
      1301534438400\,j_4^{2}\,j_5\,j_{10} +
      2742944148000\,j_3\,j_5^{2}\,j_{10} +
      1650310200000\,j_2^{2}\,j_3\,j_6\,j_{10} +
      3538216080000\,j_3\,j_4\,j_6\,j_{10} -
      1234517760000\,j_2\,j_5\,j_6\,j_{10} - 27326565000\,j_2^{3}\,j_7\,j_{10}
      - 1350253800000\,j_3^{2}\,j_7\,j_{10} +
      462944160000\,j_2\,j_4\,j_7\,j_{10} + 2314720800000\,j_6\,j_7\,j_{10},
    \end{dmath*}
    \begin{dmath*}
      0= 2296\,j_2^{10}\,j_3 - 206640\,j_2^{7}\,j_3^{3} +
      6199200\,j_2^{4}\,j_3^{5} - 61992000\,j_2\,j_3^{7} +
      146790\,j_2^{8}\,j_3\,j_4 - 8807400\,j_2^{5}\,j_3^{3}\,j_4 +
      132111000\,j_2^{2}\,j_3^{5}\,j_4 - 35779140\,j_2^{6}\,j_3\,j_4^{2} +
      1331969400\,j_2^{3}\,j_3^{3}\,j_4^{2} - 7757856000\,j_3^{5}\,j_4^{2} +
      1278861210\,j_2^{4}\,j_3\,j_4^{3} - 16894963800\,j_2\,j_3^{3}\,j_4^{3} -
      35787922200\,j_2^{2}\,j_3\,j_4^{4} + 13355466624\,j_3\,j_4^{5} +
      2789640\,j_2^{6}\,j_3^{2}\,j_5 - 167378400\,j_2^{3}\,j_3^{4}\,j_5 +
      2510676000\,j_3^{6}\,j_5 + 5556600\,j_2^{7}\,j_4\,j_5 -
      214496100\,j_2^{4}\,j_3^{2}\,j_4\,j_5 +
      1433943000\,j_2\,j_3^{4}\,j_4\,j_5 - 584809875\,j_2^{5}\,j_4^{2}\,j_5 +
      6544779750\,j_2^{2}\,j_3^{2}\,j_4^{2}\,j_5 +
      15707821320\,j_2^{3}\,j_4^{3}\,j_5 - 56479056300\,j_3^{2}\,j_4^{3}\,j_5
      - 18236693160\,j_2\,j_4^{4}\,j_5 + 37660140\,j_2^{5}\,j_3\,j_5^{2} -
      1129804200\,j_2^{2}\,j_3^{3}\,j_5^{2} -
      428269275\,j_2^{3}\,j_3\,j_4\,j_5^{2} +
      36925308000\,j_3^{3}\,j_4\,j_5^{2} +
      122126125950\,j_2\,j_3\,j_4^{2}\,j_5^{2} +
      5084118900\,j_2\,j_3^{2}\,j_5^{3} - 31744742400\,j_2^{2}\,j_4\,j_5^{3} +
      54897855300\,j_4^{2}\,j_5^{3} - 137271210300\,j_3\,j_5^{4} +
      8906940\,j_2^{7}\,j_3\,j_6 - 534416400\,j_2^{4}\,j_3^{3}\,j_6 +
      8016246000\,j_2\,j_3^{5}\,j_6 - 1079717850\,j_2^{5}\,j_3\,j_4\,j_6 +
      32391535500\,j_2^{2}\,j_3^{3}\,j_4\,j_6 +
      32992395120\,j_2^{3}\,j_3\,j_4^{2}\,j_6 -
      73249361100\,j_3^{3}\,j_4^{2}\,j_6 - 94510476000\,j_2\,j_3\,j_4^{3}\,j_6
      - 227775240\,j_2^{6}\,j_5\,j_6 + 7217853300\,j_2^{3}\,j_3^{2}\,j_5\,j_6
      - 11537883000\,j_3^{4}\,j_5\,j_6 + 9344669490\,j_2^{4}\,j_4\,j_5\,j_6 -
      102778101450\,j_2\,j_3^{2}\,j_4\,j_5\,j_6 -
      5497768080\,j_2^{2}\,j_4^{2}\,j_5\,j_6 - 44830000800\,j_4^{3}\,j_5\,j_6
      - 135755619300\,j_2^{2}\,j_3\,j_5^{2}\,j_6 +
      425249246700\,j_3\,j_4\,j_5^{2}\,j_6 + 92816170650\,j_2\,j_5^{3}\,j_6 +
      1242297000\,j_2^{4}\,j_3\,j_6^{2} - 37268910000\,j_2\,j_3^{3}\,j_6^{2} +
      118349505000\,j_2^{2}\,j_3\,j_4\,j_6^{2} -
      285039000000\,j_3\,j_4^{2}\,j_6^{2} - 12310987500\,j_2^{3}\,j_5\,j_6^{2}
      + 355781160000\,j_3^{2}\,j_5\,j_6^{2} -
      46933020000\,j_2\,j_4\,j_5\,j_6^{2} - 1893780\,j_2^{8}\,j_7 +
      113626800\,j_2^{5}\,j_3^{2}\,j_7 - 1704402000\,j_2^{2}\,j_3^{4}\,j_7 +
      283936590\,j_2^{6}\,j_4\,j_7 - 10424748600\,j_2^{3}\,j_3^{2}\,j_4\,j_7 +
      57199527000\,j_3^{4}\,j_4\,j_7 - 4568759640\,j_2^{4}\,j_4^{2}\,j_7 +
      218617659450\,j_2\,j_3^{2}\,j_4^{2}\,j_7 -
      3463683120\,j_2^{2}\,j_4^{3}\,j_7 + 8077903200\,j_4^{4}\,j_7 -
      43630650\,j_2^{4}\,j_3\,j_5\,j_7 + 1308919500\,j_2\,j_3^{3}\,j_5\,j_7 +
      109249999425\,j_2^{2}\,j_3\,j_4\,j_5\,j_7 -
      95854241700\,j_3\,j_4^{2}\,j_5\,j_7 - 454677300\,j_2^{3}\,j_5^{2}\,j_7 -
      367668598500\,j_3^{2}\,j_5^{2}\,j_7 -
      72355692150\,j_2\,j_4\,j_5^{2}\,j_7 + 1801359000\,j_2^{5}\,j_6\,j_7 -
      54040770000\,j_2^{2}\,j_3^{2}\,j_6\,j_7 -
      43461522000\,j_2^{3}\,j_4\,j_6\,j_7 -
      43226055000\,j_3^{2}\,j_4\,j_6\,j_7 +
      80450982000\,j_2\,j_4^{2}\,j_6\,j_7 -
      155003625000\,j_2\,j_3\,j_5\,j_6\,j_7 + 6200145000\,j_5^{2}\,j_6\,j_7 -
      1607445000\,j_2^{3}\,j_3\,j_7^{2} + 48223350000\,j_3^{3}\,j_7^{2} +
      17567077500\,j_2\,j_3\,j_4\,j_7^{2} - 34870500\,j_2^{6}\,j_3\,j_8 +
      2092230000\,j_2^{3}\,j_3^{3}\,j_8 - 31383450000\,j_3^{5}\,j_8 +
      8471745450\,j_2^{4}\,j_3\,j_4\,j_8 -
      254152363500\,j_2\,j_3^{3}\,j_4\,j_8 -
      574638624000\,j_2^{2}\,j_3\,j_4^{2}\,j_8 -
      292701956400\,j_3\,j_4^{3}\,j_8 - 766980900\,j_2^{5}\,j_5\,j_8 +
      23009427000\,j_2^{2}\,j_3^{2}\,j_5\,j_8 -
      9587261250\,j_2^{3}\,j_4\,j_5\,j_8 +
      871258153500\,j_3^{2}\,j_4\,j_5\,j_8 +
      379035531000\,j_2\,j_4^{2}\,j_5\,j_8 +
      554292963000\,j_2\,j_3\,j_5^{2}\,j_8 - 184144306500\,j_5^{3}\,j_8 +
      3291435000\,j_2^{3}\,j_3\,j_6\,j_8 - 98743050000\,j_3^{3}\,j_6\,j_8 +
      433780515000\,j_2\,j_3\,j_4\,j_6\,j_8 +
      226305292500\,j_2^{2}\,j_5\,j_6\,j_8 - 682015950000\,j_4\,j_5\,j_6\,j_8
      - 114817500\,j_2^{4}\,j_7\,j_8 + 3444525000\,j_2\,j_3^{2}\,j_7\,j_8 -
      2066715000\,j_2^{2}\,j_4\,j_7\,j_8 - 235605510000\,j_4^{2}\,j_7\,j_8 -
      16176510\,j_2^{7}\,j_9 - 139311900\,j_2^{4}\,j_3^{2}\,j_9 +
      18738216000\,j_2\,j_3^{4}\,j_9 + 1306674180\,j_2^{5}\,j_4\,j_9 +
      10745387100\,j_2^{2}\,j_3^{2}\,j_4\,j_9 -
      53108758080\,j_2^{3}\,j_4^{2}\,j_9 + 492668114400\,j_3^{2}\,j_4^{2}\,j_9
      + 99665264160\,j_2\,j_4^{3}\,j_9 - 10124607150\,j_2^{3}\,j_3\,j_5\,j_9 -
      145772298000\,j_3^{3}\,j_5\,j_9 - 68284263600\,j_2\,j_3\,j_4\,j_5\,j_9 +
      190468454400\,j_2^{2}\,j_5^{2}\,j_9 - 447898474800\,j_4\,j_5^{2}\,j_9 -
      17819676000\,j_2^{4}\,j_6\,j_9 - 402320520000\,j_2\,j_3^{2}\,j_6\,j_9 -
      108938844000\,j_2^{2}\,j_4\,j_6\,j_9 + 242494560000\,j_4^{2}\,j_6\,j_9 +
      31000725000\,j_2^{2}\,j_3\,j_7\,j_9 + 90935460000\,j_3\,j_4\,j_7\,j_9 +
      744017400000\,j_7^{2}\,j_9 + 1109902500\,j_2^{5}\,j_3\,j_{10} -
      33297075000\,j_2^{2}\,j_3^{3}\,j_{10} -
      47211425100\,j_2^{3}\,j_3\,j_4\,j_{10} -
      82025622000\,j_3^{3}\,j_4\,j_{10} -
      1071660618000\,j_2\,j_3\,j_4^{2}\,j_{10} +
      7065868950\,j_2^{4}\,j_5\,j_{10} +
      237534444000\,j_2\,j_3^{2}\,j_5\,j_{10} -
      93498186600\,j_2^{2}\,j_4\,j_5\,j_{10} +
      571405363200\,j_4^{2}\,j_5\,j_{10} + 1071385056000\,j_3\,j_5^{2}\,j_{10}
      + 936910800000\,j_2^{2}\,j_3\,j_6\,j_{10} -
      705438720000\,j_3\,j_4\,j_6\,j_{10} -
      744017400000\,j_2\,j_5\,j_6\,j_{10} - 13089195000\,j_2^{3}\,j_7\,j_{10}
      - 537345900000\,j_3^{2}\,j_7\,j_{10} -
      140536620000\,j_2\,j_4\,j_7\,j_{10},
    \end{dmath*}
    \begin{dmath*}
      0= -3136\,j_2^{12} + 282240\,j_2^{9}\,j_3^{2} -
      8467200\,j_2^{6}\,j_3^{4} + 84672000\,j_2^{3}\,j_3^{6} +
      281232\,j_2^{10}\,j_4 - 23319030\,j_2^{7}\,j_3^{2}\,j_4 +
      639815400\,j_2^{4}\,j_3^{4}\,j_4 - 5800599000\,j_2\,j_3^{6}\,j_4 +
      2770530\,j_2^{8}\,j_4^{2} + 226037700\,j_2^{5}\,j_3^{2}\,j_4^{2} -
      9274608000\,j_2^{2}\,j_3^{4}\,j_4^{2} - 640259640\,j_2^{6}\,j_4^{3} +
      6366498750\,j_2^{3}\,j_3^{2}\,j_4^{3} - 8138799000\,j_3^{4}\,j_4^{3} +
      20176761420\,j_2^{4}\,j_4^{4} + 363254211900\,j_2\,j_3^{2}\,j_4^{4} -
      4691010024\,j_2^{2}\,j_4^{5} - 21019554432\,j_4^{6} -
      3810240\,j_2^{8}\,j_3\,j_5 + 228614400\,j_2^{5}\,j_3^{3}\,j_5 -
      3429216000\,j_2^{2}\,j_3^{5}\,j_5 + 78287580\,j_2^{6}\,j_3\,j_4\,j_5 -
      3083856300\,j_2^{3}\,j_3^{3}\,j_4\,j_5 + 22056867000\,j_3^{5}\,j_4\,j_5
      + 11236644405\,j_2^{4}\,j_3\,j_4^{2}\,j_5 -
      178230184650\,j_2\,j_3^{3}\,j_4^{2}\,j_5 -
      277610344200\,j_2^{2}\,j_3\,j_4^{3}\,j_5 +
      181783133820\,j_3\,j_4^{4}\,j_5 - 51438240\,j_2^{7}\,j_5^{2} +
      1543147200\,j_2^{4}\,j_3^{2}\,j_5^{2} +
      2728267920\,j_2^{5}\,j_4\,j_5^{2} -
      70957597725\,j_2^{2}\,j_3^{2}\,j_4\,j_5^{2} -
      136382752485\,j_2^{3}\,j_4^{2}\,j_5^{2} -
      186183465300\,j_3^{2}\,j_4^{2}\,j_5^{2} -
      163758820230\,j_2\,j_4^{3}\,j_5^{2} - 6944162400\,j_2^{3}\,j_3\,j_5^{3}
      + 674679111750\,j_2\,j_3\,j_4\,j_5^{3} + 187492384800\,j_2^{2}\,j_5^{4}
      + 380564900100\,j_4\,j_5^{4} - 6131664\,j_2^{9}\,j_6 +
      329995890\,j_2^{6}\,j_3^{2}\,j_6 - 3244260600\,j_2^{3}\,j_3^{4}\,j_6 -
      34113555000\,j_3^{6}\,j_6 + 636743646\,j_2^{7}\,j_4\,j_6 -
      11369425410\,j_2^{4}\,j_3^{2}\,j_4\,j_6 -
      231986519100\,j_2\,j_3^{4}\,j_4\,j_6 -
      18847621944\,j_2^{5}\,j_4^{2}\,j_6 -
      192900955680\,j_2^{2}\,j_3^{2}\,j_4^{2}\,j_6 +
      92875822776\,j_2^{3}\,j_4^{3}\,j_6 + 122638780620\,j_3^{2}\,j_4^{3}\,j_6
      + 44658819936\,j_2\,j_4^{4}\,j_6 + 5036422860\,j_2^{5}\,j_3\,j_5\,j_6 -
      151092685800\,j_2^{2}\,j_3^{3}\,j_5\,j_6 -
      146045495610\,j_2^{3}\,j_3\,j_4\,j_5\,j_6 +
      363873752550\,j_3^{3}\,j_4\,j_5\,j_6 -
      797037628320\,j_2\,j_3\,j_4^{2}\,j_5\,j_6 +
      30644331480\,j_2^{4}\,j_5^{2}\,j_6 +
      2333407348125\,j_2\,j_3^{2}\,j_5^{2}\,j_6 -
      392227832115\,j_2^{2}\,j_4\,j_5^{2}\,j_6 -
      691147287360\,j_4^{2}\,j_5^{2}\,j_6 - 105898476600\,j_3\,j_5^{3}\,j_6 -
      1263729600\,j_2^{6}\,j_6^{2} + 37233189000\,j_2^{3}\,j_3^{2}\,j_6^{2} +
      20360970000\,j_3^{4}\,j_6^{2} + 2303834400\,j_2^{4}\,j_4\,j_6^{2} -
      2508605883000\,j_2\,j_3^{2}\,j_4\,j_6^{2} +
      82826841600\,j_2^{2}\,j_4^{2}\,j_6^{2} + 422227468800\,j_4^{3}\,j_6^{2}
      - 68209249500\,j_2^{2}\,j_3\,j_5\,j_6^{2} -
      419956488000\,j_3\,j_4\,j_5\,j_6^{2} +
      459086292000\,j_2\,j_5^{2}\,j_6^{2} + 48977460\,j_2^{7}\,j_3\,j_7 -
      2938647600\,j_2^{4}\,j_3^{3}\,j_7 + 44079714000\,j_2\,j_3^{5}\,j_7 -
      4568596830\,j_2^{5}\,j_3\,j_4\,j_7 +
      137057904900\,j_2^{2}\,j_3^{3}\,j_4\,j_7 -
      89936802900\,j_2^{3}\,j_3\,j_4^{2}\,j_7 -
      581967807750\,j_3^{3}\,j_4^{2}\,j_7 +
      218629053720\,j_2\,j_3\,j_4^{3}\,j_7 + 838252800\,j_2^{6}\,j_5\,j_7 -
      29500187850\,j_2^{3}\,j_3^{2}\,j_5\,j_7 +
      130578115500\,j_3^{4}\,j_5\,j_7 - 34356457800\,j_2^{4}\,j_4\,j_5\,j_7 -
      2152360817775\,j_2\,j_3^{2}\,j_4\,j_5\,j_7 +
      145612472040\,j_2^{2}\,j_4^{2}\,j_5\,j_7 +
      159765336360\,j_4^{3}\,j_5\,j_7 +
      539908626600\,j_2^{2}\,j_3\,j_5^{2}\,j_7 +
      598810004100\,j_3\,j_4\,j_5^{2}\,j_7 - 515604058200\,j_2\,j_5^{3}\,j_7 -
      42615153000\,j_2^{4}\,j_3\,j_6\,j_7 +
      1278454590000\,j_2\,j_3^{3}\,j_6\,j_7 +
      1009261134000\,j_2^{2}\,j_3\,j_4\,j_6\,j_7 -
      301261437000\,j_3\,j_4^{2}\,j_6\,j_7 +
      162244782000\,j_2^{3}\,j_5\,j_6\,j_7 +
      618062602500\,j_3^{2}\,j_5\,j_6\,j_7 +
      507240751500\,j_2\,j_4\,j_5\,j_6\,j_7 + 428652000\,j_2^{5}\,j_7^{2} -
      12859560000\,j_2^{2}\,j_3^{2}\,j_7^{2} +
      8588349000\,j_2^{3}\,j_4\,j_7^{2} - 221482957500\,j_3^{2}\,j_4\,j_7^{2}
      - 140949963000\,j_2\,j_4^{2}\,j_7^{2} + 44135280\,j_2^{8}\,j_8 -
      2648116800\,j_2^{5}\,j_3^{2}\,j_8 + 39721752000\,j_2^{2}\,j_3^{4}\,j_8 -
      8768008620\,j_2^{6}\,j_4\,j_8 + 246336863850\,j_2^{3}\,j_3^{2}\,j_4\,j_8
      + 501101842500\,j_3^{4}\,j_4\,j_8 + 406754822880\,j_2^{4}\,j_4^{2}\,j_8
      + 4176898374600\,j_2\,j_3^{2}\,j_4^{2}\,j_8 -
      67832341920\,j_2^{2}\,j_4^{3}\,j_8 - 262820013120\,j_4^{4}\,j_8 +
      15592216500\,j_2^{4}\,j_3\,j_5\,j_8 -
      467766495000\,j_2\,j_3^{3}\,j_5\,j_8 -
      687230960850\,j_2^{2}\,j_3\,j_4\,j_5\,j_8 -
      2776052922300\,j_3\,j_4^{2}\,j_5\,j_8 -
      499015225800\,j_2^{3}\,j_5^{2}\,j_8 -
      1762804559250\,j_3^{2}\,j_5^{2}\,j_8 +
      1428947418150\,j_2\,j_4\,j_5^{2}\,j_8 - 3143448000\,j_2^{5}\,j_6\,j_8 +
      94303440000\,j_2^{2}\,j_3^{2}\,j_6\,j_8 -
      396610263000\,j_2^{3}\,j_4\,j_6\,j_8 -
      200930625000\,j_3^{2}\,j_4\,j_6\,j_8 +
      324520182000\,j_2\,j_4^{2}\,j_6\,j_8 -
      6907994887500\,j_2\,j_3\,j_5\,j_6\,j_8 -
      3211675110000\,j_5^{2}\,j_6\,j_8 + 4018612500\,j_2^{3}\,j_3\,j_7\,j_8 -
      120558375000\,j_3^{3}\,j_7\,j_8 - 188071065000\,j_2\,j_3\,j_4\,j_7\,j_8
      + 773359650\,j_2^{6}\,j_3\,j_9 - 20378830500\,j_2^{3}\,j_3^{3}\,j_9 -
      84658770000\,j_3^{5}\,j_9 - 53333800380\,j_2^{4}\,j_3\,j_4\,j_9 +
      428990328900\,j_2\,j_3^{3}\,j_4\,j_9 +
      900564172200\,j_2^{2}\,j_3\,j_4^{2}\,j_9 -
      894518341920\,j_3\,j_4^{3}\,j_9 - 2537619840\,j_2^{5}\,j_5\,j_9 +
      427435699950\,j_2^{2}\,j_3^{2}\,j_5\,j_9 +
      63644718060\,j_2^{3}\,j_4\,j_5\,j_9 +
      753083389800\,j_3^{2}\,j_4\,j_5\,j_9 +
      797818124880\,j_2\,j_4^{2}\,j_5\,j_9 -
      5519741087700\,j_2\,j_3\,j_5^{2}\,j_9 - 1218700501200\,j_5^{3}\,j_9 +
      716277492000\,j_2^{3}\,j_3\,j_6\,j_9 + 1324534680000\,j_3^{3}\,j_6\,j_9
      + 5401566324000\,j_2\,j_3\,j_4\,j_6\,j_9 - 6429780000\,j_2^{4}\,j_7\,j_9
      - 110913705000\,j_2\,j_3^{2}\,j_7\,j_9 -
      58143582000\,j_2^{2}\,j_4\,j_7\,j_9 + 585293688000\,j_4^{2}\,j_7\,j_9 -
      469135800\,j_2^{7}\,j_{10} + 2125399500\,j_2^{4}\,j_3^{2}\,j_{10} +
      358460235000\,j_2\,j_3^{4}\,j_{10} + 21321354600\,j_2^{5}\,j_4\,j_{10} +
      531383044500\,j_2^{2}\,j_3^{2}\,j_4\,j_{10} +
      951209406000\,j_2^{3}\,j_4^{2}\,j_{10} +
      1236538548000\,j_3^{2}\,j_4^{2}\,j_{10} -
      673032628800\,j_2\,j_4^{3}\,j_{10} -
      325105751250\,j_2^{3}\,j_3\,j_5\,j_{10} -
      786040605000\,j_3^{3}\,j_5\,j_{10} +
      4705910055000\,j_2\,j_3\,j_4\,j_5\,j_{10} -
      1332411160500\,j_2^{2}\,j_5^{2}\,j_{10} -
      133923132000\,j_4\,j_5^{2}\,j_{10} - 336063168000\,j_2^{4}\,j_6\,j_{10}
      - 12730964400000\,j_2\,j_3^{2}\,j_6\,j_{10} +
      263069856000\,j_2^{2}\,j_4\,j_6\,j_{10} +
      1807686720000\,j_4^{2}\,j_6\,j_{10} +
      2893401000000\,j_3\,j_5\,j_6\,j_{10} +
      303807105000\,j_2^{2}\,j_3\,j_7\,j_{10} +
      3480348060000\,j_3\,j_4\,j_7\,j_{10} + 5208121800000\,j_7^{2}\,j_{10},
    \end{dmath*}
  \end{dgroup*}
  \setcounter{equation}{\themyequation}
}

\end{document}